\documentclass[a4paper,11pt,oneside,dvipsnames,reqno]{article} %

\usepackage{amsmath,amsfonts,amsthm} %

\usepackage[UKenglish]{babel}

\usepackage[T1]{fontenc}
\usepackage{newpxtext,newpxmath}
\usepackage{microtype}

\usepackage{xcolor}
\usepackage{graphicx} %
\usepackage{mathtools}

\usepackage{hyperref}
\hypersetup{colorlinks=true, citecolor=PineGreen, linkcolor=RoyalBlue, linktoc=page,urlcolor=RoyalBlue}

\usepackage{epigraph} %
\usepackage[autostyle]{csquotes} %

\usepackage[shortlabels]{enumitem}

\theoremstyle{plain}
\newtheorem{theorem}{Theorem}[section] %
\newtheorem{lemmy}[theorem]{Lemma}
\newtheorem{prop}[theorem]{Proposition}
\newtheorem{cor}[theorem]{Corollary}
\newtheorem{ass}{Assumption}[section]

\newtheorem{conv}{Convention}[section]

\theoremstyle{remark}
\newtheorem{rem}[theorem]{Remark}

\theoremstyle{definition}
\newtheorem{defn}{Definition}[section]

\numberwithin{equation}{section}

\newcommand{\ML}{\mathrm{ML}}
\newcommand{\MV}{\mathrm{MV}}
\newcommand{\vnew}{v_{\textrm{new}}}

\newcommand{\parab}{\mathrm{par}}

\newcommand{\A}{\mathbb{A}}

\newcommand{\R}{\mathbb{R}}

\newcommand{\C}{\mathbb{C}}

\newcommand{\Hb}{\mathbb{H}}

\newcommand{\Z}{\mathbb{Z}}

\newcommand{\Q}{\mathbb{Q}}

\newcommand{\Bc}{\mathcal{B}}

\newcommand{\Dc}{\mathcal{D}}

\newcommand{\Nc}{\mathcal{N}}
\newcommand{\Oc}{\mathcal{O}}
\newcommand{\Pc}{\mathcal{P}}

\newcommand{\Uc}{\mathcal{U}}

\newcommand{\gf}{\mathfrak{g}}

\newcommand{\pf}{\mathfrak{p}}
\newcommand{\af}{\mathfrak{a}}

\newcommand{\hrm}{\mathrm{h}}

\newcommand{\GL}{\operatorname{GL}}
\newcommand{\PGL}{\operatorname{PGL}}
\newcommand{\SL}{\operatorname{SL}}
\newcommand{\PSL}{\operatorname{PSL}}

\newcommand{\Mat}{\operatorname{Mat}}

\newcommand{\SO}{\operatorname{SO}}
\newcommand{\PSO}{\operatorname{PSO}}

\newcommand{\St}{\widetilde{S}}

\newcommand{\tr}{\operatorname{tr}}
\newcommand{\Tr}{\operatorname{Tr}}

\newcommand{\N}{\mathbb{N}}

\newcommand{\norm}[1]{\left\lVert#1\right\rVert}

\renewcommand{\Im}{\operatorname{Im}}

\newcommand{\Hom}{\operatorname{Hom}}

\newcommand{\sgn}{\operatorname{sgn}}

\newcommand{\Ad}{\operatorname{Ad}}
\newcommand{\ad}{\operatorname{ad}}

\newcommand{\Vol}{\operatorname{vol}}

\newcommand{\Lie}{\operatorname{Lie}}
\newcommand{\diag}{\operatorname{diag}}

\newcommand{\fin}{\operatorname{f}} %

\newcommand{\Op}{\operatorname{Op}}
\newcommand{\Proj}{\operatorname{Proj}}
\DeclareMathOperator{\dist}{dist}

\newcommand{\mat}[4]{\begin{pmatrix}
#1 & #2 \\ #3 & #4
\end{pmatrix}}

\newcommand{\abs}[1]{\lvert{#1}\rvert}

\usepackage{authblk} %

\graphicspath{ {img/} }

\usepackage[final]{todonotes}

\usepackage[style = alphabetic, backend = biber, abbreviate = true, arxiv = abs, doi = false, isbn = false, url = false, giveninits = true]{biblatex} %
\renewbibmacro{in:}{} %

\title{The orbit method in number theory through the sup-norm problem for $\GL(2)$}

\author[1]{Edgar Assing\thanks{email: assing@math.uni-bonn.de}}
\author[2]{Radu Toma\thanks{email: toma@imj-prg.fr}}

\affil[1]{Mathematisches Institut, Universität Bonn}
\affil[2]{Sorbonne Université, Université Paris Cité, CNRS, IMJ-PRG}

\begin{document}

\maketitle

\begin{abstract}
    The orbit method in its quantitative form due to Nelson and Venkatesh has played a central role in recent advances in the analytic theory of higher rank $L$-functions. 
    The goal of this note is to explain how the method can be applied to the sup-norm problem for automorphic forms on $\PGL(2)$. 
    Doing so, we prove a new hybrid bound for newforms $\varphi$ of large prime-power level $N = p^{4n}$ and large eigenvalue $\lambda$. 
    It states that $\| \varphi \|_\infty \ll_p (\lambda N)^{5/24 + \varepsilon}$, recovering the result of Iwaniec and Sarnak spectrally and improving the local bound in the depth aspect for the first time in this non-compact setting.
    We also provide an exposition of the microlocal tools used, illustrating and motivating the theory through the classical case of $\PGL(2)$, following notes and lectures of Nelson and Venkatesh.
    \end{abstract}

\tableofcontents %

\renewcommand{\epigraphsize}{\footnotesize}
\epigraph{In orbit, you’re keyed up and aware of everything going on, every little noise, anything that may have special meaning because of where you are.}{John Glenn, astronaut}

\section{Introduction}

Analysing special functions, such as Bessel or Whittaker functions, and associated integral transforms is part of the bread and butter of analytic number theorists working with automorphic forms. 
This, however, is often highly technical and can quickly become intractable. 
It is particularly the case when working with higher rank groups and, as a result, strong analytic results are sparse in these settings at the moment.

In a recent breakthrough, Nelson and Venkatesh \cite{NV} built on a quantitative form of Kirillov's orbit method and introduced a flexible microlocal calculus for representations on reductive groups. 
This is a vast collection of tools tailored for analysing periods of automorphic forms in a general setting. 
Applied to the special periods arising through GGP pairs, this has lead to ground-breaking results concerning high degree $L$-functions. 
In particular, following further developments by Nelson \cite{Nel-Un, Nel-GLn}, it was possible to prove a general subconvexity bound for $\GL(n)$ standard $L$-functions in the spectral aspect.

The strength of this new methodology is that it allows for a \textit{geometrisation} of the hard analytic problems that usually arise in the study of these periods. 
Its impact on the modern theory of automorphic forms is profound and, while most applications are currently aimed at the study of $L$-functions, many more beyond this are expected.

In this note we set out to apply this new microlocal calculus in the setting of $\PGL(2)$ to prove strong bounds on the sup-norm of Hecke--Maaß forms. 
Doing so, we give an overview of the methods of Marshall, Nelson, Venkatesh, and others, in this more classical and familiar context. 
We hope that this treatment can serve as a good introduction for analytic number theorists who wish to learn more about them.
The rank 1 case of $\PGL(2)$ does not exhibit all the subtleties and difficulties that arise more generally, but it still captures much of the essence of these ideas.

Section \ref{sec:archi-microlocalisation-big-chap} and Section \ref{sec:p-adic} form the central expository part of this paper.
The former deals with the archimedean theory, particularly as developed in \cite{Nel-Un}.
There is no new mathematical content and its purpose is to illustrate and sketch out, in the case of $\PGL(2)$, the ideas appearing in \cite{Nel-Un}, which in turn are based on \cite{NV}.
We give additional intuition and context following lectures of Nelson and Venkatesh, given in Budapest \cite{Nel-budapest} and in Oberwolfach \cite{Nel-MFO, Ven}.
This section could be used as a stepping stone or a guide to reading the sizeable but beautifully written original papers \cite{NV, Nel-Un}, which we recommend and refer to for full details.

Section \ref{sec:p-adic} offers a new exposition of the $p$-adic side of the medal that parallels the archimedean side. We focus on explaining the key ideas using examples of microlocalised vectors that previously appeared in the literature. More precisely, we consider microlocal lift vectors that appear in the work of Nelson \cite{Nel-QUE-QP} and Marshall \cite{marsh}, as well as minimal vectors that were studied by Nelson and Hu in \cite{Ne-Hu_test}, for example. An interesting consequence of our approach are certain bounds for new-vector matrix coefficients. They are qualitatively worse than bounds that are already in the literature (see \cite{Hu-Sa}), but they follow from rather soft volume bounds. In particular, no elaborate stationary phase analysis is required.

The second half of this paper deals with the application to bounding Hecke--Maaß forms in the spectral and depth aspect and contains new results.
We were inspired to look at this problem through the lens of the orbit method thanks to a talk by Farrell Brumley in Oberwolfach.
After applying the analytic machinery mentioned above together with the amplification technique, we are faced with a counting problem for matrices, which has some of the classical features appearing in the literature, but presents new challenges.
For this, we introduce some new strategies to provide strong counting results uniform in all parameters.

Due to working with non-compact spaces, we require an additional analytic input coming from the analysis of the Whittaker expansion.
This is added here as an appendix, where we present a new, powerful bound in the spectral and depth aspect, which follows from the doctoral thesis of the first author.

We defer the introduction to the quantitative orbit method to the sketch of the proof in Section \ref{sec:sketch}.
In the rest of this section, we provide context for the sup-norm problem and state our new results.

\subsection{The main theorem}

Let us now carefully state our main result. Let $\Hb=\{ x+iy\in \C\colon y>0\}$ be the upper half plane and write $\Delta=-y^2(\frac{\partial^2}{\partial x^2}+\frac{\partial^2}{\partial y^2})$ be the Laplace-Beltrami operator. The group $\SL_2(\R)$ acts on $\Hb$ via M\"obius transformations, and we define the Hecke congruence subgroups $\Gamma_H(N)$ of level $N$ by
\begin{equation}
	\Gamma_H(N) = \{ \gamma\in \SL_2(\Z)\colon N\mid \gamma_{2,1} \}.\nonumber
\end{equation}
In particular, $\Gamma_H(1)=\SL_2(\Z)$.

Given a bounded smooth function $f\colon \Gamma_H(N)\backslash \Hb\to \C$ we define the $L^{\infty}$-norm by
\begin{equation}
	\Vert f\Vert_{\infty} = \sup_{z \in \Hb}\vert f(z)\vert.\nonumber
\end{equation}
We are interested in estimates for the size of $\Vert f\Vert_{\infty}$ when $f$ is a cuspidal Hecke--Maaß newform. This size will be measured in comparison to the $L^2$-norm
\begin{equation}
	\Vert f\Vert_2^2 = \frac{1}{\Vol(\Gamma_H(N)\backslash\Hb)}\int_{\Gamma_H(N)\backslash \Hb} \vert f(x+iy)\vert^2\frac{dxdy}{y^2}.\nonumber
\end{equation}
We prove the following theorem.

\begin{theorem}\label{th:Intro}
Let $p$ be an odd prime and $n$ a positive integer. Let $N=p^{4n}$ and let $f\colon \Gamma_H(N)\backslash \Hb\to \C$ be cuspidal Hecke--Maaß newform with trivial nebentypus and Laplace-Beltrami eigenvalue $\lambda$. 
If $\lambda$ is sufficiently large, then we have
\begin{equation}
	\Vert f\Vert_{\infty} \ll_{p,\epsilon} (N\lambda)^{\frac{5}{24}+\epsilon}\cdot \Vert f\Vert_2.\nonumber
\end{equation}
\end{theorem}

Taking $N=1$ (i.e. $n=0$) recovers the classical theorem of Iwaniec and Sarnak \cite{IS95}. This is a classical version of Theorem~\ref{th:main_theorem_impro} below. The assumptions on $\lambda$ and $N$ can be relaxed with some more work.\footnote{We do not quantify the term \textit{some more work} any further, as it depends on how ambitiously one formulates the corresponding result.} We have made them mostly for convenience to balance between technical novelty and the expository nature of this paper.

\subsection{A short history of the sup-norm problem for Hecke congruence subgroups}

The sup-norm problem has received a lot of attention in the last 30 years and a detailed survey goes beyond the scope of this text. However, it seems appropriate to summarise some of the milestones which are relevant to our main theorem. We focus on results concerning Hecke--Maaß newforms on the non-compact quotients $\Gamma_H(N)\backslash \Hb$. Their sup-norm can be studied in two aspects:
\begin{itemize}
	\item \textbf{The spectral aspect:} We consider $N$ as fixed and try to prove bounds of the form
	\begin{equation}
		\Vert f\Vert_{\infty} \ll_{N,\epsilon} \lambda^{\alpha+\epsilon}\cdot \Vert f\Vert_2.\nonumber
	\end{equation}
	where $f$ travels through a sequence of Hecke--Maaß cusp forms on $\Gamma_{H}(N)\backslash \Hb$ with growing Laplace-Beltrami eigenvalue $\lambda$.
	\item \textbf{The level aspect:} We fix a spectral window $I\subseteq \R_{+}$ and try to prove bounds of the form
	\begin{equation}
		\sup_{\lambda_f\in I} \Vert f\Vert_{\infty}\ll_{I,\epsilon} N^{\beta+\epsilon}\cdot \Vert f\Vert_2
	\end{equation}
	as $N\to \infty$.\footnote{Note that the $L^2$-norm on the right-hand side is changing with $N$!} Where the supremum is taken over all Hecke--Maaß newforms $f$ with Laplace-Beltrami eigenvalue $\lambda_f$ contained in the fixed spectral window $I$.
\end{itemize}
An honest combination of these two aspects will be referred to as hybrid bound.

The local bound (or baseline bound) in the spectral aspect is $\alpha=\frac{1}{4}$ and holds in great generality. 
The first improvement was achieved by Iwaniec and Sarnak who showed that $\alpha=\frac{5}{24}$ is admissible. 
Their key new idea was to apply the amplification method from $L$-function theory to achieve this improvement for arithmetic $f$, i.e. Hecke eigenforms.
Their breakthrough inspired a great deal of research subsequently, but the exponent in the spectral aspect has not been improved upon in general, so far. 

The level aspect is more subtle. Following the historical order of things, we first discuss squarefree levels $N$. 
Here the baseline bound is arguably $\beta=\frac{1}{2}$. 
The first improvement was achieved by Blomer and Holowinsky in \cite{BoHo}. 
They prove that $\beta=\frac{1}{2}-\frac{1}{37}$ is admissible for squarefree $N$. They also combine this with (a refined version of) the result of Iwaniec and Sarnak to obtain a hybrid bound.

This was greatly improved in a series of papers by Harcos and Templier \cite{HT}. They show that $\beta = \frac{1}{3}$ is allowed in the estimate, still in the context of squarefree $N$. In \cite{Te} the argument is extended to yield the hybrid bound
\begin{equation}
	\Vert f\Vert_{\infty} \ll \lambda^{\frac{5}{24}+\epsilon}N^{\frac{1}{3}+\epsilon}\cdot \Vert f\Vert_2\nonumber
\end{equation} 
for squarefree $N$. 
This has been the state of the art for a long time. 
Only very recently the exponent $\beta=\frac{1}{4}$ has become admissible for squarefree levels through the work of \cite{KNS}.

We now turn to powerful levels $N$. 
For simplicity, we only discuss $N=N_0^2$ (i.e. squares). 
Here, the local exponent is $\beta=\frac{1}{4}$ as in the spectral aspect. 
To the best of our knowledge, this was first observed by \cite{marsh_local}. 
As shown in \cite{Sa}, one can combine this with the Iwaniec--Sarnak exponent in the spectral aspect to obtain
\begin{equation}
	\Vert f\Vert_{\infty} \ll \lambda^{\frac{5}{24}+\epsilon}N^{\frac{1}{4}+\epsilon}\cdot \Vert f\Vert_2.\nonumber
\end{equation} 
For the Hecke congruence subgroup, our result is the first to improve the local bound in the powerful level aspect and we manage to combine the best known exponents into a hybrid bound.
In the setting of compact quotients, a similar bound was achieved by Hu and Saha \cite{Hu-Sa} in the level aspect.

\section{Sketch of the proof} \label{sec:sketch}

We start by presenting an informal sketch of our approach. This serves two purposes. First and foremost, it is road map for the rest of the paper. Second, it is meant to introduce and motivate the study of quantitative microlocal analysis on groups as an important tool for questions in analytic number theory.

We focus here on the archimedean place, for which certain objects can be easier to visualise.
Where it is adequate, we point out the analogies with the non-archimedean places.

Throughout this section we will use the notation $G = \PGL(2, \R)$, $K = \PSO(2)$ and we let $D$ be the subgroup of diagonal matrices in $G$. These groups are equipped with suitably normalised Haar measures.
Since we focus on the spectral aspect, we let $\Gamma = \PSL_2(\Z)$. 

In this sketch we study the sup-norm of an $L^2$-normalised classical Hecke--Maaß cusp form $\varphi_{\lambda,\circ}$ with Laplace eigenvalue $\lambda$, which we view as a function on $G$. In particular, since we are working with an eigenfunction of all Hecke operators, there is an irreducible constituent $\pi_{\lambda}$ of the right regular representation on $L^2(\PSL_2(\Z)\backslash G)$ such that $\varphi_{\lambda,\circ}\in \pi_{\lambda}$ is the (up to scaling) unique spherical (i.e. $K$-invariant) element. We will assume that $\lambda$ is sufficiently large, so that $\pi_{\lambda}$ is tempered.

\subsection{Uniqueness, periods, and matrix coefficients} \label{sec:sketch-K-period}
It is a standard fact that the space of spherical vectors, those invariant under the action of $K$, is one-dimensional.
For any $\varphi \in \pi_{\lambda}$, this implies that its $K$-period satisfies
\begin{equation*}
	\int_{K} \varphi(k) \, dk = \langle \varphi, \varphi_{\lambda, \circ} \rangle \cdot \varphi_{\lambda,\circ}(1).
\end{equation*}
The square of the inner product can be computed as
\begin{displaymath}
    \abs{\langle \varphi, \varphi_{\lambda, \circ} \rangle}^2 = \frac{1}{\Vol(K)} \int_K \langle \pi(k) \varphi, \varphi \rangle \, dk =: Q(\varphi),
\end{displaymath}
by decomposing $\varphi = \langle \varphi, \varphi_{\lambda, \circ} \rangle \varphi_{\lambda,\circ} + \sum_{k \neq 0} \varphi_{\lambda, k}$ into $K$-types (see e.g. \cite[Chap. 2, (5.2)]{bump}) and using orthogonality.
We therefore obtain that
\begin{equation} \label{eq:intro-period}
    \left| \int_{K} \varphi(k) \, dk \right|^2 = |\varphi_{\lambda,\circ}(1)|^2 \cdot Q(\varphi).
\end{equation}
The left-hand side of \eqref{eq:intro-period} is of global or automorphic nature, while matrix coefficients and their integrals, e.g. $Q(\varphi)$, can be studied locally.

\begin{rem}
    At the non-archimedean places, the uniqueness of the spherical vector is replaced by the theory of newforms, consisting also of a multiplicity-one statement (see Section \ref{sec:new-vectors}).
    In the subconvexity problem, the analogue is provided by the theory of GGP pairs (see \cite[Sec. 2.6]{Nel-Un}).
\end{rem}

\subsection{The relative trace formula} \label{sec:intro-rel-trace-formula}

Following the work of Iwaniec and Sarnak, instead of bounding individual periods, we consider a spectral average of \eqref{eq:intro-period}, as in \eqref{eq:intro-relative-pretrace-formula} below, and use the pretrace formula to translate to a ``geometric'' setting.

Let $f$ be a smooth, compactly supported function on $G$ and let $\kappa_f$ be the corresponding automorphic kernel
\begin{displaymath}
    \kappa_f(x, y) = \sum_{\gamma \in \Gamma} f(x^{-1} \gamma y).
\end{displaymath}
By spectral expansion, we obtain that
\begin{equation} \label{eq:intro-spec-expansion}
    \kappa_f(x, y) = \int_\pi \sum_{v \in \mathcal{B}(\pi)} \pi(f)v(x) \cdot \overline{v(y)} \, d\mu_{\text{spec}}(\pi).
\end{equation}
The integral denotes the average over all spectral components of $L^2(\Gamma \backslash G)$, i.e. cusp forms, residual spectrum, and Eisenstein series, $\mathcal{B}(\pi)$ denotes an orthonormal basis for $\pi$, and $\pi(f)$ is the convolution operator defined by $f$ on the representation $\pi$.

For heuristic purposes, assume we may take the basis $\mathcal{B}(\pi)$ to consist of eigenvectors of convolutions operators and then integrate over $K$ in both variables $x$ and $y$ to obtain an average of periods as in \eqref{eq:intro-period}.
Denoting the eigenvalue of $v \in \mathcal{B}(\pi)$ by $a_f(v)$, we have
\begin{equation} \label{eq:intro-relative-pretrace-formula}
    \int_K \int_K \kappa_f(x, y) = \int_\pi \sum_{v \in \mathcal{B}(\pi)} a_f(v) \cdot \left| \int_K v \right|^2.
\end{equation}
To isolate $\pi_\lambda$ and obtain the desired bounds, we would like $a_f(v)$ to be concentrated at representations near $\pi_\lambda$ and to be non-negative.
For the latter, the standard solution is to take the test function to be a convolution of the shape $f = \omega \ast \omega^\ast$, where $\omega^\ast(g) = \overline{\omega(g^{-1})}$. 
This implies that $a_f(v) = |a_\omega(v)|^2$.
We then obtain an inequality by removing all terms but one, using non-negativity, and applying \eqref{eq:intro-period}, namely
\begin{multline} \label{eq:intro-pretrace-ineq}
    |\varphi_{\lambda,\circ}(1)|^2 \cdot \left( \sum_{v \in \mathcal{B}(\pi_\lambda)} |a_\omega(v)|^2 Q(v) \right) \\ =
    \sum_{v \in \mathcal{B}(\pi_\lambda)} |a_\omega(v)|^2 \cdot \left| \int_K v \right|^2
    \leq \int_K \int_K \kappa_f(x, y).
\end{multline}
This clearly implies a bound for the quantity of interest $|\varphi_{\lambda,\circ}(1)|$, as long as we have some lower bounds for the matrix coefficient integrals $Q(v)$ for suitable vectors $v$ and also upper bounds for the \emph{geometric side}, i.e. the right-hand side of \eqref{eq:intro-pretrace-ineq}.
Note also that we wanted $a_f(v)$ to be supported at representations near $\pi_\lambda$ so that the inequality we obtain when dropping all terms but the ones corresponding to $\pi_\lambda$ is not too lossy.

In the next sections, we sketch out how such an argument could be implemented.
We gather the desiderata:
\begin{enumerate}[a)]
    \item Basis vectors in $\mathcal{B}(\pi)$ should be, at least approximately, eigenvectors for some class of convolution operators. \label{desideratum-eigenvectors}
    \item We should have some control over matrix coefficients of the basis vectors in $\mathcal{B}(\pi)$ and their integrals. \label{desideratum-matrix-coeffs}
    \item The test function $f$ should ``concentrate'' around $\pi_\lambda$ as much as possible. \label{desideratum-test-f-concentration}
    \item The test function $f$ should also produce an automorphic kernel $\kappa_f$ for which we can bound the geometric side. \label{desideratum-automorphic-kernel}
\end{enumerate}

\subsection{Localised vectors}
The orbit method enters the picture by suggesting the existence of a practical basis $\mathcal{B}(\pi)$, satisfying the desiderata in the previous section, that also has a useful geometric interpretation.
More broadly, it provides a meaningful connection between tempered representations and coadjoint orbits together with a symplectic measure.
We focus here on the principal series representations that occur in our problem.

First, the orbit method associates to $\pi_\lambda$ a coadjoint orbit that is a one-sheeted hyperboloid
\begin{displaymath}
    \mathcal{O}_{\pi_\lambda} = \{ (a,b,c) \in \R^3 \mid a^2 + b^2 - c^2 = T^2 \}
\end{displaymath}
where $\lambda = \frac14 + T^2$.
This hyperboloid is the orbit of the element $\tau = (T, 0, 0)$ under the coadjoint action of $G$ on the dual Lie algebra $\gf^\wedge = i \gf^\ast \cong \R^3$.
We discuss this in Section \ref{sec:archi-gp-lie-algebra}.

It is useful to note that $K$ acts on $\mathcal{O}_{\pi_\lambda}$ by rotation around the $c$-axis.
Parametrising $K$ by the angle $\phi \in [0, 2\pi)$, we get the representation of the hyperboloid in polar coordinates
\begin{displaymath}
    \mathcal{O}_{\pi_\lambda} = \{ T \cdot (\cosh v \cos \phi, \cosh v \sin \phi, \sinh v) \mid \phi \in [0,2\pi), v \in \R \}.
\end{displaymath}
There is an important natural symplectic measure on $\mathcal{O}_{\pi_\lambda}$ given by the form
\begin{displaymath}
    \frac{1}{4\pi} \cdot T \cosh(v)\, dv \wedge d\phi.
\end{displaymath}

Next, we partition the hyperboloid $\mathcal{O}_{\pi_\lambda}$ into ball-like pieces $U_i$ of symplectic volume one, centred around points $\tau_i$.
Assume one of the points is $\tau = (T, 0, 0)$ in $(a,b,c)$-coordinates.

\begin{figure}[ht]
    \centering
    \includegraphics[width=0.8\textwidth]{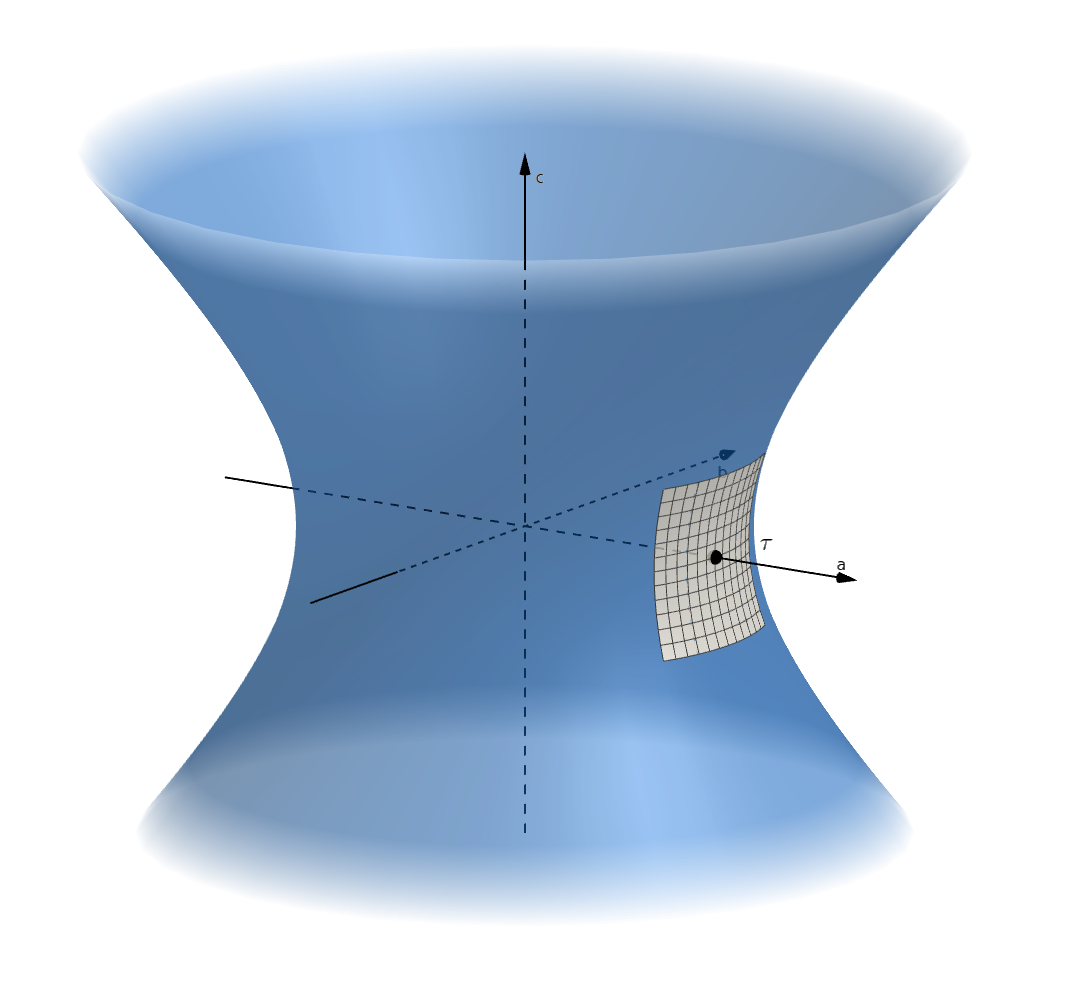}
    \caption{The region $U_\tau$ on the orbit $\Oc_{\pi_\lambda}$}
    \label{fig:localised-intro}
\end{figure}

Heuristically, we should have an orthonormal basis $\mathcal{B}(\pi)$ of vectors $v_i$ that are \emph{microlocalised} at $\tau_i$.
The postulate on the microlocalised vectors $v_i$ is that they should be approximate eigenvectors of the action of small elements of $G$, with eigenvalue determined by $\tau_i$.
We explore this idea in Section \ref{sec:archi-microloc}.
In this sketch, lets us write
\begin{displaymath}
    \pi(g) v_i \approx \chi_{\tau_i}(g) v_i,
\end{displaymath}
where $\chi_{\tau_i}$ is a certain unitary character of $G$ attached to $\tau_i$ and $g$ is a small element.
This postulate already gives a solution to desideratum \ref{desideratum-eigenvectors} in Section \ref{sec:intro-rel-trace-formula}.
Indeed, if a test function $f$ is supported on small enough elements, then
\begin{displaymath}
    \pi(f) v_i = \int_G f(g) \cdot \pi(g) v_i \approx \left( \int_G f(g) \cdot \chi_{\tau_i}(g) \right) v_i.
\end{displaymath}

As we shall see in the next section, it is useful to consider a vector $v_\tau$ microlocalised at $\tau$ with a corresponding piece of the orbit
\begin{displaymath}
    U_\tau = \{T \cdot (\cosh v \cos \phi, \cosh v \sin \phi, \sinh v) \mid \phi, v \ll T^{-1/2} \},
\end{displaymath}
which for large $T$ has volume approximately $1$ by the formula above for the symplectic measure.

\subsection{Matrix coefficients of localised vectors}\label{sec:arch_mat_scetch}
Let us therefore continue the sketch with an orthonormal basis $\mathcal{B}(\pi)$ of vectors $v_i$ microlocalised at $\tau_i \in U_i$, one of which is $v_\tau$ corresponding to $\tau = (T, 0, 0)$ and $U_\tau$ as above.
One useful fact about microlocalisation is called \emph{approximate equivariance}.
It says that $\pi(g) v_i$ should be localised at $g.\tau_i$, i.e. the image of $\tau_i$ under the coadjoint action of $g$.
Additionally, microlocalised vectors are essentially \emph{unique}, attached to their symplectic-volume-one pieces in the following way.
If $g.\tau_i \approx \tau_j$, i.e. $g.\tau_i \in U_j$, for some $j \neq i$, then equivariance says that $\pi(g)v_i$ is localised at $\tau_j$.
Uniqueness now states that $v_j$ is the unique vector localised at $\tau_j \in U_j$ and therefore that $\pi(g)v_i$ is essentially contained in the span of $v_j$.

This gives enormous control over matrix coefficients of microlocalised vectors.
Indeed, assume first that $g.\tau_i \in U_i$.
Since $U_i$ is a small neighbourhood of volume 1, one may deduce that $g$ is also a small element and therefore $\pi(g) v_i$ is approximately equal to $\chi_{\tau_i}(g) v_i$, as in the postulate.
On the other hand, if $g.\tau_i \in U_j$ for $j \neq i$, then we saw that $\pi(g) v_i$ is essentially contained in the span of $v_j$ and is thus orthogonal to $v_i$.
This allows us to easily compute inner products $\langle \pi(g) v_i, v_i \rangle$.
In particular, the preceding argument shows that these are negligible unless $g$ is small.

Zooming in further to the problem at hand, recall desideratum \ref{desideratum-matrix-coeffs}, i.e. that we want to control the integral $Q(v)$ of matrix coefficients over $K$.
The postulate defining microlocalisation states more precisely that $\chi_{\tau_j}$ is given by the exponential of the linear form on $\gf$ defined by $\tau_j$, as we explain in more detail in Section~\ref{sec:archi-microloc}.
In this sketch, we informally note that the Lie algebra of $K$ is a certain dual of the $c$-direction in the space of coadjoint orbits.
Thus, if $\tau_j = (a_j, b_j, c_j)$ in $(a,b,c)$-coordinates, then
\begin{displaymath}
    \pi(k(\phi)) v_j \approx e^{i c_j \phi} v_j,
\end{displaymath}
where $k(\phi)$ is the element of $K$ parametrised by $\phi \in [0, \pi)$ (see Section \ref{sec:archi-measures}).
This should only hold for $\phi$ small, but it already indicates that, for $v_j$ microlocalised at $t_j = (a_j, b_j, c_j)$, we have
\begin{displaymath}
    \int_K \langle \pi(k) v_j, v_j \rangle \approx \int_{[0, 2\pi)} e^{i c_j \phi} \, d\phi,
\end{displaymath}
which approximately vanishes if $c_j \neq 0$, as in orthogonality of characters.

To have useful lower bounds on $Q(v_j)$, one of the goals expressed after deriving the inequality \eqref{eq:intro-pretrace-ineq}, we should only consider vectors $v_j$ with $c_j = 0$.
This is the reason we previously singled out the vector $v_\tau$ localised at $\tau = (T, 0, 0)$.
We can compute $Q(v_\tau)$ by noting that $\langle \pi(k(\phi)) v_\tau, v_\tau \rangle \approx 1$ if $k(\phi).\tau \in U_\tau$ and we have vanishing otherwise.
Recall that $K$ acts by rotation around the $c$-axis and polar coordinates show that only $\phi \ll T^{-1/2}$ contributes.
Indeed, this is the $\phi$-width of $U_\tau$, and one of the advantages of localised vectors is precisely this possibility of \emph{computing integrals of matrix coefficients as volumes} of intersections (here $U_\tau \cap K.\tau$).
Therefore, 
\begin{equation} \label{eq:sketch-stability-lower-bd}
    Q(v_\tau) \approx T^{-1/2}.
\end{equation}
We note here the keyword \emph{stability}, which means in this context that the $K$-orbit of $\tau$ (the circle around the hyperboloid at height $c = 0$) is transitive and free under the action of $K$, which is important in performing such volume computations (see Section \ref{sec:rel-char-estimates}).

This motivates us to choose the test function $\omega$ so that $\pi_\lambda(\omega)$ acts as a projector onto $v_\tau$.
Assuming this, the inequality \eqref{eq:intro-pretrace-ineq} and the discussion above now translates into
\begin{equation} \label{eq:intro-ineq-T12}
    |\varphi_{\lambda, \circ}(1)|^2 \ll T^{1/2} \cdot \int_K \int_K \kappa_{\omega \ast \omega^\ast} (x,y).
\end{equation}

\subsection{Test functions} \label{sec:intro-test-functions}

Let us consider the shape of a test function $\omega$ for which $\pi_\lambda(\omega)$ would be the projection onto $v_\tau$.
By definition
\begin{displaymath}
    \pi(\omega) v_i = \int_G \omega(g) \cdot \pi(g) v_i.
\end{displaymath}
At least in a neighbourhood $\mathcal{U}$ of the identity, $\pi$ acts on microlocalised vectors $v_i$ by a character.
In a simpler context, if $G$ were finite, $\mathcal{U}$ a subgroup of $G$, and $\chi_{\tau_i}$ characters of $\mathcal{U}$, the choice 
\begin{equation} \label{eq:intro-archi-test-fn}
    \omega = \frac{1}{\Vol (\mathcal{U})} \vvmathbb{1}_{\mathcal{U}} \cdot \chi_{\tau}^{-1}
\end{equation}
would satisfy our requirements.

We aim to take a smoothened version of the test function $\omega$ above, but there are important conditions for $\chi_\tau$ to behave like a character of an approximate subgroup $\mathcal{U}$.
These constraints are explored in Section \ref{sec:limits} on limits to localisation and are directly related to how microlocalised vectors correspond to pieces of the orbit of symplectic volume one and not less.

\subsubsection{A simple test function}
A first, unrefined but admissible choice would be the ball $\mathcal{U} = I + O(T^{-1/2})$, where $I$ is the identity matrix.
The size $T^{-1/2}$ corresponds to the size of the arguments in the definition of $U_\tau$ above.
With this choice, the formula for $\omega$ above produces an approximate projector onto $v_\tau$ when acting on the representation $\pi_\tau$.

However, recall desideratum \ref{desideratum-test-f-concentration}, namely that for the pretrace inequality in \eqref{eq:intro-pretrace-ineq} to be efficient, we must have very few other representations $\pi$ contribute to the spectral side.
Take, for example, the representation $\pi_{\lambda'}$, with $\lambda' = \frac14 + (T + L)^2$ for some real number $L \ll T$.
We can attach to it an analogous coadjoint orbit, another one-sheeted hyperboloid, and a point $\tau' = (T+L, 0, 0)$ in $(a, b, c)$-coordinates on it.
Let $v_{\tau'}$ be a vector in $\pi_{\lambda'}$ microlocalised at $\tau'$.
We want to study when $\pi_{\lambda'}(\omega)v_{\tau'}$ is non-negligible, and we hope that this happens only for very small $L$.

Focus now on the action of $D$, which we parametrise by $d(t) = \diag(e^t, e^{-t})$.
The Lie algebra of $D$ corresponds to a certain dual in the space of coadjoint orbits of the $a$-direction.
Thus, as we explain more in Section~\ref{sec:archi-microloc}, the postulate on localisation states in this case that
\begin{displaymath}
    \pi(d(t)) v_\tau \approx e^{itT} v_\tau, \quad \pi(d(t)) v_{\tau'} \approx e^{it(T+L)} v_{\tau'},
\end{displaymath}
where the exponents give the pairing of $\tau$ and $\tau'$ with the Lie algebra element corresponding to $d(t)$.
Considering the Iwasawa decomposition $G = NDK$, where $N$ is the subgroup of unipotent upper-triangular matrices, let us informally compute the $D$-part of the integral defining the convolution operator.
Computing at the level of the Lie algebra, the intersection $\mathcal{U} \cap D$ is captured by the bound $t \ll T^{-1/2}$, so $\pi(\omega) v_{\tau'}$ as an integral partly consists of
\begin{displaymath}
    T^{1/2} \int_{t \ll T^{-1/2}} \chi_{\tau}^{-1}(d(t)) \cdot \pi(d(t)) v_{\tau'} \approx 
    T^{1/2} \int_{t \ll T^{-1/2}} e^{itL} \approx
    \frac{2T^{1/2}}{L} \sin\left( \frac{L}{T^{1/2}} \right).
\end{displaymath}
We conclude that representations $\pi_{\lambda'}$ with $L$ as large as $T^{1/2}$ can contribute significantly to the spectral average.

In effect, this choice of $\mathcal{U}$ corresponds through dualising to a ball around $\tau$ of radius $T^{1/2}$.
Informally, if the coadjoint orbit of a representation intersects this ball, then it contributes to the spectral average.
By the observations above, we must thus reduce the size of this ball in the $a$-direction.

\subsubsection{A refined test function} \label{sec:intro-refined-test-fn}
The issue with the choice made above has to do with the subgroup $D$ of diagonal matrices, which is dual to the $a$-direction.
This is no coincidence, since $D$ is the stabiliser $G_\tau$ of the point $\tau$ under the coadjoint action (we will discuss this action in Section~\ref{sec:archi-microlocalisation-big-chap}).
However, as we shall see in Section~\ref{sec:limits}, localisation in the direction of the stabiliser imposes much weaker restrictions.
This is one of the key insights that allowed Nelson \cite{Nel-Un} to refine the microlocal calculus of \cite{NV} and prove subconvexity.

\begin{figure}[ht]
    \centering
    \includegraphics[width=0.8\textwidth]{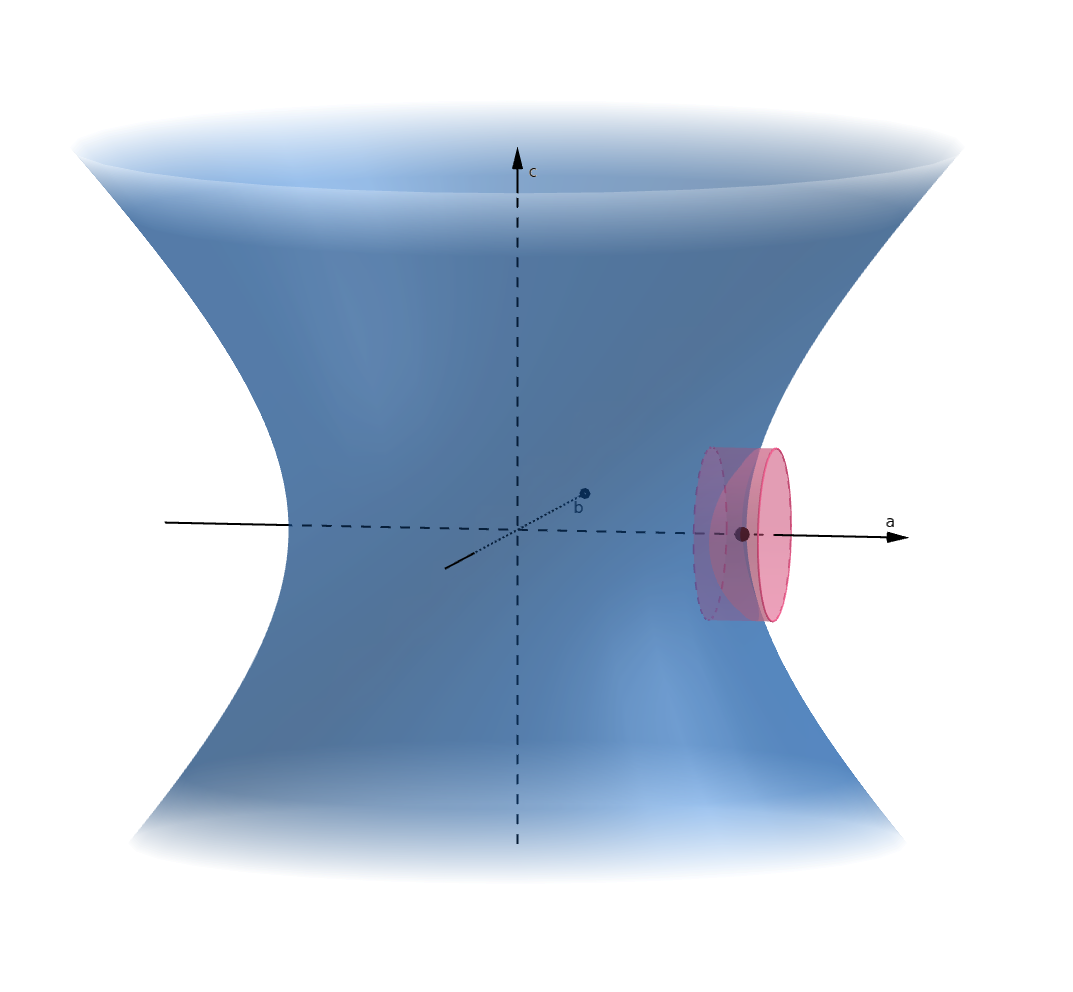}
    \caption{The coin-shaped region (in pink) containing the support of $\psi^\tau_{1/2, 1}$ for $\tau = TA$ and $\hrm = T^{-1}$}
    \label{fig:coin-shaped}
\end{figure}

The result is that we can enlarge $\mathcal{U}$ to be
\begin{displaymath}
    \mathcal{U} = (G_\tau + O(T^{-1/2})) \cap (I + O(T^{-\varepsilon})).
\end{displaymath}
We therefore extend the domain in the $D$-direction, giving now a volume
\begin{displaymath}
    \Vol(\mathcal{U} \cap D) \approx T^{-\varepsilon}.
\end{displaymath}
In the dual space of coadjoint orbits, this choice of $\mathcal{U}$ corresponds to squishing the ball of radius $T^{1/2}$ around $\tau$ in the $a$-direction down to a width of size $T^{\varepsilon}$.
Nelson calls this a \emph{coin-shaped region}, and we depict it in Figure~\ref{fig:coin-shaped}.
The number of hyperboloids intersecting it is much smaller.

Indeed, the same computation as above now shows the $D$-contribution in the convolution operator to be
\begin{displaymath}
    T^{\varepsilon} \int_{t \ll T^{-\varepsilon}} e^{itL} \approx
    \frac{2T^{\varepsilon}}{L} \sin\left( \frac{L}{T^{\varepsilon}} \right) \ll \frac{T^{\varepsilon}}{L}.
\end{displaymath}
This quantity becomes small in the $T \rightarrow \infty$ limit once $L \gg T^{\varepsilon'}$ for some $\varepsilon' > \varepsilon$.
Therefore, the spectral average is now over a very short family with parameters $T + O(T^\varepsilon)$, essentially the same as in \cite{IS95}.

\begin{rem}
    At the non-archimedean places, the implementation of these ideas is simpler.
    We may work with actual subgroups $\mathcal{U}$ and characters, and the operators we obtain this way are true (as opposed to approximate) projectors.
    See Section \ref{sec:p-adic} for more details.
\end{rem}

\subsection{The baseline bound} \label{sec:sketch_baseline}
We turn our attention back to the bound \eqref{eq:intro-ineq-T12} and aim, for now, at proving the baseline bound $\lvert \varphi_{\lambda, \circ}(1) \rvert \ll T^{1/2}$.
This amounts to showing that
\begin{displaymath}
    \int_K \int_K \kappa_{\omega \ast \omega^\ast} (x,y) \ll T^{1/2}.
\end{displaymath}
To do so we insert the definition of $\kappa_{\omega \ast \omega^\ast}(x,y)$ and apply the triangle inequality to obtain
\begin{displaymath}
    \left \lvert \int_K \int_K \kappa_{\omega \ast \omega^\ast} (x,y) \right \rvert \leq \sum_{\gamma\in \Gamma} \int_K\int_K \lvert [\omega \ast \omega^\ast](x^{-1}\gamma y)\rvert \,dx \,dy.
\end{displaymath}
The main contribution is expected to arise from $\gamma\in \Gamma\cap K$. 
For such elements, a simple change of variables yields
\begin{displaymath}
     \int_K\int_K \lvert [\omega \ast \omega^\ast](x^{-1}\gamma y)\rvert \, dx\, dy = \Vol(K) \int_K \lvert [\omega \ast \omega^\ast]( y)\rvert \,dy.
\end{displaymath}

To evaluate this we observe that, if $\mathcal{U}$ were an actual subgroup, then $\omega^\ast = \omega$ and $\omega \ast \omega^\ast = \omega$.
In our case, we have a good approximation to these identities.
Next, to compute $\Vol(\mathcal{U})$, we work at the level of the Lie algebra, where we have three orthogonal basis vectors
\begin{displaymath}
    \mat{1}{0}{0}{-1},\; \mat{0}{1}{1}{0},\; \mat{0}{-1}{1}{0}.
\end{displaymath}
The first of these generates the $D$-direction and therefore contributes by $T^{-\varepsilon}$.
The last two contribute each by a factor of $T^{-1/2}$, showing that $\Vol(\mathcal{U}) \approx T^{-1-\varepsilon}$.
Since $\lvert \omega \rvert  = \frac{1}{\Vol(\mathcal{U})} \vvmathbb{1}_{\mathcal{U}}$, we obtain that
\begin{displaymath}
    \int_K \lvert [\omega \ast \omega^\ast]( y) \rvert \,dy \approx 
    \int_K \lvert \omega( y) \rvert \,dy =
    \frac{\Vol(K\cap \mathcal{U})}{\Vol(\mathcal{U})} \ll T^{1/2+\varepsilon}.
\end{displaymath}

Finally, note that $\abs{\Gamma \cap K} \ll 1$.
For the elements $\gamma \in \Gamma \setminus K$, we again consider the size and support of $\omega$, so that
\begin{displaymath}
    \int_K \int_K \kappa_{\omega \ast \omega^\ast} (x,y) \ll T^{\frac{1}{2} + \varepsilon}+ \frac{1}{\Vol(\mathcal{U})} 
    \sum_{\gamma\in \Gamma\setminus K} \Vol(\{ (x,y)\in K\times K\colon x^{-1}\gamma y \in \mathcal{U} \}).
\end{displaymath}
At this point, notice that $K\mathcal{U}K \subset K + O(T^{-\varepsilon})$, since $K$ is a bounded subgroup.
Since $\Gamma$ is discrete (and \enquote{fixed}), we could choose the implicit constants in the definition of $\mathcal{U}$ to have $\gamma \in K \mathcal{U} K$ only if $\gamma \in K$.
In this way, the second term above disappears and the baseline bound follows.

\subsection{Amplification}\label{sec:sketch_amp}

To improve upon the baseline bound, Iwaniec and Sarnak \cite{IS95} use the amplification technique.
The same idea was successful for proving subconvexity using the orbit method as in \cite{Nel-Un} and \cite{marsh}, and it applies to our set-up in a straight-forward way.

At an adelic level, the amplification method consists in choosing a certain test function at the (unramified) finite places to further increase the contribution of $\pi_\lambda$ in the spectral expansion.
Thus, we would obtain an inequality analogous to \eqref{eq:intro-pretrace-ineq}, with $\omega = \omega_\infty \cdot \omega_f$, where $\omega_\infty$ is given by our previous choice \eqref{eq:intro-archi-test-fn}.
We explain the standard choice of $\omega_f$ in Section \ref{sec:amp}.

Classically, this amounts to applying a linear combination of Hecke operators to the spectral expansion \eqref{eq:intro-spec-expansion}.
For this purpose, define
\begin{displaymath}
    R(n) \coloneqq \{ \gamma \in \mathcal{M}_2(\Z) \mid \det(\gamma) = n \}
\end{displaymath}
and let
\begin{displaymath}
    T_n(\phi)(x) = \sum_{\gamma \in R(1) \backslash R(n)} \phi(\gamma x)
\end{displaymath}
be the associated Hecke operator.
By letting it act on the $y$ coordinate of the automorphic kernel, defined as $\kappa_f(x,y) = \sum_{\gamma \in R(1)} f(x^{-1} \gamma y)$, we may fold the sums over $R(1)$ and $R(1)\backslash R(n)$ together into a sum over $R(n)$, and we obtain
\begin{displaymath}
    \sum_{\gamma \in R(n)} [\omega \ast \omega^\ast](x^{-1} \gamma y) = T_n(\kappa_{\omega \ast \omega^\ast})(x,y) = \int_\pi \lambda_n(\pi) \sum_{v \in \mathcal{B}(\pi)} \pi(f)v(x) \cdot \overline{v(y)},
\end{displaymath}
where $\lambda_n(\pi)$ is the eigenvalue of $T_n$ on the representation $\pi$.
This is a more general version of Equation (1.4) in \cite{IS95}.
To achieve positivity and obtain an inequality as in \eqref{eq:intro-pretrace-ineq}, an \emph{amplified pretrace inequality}, we take a convolution (as in the archimedean test function) of a weighted average of Hecke operators over all unramified primes of size $p \asymp X$, for some parameter $X > 0$ to be chosen later.

Skipping technical details, this method boosts the contribution of $\pi_\lambda$ by $X^2$ at the cost of having on the geometric side an average over the ``denser" sets $R(n)$ for $n$ of size given by some power of $X$.
The overall shape of the problem is
\begin{displaymath}
    \abs{\varphi_{\lambda,\circ}(1)}^2 \ll \frac{T^{1/2}}{X} + \sum_{n} \sum_{\gamma \in R(n)} \frac{1}{ X^2 \sqrt{n}} \iint \omega \ast \omega^\ast (x^{-1} \gamma y).
\end{displaymath}
Taking $X = T^{\delta}$ is the goal, as long as we can bound the average over the $R(n)$ when $X$ is this large.

Thankfully, using number-theoretic and geometric methods, we have a rather good understanding of the number of contributing matrices $\gamma$.
However, simply plugging in a trivial bound
\begin{displaymath}
    \iint \omega \ast \omega^\ast (x^{-1} \gamma y) \ll T^{1/2}
\end{displaymath}
would make our gambit fail.
A better bound involves the distance $d(\gamma)$ from $\gamma$ to the subgroup $K$, which is simply defined as the hyperbolic distance between $\gamma \cdot i$ and $i$ in the upper half-plane.

As observed before, the main term on the geometric side should come from $\gamma \in K$ or, more generally, those with very small $d(\gamma)$, close to zero.
There are only a few of these even in the denser lattice $R(n)$, although proving strong uniform bounds in a general setting is difficult and discussed in Section \ref{sec:amp}.

To win, we must therefore save a power of $T$ when $d(\gamma)$ grows larger.
If $d(\gamma) \neq 0$, we use a bound of the shape
\begin{equation} \label{eq:intro-decay}
    \int_K\int_K \lvert [\omega \ast \omega^\ast](x^{-1}\gamma y)\rvert \,dx \,dy \ll 
    \frac{T^{1/2 - \delta}}{d(\gamma)}.
\end{equation}
Strong counting results and optimising the choice of $X$ would then give the sub-baseline sup-norm bound.

\subsection{Transversality} \label{sec:intro-transv}
To sketch the argument for the bound \eqref{eq:intro-decay}, we again approximate $\abs{\omega \ast \omega^\ast}$ by $\Vol(\mathcal{U})^{-1} \vvmathbb{1}_{\mathcal{U}}$.
Inserting the estimate on the volume of $\mathcal{U}$, we thus work with the double integral
\begin{displaymath}
    T \int_K \int_K \vvmathbb{1}_\mathcal{U} (x^{-1} \gamma y) \, dx \, dy.
\end{displaymath}

Consider first the inner integral in the variable $x$.
Clearly, if $\gamma y \notin K \mathcal{U}$, then the integrand is always zero.
Therefore, let $\gamma y = k \cdot u$ for $k \in K$ and $u \in \mathcal{U}$.
A change of variables allows us to assume that $k = 1$ and, since $\mathcal{U}$ is an approximate subgroup, we can also assume that $u = 1$.
We then have
\begin{displaymath}
    \int_K \vvmathbb{1}_\mathcal{U} (x^{-1} \gamma y) \, dx \approx \Vol(\mathcal{U} \cap K) \approx T^{-1/2}.
\end{displaymath}

\begin{figure}[ht]
    \centering
    \includegraphics[width=\textwidth]{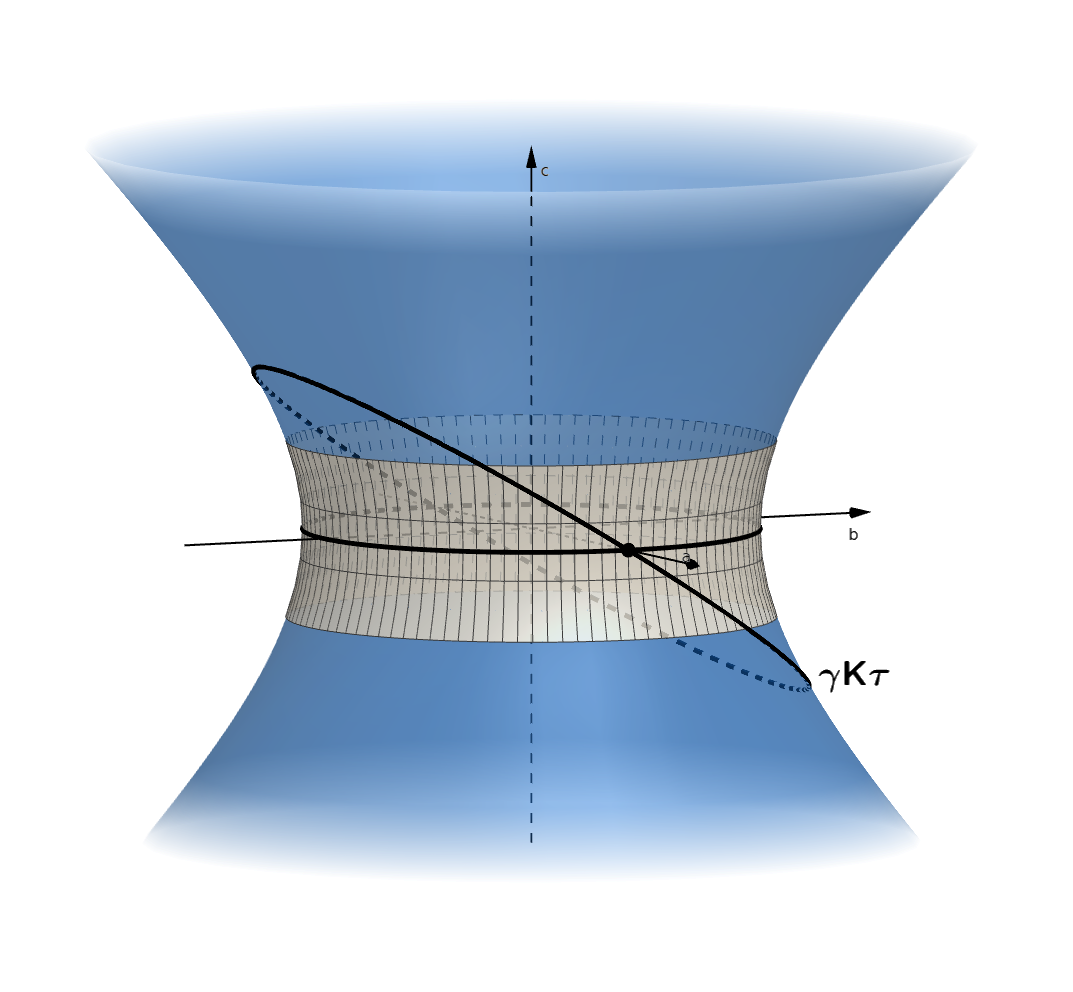}
    \caption{The circular band $K \mathcal{U} \tau$ and the tilted circle $\gamma K \tau$}
    \label{fig:band-ring}
\end{figure}

It remains to show now that
\begin{displaymath}
    \Vol \{ y \in K \mid \gamma y \in K \mathcal{U} \} \ll \frac{T^{-\delta}}{d(\gamma)}.
\end{displaymath}
By invariance of the measure on $K$ and because the distance function is bi-$K$ invariant, we can replace $\gamma$ by $k_1 \gamma k_2$ for any $k_1, k_2 \in K$.
The Cartan decomposition now allows us to assume that $\gamma \in D$.

By stability, we can measure volumes at the level of the coadjoint orbit.
Being careful to normalise measures properly (note that the circle passing through $\tau$ has length $2 \pi T$, which is $T$ times the volume of $K$), we may compute
\begin{displaymath}
    \Vol \{ y \in K \mid \gamma y \in K \mathcal{U} \} = \frac{1}{T} \Vol (K \gamma. \tau \cap K\mathcal{U}.\tau).
\end{displaymath}
Notice that the image of $\tau$ under the action of $\mathcal{U}$ lies in a ball of radius $T^{1/2}$ around $\tau$, by definition.
Since $K$ acts by rotations around the axis of the hyperboloid, the set $K\mathcal{U} \tau$ is a circular band of width $T^{1/2}$.

On the other hand $\gamma K \tau$ is the image under $\gamma$ of the circle passing through $\tau$.
For $\gamma \ll 1$, this is also a tilted circle, and we are left with computing the length of its intersection with the circular band (see Figure~\ref{fig:band-ring}).
Assuming that $\gamma = \diag(e^t, e^{-t})$ with $t \ll 1$ small, the angle between the tilted circle and the circle at the centre of the band is approximately $t$, as we compute in Remark~\ref{rem:tilted}.
Ignoring curvature, we can estimate the length of the intersection by the size of the hypotenuse in a right triangle with a leg of size $T^{1/2}$ and an opposite angle of size $t$.
For small $t$, this length is roughly $T^{1/2}/t$.
Since $d(\gamma) \approx t$, we put everything together and obtain the bound~\eqref{eq:intro-decay} with $\delta = 1/2$.

\section{Archimedean microlocalisation} \label{sec:archi-microlocalisation-big-chap}

This section presents the ideas behind the Nelson-Venkatesh orbit method, providing motivation, intuition, and a certain amount of details.
For the most precise and complete account of the statements below, we refer to the original papers \cite{NV, Nel-Un}.
We roughly follow \cite{Nel-Un} and make most ingredients and steps as explicit as possible for the group $\PGL(2, \R)$.

\subsection{The group and its Lie algebra} \label{sec:archi-gp-lie-algebra}
Let $G = \PGL(2, \R)$ and let $\gf = \mathfrak{sl}_2(\R)$ be the Lie algebra of $G$.
The latter consists of real two-by-two matrices with vanishing trace and is equipped with the adjoint action of $G$.
More precisely, $g \in G$ acts on $X \in \gf$ by conjugation, i.e. $\Ad_g X = gXg^{-1}$.
The Lie algebra also acts on itself by $\ad_X Y = [X,Y] = XY - YX$.

\subsubsection{A choice of basis and measures} \label{sec:archi-measures}
In this paper, we mainly work with the basis for $\gf$ consisting of
\begin{displaymath}
    A = \mat{1}{0}{0}{-1},\; B = \mat{0}{1}{1}{0},\; C = \mat{0}{-1}{1}{0}.
\end{displaymath}
Given this choice, we now relate measures on the Lie algebra and the group, carefully computing the correct normalisations.
To make it easy to compare to common normalisations, we compute measures in the Iwasawa decomposition as well.

We start by recalling the general description of a Haar measure on a Lie group using a basis for the Lie algebra as canonical coordinates through the exponential map.
First, we interpret a given (left) Haar measure $dg$ as a top-degree differential form (see \cite[Chap.\ I.1, Section 2]{helgason-gga} for more details).
Now take a top-degree form $dX$ on the vector space $\gf$, which can also be interpreted as a measure, and assume that $(dg)_I = dX$, where $I$ is the identity element in $G$.

With these assumptions, we have that the pullback of $dg$ under the exponential map satisfies
\begin{displaymath}
    \exp^\ast (dg) = \det \left( \frac{1 - e^{-\ad X}}{\ad X} \right) dX
\end{displaymath}
on neighbourhoods where the exponential map is a diffeomorphism (see \cite[Chap.\ I.1, Thm.\ 1.14]{helgason-gga}).
In the expression above, we interpret $\ad X$ as a matrix operating on $\gf$ and plug it into the corresponding Taylor series, as we exemplify below.
Thus, the function called $j(X)$ or $\operatorname{jac}(X)$ in \cite[Sec.\ 2.1]{NV} or \cite[Thm.\ 7.1]{Nel-Un} has the form
\begin{displaymath}
    j(X) = \det \left( \frac{1 - e^{-\ad X}}{\ad X} \right).
\end{displaymath}

We now make this computation concrete.
Let us take the measure $dX = dA \wedge dB \wedge dC$ on $\gf$ and fix the Haar measure $dg$ on $G$ by setting $(dg)_I = dX$.

Let $X = aA + bB + cC$. It is easy to compute that
\begin{displaymath}
    \ad X = \begin{pmatrix}
        0 & -2c & 2b \\
        2c & 0 & -2a \\
        2b & -2a & 0
    \end{pmatrix}
\end{displaymath}
in our basis.
This is a diagonalizable matrix, with eigenvalues 
\begin{displaymath}
    \pm 2 \sqrt{a^2 + b^2 - c^2}
\end{displaymath} and $0$.
Since we are taking the determinant for computing $j(X)$, we can assume a diagonal form when plugging into the Taylor series and we obtain that
\begin{displaymath}
    j(X) =
    \frac{1 - e^{- 2 \sqrt{a^2 + b^2 - c^2}}}{2 \sqrt{a^2 + b^2 - c^2}} \cdot \frac{1 - e^{2 \sqrt{a^2 + b^2 - c^2}}}{-2 \sqrt{a^2 + b^2 - c^2}}
    = \frac{\sinh^2(\sqrt{a^2 + b^2 - c^2})}{a^2 + b^2 - c^2}.
\end{displaymath}

We now compute the normalisation of the Haar measure defined above in the perhaps more familiar Iwasawa coordinates.
For technical reasons, it is useful to work here with the connected group $\PSL_2(\R)$, which we view as an index-two subgroup of $G$.
We now use the standard Iwasawa decomposition
\begin{displaymath}
    \PSL_2(\R) = N D^0 K,
\end{displaymath}
where $N$ is the subgroup of unipotent upper-triangular matrices, $D^0$ is the subgroup of diagonal matrices, and $K$ is the projective special orthogonal group.
When following the general construction, say in \cite[Chap.\ I.5, Sec.\ 1]{helgason-gga}, this decomposition corresponds to taking the character $A \mapsto 2$ as positive root.

Using the exponential map, these subgroups have canonical coordinates
\begin{align*}
    N &= \left\{ \exp(x (B-C)/2) = \begin{pmatrix}
        1 & x\\
         & 1
    \end{pmatrix} \mid x \in \R \right\};\\
    D^0 &=  \left\{ \exp(t A) = \begin{pmatrix}
        e^t & \\
         & e^{-t}
    \end{pmatrix} \mid t \in \R \right\};\\
    K &=  \left\{ \exp(\phi C) = \begin{pmatrix}
        \cos \phi & -\sin \phi \\
        \sin \phi & \cos \phi
    \end{pmatrix} \mid \phi \in [0, \pi) \right\}.
\end{align*}
On each of them, we take the invariant measures induced by the differential forms $d(B-C)/2$, $dA$, and $dC$, respectively.
Notice that there is no Jacobian factor between the measures on the subgroups and those on their Lie algebras, meaning that the corresponding function $j$ as above is the constant $1$, since these are one-parameter groups. 
In any case, this gives a measure on $N \times D^0 \times K$.

As in \cite[Chap.\ I, Prop.\ 5.1]{helgason-gga}\footnote{Our Iwasawa decomposition is in reversed order when compared to \cite{helgason-gga}. Thus, one needs to apply an inversion in the cited formula, which does not introduce any additional factors since $\PSL_2(\R)$ and the subgroups above are unimodular.}, we can now define a Haar measure $\mu$ on $\PSL_2(\R)$ such that
\begin{multline*}
    \int_G f(g) d\mu = \int_{x \in \R} \int_{t \in \R} \int_{\phi \in [0, \pi)} f(e^{x (B-C)/2} \cdot e^{t A} \cdot e^{\phi C}) e^{-2t} \, dx\, dt\, d\phi.
\end{multline*}
It is more common to parametrise $D^0$ using matrices
\begin{displaymath}
    \begin{pmatrix}
        \sqrt{y} & \\
        & 1/\sqrt{y}
    \end{pmatrix},
\end{displaymath}
thus making a change of variables $e^{2t} = y$.
We then have the integration formula
\begin{multline*}
        \int_G f(g) d\mu \\
        = \frac12 \int_{x \in \R} \int_{y \in \R_{>0}} \int_{\phi \in [0, \pi)}
        f\left(
        \begin{pmatrix}
        1 & x\\
         & 1
        \end{pmatrix}
        \begin{pmatrix}
        \sqrt{y} & \\
        & 1/\sqrt{y}
        \end{pmatrix}
        \begin{pmatrix}
        \cos \phi & -\sin \phi \\
        \sin \phi & \cos \phi
        \end{pmatrix}
        \right) \\
        \frac{dx\, dy\, d\phi}{y^2}.
\end{multline*}
As a differential form, we have
\begin{displaymath}
    (d\mu)_I = (d(B-C)/2) \wedge dA \wedge dC = -\frac{1}{2} dA \wedge dB \wedge dC.
\end{displaymath}

We now extend the measure $\mu$ to a Haar measure on $\PGL_2(\R)$ by invariance ($\PSL_2(\R)$ is a normal subgroup of index $2$). 
By the previous computations, up to orientation, the Haar measure $dg$ defined above is now given by $\frac{dx\, dy\, d\phi}{y^2}$ in Iwasawa coordinates. 
At least on $\PSL_2(\R)$, this is the same as the push-forward under the projection map of the measure on $\GL_2^+(\R)$ defined in \cite[Chap.\ 2, (1.20)]{bump} restricted to $\SL_2(\R)$, divided by $2$ (since $\SO(2)$ gets folded in half).
For future reference, we note that the volume of $K$ is $\pi$.

\subsubsection{The coadjoint action}
Let $\gf^\ast$ be the real dual space of $\gf$. Using the trace pairing
\begin{displaymath}
    (X, Y) \mapsto \tr (X \cdot Y),
\end{displaymath}
we identify $\gf^\ast$ with $\gf$ and we have the dual basis
\begin{displaymath}
    A^\ast = \frac12 A,\; B^\ast = \frac12 B,\; C^\ast = -\frac12 C.
\end{displaymath}
We generally denote the value of a functional $\tau \in \gf^\ast$ on a vector $X \in \gf$ by $\langle \tau, X \rangle$. This pairing induces a representation dual to $\Ad$, called the \emph{coadjoint representation} $\Ad^\ast$, defined by
\begin{displaymath}
    \langle \Ad^\ast_g \tau, X \rangle = \langle \tau, \Ad_{g^{-1}} X \rangle.
\end{displaymath}
In its infinitesimal version, we have
\begin{displaymath}
    \langle \ad^\ast_X \tau, Y \rangle = - \langle \tau, [X, Y] \rangle.
\end{displaymath}
It is straightforward to compute using the trace pairing that, in the identification $\gf^\ast \cong \gf$ above, the coadjoint action reduces to
\begin{displaymath}
    \Ad_g^\ast \tau = g \tau g^{-1}
\end{displaymath}
and the infinitesimal action to
\begin{displaymath}
    \ad^\ast_X \tau = [X, \tau].
\end{displaymath}

Of course, this coadjoint representation extends to the dual $\gf_\C^\ast$ of the complexification $\gf_\C = \gf \otimes \R$.
The same applies to the imaginary dual $i \gf^\ast = \Hom(\gf, i\R)$, which we view as real form of $\gf_\C^\ast$ and also identify with the Pontryagin dual $\gf^\wedge$ by applying the exponential map.

Elements of $\gf_\C^\ast$ whose orbits under the coadjoint action have maximal dimension are called regular.
Equivalently, they are the elements whose centraliser, i.e. the stabiliser under the coadjoint action, has minimal dimension.
The subset of $\gf_\C^\ast$ of regular elements is denoted by $(\gf_\C^\ast)_\textrm{reg}$ and similarly for other sets.
Many statements (e.g. in \cite{NV})  are restricted to regular elements and this notion becomes important in higher rank.
In our case, it can be easily checked that regular simply means non-zero, i.e. the only element that is not regular is $0 \in \gf_C^\ast$.

\subsubsection{Coadjoint orbits} \label{sec:orbits}
Let us compute some examples of orbits under the coadjoint action explicitly.
For the orbit method, we are mostly interested in orbits of elements in $i \gf^\ast$, but for this section we simply ignore the imaginary unit and work in $\gf^\ast$ directly.

The prototype we focus on in this paper is the orbit of $2 A^\ast$, which is identified with the matrix $A$ as above.
More generally, let us denote the orbit of $X$ by $\Oc_X$.

We write an orbit element $\tau = a A + b B + c C \in \Oc_A$ in our chosen basis and notice that, since $G$ acts by conjugation, the determinant of $\tau$ is an invariant.
More precisely,
\begin{displaymath}
    -\det(\tau) = a^2 + b^2 - c^2 = -\det(A) = 1
\end{displaymath}
for any $\tau \in \Oc_A$. 
The orbit is therefore contained in a one-sheeted hyperboloid.

\begin{figure}[ht]
    \centering
    \includegraphics[width=0.8\textwidth]{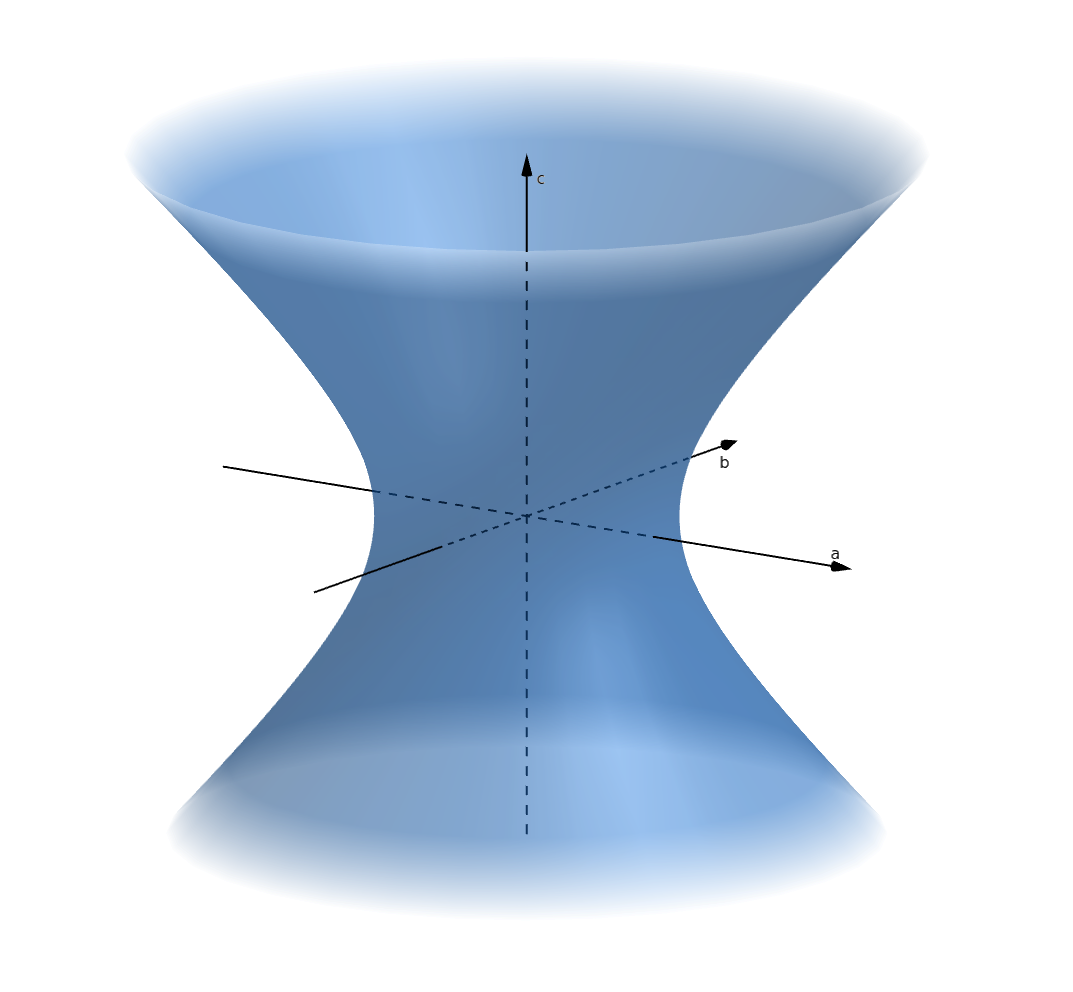}
    \caption{The orbit $\Oc_A$}
    \label{fig:orbit-OA}
\end{figure}

In fact, the hyperboloid describes the full orbit $\Oc_A$.
To see this, we parametrise the hyperboloid by the coordinates
\begin{equation} \label{eq:hyperboloid-coords}
    (a, b, c) = (\cosh v \cos \phi, \cosh v \sin \phi, \sinh v) =: \tau_{v, \phi}
\end{equation}
for $v \in \R$ and $\phi \in [0, 2 \pi)$.
Notice that we can use the same coordinates for the scaled orbits $\Oc_{T A}$ for $T > 0$, for which elements can be written as $T \tau_{v, \phi}$.

Using trigonometric identities, we compute that
\begin{displaymath}
    g_{v, \phi} = \begin{pmatrix}
        \cosh v \cos \phi - \sinh v \sin \phi && \sinh v \cos \phi - \cosh v \sin \phi \\
        \sinh v \cos \phi + \cosh v \sin \phi && \cosh v \cos \phi + \sinh v \sin \phi
    \end{pmatrix}
\end{displaymath}
operates by conjugating $A = \tau_{0, 0}$ to $\tau_{2v, 2\phi}$.
We remark that $g_{v, 0} = \exp(v \cdot B)$, and  $g_{0, \phi} = \exp(\phi \cdot C)$.
Notice also that the diagonal torus $D$ inside $G$ is the stabiliser of $A$.
Thus, putting everything together, we can also deduce a decomposition of $G$ as
\begin{displaymath}
    G = \{ g_{v, \phi} \mid v \in \R, \phi \in [0, 2 \pi) \} \cdot D,
\end{displaymath}
and we can write the orbit as the homogeneous space $\Oc_A = G / D$.

Another representative coadjoint orbit is that of $C = 2C^\ast$. 
The determinant gives us the equation
\begin{displaymath}
    a^2 + b^2 - c^2 = -1
\end{displaymath}
and, as above, one can prove that $\Oc_C$ is precisely a double-sheeted hyperboloid. 

\begin{figure}[ht]
    \centering
    \includegraphics[width=0.8\textwidth]{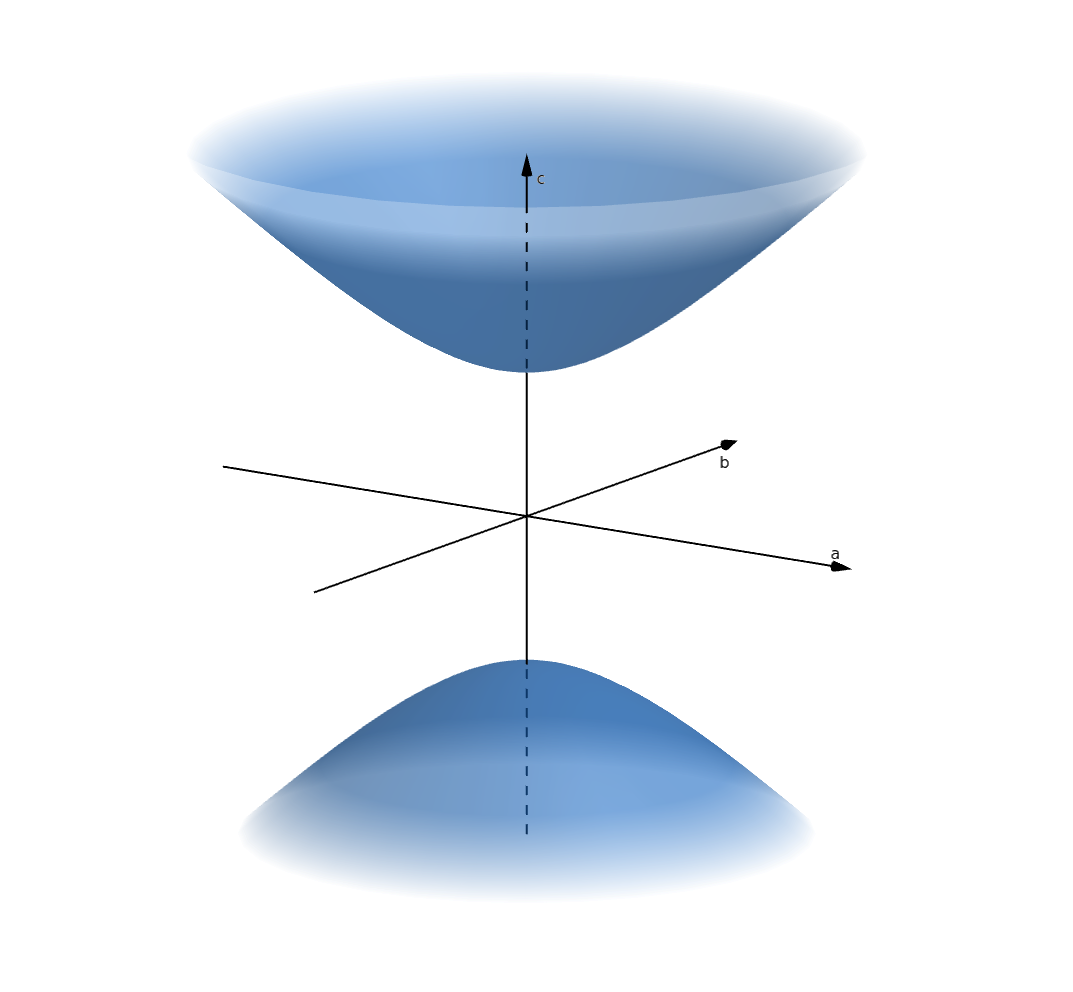}
    \caption{The orbit $\Oc_C$}
    \label{fig:orbit-OC}
\end{figure}

By the canonical rational form, we can classify all coadjoint orbits. 
We have the family of one-sheeted hyperboloids $\Oc_{T A}$ for $T > 0$, the double-sheeted hyperboloids $\Oc_{T C}$ for $T > 0$, the conic orbit $\Oc_{B+C}$ and the singular orbit $\Oc_{0}$. 
We call the latter two degenerate, since the closure of $\Oc_{B+C}$ is the union $\Oc_{B+C} \cup \Oc_0$.
This observation is mirrored in the fact that these orbits are glued together in the geometric invariant theory quotient below.

\begin{figure}[ht]
    \centering
    \includegraphics[width=0.8\textwidth]{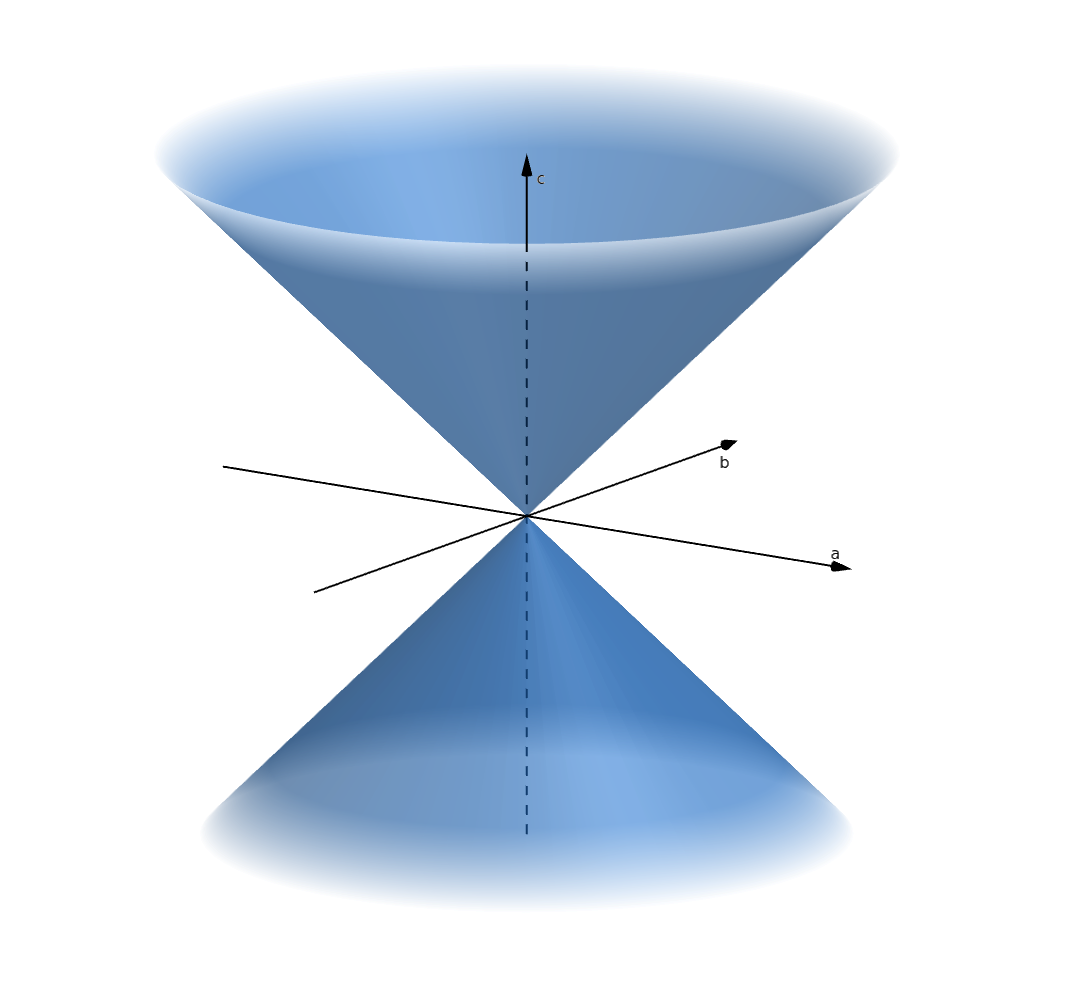}
    \caption{The union of orbits $\Oc_{B+C} \cup \Oc_0$, called the nilcone}
    \label{fig:orbit-nilcone}
\end{figure}

An essential feature of coadjoint orbits is their symplectic structure.
To define it, let $\tau \in \gf^\ast$ lie in some orbit $\Oc$ and observe that the tangent space $T_\tau \Oc$ is spanned by $\ad_X^\ast \tau$ for $X \in \gf$.
Explicitly, for instance, the curve $\Ad_{\exp tB} \tau$ corresponds to $\ad_B^\ast \tau$.
The symplectic form at $\tau$ is now defined as
\begin{displaymath}
    \sigma_\tau(\ad_X^\ast, \ad_Y^\ast) = \langle \tau, [X, Y] \rangle.
\end{displaymath}
Its existence implies that coadjoint orbits are even-dimensional, yet this is clear for our group $G = \PGL_2(\R)$ by the classification of orbits above.
Since our non-degenerate orbits are $2$-dimensional, the symplectic form already defines a volume form,  using the normalisation
\begin{displaymath}
    \omega_{\Oc} = \frac{\sigma}{4 \pi}.
\end{displaymath}
We refer to \cite[Sec. 6.1]{NV} and \cite{Kir} for more details in general.

Let us compute the symplectic form on $\Oc_{T A}$, our orbit of choice, in the coordinates given in \eqref{eq:hyperboloid-coords}.
Let $\tau \in \R$ and $\phi \in [0, 2 \pi)$.
By direct computation using the conjugation action on $T \tau_{v,\phi}$ and the curves given by $\exp(t B)$ and $\exp(t C)$, we find that
\begin{align*}
    \ad^\ast_B \tau_{v, \phi} &= 
    2 \left( \cos \phi \cdot \frac{\partial}{\partial v} - \sin \phi \tanh v \cdot \frac{\partial}{\partial \phi} \right) \\
    \ad^\ast_C \tau_{v, \phi} &=
    2 \cdot \frac{\partial}{\partial \phi}.
\end{align*}
Now, since there is only one top degree form up to scalars, there is $f(T \tau_{v, \phi})$ such that
\begin{displaymath}
    \sigma_{T \tau_{v, \phi}} (\ad^\ast_B, \ad^\ast_C) = f(T \tau_{v, \phi}) \cdot dv \wedge d\phi.
\end{displaymath}
Using the definition, we now easily compute that $f(T \tau_{v, \phi}) = T \cosh v$ and thus
\begin{equation} \label{eq:symplectic-meas}
    \omega_{\Oc_{T A}} = \frac{T \cosh v}{4 \pi} dv \wedge d\phi.
\end{equation}

\begin{rem} \label{rem:tilted}
    At this point, we can also observe the action of $\exp(tA)$ on the circle $T \tau_{0, u}$, as used in Section \ref{sec:intro-transv}.
    Taking $T = 1$ for simplicity, we have 
    \begin{displaymath}
        \exp(tA) \cdot \tau_{0, u} = \cos(u) A  + \cosh(2t) \sin(u) B - \sinh(2t) \sin(u) C.
    \end{displaymath}
\end{rem}

For completeness, we also describe the symplectic form on the other orbits.
We parametrise $\Oc_{TC}$ by the coordinates
\begin{displaymath}
    (a,b,c) = T \cdot (\sinh v \cos \phi, \sinh v \sin \phi, \pm \cosh v)
\end{displaymath}
where $v \in [0, \infty)$ and $\phi \in [0, 2 \pi)$.
Then the volume form is given by
\begin{displaymath}
    \omega_{\Oc_{TC}} = \frac{T \sinh v}{4 \pi} dv \wedge d\phi .
\end{displaymath}

We finally parametrise the union $\Nc := \Oc_{B-C} \cup {0}$, also called the nilcone, by
\begin{displaymath}
    (a,b,c) = (v \cos \phi, v \sin \phi, v)
\end{displaymath}
with $v \in \R$ and $\phi \in [0, 2\pi)$.
The volume form is
\begin{displaymath}
    \omega_{\Nc} = \frac{1}{2 \pi} dv \wedge d\phi. 
\end{displaymath}

\subsubsection{The GIT quotient}
The papers \cite{NV} and \cite{Nel-Un}, among others, use the language of geometric invariant theory (GIT) to put structure on the space of coadjoint orbits and make it precise.
It also allows identifying the regular orbits and the degenerate ones.
While this is not quite necessary for the rest of the present paper, we spend a few paragraphs to work out some details of this theory.
This could aid those who wish to fully understand \cite{NV} (the general picture is given in \cite[Sec.\ 9.1-2]{NV}), but we suggest to the reader to skip this section on a first reading.

In our case, we define the ring $R_\C = \C[A, B, C]$ and identify $\gf_\C^\ast$ with the complex points of the spectrum of $R_\C$, that is $\Hom_\C(R_\C, \C)$, in the usual way.
We are particularly interested in the subspace $\gf^\wedge = i \gf^\ast$.
It can be interpreted as the real points of the spectrum of
\begin{displaymath}
    R = \{ f \in R_\C \mid f(A,B,C) = \overline{f(-A, -B, -C)} \}.
\end{displaymath}
The decomposition $R_\C = R + iR$ allows us to recover the inclusion $i\gf^\ast \subset \gf_\C^\ast$ in this algebraic interpretation.
The usual action of complex conjugation on $\gf_\C^\ast$ can be also interpreted directly on $\Hom_\C(R_\C, \C)$, as well as the action of scalars in $\C^\times$.
This also recovers the subset $i \gf^\ast$ as the locus of points $\lambda$ for which $\overline{\lambda} = -\lambda$.

The algebraic group $\mathbf{G} = \mathbf{PGL(2)}$ over $\C$ acts on the spectrum of $R_\C$ by the coadjoint representation.
We denote by $[\gf_\C^\ast]$ the complex points of the GIT quotient of $\gf_\C^\ast$ by this action.
The latter is defined as the spectrum of $R_\C^\mathbf{G}$, the ring of $\mathbf{G}$-invariant polynomials.
Over the algebraically closed field $\C$ this is easy to compute.
Practically, we are looking at polynomials on the space of $2 \times 2$ traceless complex matrices that are invariant under conjugation by complex matrices.
The Jordan normal form implies that this algebra is generated by the coefficients of the characteristic polynomial.
The vanishing of the trace now implies that $R_\C^\mathbf{G}$ is generated by the determinant.

In the same way as above, the ring $R^\mathbf{G}:= R \cap R_\C^\mathbf{G}$ gives a subset $[i \gf^\ast]$ of $[\gf_\C^\ast]$.
In the interpretation of points as homomorphisms, complex conjugation and the scaling action descend to $[\gf_\C^\ast]$ and we can again see that
\begin{displaymath}
    [i \gf^\ast] = \{ \lambda \in [\gf_\C^\ast] \mid \overline{\lambda} = - \lambda \}.
\end{displaymath}
We also have natural maps $\gf_\C^\ast \longrightarrow [\gf_\C^\ast]$ and $i \gf^\ast \longrightarrow [i \gf^\ast]$, denoted by $\tau \mapsto [\tau]$.

Let us consider a few examples to understand this construction.
Firstly, notice that all elements in $i \gf^\ast$ with determinant $0$ are now made equivalent.
In particular, the distinct degenerate orbits of $0$ and $i(B+C)$ both project to $[0]$.
This is, actually, a desired effect of the GIT quotient.

Secondly, elements in the real part $\gf$ project to $[i \gf^\ast]$, since real matrices and $i$ times real matrices all have real determinants.
From the point of view of representation theory, this is natural.
We return to this idea when discussing infinitesimal characters.

Another way to parametrise $[\gf_\C^\ast]$ is using Cartan subalgebras. 
The standard choice here is the subalgebra $\af_\C = \C \cdot A$ of diagonal matrices.
Its corresponding Weyl group $W_\C$ has two elements and the non-trivial one acts by $A \mapsto -A$.

We choose $2A^\ast$, the homomorphism sending $a A \mapsto 2 a$, to be the positive root in $\af_\C^\ast$.
This induces a decomposition $\gf_\C = \mathfrak{n}^{-} \oplus \af_\C \oplus \mathfrak{n}^{+}$, where
\begin{displaymath}
    \mathfrak{n}^{-} = \C \cdot \begin{pmatrix}
        0 & 0 \\ 1 & 0
    \end{pmatrix}, \,
    \mathfrak{n}^{+} = \C \cdot \begin{pmatrix}
        0 & 1 \\ 0 & 0
    \end{pmatrix}
\end{displaymath}
are the corresponding negative and positive root spaces, respectively.
For the following discussion, we also denote the matrices above, spanning the root spaces, by $f$ and $e$, respectively (this is fairly standard notation).
We can now directly identify $\af_\C^\ast$ with a subspace of $\gf_\C^\ast$, extending linear forms by zero on the direct complement.
It is a fact (see Chevalley's theorem \cite[(9.2)]{NV}) that composing this inclusion with the projection to the GIT quotient gives an isomorphism
\begin{displaymath}
    \af_\C^\ast / W_\C \cong [\gf_\C^\ast].
\end{displaymath}

\subsubsection{The infinitesimal character}
To any irreducible representation $\pi$ of $G$ we can attach an infinitesimal character $\lambda_\pi \in [\gf_\C^\ast]$ that describes the action of the centre of the universal enveloping algebra of $\gf$ on $\pi$. 
We make this explicit in our case, following \cite[Sec. 9.4-5]{NV}.
This infinitesimal character then picks out the coadjoint orbit we attach to the representation $\pi$.

The centre of the universal enveloping algebra of our $\gf$ is generated by the Casimir element $\Omega_C$.
It is defined here as (see the definition of $\Delta$ in \cite[Chap.\ 2, (2.20)]{bump})
\begin{equation} \label{eq:def-casimir}
    -4 \cdot \Omega_C = A^2 + B^2 - C^2.
\end{equation}
Standard computations show that
\begin{displaymath}
    A^2 + B^2 - C^2 = A^2 + 2A + 4fe \in (A^2 + 2A) \oplus \mathfrak{n}^- \mathfrak{n}^+.
\end{displaymath}
Thus, in the decomposition of \cite[(9.4)]{NV}, the projection of $-4 \Omega_C$ to the universal enveloping algebra of $\af$ (by commutativity, this is simply the symmetric algebra) is $A^2 + 2A$.

Now let $\rho = (2A^\ast)/2 = A^\ast$ be the half-sum of the positive roots.
The Harish-Chandra homomorphism is given by composing the projection above with the algebra isomorphism extending $X \mapsto X - \rho(X) \cdot 1$ on $X \in \af$.
Therefore, the operator $-4 \Omega_C$ corresponds to
\begin{displaymath}
    (A - \rho(A))^2 + 2(A - \rho(A)) = A^2 - 1.
\end{displaymath}
Such a polynomial now provides a Weyl-invariant function on $\af_\C^\ast$ by evaluation.
In the case above we obtain the map
\begin{displaymath}
    c A^\ast \mapsto \langle c A^\ast, A \rangle^2 - 1 = c^2 - 1.
\end{displaymath}
The infinitesimal character $\lambda_\pi$ is given by the element $c A^\ast \in \af_\C^\ast / W_\C$ such that $-4 \Omega_C$ acts on $\pi$ by the scalar $c^2 - 1$. 
By Chevalley's isomorphism mentioned above, we view $\lambda_\pi$ as an element of $[\gf_\C^\ast]$.

As in \cite[Sec.\ 9.5]{NV}, one can show that the infinitesimal character of a unitary representation must lie in $[i \gf^\ast]$.
We can check this in our case by hand when discussing the classification of representations for $G$.

\subsubsection{The Kirillov formula}

Let $\pi$ be an irreducible unitary representation of $G$.
For any smooth and compactly supported function $f$ on $G$ we now obtain operators
\begin{displaymath}
    \pi(f) = \int_G \pi(g)f(g)\, dg.
\end{displaymath}
These are trace class, and Harish-Chandra showed that there exists a locally $L^1$-function $\chi_\pi: G \to \C$ for which the trace of $\pi(f)$ is, by a slight abuse of notation,
\begin{displaymath}
    \chi_\pi(f) = \int_G \chi_\pi(g) f(g)\, dg.
\end{displaymath}
We call $\chi_\pi$, as a distribution on $G$, the (generalised or distributional) character of $\pi$.
See \cite[Sec.\ 6]{NV} for more details.

Kirillov's formula essentially expresses the character $\chi_\pi$ as the Fourier transform of the symplectic measure on a coadjoint orbit attached to $\pi$.
For this we assume that $\pi$ is tempered.\footnote{We remark that one should not expect the orbit method to be very useful for non-tempered representations. See \cite[Sec.\ 8.4]{Kir}.}
In our case, this means that $\pi$ is not in the complementary series, as we explain in the next section.

Let $\lambda_\pi \in [\gf^\wedge]$ be the infinitesimal character of $\pi$. 
As in \cite[Sec.\ 6]{NV} and \cite[Thm.\ 7.1]{Nel-Un}, there exists a unique finite union of coadjoint orbits\footnote{The references often call this union a \emph{multiorbit}. 
For the group and almost all representations that interest us, this is a single orbit.} $\Oc_\pi$, made up of regular elements of $\gf^\wedge$ and contained in the preimage of $\lambda_\pi$, such that
\begin{equation} \label{eq:kirillov}
    \chi_\pi(\exp X) \sqrt{j(X)} = \int_{\xi \in \Oc_\pi} e^{\langle \xi, X \rangle} \, \omega_{\Oc_\pi}
\end{equation}
for $X \subset \gf$ is a small enough neighbourhood of $0$.

One way to understand this formula intuitively is by thinking in terms of microlocalised vectors.
We explore this idea below.

\subsubsection{Classification of representations}
We now briefly recall the classification of the relevant class of representations of $G$. 
For this we specialise \cite[Thm.\ 2.6.7]{bump} to our case, that is, assuming a trivial central character.

\begin{prop}[Classification of representations] \label{prop:classif-rep}
    Every infinite dimensional unitary admissible irreducible representation of $G$ is isomorphic to one of the following:
    \begin{itemize}
        \item the principal series representation $\Pc(s)$, where $s \in i \R$;
        \item the complementary series representation $\Pc(s)$, where \begin{displaymath}
            s \in \left( -\frac12, \frac12 \right) \setminus \{0\};
        \end{displaymath}
        \item the discrete series representation $\mathcal{D}(k)$, for $k$ an even positive integer.
    \end{itemize}
    In this list, the only equivalence is $\Pc(s) \cong \Pc(-s)$. 
    The above representations are tempered except for the complementary series.
    The Casimir operator $\Omega_C$ acts on principal and complementary series representations $\Pc(s)$ by
    \begin{displaymath}
        \lambda = \frac14 - s^2
    \end{displaymath}
    and on discrete series representations $\Dc(k)$ by
    \begin{displaymath}
        \lambda = \frac14 - \left(\frac{k-1}{2}\right)^2.
    \end{displaymath}
\end{prop}

For concreteness, we remark that the principal series $\Pc(s)$ are representations induced by the character
\begin{displaymath}
    \rho \left( \mat{y_1}{x}{0}{y_2}\right) = \left| \frac{y_1}{y_2} \right|^{s}
\end{displaymath}
on the Borel subgroup of upper triangular matrices.
For more details, we refer to \cite[Sec.\ 2.4-2.6]{bump}.%

The Casimir eigenvalue shows that the infinitesimal character of $\Pc(s)$ corresponds to $2 s A^\ast$ and that of $\Dc(k)$ to $(k-1) A^\ast$.
As already mentioned, these lie in the real form $[\gf^\wedge]$, since $2 s A^\ast$ and $(k-1) A^\ast$ are contained in $\gf^\ast \cup i \gf^\ast$ for all relevant parameters.

For the Kirillov formula, we directly see that the orbit attached to $\Pc(iT)$ for $T > 0$ is precisely the $G$-orbit of $2 i T A^\ast$, made up of regular elements of $\gf^\wedge$. 
For the discrete series $\Dc(k)$, note that the determinant of $(k-1) A^\ast$, identified with $\tfrac{1}{2}(k-1) A$, is $-(k-1)^2 / 4$.
The element $\tfrac12 (k-1) i C \in \gf^\wedge$ has the same determinant and so projects to the infinitesimal character of $\Dc(k)$, by the construction of the GIT quotient.
Thus, the orbit we attach is that of $\tfrac12 (k-1) i C$.

\subsection{Microlocalised vectors} \label{sec:archi-microloc}

The following discussion is heavily inspired by lectures and notes of Paul Nelson (\cite{Nel-loc}, \cite{Nel-budapest}) and of Akshay Venkatesh (\cite{Ven}).

Let $\pi$ be a unitary representation of $G$. The Lie algebra $\gf$ acts on smooth vectors of $\pi$, so let us naively consider how an eigenvector $v \in \pi$ for this action might look like. 
Clearly, such a vector would give rise to a complex linear functional $\lambda$ on $\gf$ such that
\begin{displaymath}
    Xv = \lambda(X) \cdot v
\end{displaymath}
for $X \in \gf$. Since $\pi$ is unitary, we have that
\begin{displaymath}
    \langle Xv, v \rangle_\pi = \langle v, -X v \rangle_\pi
\end{displaymath}
and thus $\lambda$ must have values in $i \R$. Therefore,
\begin{displaymath}
    \lambda(X) = \langle \tau, X \rangle
\end{displaymath}
for $\tau$ in the imaginary dual $i \gf^\ast = \gf^\wedge$.

Such eigenvectors would heuristically explain the Kirillov formula \eqref{eq:kirillov}.
Indeed, say we had a basis of $\pi$ consisting of eigenvectors $v_j$ associated to $\lambda_j \in \gf^\wedge$.
Formally, expanding the trace of $\exp(X)$ in terms of this basis and ignoring the $j(X)$ factor, we would have that
\begin{equation} \label{eq:kirillov-microlocal-explanation}
    \chi_\pi(\exp(X)) \approx \sum_j \exp \langle i \lambda_j, X \rangle \approx \int_{\xi \in \Oc_\pi} e^{\langle i \xi, X \rangle}.
\end{equation}
Here, the vectors $v_j$ should correspond to a partition of $\Oc_\pi$ into pieces of unit symplectic area for the last approximation sign to make sense.
Though only a heuristic, this is an essential idea in the present application of the orbit method.

Unfortunately, such eigenvectors as above do not exist in general (see Section \ref{sec:limits} below).
However, we could ask for a reasonable approximation of one. 
This leads to the concept of microlocalised vectors. 
To maximise usefulness, we consider these approximations at several orders of differentiation, that is, we consider the action of the universal enveloping algebra.  %

\begin{defn}\label{def:microlocal}
    Let $\pi$ be an irreducible unitary representation of $G$ and let $\tau \in \gf^\wedge$. 
    For $N \in \Z_{>0}$ and $R > 0$, we say that $v \in \pi$ is \emph{$R$-localised at $\tau$} to order $N$ if
    \begin{displaymath}
        \norm{ (X_1 - \langle \tau, X_1 \rangle) \cdots (X_k - \langle \tau, X_k \rangle) v } \leq R^k \norm{v}
    \end{displaymath}
    for any collection $X_1, \ldots, X_k$ of elements in the basis of $\gf$ and any $k \leq N$.
\end{defn}

\begin{rem}
    Localisation is almost $G$-equivariant. 
    More precisely, let $v$ be $R$-localised at $\tau$ and consider the shifted vector $\pi(g)v$ for $g \in G$. One computes that
    \begin{displaymath}
        (X - \langle \tau, X \rangle) \pi(g) v = \pi(g) (\Ad_{g^{-1}} X - \langle \Ad_{g}^\ast (\tau), \Ad_{g^{-1}} X \rangle) v 
    \end{displaymath}
    The norm of this vector is bounded by assumption by $C_g R \norm{v}$, where $C_g > 0$ is a constant depending on $\norm{\Ad_{g}}$. 
    Therefore, $\pi(g)v$ is localised at $\Ad_{g}^\ast \tau = g.\tau$, possibly to a lesser extent.
\end{rem}

\subsubsection{Limits to localisation} \label{sec:limits}

There are limits to how much a vector can be localised, meaning how small the parameter $R$ in Definition \ref{def:microlocal} can be. 
Broadly speaking, these constraints are a result of the non-commutativity of differential operators, given in this paper by Lie algebra vectors.
Indeed, they are a manifestation of the Heisenberg uncertainty principle. 
We give here two perspectives on these constraints.

First, we consider how small group elements act on localised vectors. 
Let $v$ be $R$-localised at $\tau$, to some order. 
In the next few paragraphs we perform formal computations, without making approximations precise, and we assume that $v$ is an analytic vector (see e.g. \cite[Sec.\ VI.1]{lang}).%

We rewrite the localisation condition simply as
\begin{displaymath}
    \frac1R Xv \approx \langle \tau, \frac1R X \rangle v,
\end{displaymath}
for any $X$ in the basis of $\gf$.
We thus have
\begin{displaymath}
    Xv \approx \langle \tau, X \rangle v,
\end{displaymath}
for any $X = O(R^{-1})$, the last condition meaning that $X$ lies in a ball $\mathcal{D}$ of radius $O(R^{-1})$ in a fixed norm, e.g. the Frobenius norm $X \mapsto \tr(X^t X)$.

Exponentiating, we have a subset $\exp(\mathcal{D})$ of $G$ that acts on vectors by
\begin{displaymath}
    \pi(e^X) v \approx e^{\langle \tau, X \rangle} v.
\end{displaymath}
In the simpler situation of finite groups, if $\exp(\mathcal{D})$ were a subgroup of $G$ and
\begin{displaymath}
    e^X \mapsto e^{\langle \tau, X \rangle}
\end{displaymath}
a character thereof, integrating against this character would define a projection onto the space of vectors localised at $\tau$, thus producing such vectors.
We now attempt to approximate this process and deduce some necessary conditions on $R$.

For this we recall the Baker--Campbell--Hausdorff formula (see \cite[Thm.\ 9.14]{Nel-Un}), which asserts that
\begin{equation} \label{eq:BCH}
    e^X e^Y = e^{X \ast Y},
\end{equation}
where
\begin{align*}
     X \ast Y &= X + Y + \frac12 [X,Y] + \frac{1}{12} \left( \left[X,[X,Y]\right] - \left[Y,[X,Y]\right] \right) + \ldots \\
     &= X + Y + O(\norm{X} \cdot \norm{Y}),
\end{align*}
for small enough elements $X, Y \in \gf$.
To have an approximate character of $\exp(\mathcal{D})$, we would need
\begin{displaymath}
    \exp(\langle \tau, X \ast Y \rangle) \approx \exp(\langle \tau, X \rangle) \cdot \exp(\langle \tau, Y \rangle).
\end{displaymath}
By the formula above, this would require
\begin{displaymath}
    \langle \tau, X \ast Y - X - Y \rangle = o(1).
\end{displaymath}
Now assume that $\tau$ has size $T$ in the dual norm induced by the trace pairing from the Frobenius norm on $\gf$. 
The required bound would now follow from Cauchy--Schwarz if 
\begin{displaymath}
    \norm{X} \cdot \norm{Y} \cdot T = o(1),
\end{displaymath}
which would be the case for
\begin{equation} \label{eq:limit-to-loc}
    R \gg T^{1/2+\varepsilon}.
\end{equation}

A more careful inspection of the Baker--Campbell--Hausdorff formula reveals that one should be able to enlarge the domain $\exp(\mathcal{D})$. 
Indeed, using the explicit description of $X \ast Y$ in terms of iterated commutators and the formula
\begin{displaymath}
    \langle [\tau, X], Y \rangle = \langle \tau, [X,Y] \rangle
\end{displaymath}
for the infinitesimal coadjoint representation, we find that
\begin{displaymath}
    \exp(\langle \tau, X \ast Y\rangle) = \exp(\langle \tau, X\rangle) \cdot \exp(\langle \tau, Y \rangle)
\end{displaymath}
for $X$ and $Y$ in the centraliser of $\tau$, which we denote by
\begin{displaymath}
    \gf_\tau = \{ X: [X, \tau] = 0 \}.
\end{displaymath}
Thus, the limits to localisation in the direction of $\gf_\tau$ should be weaker than \eqref{eq:limit-to-loc}, perhaps allowing $R \gg T^{\varepsilon}$. 
If we denote by $\gf_\tau^\flat$ the orthogonal complement of $\gf_\tau$ in $\gf$ (with respect to the Frobenius inner product, although any fixed inner product is fine), we expect to produce localised vectors for the enlarged domain
\begin{displaymath}
    \mathcal{D} = \{ (\xi, \xi') \in \gf_\tau \oplus \gf_\tau^\flat \mid \xi = o(1),\, \xi' = o(T^{-1/2}) \}.
\end{displaymath}
The quantity $T^{-1/2}$ is often called the \emph{Planck scale} (see Sec. 1.7 or 4.6 in \cite{NV}), and the term is also used to describe the smallest scale at which we can expect localisation.

It might be insightful to consider a concrete example, where we can compute such limits to localisation in a more direct way. 
We consider our favourite type of element
\begin{displaymath}
    \tau = i \mat{T}{0}{0}{-T} = 2iT \cdot A^\ast,
\end{displaymath}
for some $T > 0$. 
Let $v$ be a unit vector in an irreducible representation of $G$ and assume that $v$ is $R$-localised at $\tau$. 
We will now compute directly what the limits to localisation are using coordinates on $\gf$, showing that the direction $\gf_\tau$ enjoys weaker conditions.

\begin{rem}
    Having this example at hand, observe that, under the reasonable and desirable assumption that $R = o(T)$, localisation already puts some restrictions on how the Casimir operator acts. 
    Indeed, using the expression \eqref{eq:def-casimir}, we compute the rough estimate
    \begin{equation} \label{eq:casimir-eigval-est}
        \Omega_C \cdot v = (T^2 + o(T^2)) \cdot v.
    \end{equation}
    Comparing this with the interpretation of localised vectors as unit area pieces of the orbit $\Oc_\pi$ from the Kirillov formula confirms the correctness of this assignment.
    Namely, if $\pi$ is the principal series representation $\Pc(iT)$, then $\Oc_\pi = \Oc_{i T A}$ and the localisation of $v$ at $s A$ would imply $s \approx iT$ and vice versa.
\end{rem}

Now let 
\begin{displaymath}
    a:= \norm{Av},\; b := \norm{Bv},\; c := \norm{Cv},
\end{displaymath}
and note that $b, c \leq R$ since $\tau$ evaluates to $0$ on $B$ and $C$.
As observed before, limits to localisation stem from non-commutativity, and here we exhibit this through the relation
\begin{displaymath}
    [B, C] = 2A.
\end{displaymath}
We can write
\begin{align*}
    2 \langle A v , v \rangle &= \langle BC v, v \rangle - \langle CB v, v \rangle \\
    &= - \langle Cv, Bv \rangle + \langle Bv, Cv \rangle \leq 2 bc,
\end{align*}
using the Cauchy--Schwarz inequality. 
On the other hand
\begin{displaymath}
    \langle A v , v \rangle = 2iT + O(R)
\end{displaymath}
by localisation. 
Thus, recalling that $R = o(T)$, we must have
\begin{equation} \label{eq:limit-to-loc-volume}
    bc \gg T.
\end{equation}
This again suggests the bound $R \gg T^{1/2}$.

Notice however that the strength of localisation in the direction $A$ is not limited through the argument above.
Indeed, comparing with the heuristic based on the Baker--Campbell--Hausdorff formula, we have
\begin{displaymath}
    \gf_\tau = \R \cdot A
\end{displaymath}
in this example.

Let us provide yet another, dual perspective on these limits to localisation, now in connection with the Kirillov formula.
First, we observe how our computations are related to the symplectic form on the coadjoint orbit.
One can check, formally or by our explicit description of the orbits, that the tangent space of $\Oc_\pi$ at $\tau$ is spanned by $\ad^\ast B$ and $\ad^\ast C$.
Recall now that the symplectic form at $\tau$ is therefore essentially given by the expression
\begin{displaymath}
    \langle \tau/i, [B, C] \rangle, 
\end{displaymath}
which we recognise in the computations above.
For instance 
\begin{displaymath}
    bc \geq |2 \langle \tau, A \rangle| = |\langle \tau, [B,C] \rangle|.
\end{displaymath}
Thus, at least intuitively, the symplectic form controls the strength of localisation.
One can make this more concrete as follows.

The more general guiding principle is that a localised vector corresponds to a piece of the orbit of unit symplectic area that contains $\tau$, as we noted when discussing the Kirillov formula in \eqref{eq:kirillov-microlocal-explanation}.\footnote{The geometry of the piece of the orbit must still satisfy certain constraints. Nelson and Venkatesh only work with ``balls'' or squares of unit volume. For more details, the interested reader may consult the highly influential paper of Fefferman on the uncertainty principle \cite{fefferman}.}
As in Section \ref{sec:sketch}, using the hyperboloid coordinates in \eqref{eq:hyperboloid-coords}, let
\begin{displaymath}
    U_\tau = \{ T \tau_{v, \phi} \mid |v| \leq \tfrac{1}{b'}, \, |\phi| \leq \tfrac{1}{c'} \} \subset \Oc_\pi,
\end{displaymath}
which is centred around $\tau = T \tau_{0,0}$ (see Figure \ref{fig:localised-intro}).
We would like this set to correspond to essentially \emph{one} vector microlocalised at $\tau$.
By the formula for the symplectic form, we compute that 
\begin{displaymath}
    \int_{U_\tau} \omega_{\Oc_\pi} \approx T \cdot \frac{1}{c'} \cdot \sinh \frac{1}{b'} \approx T (b'c')^{-1},
\end{displaymath}
assuming that $b'$ is large.
From the explicit computations of $\ad_B^\ast$ and $\ad_C^\ast$, one can recognise that $b'$ and $c'$ here correspond to $b$ and $c$ above, respectively.
Asking for the volume of $U_\tau$ to be 1 now implies the same conclusion.

\subsubsection{Orthogonality of localised vectors} \label{sec:orth}

Part of the usefulness of microlocalised vectors stems from orthogonality properties, which we indicate here at a formal level. 
Let $v_1$ and $v_2$ be unit vectors $R$-localised at $\tau_1$ and $\tau_2$ and assume that $\tau_1 \neq \tau_2$. 
If $X \in \gf$ is a basis element, then
\begin{displaymath}
    0 = \langle X v_1, v_2 \rangle + \langle v_1, X v_2 \rangle = \langle \tau_1 - \tau_2, X \rangle \cdot \langle v_1, v_2 \rangle + O(R).
\end{displaymath}
We deduce a bound of the form
\begin{displaymath}
    \langle v_1, v_2 \rangle = O\left( \frac{R}{\norm{\tau_1 - \tau_2}} \right).
\end{displaymath}
The bound can be improved by requiring localisation to higher orders.

More generally, let $l$ be a continuous linear functional on a representation of $G$ and assume that $l$ is invariant under a subgroup $H \leq G$.
Up to technicalities (being only densely defined), the main example to keep in mind is the $K$-period considered in Section \ref{sec:sketch-K-period}.
We can then find a vector $w$ such that
\begin{displaymath}
    l(v) = \langle v, w \rangle.
\end{displaymath}
The invariance can be translated into the property that $Xw = 0$ for $X \in \Lie(H)$.
Now let $v$ be $R$-localised at $\tau$. The discussion above shows that $l(v)$ is small unless $\tau(X) \approx 0$ for $X \in \Lie(H)$.
From a dual perspective, $l(v)$ is small if $v$ is localised away from the orthogonal complement of $\Lie(H)^\wedge$ inside $\gf^\wedge$.

\subsubsection{Examples}
There are some microlocalised vectors that are very familiar to number theorists.
These are namely the holomorphic modular forms.
When lifted to the group, these give rise to weight-$k$ unit vectors $v_k$ in a discrete series representation $\Dc(k)$, where $k$ is the weight of the modular form (see e.g.\ \cite[Sec.\ 3.2]{bump}).
We now show that they are $O(\sqrt{k})$-localised at $i(k-1)C^\ast \in \Oc_{\Dc(k)}$.
Notice here that the norm of this point is approximately $k$, so that by Section \ref{sec:limits} these vectors are maximally localised.

For the action of $\gf_\C$ on weight vectors, we refer to \cite[Prop.\ 2.5.4]{bump}. 
Being a weight $k$ vector means that
\begin{displaymath}
    C v_k = ik v_k.
\end{displaymath}
From the classification of representations, we recall that
\begin{displaymath}
    -4\Omega_C v_k = (A^2 + B^2 -C^2) v_k = \left((k-1)^2 - 1\right) v_k.
\end{displaymath}
Using that $\langle A v_k, A v_k \rangle = -\langle v_k, A^2 v_k \rangle$, we get the bound
\begin{align*}
    \norm{Av_k}^2 + \norm{Bv_k}^2 &\leq \norm{v_k} \cdot \norm{(A^2 + B^2)v_k} \\
    &= \norm{v_k} \cdot \norm{(-4\Omega_C + C^2)v_k} = 2k \norm{v_k}^2.
\end{align*}
The claim now follows. 

As observed before, localisation is stronger in the $C$-direction, which corresponds to the centraliser.
It should also be pointed out that a linear combination of vectors of weight $k + O(1)$, would have the same localisation property, suggesting that this type of analysis is not sufficient for isolating \emph{a single} vector.

Another example comes from the problem of quantum ergodicity and is given by the microlocal lift of Shnirelman.
This is a vector in a principal series representation $\pi = \Pc(iT)$.
This observation is made in \cite[Sec.~7]{Nel-QV-III}.

For this we introduce the raising and lowering operators
\begin{displaymath}
    R = \frac12(A + iB)\, \text{ and } \, L = \frac12 (A - iB)
\end{displaymath}
as in \cite[Sec.\ 2.2, (2.23)]{bump}.
We recall that $v_k$ is a weight $k$ vector in $\pi$, then $Rv_k$ and $Lv_k$ have weights $k+2$ and $k-2$, respectively.

We follow the construction of Wolpert, as presented in \cite{einsiedler-aque}.
Let $v_0$ be the unit spherical vector in $\pi$, that is, of weight $0$.
Define inductively 
\begin{displaymath}
    v_{n+1} = \frac{1}{iT + 1/2 + n} R v_{n}
\end{displaymath}
for $n \geq 0$, and similarly
\begin{displaymath}
    v_{n-1} = \frac{1}{iT + 1/2 - n} L v_n
\end{displaymath}
for $n \leq 0$.
We then put
\begin{displaymath}
    \psi = \sum_{|n| \leq T^{1/2}} v_n.
\end{displaymath}
It is a fact that the $v_n$ are unit vectors and thus, by orthogonality of $K$-types, the norm of $\psi$ is roughly $T^{1/4}$.
The standard microlocal lift is given by $\psi$ normalised.

We now claim that $\psi$ is $O(T^{3/4})$-localised at $\tau = 2iTA^\ast$.
This is not as strong as it can be, and we remedy this below by taking a smooth average of the $v_n$.

First,
\begin{displaymath}
    C \psi = \sum_{|n| \leq T^{1/2}} C v_n = \sum_{|n| \leq T^{1/2}} 2in v_n \ll \sqrt{T} \norm{\psi}
\end{displaymath}
since the support of the average is of size $\sqrt{T}$.
On the other hand, we can compute that
\begin{displaymath}
    R \psi = \sum_{|n| \leq T^{1/2}} (iT + 1/2 + n) v_{n+1} = iT \psi + O(T)
\end{displaymath}
since the marginal terms, say $(iT + 1/2 + n) v_{n+1}$, have norm $O(T)$.
Analogously,
\begin{displaymath}
    L \psi = iT \psi + O(T).
\end{displaymath}
These estimates now imply the claim, recalling that $\norm{\psi} = T^{1/4}$.

To obtain better localisation, we can define
\begin{displaymath}
    \tilde{\psi} = \sum_n f\left( \frac{n}{T^{\frac12}} \right) v_n,
\end{displaymath}
for some fixed, smooth and compactly supported function $f$ on $\R$.
Again, we have $\lVert \tilde{\psi} \rVert \asymp T^{1/4}$.

The localisation in the $C$ direction now follows in the same way.
Concerning the stronger localisation in the $R$ and $L$ directions, for example when considering $R\tilde{\psi} - iT \tilde{\psi}$, we bound
\begin{displaymath}
    \left| f\left( \frac{n+1}{T^{\frac12}} \right) - f\left( \frac{n}{T^{\frac12}} \right) \right| \ll_f T^{-1/2},
\end{displaymath}
by the mean value theorem.
By bounding the support of the average, we get that
\begin{displaymath}
    R\tilde{\psi} - iT \tilde{\psi} \ll_f \sqrt{\sum_n T^2 \cdot T^{-1}} \ll T^{3/4}.
\end{displaymath}
This and similar computations now show that $\tilde{\psi}$ is $O(T^{1/2})$-localised at $\tau$.
In fact, the smooth average gives us localisation to arbitrary order.

\subsection{The Op-calculus}
We now turn to the issue of producing microlocalised vectors in a general way.
We would also like to decompose vectors into microlocalised ones, following the orbit method heuristic where a basis of localised vectors formally indicates the shape of the Kirillov formula, as in \eqref{eq:kirillov-microlocal-explanation}.

The commonly used device for these purposes is an operator calculus.
It was first developed in a quantitative form in \cite{NV} and then refined in \cite{Nel-Un}.
Many of its basic properties are described very efficiently in \cite[Sec.~6.7]{Nel-QV-III}.

Recall that a representation $\pi$ should informally have a basis of vectors $v_i$, microlocalised at points $\lambda_i \in \Oc_\pi$.
These vectors correspond to pieces $B_i \subset \Oc_\pi$ of unit volume in the symplectic measure, such that $\lambda_i \in B_i$.

Now, to suitable functions $a$ on $\Oc_\pi$, or more generally on $\gf^\wedge$, we would like to assign an operator $\Op(a)$, an endomorphism of $\pi$, such that
\begin{equation} \label{eq:op-calc-desire}
    \Op(a) \approx \sum_i a(\lambda_i) \cdot \Proj_{v_i},
\end{equation}
where $\Proj_{v_i}$ is the orthogonal projection onto $v_i$.
Producing a microlocalised vector would then amount to taking $a$ to be a bump function on some $B_i$.

One of the most important properties of these operators is that they formalise the orbit method heuristic by rephrasing the Kirillov formula \eqref{eq:kirillov} as
\begin{displaymath}
    \tr(\Op a) \approx \int_{\Oc_\pi} a.
\end{displaymath}
Other crucial properties are approximate $G$-equivariance and approximate multiplicativity, that is,
\begin{displaymath}
    \Op (a) \cdot \Op(b) \approx \Op(a \cdot b).
\end{displaymath}
The meaning and precise form of these properties are described in the following sections.
However, we note already that the approximation signs above become accurate only for functions $a$ and $b$ that are roughly constant at Planck scale.
This covers the case of bump functions on the pieces $B_i$ mentioned above.

\subsubsection{Operators}
We selected the Haar measure $dX = dA \wedge dB \wedge dC$ on $\gf$ in Section \ref{sec:archi-gp-lie-algebra}, and we choose a dual Haar measure $d\xi$ on $\gf^\wedge$.
More precisely, if $a$ and $\phi$ are Schwartz functions on $\gf^\wedge$ and $\gf$, respectively, then the Fourier transforms
\begin{displaymath}
    a^\vee(X) = \int_{\xi \in \gf^\wedge} a(\xi) e^{-\langle \xi, X \rangle} \, d\xi, \quad \phi^\wedge(\xi) = \int_{X \in \gf} \phi(X) e^{\langle \xi, X \rangle} \, dX,
\end{displaymath}
should be mutually inverse.

Now consider the example of a linear function $a(\xi) = \langle \xi, X \rangle$, for some fixed $X \in \gf$, e.g. one of the basis elements. 
For $v$ localised at $\tau$, the desideratum for the operator calculus is then that $\Op(a) v \approx \langle \tau, X \rangle v = a(\tau) v$.
From the definition of localisation, this suggests that $\Op(a)$ should be the differential operator defined by $X$.

However, to make the approximation sign more meaningful, we perform the following renormalisation, similar to the beginning of Section \ref{sec:limits} on limits to localisation.
Assume that $v$ is $R$-localised at $\tau$, and, to have errors in approximations tend to $0$ rather than just be smaller in order of magnitude, we introduce the parameter $\hrm = 1/R^2$.
Indeed, as in the discussion in Section \ref{sec:limits}, we can think of $\hrm$ as a very small quantity.
If we define $a_\hrm(\xi) = a(\hrm \, \xi)$, then $\Op(a_\hrm)$ should act as $\hrm \cdot X$, and we assume this for our discussion.
Microlocalisation now means that
\begin{displaymath}
    \Op(a_\hrm) v = \hrm X \cdot v = a_\hrm(\tau) \cdot v + O(\hrm^{1/2}).
\end{displaymath}
In this normalisation, the error term can be thought of as $o(1)$, as we will generally consider the limiting behaviour as $R \rightarrow \infty$.
It is therefore practical to work with the quantity $\hrm$, reminiscent of the Planck constant, and this is the common way of phrasing the microlocal calculus.

\begin{rem} \label{rem:renormalisation}
    We stress that, from now on, we work with a renormalised picture.
    We will generally have symbols $a$ concentrated around some \emph{bounded} element $\tau \in \gf^\wedge$, but $\Op(a_\hrm)$ will pick up vectors localised at $\hrm^{-1} \tau$.
    Thus, in practice, instead of working geometrically on coadjoint orbits getting bigger and bigger, we can ``fix'' the orbit, e.g. a hyperboloid, and consider different representations by varying the parameter $\hrm$.
\end{rem}

Recalling the standard properties of the Fourier transform, where differentiation corresponds to multiplication by a linear function on the dual side, we are now led to consider the following pseudo-differential operators.
First, let $\chi$ be a nice cutoff function on $G$, as in \cite[Sec. 8.2]{Nel-Un}, that is a smooth, symmetric bump function in a small neighbourhood of the identity. 
For $\hrm > 0$ and $a$ a Schwartz function on $\gf^\wedge$, define the operator $\Op_\hrm(a)$ on a unitary representation $\pi$ of $G$ by
\begin{displaymath}
    \Op_\hrm(a) v = \int_{X \in \gf} \chi(X) a_\hrm^\vee (X) \pi(e^X) v \, dX,
\end{displaymath}
for $v \in \pi$, and we often write $\Op_\hrm(a:\pi)$ to make the representation explicit.
We include the truncation $\chi$ in this definition because of the use of the exponential map, but the operators we obtain are not very sensitive to the precise choice (see \cite[Sec. 5.4]{NV}).
This allows us to describe $\Op_\hrm$ as a convolution operator at the group level, which will thus guide our choice of test function in the spectral expansion \eqref{eq:intro-spec-expansion} in the introduction.
Indeed, we have that
\begin{displaymath}
    \Op_\hrm(a:\pi) = \pi(\widetilde{\Op}_\hrm(a)) = \int_G \widetilde{\Op}_\hrm(a)(g) \pi(g) \, dg,
\end{displaymath}
as a distribution, where
\begin{displaymath}
    \widetilde{\Op}_\hrm(a)(\exp(X)) = j(X)^{-1} \chi(X) a_\hrm^\vee(X).
\end{displaymath}
These are the type of functions that we will use in the pretrace formula (see Section \ref{sec:def-test-function-archi}).

We make a few easy observations (see \cite{NV}).
First, if $a \equiv 1$ is the constant function, then $\Op_\hrm(a)$ is the identity operator.
Second, if $a$ is a real function, then $\Op_\hrm(a)$ is self-adjoint.
The latter is a consequence of the more general fact that the adjoint $\Op_\hrm(a)^*$ is given by $\Op_\hrm(\overline{a})$, by our assumption that $\chi(-X) = \chi(X)$.
Finally, ignoring the truncation function $\chi$, if $v$ is localised at $\tau$, then we have the desired action
\begin{displaymath}
    \Op_\hrm(a) v \approx \int_{X \in \gf} a_\hrm^\vee (X) \pi(e^X) v \, dX \approx \int a_\hrm^\vee (X) e^{\langle \tau, X\rangle} v = a_\hrm(\tau) v.
\end{displaymath}

\subsubsection{Symbols and star products}
The main properties of the Op-calculus are given by approximate identities and their accuracy depends heavily on the functions $a$ that we consider, commonly called \emph{symbols}.
The sizeable Section 9 of \cite{Nel-Un} discusses the relevant classes of symbols.
The calculus developed there and in \cite{NV} allows for varying degrees of regularity and employs asymptotic notation from nonstandard analysis (see the remarks in \cite[Sec. 9.3]{Nel-Un}).
We simplify this by focusing on the types of symbols used in our applications.

Fix a Schwartz function $\psi$ on $\gf^\wedge$ with values in the interval $[0, 1]$.
For $\delta \geq 0$ and $\tau \in \gf^\wedge$ define
\begin{displaymath}
    \psi^\tau_\delta(\xi) \colon = \psi \left( \frac{\xi - \tau}{\hrm^\delta} \right).
\end{displaymath}
Thus, $\psi^\tau_\delta$ can be interpreted as a bump function on a ball of radius~$h^\delta$ around~$\tau$.
One can check that such a function belongs to the class $S_\delta^{-\infty}$ as defined in \cite[Sec. 9.3]{Nel-Un}. 
We do not define the latter precisely, but encourage the reader to think of functions of the form $\psi_\delta^\xi$ for some bounded $\xi \in \gf^\wedge$ and $\phi$ in some fixed bounded subset of the Schwartz functions.
We also often multiply such classes by scalars, such as $\hrm$ and powers thereof, analogously to subsets of $\R$-vector spaces.

To motivate the following construction, we recall from Section \ref{sec:sketch} that our strategy involves using the operator $\pi(f)$ (see \eqref{eq:intro-spec-expansion}) for a test function~$f$ that is a convolution of two functions.
We shall specify the latter as $\widetilde{\Op}_\hrm(a)$ for some well-chosen symbols $a$.
Applying one of our main tools, the Kirillov formula, requires writing
\begin{displaymath}
    \Op_\hrm(a) \cdot \Op_\hrm(b) \overset{!}{=} \Op_\hrm(a \star_\hrm b)
\end{displaymath}
as a single operator, for some operation $a \star_\hrm b$, or
\begin{displaymath}
    \widetilde{\Op}_\hrm(a) \ast \widetilde{\Op}_\hrm(b) \overset{!}{=} \widetilde{\Op}_\hrm(a \star_\hrm b),
\end{displaymath}
equivalently.\footnote{To be precise, the way we define $a \star_\hrm b$ gives that $\Op_\hrm(a) \cdot \Op_\hrm(b) = \Op'_\hrm(a \star_\hrm b)$, where $\Op'$ is defined with respect to a different truncation function $\chi'$ (see \cite[Sec. 8.7]{Nel-Un}). However, we noted that this only introduces negligible errors (see e.g. \cite[(6.8)]{Nel-QV-III}).}

To make this explicit, write $e^{X \ast Y}= e^X \cdot e^Y $, at least for small $X, Y \in \gf$, as in the BCH formula \eqref{eq:BCH}. 
Ignoring truncation, we formally have
\begin{align*}
    \Op_\hrm(a) \cdot \Op_\hrm(b) &= \int_X \int_Y a^\vee_\hrm (X) b^\vee_\hrm (Y) \pi(e^X) \pi(e^Y) \\
    &= \int_Z \int_{X \ast Y = Z} a^\vee_\hrm (X) b^\vee_\hrm (Y) \pi(e^Z) = \Op_\hrm(a \star_\hrm b).
\end{align*}
We use Fourier inversion to deduce that
\begin{displaymath}
    a \star_\hrm b(\tau) = \int_X \int_Y a^\vee_\hrm (X) b^\vee_\hrm (Y) e^{\langle \tau/\hrm, X \ast Y \rangle}. 
\end{displaymath}
Now, if $\gf$ were commutative, $X \ast Y$ would be equal to $X + Y$ and the above would reduce to standard convolution, meaning that $a \ast_\hrm b = a \cdot b$.
However, we recall that
\begin{displaymath}
    X \ast Y = X + Y + \frac12 [X,Y] + \ldots.
\end{displaymath}
By expanding $e^{\langle \tau, X \ast Y - X - Y \rangle}$ into a power series, we obtain
\begin{displaymath}
    a \star_\hrm b(\tau) = \int_X \int_Y a^\vee_\hrm (X) b^\vee_\hrm (Y) e^{\langle \tau/\hrm, X\rangle} e^{\langle \tau/\hrm, Y\rangle} (1 + \tfrac12 \langle \tau/\hrm, [X,Y] \rangle + \ldots).
\end{displaymath}
The term corresponding to $1$ in the expansion gives precisely $a(\tau) \cdot b (\tau)$.

We sketch the computation for the term corresponding to $\frac12 \langle \tau, [X,Y] \rangle$.
To compute the integrals, we can parametrise them using our basis $(A,B,C)$ for $\gf$.
We then write $[X, Y]$ using this basis for all $X, Y \in \{A, B, C\}$.
For example, $[B,C] = 2 A$ and also $[C, B] = -2 A$.
Thus, along $X = t B$ and $Y = u C$, we have $\frac12 \langle \tau, [X,Y] \rangle = \tau(A) t u$.
By standard properties of Fourier transforms (here simply by partial integration), we thus arrive at
\begin{align}
    &\int_X \int_Y a^\vee_\hrm (X) b^\vee_\hrm (Y) e^{\langle \tau/\hrm, X\rangle} e^{\langle \tau/\hrm, Y\rangle} \tfrac12 \langle \tau/\hrm, [X,Y] \rangle \nonumber \\
    &= \frac{1}{\hrm} \tau(A) \cdot \left(\partial_B a_\hrm \cdot \partial_C b_\hrm - \partial_C a_\hrm \cdot \partial_B b_\hrm\right) (\tau/\hrm) + \ldots  \nonumber \\
    &= \hrm \tau(A) \cdot \left(\partial_B a \cdot \partial_C b - \partial_C a \cdot \partial_B b\right) (\tau) + \ldots \label{eq:poisson-bracket}
\end{align}
where we continue with the other possibilities for $X$ and $Y$.
This last expression, excluding $\hrm$, is often called the Poisson bracket $\{a, b\}$.

In short, we have
\begin{displaymath}
    a \star_\hrm b = a \cdot b + \hrm \{a, b\} + \ldots,
\end{displaymath}
where the omitted terms involve higher order derivatives.
A generalised version of this computation is the content of \cite[Thm. 4.5]{NV} or \cite[Thm. 9.1]{Nel-Un}.

Consider now the case of $a$ and $b$ of the form $\psi^\tau_\delta \in S_\delta^{-\infty}$ as above.
Differentiating each function once, we observe that $\{a, b\}$ lies in $\hrm^{-2\delta} S_\delta^{-\infty}$ and therefore
\begin{equation} \label{eq:star-prod-approx}
    a \star_\hrm b = a \cdot b + \hrm^{1 - 2 \delta} S_\delta^{\-\infty}.
\end{equation}
This is a special case of Theorem 9.5 in \cite{Nel-Un}.
For this statement to be useful, we must assume that $\delta \in [0, 1/2)$.
Upon letting $\hrm = T^{-1}$, this limitation of the calculus with symbols corresponds exactly to the limits to localisation of Section \ref{sec:limits}.
Put another way, we must work with bump functions on balls of radius at least $T^{-1/2}$ (recall Remark \ref{rem:renormalisation} and note that $T^{-1/2} = T^{1/2} \cdot \hrm$).
Under such conditions, we also have an informal approximation
\begin{equation} \label{eq:op-prod-approx}
    \Op_\hrm(a) \cdot \Op_\hrm(b) \approx \Op_\hrm(a \cdot b).
\end{equation}
For brevity, we refer to Section 10 of \cite{Nel-Un} for more details on this last statement.

If we let $a = \psi_{1/2-\varepsilon}^\tau$, then \eqref{eq:op-prod-approx} implies that $\Op_h(a)$ is an approximate idempotent.
Thus, it is an approximate orthogonal projector and the image consists of microlocalised vectors, as desired in \eqref{eq:op-calc-desire}.
For more details on this intuition, see \cite[Rem.~6.6]{NV}.

The rescaled version of Kirillov's formula is given in \cite[Thm. 10.11]{Nel-Un}.
For tempered $\pi$ and $a \in S_\delta^{-\infty}$, with $\delta < 1$, a simplified version states that, for any fixed $J \in \Z_{\geq 1}$, there exist differential operators $\Dc_j$ for $j \in \{1, \ldots, J\}$ such that
\begin{equation} \label{eq:kirillov-explicit}
    \hrm^2 \tr \Op_\hrm(a)  = \int_{\hrm \Oc_\pi} a + \sum_{j = 1}^{J-1} \hrm^j \int_{\hrm \Oc_\pi} \Dc_j a + O(\hrm^{(1-\delta)J}).
\end{equation}
Here, the $2$ in the exponent of $\hrm$ on the left-hand side is the dimension of the orbit $\Oc_\pi$.
Notice that, although our calculus generally requires $\delta<1/2$, Kirillov's formula is more flexible and this is very important.

Indeed, one of the main new tools in \cite{Nel-Un} is an extension of \eqref{eq:star-prod-approx} and \eqref{eq:op-prod-approx} that allows symbols $a$ and $b$ localising at some $\tau \in \gf^\wedge$ to have less regularity in the direction of the centraliser $\gf_\tau$.
This implements the observations in Section \ref{sec:limits}.

\subsubsection{Tailored coordinates and refined symbols} \label{sec:tailored-coords}
To talk about stronger or weaker localisation at $\tau \in \gf^\wedge$ along different directions, we generally need to introduce new coordinates, as in \cite[Sec. 9.4]{Nel-Un}.
If $\tau$ is regular, meaning non-zero in our case, then its centraliser $\gf_\tau \subset \gf$ is a subalgebra of minimal dimension, namely $1$.
For example, if $\tau \in i\R A^\ast \setminus 0$, then $\gf_\tau = \R A$.

We now define $\gf_\tau^\flat \subset \gf$ to be the orthogonal complement of $\gf_\tau$, using the fixed inner product, e.g. the Frobenius inner product.
Dualising, this gives a decomposition
\begin{displaymath}
    \gf^\wedge = \gf_\tau^\perp \oplus \gf_\tau^{\perp \flat},
\end{displaymath}
where $\gf_\tau^\perp$ are the elements that send $\gf_\tau$ to $0$.
If $\tau = iA^\ast$, then $\gf_\tau^\perp$ is the span of $iB^\ast$ and $iC^\ast$, while $\gf_\tau^{\perp \flat} = i\R A^\ast$.
Notice also that $\gf_\tau^\perp$ is isomorphic to the tangent space at $\tau$ of the orbit of $\tau$.

These decompositions give coordinates 
\begin{displaymath}
    X = (X', X'') \in \gf  = \gf_\tau^\flat \oplus \gf_\tau
\end{displaymath}
and
\begin{displaymath}
    \xi = (\xi', \xi'') \in \gf^\wedge = \gf_\tau^\perp \oplus \gf_\tau^{\perp \flat},
\end{displaymath}
tailored to our chosen element $\tau$.
We have $\langle \xi, X \rangle = \langle \xi', X' \rangle + \langle \xi'', X'' \rangle$.

Now let $0 < \delta' < \delta'' < 2 \delta' < 1$ and, for a Schwartz function $\psi$ on $\gf^\wedge$, define
\begin{displaymath}
    \psi_{\delta', \delta''}^\tau \colon = \psi \left( \frac{\xi'}{\hrm^{\delta'}}, \frac{\xi'' - \tau}{\hrm^{\delta''}} \right) 
\end{displaymath}
using $\tau$-tailored coordinates.
This is the prototype of a function in the class $S^{\tau}_{\delta', \delta''}$ of \cite[Sec. 9.4.3]{Nel-Un}.
It is a bump on a \enquote{cylinder} or rectangle around $\tau$, with width $\hrm^{\delta'}$ in the $\gf_\tau^\perp$-direction and width $\hrm^{\delta''}$ in the $\gf_\tau^{\perp \flat}$-direction, i.e. $\tau$-direction.
Our goal is to take $\delta'$ close to $1/2$ and $\delta''$ close to $1$, in which case we have a bump on a \emph{coin-shaped region}, as described by Nelson.
We also recall again the renormalisation we are making, as opposed to the introduction, e.g. Section \ref{sec:intro-refined-test-fn}, and refer to Remark \ref{rem:renormalisation}.

The analogue of \eqref{eq:star-prod-approx} for symbols in $S^{\tau}_{\delta', \delta''}$ is given in Theorem 9.12 of \cite{Nel-Un}.
In particular, it states that
\begin{equation}
    a \star_\hrm b = a \cdot b + \hrm^{1 - 2\delta'} S_{\delta', \delta''}^{\tau},
\end{equation}
for $a, b \in S^\tau_{\delta', \delta''}$.
The $\Op$-composition formula \eqref{eq:op-prod-approx} also holds analogously by \cite[Thm. 10.9]{Nel-Un}.

We sketch the idea for the proof, which is based on an analysis of the BCH formula, in the case $\tau = iA^\ast$, essentially reiterating some observations from Section \ref{sec:limits}.
Here we can work directly with the formula for the Poisson bracket given in \eqref{eq:poisson-bracket}, since the coordinates $(A, B, C)$ are tailored for $\tau$, as noted above.
Since $\tau$ vanishes on $B$ and $C$, we have
\begin{displaymath}
    \{a, b\} = \tau(A) (\partial_B a \cdot \partial_C b - \partial_C a \cdot \partial_B b)(\tau)
\end{displaymath}
with no additional terms.
We are done by noting that derivatives in the $B$ and $C$ directions only introduce factors of $\hrm^{\delta'}$.

\subsubsection{Relative character estimates} \label{sec:rel-char-estimates}

Turning back to our application, we indicate now that our strategy involves applying the relative pretrace inequality \eqref{eq:intro-pretrace-ineq} with $\omega = \widetilde{\Op}_\hrm(a)$ where $a \in S_{\delta', \delta''}^{\tau}$.
Rewriting the inequality without microlocalised vectors (thus making it more precise), observe that we require a lower bound for
\begin{displaymath}
    \sum_{v \in \Bc(\pi_\lambda)} Q(\Op_\hrm(a) v),
\end{displaymath}
where
\begin{displaymath}
    Q(v) = \int_K \langle \pi(k) v, v \rangle \, dk.
\end{displaymath}

To apply Kirillov's formula, the expression above needs to be rewritten as a trace.
We record here the formal computations that are described in more detail in \cite[Sec. 12.2, Sec. 12.3]{Nel-Un}.
For this, write $T = \Op_\hrm(a)$ for short and decompose $Tv$ as 
\begin{displaymath}
    Tv = \sum_{v' \in \Bc(\pi_\lambda)} \langle Tv, v' \rangle v' = \sum_{v'} \langle v, T^\ast v' \rangle v'.
\end{displaymath}
Then
\begin{align*}
    \sum_v Q(Tv) &= 
    \sum_v \int \langle \pi(k) Tv, \sum_{v'} \langle v, T^\ast v' \rangle v' \rangle \\
    &= \sum_v \sum_{v'} \overline{\langle v, T^\ast v'\rangle} \int \langle \pi(k) Tv, v' \rangle \\
    &= \sum_{v'} \sum_v \langle T^\ast v', v\rangle \int \langle v, T^\ast \pi(k)^\ast v' \rangle.
\end{align*}
Applying the same decomposition argument on $T^\ast v'$ backwards, we get
\begin{align}
    \sum_v Q(Tv) &= 
    \sum_{v'} \int \langle T^\ast v', T^\ast \pi(h)^\ast v' \rangle \nonumber \\
    &= \sum_{v'} \int_K \langle \pi(k) T T^\ast v', v' \rangle \, dk = \colon \mathcal{H}(T T^\ast). \label{eq:H-T-Tstar}
\end{align}
This last expression is the relative trace of the operator $T T^\ast$.

To apply the Kirillov formula to relative traces, one introduces the notion of relative coadjoint orbits and stability.
The latter notion is particularly important in the subconvexity problem, where it is related to the conductor dropping phenomenon (see \cite[Sec. 15.4]{NV}).

For this, we generally first define a map $\gf^\ast \rightarrow \Lie(K)^\ast$, denoted $\xi \mapsto \xi_K$, by restriction of linear forms. 
Since $\mathfrak{k} \colon = \Lie(K) = \R C$, this is simply the projection onto the $c$-coordinate used in the discussion of orbits.

As in geometric invariant theory, we say that $\xi \in \gf^\ast$ is \emph{stable} with respect to $K$, which acts as usual on $\gf^\ast$, if it has finite stabiliser in $K$ and a closed $K$-orbit.
We now recall that $K$ acts by rotation around the $C$-axis (see Section \ref{sec:orbits}).
Thus, in our case, all elements $\xi \in \gf^\ast \setminus \R C$ are stable.
We can obviously extend this discussion to $\xi \in \gf^\wedge$ by the usual isomorphism $\gf^\wedge \cong \gf^\ast$.

More generally, one can talk about pairs $(\lambda, \mu) \in [\gf^\wedge] \times [\mathfrak{k}^\wedge]$.
These are called stable if the eigenvalues of $\lambda$ and those of $\mu$ are disjoint.
Here, if $\lambda$ is the class of $2TA^\ast$ for $T \in i\R$, the eigenvalues (independent of representative) are $\{-T, T\}$, i.e. those of $TA$.
We also view $\mathfrak{k}^\wedge$ inside $\gf^\wedge$ to use the same definition for $\mu$, which we can therefore construe as $TC$ for some $T \in i\R$.

Now take $\pi$ and $\sigma$ tempered irreducible unitary representations of $G$ and $K$.
Let $\Oc_\pi$ and $\Oc_\sigma \subset \mathfrak{k}^\wedge \subset \gf^\wedge$ be the corresponding orbits (or multiorbits).
Since $K$ is abelian, in our case $\sigma$ is a character $\exp(\phi C) \mapsto e^{ik\phi}$ for some $k \in \Z$.
One can easily check that, following an analogous setup as for $G$, the orbit $\Oc_\sigma = \{ i\tfrac{k}{2} C\}$ is a singleton.

Define now the \emph{relative coadjoint orbit}
\begin{displaymath}
    \Oc_{\pi, \sigma} = \{ \xi \in \Oc_\pi \colon \xi_K \in \Oc_\sigma \}.
\end{displaymath}
As in \cite[Sec. 11.2]{Nel-Un}, if $(\lambda_\pi, \lambda_\sigma)$ is stable, then $\Oc_{\pi, \sigma}$ is either empty or a $K$-torsor, meaning a closed $K$-invariant subset on which $K$ acts simply-transitively.

Take for example $\pi = \mathcal{P}(iT)$ and $\sigma$ the character $\exp(\phi C) \mapsto e^{ik\phi}$.
Then the eigenvalues of $\lambda_\pi$ are $\{-iT, iT\}$ and the eigenvalues of $\lambda_\sigma$ are $\{ -k/2, k/2\}$.
Thus, if $T > 0$, this pair is always stable.
This corresponds to the obvious visual feature of one-sheeted hyperboloids: the level sets of the projection to the $c$-coordinate (the height coordinate) are always circles and, in fact, $K$-torsors.

\begin{figure}[ht]
    \centering
    \includegraphics[width=\textwidth/2]{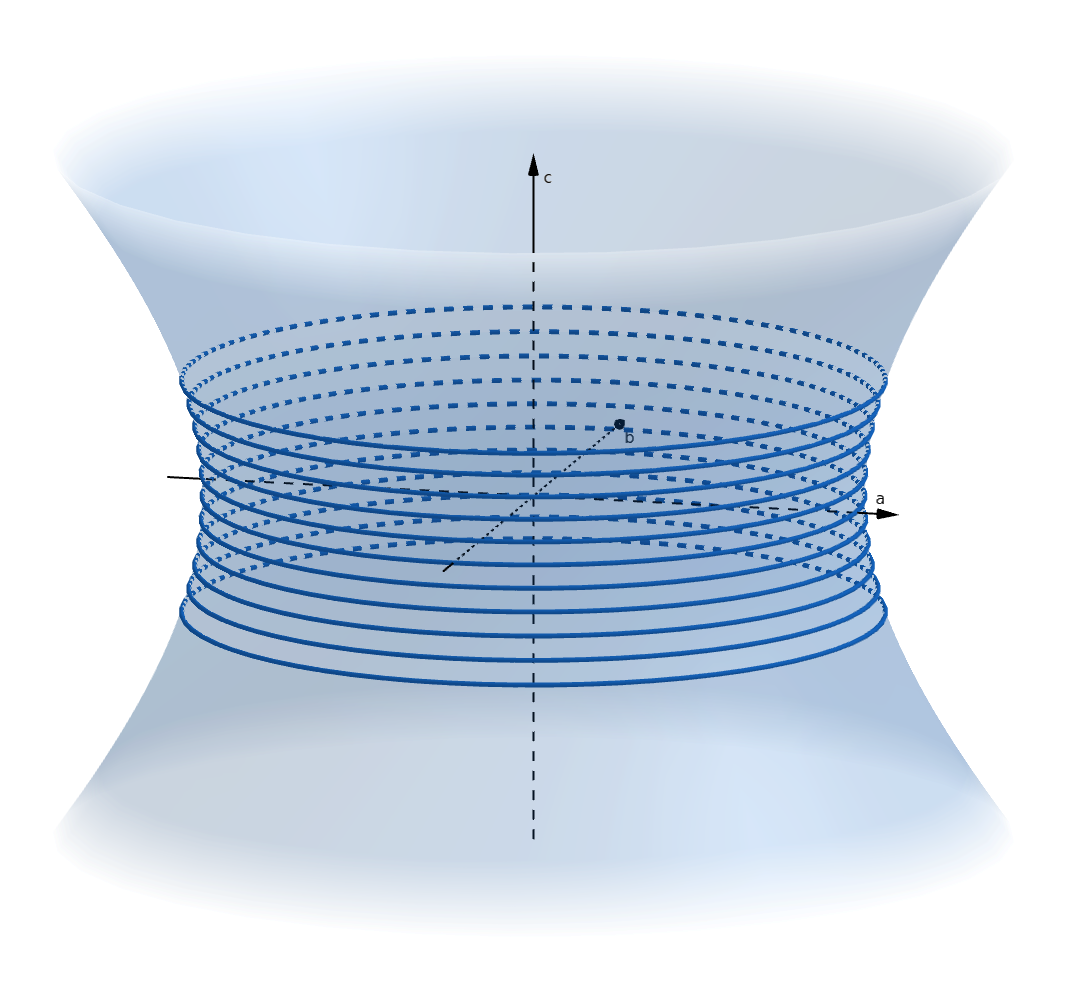}
    \caption{The relative coadjoint orbits for $\pi = \mathcal{P}(iT)$ and the characters $e^{ik\phi}$ on $K$}
    \label{fig:level-lines}
\end{figure}

We have a contrasting picture if we take $\pi = \mathcal{D}(k + 1)$.
Here, the eigenvalues coincide, and this corresponds to the fact that $i\tfrac{k}{2} C \in \Oc_\pi$ is fixed by $K$.
Even if we take a different $\sigma$ to have stability, it is visually clear that there are relative coadjoint orbits which are $K$-torsors, but also some which are empty, since $\Oc_\pi$ is a double-sheeted hyperboloid.

Let us now assume that $\Oc_{\pi, \sigma}$ is a $K$-torsor and equip it with the pushforward $\mu_K$ of our fixed Haar measure on $K$.
Let $\sigma$ be the character $\exp(\phi C) \mapsto e^{il\phi}$ on $K$ and, for $\delta < 1/2$, let $a \in S_\delta^{-\infty}$ be supported on stable pairs.
Then one of the main local results of \cite{NV}, as stated in \cite[Thm. 12.2]{Nel-Un}, is that
\begin{equation} \label{eq:relative-kirillov}
    \sum_{v \in \mathcal{B}(\pi)} \int_K \langle \pi(k) \Op_\hrm(a) v, v \rangle \cdot \overline{\sigma(k)} \, dk
    = \int_{\hrm \Oc_{\pi, \sigma}} a \, d\mu_K + O(\hrm^{1- 2\delta}).
\end{equation}
Notice that, for $l = 0$, we obtain a formula for $\mathcal{H}(\Op_\hrm(a))$.
Looking towards computing the quantity in \eqref{eq:H-T-Tstar}, we can use the composition rule \eqref{eq:op-prod-approx} to obtain the analogous formula 
\begin{equation} \label{eq:relative-kir-triv-sigma}
    \mathcal{H}(\Op_\hrm(a) \Op_\hrm(a)^\ast) \approx \int_{\hrm \Oc_{\pi, \sigma}} \overline{a(\xi)} \cdot a(\xi) \, d\mu_K(\xi).
\end{equation}

\begin{rem}
    This is part of the orbit method philosophy, where restricting a representation to a subgroup corresponds to disintegration along relative adjoint orbits (see \cite[Sec. 1.9]{NV}).
    For example, the more general estimate for any character $\sigma$ of $K$ would be useful for studying the sup-norm problem for non-spherical Maaß forms (see e.g. \cite{BHMM-beyond}).
\end{rem}

The motivation behind \eqref{eq:relative-kirillov} is sketched in \cite[(19.2), (19.3)]{NV}, and we provide here similar computations.
First, take the left-hand side of the equation, and write the integral over $K$ using its Lie algebra $\mathfrak{k}$ and ignoring the Jacobian factor, giving
\begin{displaymath}
    \sum_{v \in \mathcal{B}(\pi)} \int_{Y \in \mathfrak{k}} \langle \pi(e^Y) \Op_\hrm(a) v, v \rangle \cdot \overline{\sigma(e^Y)}.
\end{displaymath}
Insert the definition of $\Op_\hrm(a)$, ignoring the truncation, to have
\begin{displaymath}
    \sum_{v \in \mathcal{B}(\pi)} \int_{Y \in \mathfrak{k}} \int_{X \in \gf} a_\hrm^\vee(X) \langle \pi(e^Y \cdot e^X) v, v \rangle \cdot \overline{\sigma(e^Y)}.
\end{displaymath}
Put the sum over $v$ inside the integrals to obtain
\begin{displaymath}
    \int_{Y \in \mathfrak{k}} \int_{X \in \gf} a_\hrm^\vee(X) \chi_\pi(e^Y \cdot e^X) \cdot \overline{\sigma(e^Y)}.
\end{displaymath}
Assuming that $e^X e^Y \approx e^{X + Y}$ and ignoring the Jacobian again, the Kirillov formula \eqref{eq:kirillov} gives
\begin{displaymath}
    \int_{Y \in \mathfrak{k}} \int_{X \in \gf} a_\hrm^\vee(X) \int_{\tau \in \Oc_\pi} e^{\langle \tau, X \rangle} e^{\langle \tau, Y \rangle} \cdot \overline{\sigma(e^Y)}.
\end{displaymath}
Move the integrals over $X$ and $Y$ inside and simply apply Fourier inversion to compute the $X$-integral and obtain 
\begin{displaymath}
    \int_{\tau \in \hrm \Oc_\pi} a(\tau) \int_{Y \in \mathfrak{k}} e^{\langle \tau, Y \rangle} \cdot \overline{\sigma(e^Y)}.
\end{displaymath}

At this point, one disintegrates the symplectic measure on $\Oc_\pi$ along the projection to the $c$-coordinate.
As noticed above, if $\pi$ is a principal series representation, for instance, then all the fibres of this projection are circles.
They become equipped with measures that can be normalised to be the pushforward of the measure on $K$.
Finally, we use character orthogonality through the integral over $\mathfrak{k}$, since $e^{\langle \tau, Y \rangle}$ gives a character of $\mathfrak{k}$ when fixing the $c$-coordinate of $\tau$.

\subsection{Test functions}
We are now ready to define a test function for our relative trace formula application, following \cite[Sec. 14]{Nel-Un}, where we choose the representation $\sigma$ of $K$ to be the trivial character.
Let $\pi = \mathcal{P}(iT)$ for some large $T \in \R_{>0}$ and let $\hrm = T^{-1}$ and $\tau = 2iA^\ast \in \gf^\wedge$.
Recall here that our usual $(a,b,c)$-coordinates are $\tau$-tailored coordinates, where the $a$-coordinate is the $\xi''$-coordinate. 

\subsubsection{The lower bound}
Let $0<\delta' < \delta'' < 2\delta' < 1$, with $\delta' \approx 1/2$ and $\delta'' \approx 1$, and let $a = \psi^\tau_{\delta', \delta''}$ for a smooth bump function $\psi$ with values in the interval $[0,1]$.
More precisely, assume that $a(\tau + \xi) \neq 1$ implies that $\xi' \gg \hrm^{\delta'}$ or $\xi'' \gg \hrm^{\delta''}$, and also that $a(\tau + \xi) \neq 0$ implies that $\xi' \ll \hrm^{\delta'}$ and $\xi'' \ll \hrm^{\delta''}$.
Thus, $a$ is a bump on a coin-shaped region, of width $\hrm^{\delta'}$ in the $B$- and $C$-directions, and width $\hrm^{\delta''}$ in the $A$-direction.

In an ideal world, we would apply \eqref{eq:relative-kir-triv-sigma} to $a$.
Recall that 
\begin{displaymath}
    \Oc_{\pi, \sigma} = \{ T(\cos \phi, \sin \phi, 0) \mid \phi \in [0, 2\pi] \}
\end{displaymath}
(since $\sigma$ here is the trivial representation), with the obvious measure coming from $\PSO(2)$.
We would then have
\begin{displaymath}
    \mathcal{H}(\Op_\hrm(a) \Op_\hrm(a)^\ast) \approx \int_{\hrm \Oc_{\pi, \sigma}} a(\xi)^2 \, d\mu_K(\xi).
\end{displaymath}
Estimating this integral boils down to computing the intersection of $\hrm \Oc_{\pi, \sigma}$ with the coin-shaped region that forms the support of $a$.

Notice that, if $\xi = (\cos \phi, \sin \phi, 0) \in \hrm \Oc_{\pi, \sigma}$, then $\xi' = (0, \sin \phi, 0)$ and $\xi'' = (\cos \phi, 0, 0)$.
Recall that $a(\xi) = \psi^\tau_{\delta', \delta''}(\xi) = 1$ if $\xi' \ll \hrm^{\delta'}$ and $\xi'' - 1 \ll \hrm^{\delta''}$.
For $\hrm$ small, we have 
\begin{displaymath} 
    \abs{\xi'} \asymp \phi, \quad \abs{\xi'' - 1} \asymp \phi^2.
\end{displaymath}
By our assumptions on $\delta'$ and $\delta''$, we obtain a lower bound on the integral above of the shape 
\begin{displaymath}
    \mathcal{H}(\Op_\hrm(a) \Op_\hrm(a)^\ast) \gg \hrm^{\delta'} \gg \hrm^{1/2}.
\end{displaymath}
In fact, analogous computations also give an upper bound of $\hrm^{\delta'}$ by considering the support of $a$.

However, this argument is not valid as it is, since \eqref{eq:relative-kir-triv-sigma} only applies to $a \in S^{-\infty}_\delta$ with $\delta < 1/2$.
In \cite[Sec. 14]{Nel-Un}, this is handled by defining the symbol $a$ to first be a bump function on the larger region $\{\tau + \xi \mid \xi \ll \hrm^{\delta'} \}$.
This region can be partitioned into the coin-shaped region $\{\tau + \xi \mid \xi' \ll \hrm^{\delta'}, \xi'' \ll \hrm^{\delta''} \}$, on which we take a bump function $a'$ exactly as before, and the complement $\{\tau + \xi \mid \xi' \ll \hrm^{\delta'}, \hrm^{\delta''} \ll \xi'' \ll \hrm^{\delta'} \}$ with corresponding bump function $a''$.
We can then assume that $a = a' + a''$.

The arguments above are now valid for $a$.
To deduce the advertised bound for $\mathcal{H}(\Op_\hrm(a') \Op_\hrm(a')^\ast)$, we show that $\Op_\hrm(a'')$ is negligible on $\pi$.
For this, the crucial observation is that the support of $a''$ is disjoint from $\hrm \Oc_{\pi}$.
In that case, we are done using the Kirillov formula \eqref{eq:kirillov-explicit}, expanding to arbitrary degree.

To see the statement about the support of $a''$, let $\tau + \xi \in G \cdot \tau$ and $\xi \ll \hrm^{\delta'}$.
Thus, $\tau + \xi$ lies in the intersection of $\hrm \Oc_{\pi}$ and the support of $a$.
Using our parametrisation of the orbit, if $\tau + \xi = \tau_{v, \phi}$, then we can approximate 
\begin{displaymath}
    \tau + \xi \approx (1 + v^2/2 - \phi^2/2, \phi, v)
\end{displaymath} 
by Taylor expansion and the fact that $v$ and $\phi$ are small, since we assume that $\xi$ is small.
It follows that $v, \phi \ll \hrm^{\delta'}$ and, therefore, the $a$-coordinate of $\xi$, approximated by $ v^2/2 - \phi^2/2$ is bounded by $\hrm^{2 \delta'}$ and we are done.
The same argument, using Taylor expansion, works more generally to understand the curvature of orbits in tailored coordinates, as in \cite[Lemma 13.10]{Nel-Un}.

\subsubsection{Definition of the test function} \label{sec:def-test-function-archi}
The conclusion of our arguments is as follows.
Let $f = \widetilde{\Op}_\hrm(a')$, meaning that $\pi(f) = \Op_\hrm(a')$.
Then we have
\begin{equation} \label{eq:archi-lower-bd-spec-side}
    \sum_{v \in \Bc(\pi_\lambda)} Q(\pi(f) v) = \mathcal{H}(\Op_\hrm(a') \Op_\hrm(a')^\ast) \gg \hrm^{1/2} = T^{-1/2}.
\end{equation}
This is the precise version of \eqref{eq:sketch-stability-lower-bd} in our sketch argument and can be found as a special case of \cite[(4.11)]{Nel-Un}.

We recall that there is a truncation function $\chi$ in the definition of $\widetilde{\Op}_\hrm(a')$.
This can be replaced, with negligible error, so that $f$ is truncated to the essential support of $(a'_\hrm)^\vee$.
We assume this from now on, referring to Section 14.6 in \cite{Nel-Un} for details, and study the support and size of $f$ below.

\subsubsection{The upper bound}\label{tub}

We start by making the test function more explicit.
We have
\begin{displaymath}
    f(\exp X) = j(X)^{-1} \chi(X) (a'_\hrm)^\vee(X),
\end{displaymath}
where $\chi(X)$ is a truncation function and $(a'_\hrm)^\vee$ is the Fourier transform of $a'_\hrm$.
Recall now from Section \ref{sec:tailored-coords} that we have $\tau$-tailored coordinates $(X', X'')$ on $\gf$.
For $\tau = 2 i A^\ast$, our choice above, $X''$ is the $A$-direction and $X'$ the $B$- and $C$-direction.
Since $a'(\xi) = \psi^\tau_{\delta', \delta''}$, by standard properties of the Fourier transform, we see that
\begin{displaymath}
    \abs{(a'_\hrm)^\vee (X)} = \hrm^{-3} \cdot \hrm^{2 \cdot \delta'} \cdot \hrm^{\delta''} \cdot \psi^\vee \left( \frac{X'}{\hrm^{1-\delta'}}, \frac{X''}{\hrm^{1-\delta''}} \right).
\end{displaymath}
Here, the first factor comes from rescaling $a'$ to $a'_\hrm$.

It follows that the essential support of $(a'_\hrm)^\vee$ is $X' \ll \hrm^{1-\delta'}$ and $X'' \ll \hrm^{1-\delta''}$.
Exponentiating and recalling that $\delta' \approx 1/2$ and $\delta'' \approx 1$, this is the same as the set $\mathcal{U}$ in the introduction, Section \ref{sec:intro-refined-test-fn}.
A priori, $(a'_\hrm)^\vee$ is certainly not compactly supported, but as mentioned before, we can assume that $f$ is truncated to this essential support.

Furthermore, it is clear from its abstract or explicit form (computed in Section \ref{sec:archi-gp-lie-algebra}) that $j(X) \gg 1$ for small $X$.
Therefore,  we conclude that
\begin{displaymath}
    f(X) \ll \hrm^{-1 -\varepsilon} = T^{1 + \varepsilon}.
\end{displaymath}
Taking into account the oscillation of $f$ described by $\tau$, this is essentially the same test function $\omega$ as in Section \ref{sec:intro-refined-test-fn}.

Continuing to follow Section \ref{sec:sketch}, we are left with bounding certain integrals of $f$ from above.
More precisely, we concentrate on the bound \eqref{eq:intro-decay} that states
\begin{displaymath}
    \int_K \int_K [f \ast f^\ast](x^{-1} \gamma y) \, dx \, dy \ll \frac{T^{1/2 - \delta}}{d(\gamma)}
\end{displaymath}
for any $\gamma \in G \setminus K$.
The proof in our case follows the sketch given in Section \ref{sec:intro-transv} quite closely, but we remark that it overly simplifies the more general problem tackled by Nelson in \cite{Nel-Un}.
In higher rank, the proof is considerably more intricate and interesting and relies on reductions to some Lie-algebraic results (see Sections 16 and 17 in \cite{Nel-Un}) that are trivial in our situation.

First, notice that the support of $f \ast f^\ast$ is essentially the same as that of $f$.
This is because, if $X = (X', X'')$ and $Y = (Y', Y'')$ satisfy that $X', Y' \ll \hrm^{1-\delta'}$ and $X'', Y'' \ll \hrm^{1-\delta''}$, then so do $-X$ and $X+Y$.
We call this support $\mathcal{U}$ and note that it is of the shape
\begin{displaymath}
    (D + O(T^{-1/2-\varepsilon})) \cap (1 + O(T^{-\varepsilon})).
\end{displaymath}

Using this and the bound above, we can follow the first few paragraphs of Section \ref{sec:intro-transv} and reduce the problem to showing that
\begin{displaymath}
    \Vol \{ y \in K \mid \gamma y \in K \mathcal{U} \} \ll T^{-\delta}/d(\gamma)
\end{displaymath}
for some $\delta > 0$.
This is done more generally and precisely in Sections 15.3-5 in \cite{Nel-Un}.
As mentioned in Section \ref{sec:intro-transv}, we are also allowed to assume that $\gamma \in D$, and this reduction is done in Section 16.2 in \cite{Nel-Un}.

We again pass to the coadjoint orbit by noting that $\gamma y \in K \mathcal{U}$ implies that $\dist(\gamma y \tau, K \tau) := \min_{k \in K} \norm{\gamma y \tau - k \tau} \ll T^{-1/2 + \varepsilon}$.
Instead of an ad-hoc argument with angles and the size of hypotenuses, we now follow Nelson by reducing the volume estimate to the size estimate
\begin{equation} \label{eq:size-estimate}
    \norm{y - 1} \cdot \norm{\Ad(\gamma) - 1} \ll \dist(\gamma y \tau, K \tau),
\end{equation}
as in \cite[Theorem 16.2]{Nel-Un}.
This separates the two variables $y$ and $\gamma$ and allows us to compute the integral over $y$ directly.

Before that, we only need the remark that the operator norm of $\Ad(\gamma) - 1$ is roughly equal to $d(\gamma)$ for $\gamma \in D$.
Indeed, if $\gamma = \diag(e^t, e^{-t})$, then $\Ad(\gamma)$ acting on $\gf$ has eigenvalues $1$, $e^{2t}$, and $e^{-2t}$.
Therefore, $\norm{\Ad(\gamma) - 1} \asymp \abs{2t}$ for small $t$ (we can assume $t$ is small since $f$ is supported on a small neighbourhood).
On the other hand, it is well-known that the hyperbolic distance between $\gamma i = e^{2t} i$ and $i$, by definition $d(\gamma)$, is $\abs{2t}$. 

Therefore, we are computing the volume of the set of $y \in K$ such that
\begin{displaymath}
    \norm{y - 1} \cdot d(\gamma) \asymp \norm{y - 1} \cdot \norm{\Ad(\gamma) - 1} \ll \dist(\gamma y \tau, K \tau) \ll T^{-1/2+\varepsilon}.
\end{displaymath}
This volume is bounded by $T^{-1/2}/d(\gamma)$, since $K$ is one-dimensional, and this is the desired bound.

We now give some indications regarding the proof of \eqref{eq:size-estimate}, following Section 16.7 in \cite{Nel-Un}.
Take $k \in K$ and note that 
\begin{displaymath}
    \norm{\gamma y \tau - k \tau} = \norm{\gamma y \gamma^{-1} \tau - k \tau} = \norm{\Ad(\gamma)(y) \tau - k \tau},
\end{displaymath}
since $\gamma$ lies in the stabiliser of $\tau$.
Since $\gamma \in \exp(\R A)$ and $y \in \exp(\R C)$, we must have that the logarithm of $\Ad(\gamma)y$ (we may assume here that all elements are small enough) lies in the subspace $\gf_\tau^\flat = \R B + \R C$.
The space $\gf_\tau^\flat$ is a complement (in fact, the orthogonal complement) to $\gf_\tau$.
This observation can be used to show that the exponential map on small elements of $\gf_\tau^\flat$ preserves distances when considering the action on $\tau$ (see Section 16.6 in \cite{Nel-Un} -- there, one should consider $W = \gf_\tau^\flat$).
More precisely, we have
\begin{align*}
    \norm{\Ad(\gamma)(y) \tau - k \tau} &\asymp \norm{\Ad(\gamma)(\log y) \tau - \log(k) \tau} \\
    &= \norm{[\Ad(\gamma)(\log y), \tau] - [\log(k), \tau]},
\end{align*}
essentially by plugging in the first two terms in the series expansion of the exponential.
The square brackets denote that we now think of all elements in the space $\gf_\C$.

Now, $\log(k) \in \R C$ and therefore $[\log(k), \tau] \in i \R B$, since we identify $\tau$ as usual with $i A$.
If we let $p_\tau$ be the orthogonal projection inside $\gf^\wedge$ onto the space spanned by $iA$ and $iC$, we have that
\begin{displaymath}
    \norm{[\Ad(\gamma)(\log y), \tau] - [\log(k), \tau]} \geq \norm{p_\tau([\Ad(\gamma)(\log y), \tau])}.
\end{displaymath}

Next, $\log y$ is of the form $\alpha \cdot C$, where $\alpha = \norm{\log(y)}/\norm{C}$.
For small $y$, we approximate $\log y \asymp \norm{y - 1}$.
By linearity,
\begin{displaymath}
    \norm{p_\tau([\Ad(\gamma)(\log y), \tau])} \asymp \norm{y - 1} \cdot \norm{p_\tau([\Ad(\gamma)C , \tau])}.
\end{displaymath}
To prove that the last factor is approximated by $\norm{\Ad(\gamma) - 1}$, Nelson uses a version of the implicit function theorem and a crucial Lie-theoretic result, namely Theorem 16.3 in \cite{Nel-Un}.
The latter states in our case that $[A, [C, \tau]] \notin [\R \cdot C, \tau]$.
For us, this is a trivial computation, yet in more general situations, especially in higher rank, this is an interesting problem, and it is solved in Section 17 of \cite{Nel-Un}.

Since $[A, [C, \tau]] = [[A, C], \tau]$ and $[C, A] = 2B$, the Lie-theoretic result states that $p_\tau([\Ad(\gamma)(\log y), \tau])$ is not zero, and this is an algebraic way of describing transversality as in Figure \ref{fig:band-ring}.
We can directly compute that
\begin{align*}
    p_\tau([\Ad(\diag(e^t, e^{-t})C , \tau]) &= i p_\tau([-\sinh(2t) B + \cosh(2t) C, A]) \\
    &= -2i \sinh(2t) C 
\end{align*}
This is essentially the same computation as in Remark \ref{rem:tilted}, used to quantify transversality in the introductory sketch in the last part of Section \ref{sec:sketch}.
As above, $\norm{\Ad(\gamma) - 1} \asymp \abs{2t} \asymp \sinh(2t)$, and we are done.

\subsection{The main archimedean estimates}

All the considerations above can now be condensed into the following theorem, a special version of \cite[Thm.\ 4.2]{Nel-Un}.

\begin{theorem} \label{thm:archi-main-estimate}
    Let $T$ be a large positive real number and let $\pi = \mathcal{P}(iT)$ be a principal series representation.
    For any $\varepsilon > 0$, there exists a test function $f \in C_c^\infty(G)$, supported on a small neighbourhood of the identity, such that:
    \begin{enumerate}
        \item for any orthonormal basis $\mathcal{B}(\pi)$ of $\pi$, we have
            \begin{displaymath}
                \sum_{v \in \Bc(\pi_\lambda)} Q(\pi(f) v) \gg T^{-1/2};
            \end{displaymath}
        \item we have
            \begin{displaymath}
                \int_K \abs{f \ast f^\ast} \ll T^{1/2 + \varepsilon};
            \end{displaymath}
        \item for any $\gamma \in G \setminus K$, we have
            \begin{displaymath}
                \int_K \int_K \abs{f \ast f^\ast(x^{-1} \gamma y)} \, dx \, dy \ll \frac{T^{\varepsilon}}{d(\gamma)}.
            \end{displaymath}
    \end{enumerate}
\end{theorem}

We can compare these estimates to those provided by Iwaniec and Sarnak in Lemma 1.1 of \cite{IS95}.
First, notice that the latter provides a function $k$ that is bi-$K$ invariant.
As such, it should be compared to the function
\begin{displaymath}
    g \mapsto \int_K \int_K f \ast f^\ast (x^{-1} g y) \, dx \, dy.
\end{displaymath}
Scaling $k$ to $T^{-1/2} k$, we now have similar behaviour on the spectral side.
Indeed, as sketched out in Section \ref{sec:sketch}, both test functions pick out a family of representations with spectral parameter in an $O(T^{\varepsilon})$-ball around $T$.
The contribution of our automorphic form of interest, $\phi_\lambda$, is at least $T^{-1/2}$ (see \cite[(1.6)]{IS95}).

On the geometric side, the trivial bound (i.e. part (2) in Theorem \ref{thm:archi-main-estimate}) is essentially given by $T^{1/2}$ in both cases.
However, Iwaniec and Sarnak manage to prove stronger decay properties.
Indeed, in their notation we have $u(\gamma, i) \asymp d(\gamma, i)^2$ for small $d(\gamma, i) = d(\gamma)$ (see Section \ref{sec:analysing-supp-condition} for more details), and their bound is of the shape
\begin{displaymath}
    k(\gamma) \ll d(\gamma)^{-1/2}, \label{eq:comp_IS_bounds}
\end{displaymath}
a big improvement over the bound achieved above.

Of course, the reason for having a weaker bound is the fact that we put absolute values inside the two $K$-integrals.
Since our test function oscillates, there is room for improvement, as mentioned in Remark 4.3 of \cite{Nel-Un}.
We will observe a similar phenomenon in the $p$-adic case. See Table~\ref{table:Whatnot} below.

However, Nelson observes \cite[Rem. 1.4]{Nel-Un} that optimising the entire argument would give the same bound as Iwaniec and Sarnak's, which is indeed the case.
More precisely, the optimisation can be performed in the counting methods, as we discuss in Remark \ref{rem:com_counting_IS} and the sections following it.

\section{Non-archimedean localisation}\label{sec:p-adic}

In the $p$-adic setting, versions of the orbit method were already developed very early on by Howe in \cite{How1} or \cite{How2}, for example. 
Modern applications of ideas along these lines have been given in \cite{Hu-Nel-Sa, Nel-QUE-QP, Ne-Hu_test, Hu-Nel, marsh}, among others.

We take a slightly different approach than in Section \ref{sec:archi-microlocalisation-big-chap} and work mostly at the group level instead of the Lie algebra.
One reason for doing so is that microlocalisation can be achieved in a much stronger sense when doing $p$-adic analysis.
It can be stated in simpler terms that are perfectly adequate for practical applications, and we prefer to do so in this paper.
Besides, twisting the $p$-adic picture to fit its much more pleasantly visual real counterpart seemed to offer more confusion than helpful intuition.

To further simplify computations, we eventually assume that the representations we are considering have a level with exponent divisible by $4$. Similar assumptions are made in \cite{marsh, Hu-Nel}. Treating general levels would require us to consider more cases and greatly complicate the notation. While for supercuspidal representations all the ingredients we need are developed in \cite{Ne-Hu_test}, the case of principal series would require some additional computations. Note that our method will fail for representations with small level such as the Steinberg representation (and quadratic twists thereof) or principal series with little ramification.

The main goal of this section is thus to complement the discussion Section~\ref{sec:archi-microlocalisation-big-chap} with the natural $p$-adic analogues. However, we are doing so with our application to the sup-norm problem in mind. To make the navigation of this section easier we will now present a little summary of the key points in each section:
\begin{itemize}
	\item In Section~\ref{glam} we introduce relevant groups and measures over non-archimedean local fields. We also briefly talk about the Lie algebra.
	\item Section~\ref{bn-ar} contains a crash course in the theory of admissible smooth representations of $\GL_2$ over non-archimedean fields. The main goal is to provide the language necessary for our discussion of $p$-adic microlocalisation. On the way we provide results that are indispensable for our application to the sup-norm problem. Most notably, Convention~\ref{conv:newvec} contains the exact local description of the \textit{translated newform} that plays a key role in the global application. Furthermore, Remark~\ref{rem:conductor_class} classifies the representations that can appear as local constituents of (cuspidal) automorphic representations of level $p^{4n}$ and with trivial central character. Finally, in Section~\ref{mcur} we discuss properties of truncated matrix coefficients that will be used when constructing the test function in the global application. 
	\item The conceptual core of this chapter is contained in Section~\ref{mdc}. Indeed, Definition~\ref{def:loc_vec} is a precise definition of what we mean by $p$-adic localisation. We then discuss several basic consequences of this definition. Section~\ref{ltl} is analogous to Section~\ref{sec:limits} and answers the question: \textit{What is the smallest scale at which localised vectors can exist?} Finally, in Section~\ref{ooo} we establish orthogonality properties of localised vectors. The corresponding archimedean theory was contained in Section~\ref{sec:orth}. The last part of this section links back to the discussion in Section~\ref{mcur}. Here we formulate a heuristic expectation, see \eqref{eq:matrix_coeff_exp.} and Remark~\ref{rem:about_f1f2},  concerning the shape of matrix coefficients associated to localised vectors.
	\item The abstract notion of localised vectors will be made concrete in Section~\ref{emv}, where we will present two important examples, namely \textit{minimal vectors} and \textit{microlocal lift vectors}. We describe their microlocal support (see Definition~\ref{def:cusp_type} and Definition~\ref{def:chi-type} respectively) and compute (the support of) their matrix coefficients (see Lemma~\ref{lm:mat_coeff_super} and Lemma~\ref{lm:mat_coeff_princ} respectively). The latter will be the basis for the geometric estimates that drive our global application. Further, we compute representations of these localised vectors in the Kirillov model and relate them to the generalised new-vector $v^{\circ}$ (see Corollary~\ref{cor:decomp_newform_super} and Corollary~\ref{cor:newform_exp_princ} respectively). This will be needed in order to control the spectral side of the relative trace formula, which appears in our application.
	\item In Section~\ref{sec:loc_perio_fin} we use the theory developed so far, in particular the connection between our examples of localised vectors with generalised new-vectors, to estimate an important local period. Indeed, the local period $Q$ defined in \eqref{eq:loc_per_def} appears as a constituent in the factorization of the global period from \eqref{Q_as_inner_prod}. The most important result of this section is Lemma~\ref{lm:summary_Q}, which is the $p$-adic counterpart to the first statement of Theorem~\ref{thm:archi-main-estimate}.
	\item Finally, in Section~\ref{anavs} we establish $p$-adic volume bounds analogous to the archimedean estimates from Section~\ref{tub}. The main technical result is Theorem~\ref{th:local_vol_bound} below. These volume bounds can be used to estimate the local orbital integrals defined in \eqref{eq:def_Igamma}. These appear on the geometric side of the relative trace formula from Proposition~\ref{pr:pre_amplified_guy}. In practice, we will use the bound stated in Corollary~\ref{cor:vol_bound_p-adig}, which can be seen as a $p$-adic version of the second and third part of Theorem~\ref{thm:archi-main-estimate}. Let us stress that the key point of this section is that we can control the non-archimedean orbital integral using relatively soft geometric ideas. We end this chapter with Section~\ref{sec:p-adic_relmat}, where we compare our bounds to the results obtained by Saha and Hu in \cite{Hu-Sa}. The latter bounds are slightly stronger, but their proof relies on hard computations using the $p$-adic method of stationary phase.\footnote{We have been informed by Paul Nelson in private communication that our soft geometric tools can be pushed further to fully recover the bounds from Saha and Hu.}
\end{itemize}

All the computations here are local in nature, and it will be notationally convenient to work directly over some non-archimedean local field $F$ of characteristic $0$. 
We let $\mathcal{O}$ be the valuation ring of $F$ and we denote its unique maximal ideal by $\mathfrak{p}$. 
The residual field is defined as $\mathfrak{k}=\mathcal{O}/\mathfrak{p}$ and we assume that $\operatorname{char}(\mathfrak{k})\neq 2.$ 
Further we put $q= \# \mathfrak{k}$, choose a uniformiser $\varpi$, let $\nu$ be the corresponding valuation, and normalise the absolute value on $F$ by $\vert\varpi\vert = q^{-1}$.
Of course, in our setting, the relevant field is $F = \Q_p$, where $q = p$ and we assume $p > 2$.

The local zeta-factor is given by
\begin{equation*}
    \zeta_q(s)=(1-q^{-s})^{-1}.
\end{equation*}
If we equip $(F,+)$ with a Haar measure normalised $dy$ such that $\Vol(\mathcal{O})=1$, then we have
\begin{equation*}
    \int_{\mathcal{O}^{\times}} \frac{dy}{\vert y\vert} = \zeta_q(1)^{-1}.
\end{equation*}
Recall that $\frac{dy}{\vert y\vert}$ is a Haar measure for $F^{\times}$. 
Finally, we fix an unramified additive character $\psi\colon F\to S^1$.

\subsection{Groups, Lie algebras and measures}\label{glam}

Departing slightly from the notation in the previous part of this paper, we work with the group $G = \GL_2(F)$. We will start by introducing relevant subgroups, measures as well as the $p$-adic version of the Lie algebra.

\subsubsection{Subgroups and measures}

The diagonal torus in $G$ is
\begin{equation*}
    D=\{\diag(y_1,y_2)\colon y_1,y_2\in F^{\times}\}.
\end{equation*}
Note that $D\cong F^{\times}\times F^{\times}$ and we obtain a Haar measure by setting
\begin{equation*}
    \int_{D} f(t) \, dt = \zeta_q(1)^2\cdot \int_{F^{\times}} \int_{F^{\times}} f(\diag(y_1, y_2)) \, \frac{dy_1\, dy_2}{\vert y_1y_2\vert}.
\end{equation*}
Denote $a(y) = \diag(y,1)$ and $a(y_1, y_2) = \diag(y_1, y_2)$.

Next we define
\begin{equation*}
    N=\left\{ n(x) = \left(\begin{matrix} 1 & x \\0&1\end{matrix}\right)\colon x\in F\right\}.
\end{equation*}
Since $N\cong F$, it inherits the Haar measure from $F$. 
Similarly, we can lift the character $\psi$ from $F$ to $N$ by setting 
\begin{displaymath}
    \boldsymbol{\psi}^{(a)}(n(x)) = \psi(a\cdot x). 
\end{displaymath}
The standard Borel subgroup of $G$ is given by $B=DN$. %

Finally, let $K_0 = \GL_2(\mathcal{O})$.\footnote{In the global setting we will usually write $K_p=\GL_2(\Z_p)$. However, in this setting the notation $K_0$ will be more convenient.} 
For $n>0$, we will also encounter the Hecke congruence subgroups
\begin{equation}
    K_H(n)=\{ k \in K_0\colon k_{2,1}\in \mathfrak{p}^n\}, \label{eq:def_KH}
\end{equation}
where $k_{2,1}$ is the lower left entry.
We equip $K_0$ with the probability Haar measure $\mu_{K_0}$.

Using the Iwasawa decomposition we find that the Haar measure $\mu$ of $G$ can be written as
\begin{equation*}
    \int_{G} f(g) \, d\mu(g) = \zeta_q(1)^2 \cdot \int_{F}\int_{F^{\times} \times F^\times} \int_{K_0} f(n(x)a(y_1, y_2)k) \, d\mu_{K_0}(k) \, \frac{dy_1 \, dy_2}{\vert y_1 \vert^2} \,dx.
\end{equation*}

\subsubsection{The Lie algebra}

The Lie-algebra of $G(F)$ is given by $\mathfrak{g}' = \Mat_{2\times 2}(F)$ and the adjoint action is simply conjugation $\Ad_g(X)=gXg^{-1}$. 
Neighbourhoods of $0\in \mathfrak{g}'$ are given by $$L(n)=\Mat_{2\times 2}(\mathfrak{p}^n) \text{ for } n \in \Z$$ 
Note that if $\Vert X\Vert = \max_{i,j}\vert X_{i,j}\vert$ is the maximum norm on $\mathfrak{g}'$, then 
\begin{displaymath}
    L(n) = \{ X\in \mathfrak{g}'\colon \Vert X\Vert \leq q^{-n}\}
\end{displaymath} 
is nothing but a ball of radius $q^{-n}$. We also set 
\begin{displaymath}
    K(n)=1+L(n)\subseteq K_0
\end{displaymath} 
for $n\in \N$. As for the exponential map, we define $e\colon L(1) \to K(1)\subseteq G(F)$ by
\begin{equation*}
    e(X)=I_2+X \text{ for }X\in L(1),
\end{equation*}
which is a good approximation in practice.

\begin{rem}
Let $\alpha_F$ be the smallest integer such that $\vert \varpi^{\alpha_F(p-1)}\vert<\vert p\vert$, where $p\mid q$ is the characteristic of $\mathfrak{k}$. 
For instance, if $F = \Q_p$ and $p > 2$, then $\alpha_F = 1$.
Then we can define 
\begin{equation*}
    \exp\colon L(\alpha_F) \to K(\alpha_F) \text{ and }\log\colon K(\alpha_F)\to L(\alpha_F)
\end{equation*}
via the standard Taylor expansions. These maps will define homeomorphisms and satisfy the Baker--Campbell--Hausdorff formula 
\begin{equation}
	\exp(X)\exp(Y) = \exp\left( X+Y+\frac{1}{2}[X,Y]+O(X^2Y,XY^2)\right)\nonumber
\end{equation}
see \cite[Definition~3.28]{Ne-Hu_test}. One compares this to \eqref{eq:BCH}.
\end{rem}

\begin{rem}
If we equip $\mathfrak{g}'\cong F^4$ with the obvious measure, then we can use the approximate exponential function $e$ to pull this measure back to $K(1)\subseteq K_0$. More precisely, we get
\begin{equation*}
    \int_{K(1)}f(k) \, d\mu_{K_0}(k) = \zeta_q(1)\zeta_q(2)\int_{L(1)}f(1+X) \, dX.
\end{equation*}
To see that this is the correct normalisation we note that $\Vol(L(1))=q^{-4}$. On the other hand, since $K_0/K(1)=\GL(\mathfrak{k})$, we have
\begin{equation*}
    \Vol(K(1)) = [\sharp\GL(\mathfrak{k})]^{-1} = q^{-4}\zeta_q(1)\zeta_q(2).
\end{equation*}
\end{rem}

As in the archimedean case we identify $\mathfrak{g}'$ with its dual via the trace pairing $(X,Y)=\tr(X\cdot Y)$. 
More precisely, let $(\mathfrak{g}')^{\wedge}$ be the Pontryagin dual of $\mathfrak{g}'$. 
The identification $\Psi\colon \mathfrak{g}'\to (\mathfrak{g}')^{\wedge}$ is defined by $$\Psi(Y)X=\psi((Y,X))=\psi(\tr(Y\cdot X)).$$ Given a lattice $\Lambda\subseteq \mathfrak{g}'$ we associate the usual dual lattice and denote it by $\Lambda^{\wedge}$. Using $\Psi$ we can identity $\Lambda^{\wedge}$ with a lattice in $\mathfrak{g}'$. We set $\Lambda^{\ast} = \Psi^{-1}(\Lambda^{\wedge})$. Spelling out this identification concretely shows that
\begin{equation}
	\Lambda^{\ast} = \{ x\in \mathfrak{g}'\colon \tr(xy)\in \mathcal{O}\text{ for all }y\in \Lambda\} \label{eq:dual_lat}
\end{equation}
It is now easily checked that $L(m)^{\ast} =  L(-m)$.

\subsection{Background on non-archimedean representations}\label{bn-ar}

We will now recall some important results concerning the representation theory of $\GL(2)$ over $F$. Focusing on the theory of smooth admissible representations of $G$ with trivial central character.\footnote{Such representations factor through the canonical projection $G \to \PGL = Z\backslash G$ and all representations of $\bar{G}$ arise this way.} We will start by recalling the main classification result for such representations in Section~\ref{sec:classification}. Next, in Section~\ref{sec:new-vectors} we will introduce the theory of new-vectors. This is a local version of the classical Atkin-Lehner theory. Section~\ref{wkm} contains some background on Whittaker and Kirillov models. Finally, Section~\ref{mcur} contains a short discussion of (truncated) matrix coefficients. The results obtained here will later be essential for the definition of the global test function in Section~\ref{sec:test_fct} below.

\subsubsection{The classification}\label{sec:classification}

Let $\pi$ be an irreducible unitary smooth representation of $G(F)$ with trivial central character. 
If $\pi$ is infinite dimensional, it is isomorphic to one of the following:
\begin{enumerate}
    \item an irreducible principal series representation $\chi \boxplus \chi^{-1}$ for a unitary character $\chi\colon F^{\times}\to \C^{\times}$;
    \item a quadratic twist of a spherical irreducible complementary series representation $\omega\vert\cdot\vert^{\alpha}\boxplus \omega\vert\cdot\vert^{-\alpha}$, where $\alpha\in (-\frac{1}{2},\frac{1}{2})$ and $\omega\colon F^{\times}\to\C^{\times}$ is a quadratic character;
    \item a quadratic twist of the Steinberg representation: $\omega\cdot \textrm{St}$, where $\omega\colon F^{\times}\to\C^{\times}$ is a quadratic character (these representations are sometimes called special representations);
    \item a supercuspidal representation with trivial central character. These are parametrised by conjugacy classes of cuspidal types $(\mathfrak{A},J,\Lambda)$, where $\mathfrak{A}\subseteq \mathfrak{g}'$ is a hereditary order, $J$ is a subgroup of the centraliser of $\mathfrak{A}$ and $\Lambda$ is a smooth representation of $J$. Given such a cuspidal type we obtain a supercuspidal representation by compact induction: $\textrm{c-Ind}_J^{G}(\Lambda)$.
\end{enumerate}

This general classification can be found in \cite[Section~9.11]{bushnell}. 
Supercuspidal representations are then discussed in more detail in \cite[Section~15]{bushnell}. 
Alternatively, properties of parabolically induced representations of $G$ are discussed in \cite[Theorem~4.5.1, 4.5.2 and 4.5.3]{bump}.

\begin{rem}
Since we are assuming that $q$ is odd all supercuspidal representations are dihedral and can be constructed from the Weil representation. 
The construction is carried out in \cite[Theorem~4.8.6]{bump}, for example. 
Exhaustion then follows from \cite[Proposition~4.9.3]{bump} and the local Langlands correspondence.
\end{rem}

\begin{rem}
In \cite{Hu-sup} representations of $\PGL_n(F)$ with generic induction datum are considered. 
See \cite[Definition~3.12 and~3.15]{Hu-sup} for precise definitions. 
In our case, i.e. $n=2$ and $q$ odd, all supercuspidal representations have generic induction datum. 
This follows from the discussion in \cite[Section~19.2]{bushnell}.
Also, principal series $\chi\boxplus\chi^{-1}$ where $\chi\vert_{\mathcal{O}^{\times}}$ is not quadratic have a generic induction datum.
\end{rem}

Spherical, or unramified, representations are either of type (1) with unramified $\chi$ or of type (2) with unramified $\omega$.
These feature a unique spherical, i.e. $K_0$-fixed, vector and are particularly easy to understand. 
Their theory will play a role for the amplification method, see Section~\ref{sec:amp} below. 
For the moment, we postpone their discussion and study ramified representations in what follows.

\subsubsection{New-vector theory} \label{sec:new-vectors}

The \emph{level} of a representation $\pi$ of $G$ is defined by
\begin{equation}
    a(\pi) = \min\{ m\colon \pi^{K_H(m)}\neq \{ 0\}\}.\nonumber
\end{equation}
We call the integer $a(\pi)$ the exponent conductor of $\pi$. 
It is a consequence of classical Atkin-Lehner theory, see for example \cite{cass, ralf} in the local setting, that the space $\pi^{K_H(a(\pi))}$ is one dimensional. 
Up to a choice of normalisation, a non-zero element of the latter is referred to as the \emph{new-vector}. 
For the types of representations listed above the level and the new-vector can be determined quite explicitly.

For $m,m'\in \Z$ with $n=m+m'>0$ it will be convenient to consider the generalised subgroups
\begin{equation}
    K_H(m,{m'}) = \left\{\gamma= \left(\begin{matrix} a & \varpi^{m'}b\\ \varpi^mc& d\end{matrix}\right)\colon a,b,c,d\in \mathcal{O} \text{ and }\det(\gamma)\in \mathcal{O}^{\times}\right\}.\label{eq:defHecke}
\end{equation}
All these are conjugate to the usual Hecke congruence subgroups $K_H(n)$. Indeed, we have
\begin{equation}
    a(\varpi^{-m'})K_H(m,{m'})a(\varpi^{m'}) = K_H(m+m').\nonumber
\end{equation}
Following P.\ Nelson, see \cite{Nel-QUE-QP}, we call a non-zero vector $v\in \pi$ a \emph{generalised new-vector} if it is $K_H(m,{m'})$-invariant and $m+m'=a(\pi)$. 
All of these are translates of the usual new-vector. 
Indeed, if $\vnew \in \pi^{K_H({a(\pi)})}$ is the standard new-vector, then $\pi(a(\varpi^m))\vnew$ is invariant under~$K_H({a(\pi)-m},m)$ and thus a generalised new-vector. 
It is clear that all generalised new-vectors are of this form. We finally caution that, if $n=2m$, then we almost always work with the $K_H(m,m)$-invariant new-vector. For later reference we introduce the following convention.

\begin{conv}\label{conv:newvec}
If $n=2m$, then we write $v^{\circ}$ for the $K_H(m,m)$-invariant new-vector. This is often referred to as simply the new-vector.	
\end{conv}

We do give a detailed discussion of the general new-vector theory. Instead, we will discuss the special case where $\pi = \chi\boxplus \chi^{-1}$ for $\chi\colon F^{\times}\to S^1$ with $a(\chi)=m$. In this case the generalized new-vector $v^{\circ}$ (i.e. the unique $K_H(m,m)$-invariant vector) can be constructed rather explicitly, and we hope that this provides some useful intuition. Let us write $\chi_B\colon B\to S^1$ for the character given by
\begin{equation}
	\chi_B\left(\left(\begin{matrix} t_1 & x\\ 0 &t_2\end{matrix}\right)\right) = \chi(t_1/t_2). \label{eq:chi_BBBB}
\end{equation}
Thus, the action of $\pi$ is given quite explicitly on functions $f\colon K_0\to \C$ satisfying 
\begin{equation*}
	f(bk) = \chi_B(b)f(k) \text{ for }k\in K_0 \text{ and }b\in B(\mathcal{O}).\nonumber
\end{equation*}
This is called the compact model of $\pi$, and we will denote it by $I(\chi_B)$. 

Note that the generalised new-vector $v^\circ$ is per definition $K_H(m,m)=D(\mathcal{O})K(m)$ invariant. Thus, in order to understand $v^\circ$ in the compact model, we use the double coset decomposition
\begin{equation}
	K_0=\bigsqcup_{i=1}^{m}B(\mathcal{O})\left(\begin{matrix} 1 & 0 \\ \varpi^i &1\end{matrix}\right)K_H(m,m) \sqcup\bigsqcup_{j=0}^m B(\mathcal{O})\left(\begin{matrix} 0 & -1 \\ 1 &\varpi^j\end{matrix}\right)K_H(m,m). \label{some_double_coset}
\end{equation}
This is analogous to the corresponding decomposition for $K_H(m)$ given in  \cite[Lemma~2.1.1]{ralf}. For completeness, we will sketch the proof. Let $k=\begin{psmallmatrix} a & b \\ c & d \end{psmallmatrix}\in K$ and define $(i,j) = (v(c), v(d))$. If $0<i$, then $j=0$ and we can simply write
\begin{equation}
	k =\left(\begin{matrix} \frac{\det(k)\varpi^i}{c} & b\\ 0 & d\end{matrix}\right)\left(\begin{matrix} 1 & 0 \\ \varpi^i &1\end{matrix}\right)\left(\begin{matrix}\frac{c}{\varpi^i d} & 0 \\ 0 & 1 \end{matrix}\right).\nonumber
\end{equation}
This accounts for the first part of the decomposition. From now on let us assume that $i=0$. In this case we can simply write
\begin{equation}
	k = \left(\begin{matrix} \frac{\det(k)\varpi^j}{d} & a\\ 0 & c\end{matrix}\right)\left(\begin{matrix} 0 & -1 \\ 1 & \varpi^j \end{matrix}\right)\left(\begin{matrix}1 & 0 \\ 0 & \frac{d}{\varpi^jc} \end{matrix}\right)\nonumber
\end{equation}
accounting for the second part of the decomposition. Disjointness is easily checked by computing the valuations of the entries in the bottom row of each cell.

With this decomposition at hand we can compute an exact expression for the generalised new-vector in the induced model:
\begin{lemmy}\label{lm:newvec_princ_ind}
	In the compact model, the generalised new-vector $v^\circ$ is given by
	\begin{equation*}
		f_{\circ}(k)= \begin{cases}
			\chi(t_1/t_2) &\text{ if }k\in \left(\begin{matrix} t_1 & \mathcal{O} \\ 0 & t_2\end{matrix}\right)\left(\begin{matrix} 0 & -1 \\1 & 1 \end{matrix}\right)K_H(m,m), \\
			0 &\text{ else.}
		\end{cases}
	\end{equation*}
\end{lemmy}
\begin{proof}
	We will follow the argument from \cite[p.9]{ralf} and deduce this from the double coset decomposition \eqref{some_double_coset}. 
	
	Let $\gamma_i = \left(\begin{matrix} 1 & 0 \\ \varpi^i & 1 \end{matrix}\right)$ and $\tilde{\gamma}_j =  \left(\begin{matrix} 0 &-1 \\ 1 & \varpi^j \end{matrix}\right).$ Recall that we are trying to construct a $K_H(m,m)$-invariant element in the compact model for $\pi= \chi \boxplus \chi^{-1}$. In view of \eqref{some_double_coset}, we can attempt to define functions $f_{\gamma}\colon K_0 \to \C$ via
	\begin{equation}
		f_{\gamma}(k) = \begin{cases}
			\chi_B(b) &\text{ if }k=b\gamma k_1 \text{ for }b\in \mathcal{O} \text{ and }k_1\in K_H(m,m),\\
			0&\text{ else.}
		\end{cases} \nonumber
	\end{equation}
	Here, $\gamma = \gamma_i$ for $i=1,\ldots,m$ or $\gamma=\tilde{\gamma}_j$ for $j=0,\ldots, m$. For this to be well-defined, we need to check that $\chi_B$ is trivial on $B(\mathcal{O}) \cap \gamma K_H(m,m) \gamma^{-1}$. We do so by considering different cases:
	\begin{itemize}
		\item If $\gamma = \gamma_i$ for $1\leq i<m$, then we can take $k = \left(\begin{matrix} 1 & 0 \\ c\varpi^m & 1+\varpi^{m-i}c\end{matrix}\right)\in K_H(m,m)$, for $c\in \mathcal{O}^{\times}$. We compute that
		\begin{equation}
			\chi_B(\gamma_i k\gamma_i^{-1}) = \chi_B\left(\left(\begin{matrix}  1 & 0 \\ 0 & 1+\varpi^{m-i}c\end{matrix}\right)\right) = \chi^{-1}(1+c\varpi^{m-i}).\nonumber
		\end{equation}
		Since the conductor of $\chi$ is $m$ we can find $c$ with $\chi_B((\gamma_i k\gamma_i^{-1})\neq 1$. We conclude that $f_{\gamma_i}$ is not well-defined for $i=1,\ldots,m-1$.
		\item If $\gamma = \gamma_m$, then we observe that $\gamma_m K_H(m,m)\gamma_m^{-1} = K_H(m,m)$. However, $\chi_B$ is a non-trivial character of $B(\mathcal{O})$. Thus $f_{\gamma_m}$ is not well-defined.
		\item If $\gamma = \tilde{\gamma}_j$ for $1\leq j \leq m$, then we can take $k= \left(\begin{matrix} 1 & b\varpi^m\\ 0 & 1+b\varpi^{m-j}\end{matrix}\right)\in K_H(m,m)$. For this choice we check that
		\begin{equation}
			\chi_B(\tilde{\gamma}_j k\tilde{\gamma}_j^{-1}) = \chi_B\left(\left(\begin{matrix}  1 & 0 \\ 0 & 1+\varpi^{m-j}b\end{matrix}\right)\right) = \chi^{-1}(1+b\varpi^{m-j}).\nonumber
		\end{equation}
		Again we can arrange $b\in \mathcal{O}^{\times}$ such that $\chi_B(\tilde{\gamma}_j k\tilde{\gamma}_j^{-1})  \neq 1$. This shows that $f_{\tilde{\gamma}_j}$ is not well-defined.
		\item Finally, if $\gamma = \tilde{\gamma}_0$, then	we compute that
		\begin{equation}
			B(\mathcal{O}) \cap \tilde{\gamma}_0K_H(m,m)\tilde{\gamma}_0^{-1} = \left\{ \left(\begin{matrix} a & x \\ 0 & d \end{matrix}\right)\colon x\in \mathfrak{p}^m,\, a,d\in \mathcal{O}^{\times}\text{ and }a\equiv d \text{ mod }\mathfrak{p}^m\right\}.\nonumber
		\end{equation}
		Because $\chi$ has conductor $m$, we see that $\chi_B$ is trivial on this subgroup. This ensures that $f_{\tilde{\gamma}_0}$ is well-defined.
	\end{itemize}
	Thus we have seen that $f_{\circ} = f_{\tilde{\gamma}_0}$ is (up to multiplication by constants) the only possible element in the compact model for $\pi=\chi\boxplus \chi^{-1}$, which is right $K_H(m,m)$ invariant. This completes the proof. 
\end{proof}

Let us return to our general discussion where $\pi$ is an arbitrary irreducible admissible smooth representation of $G(F)$. We call $\pi$ \emph{twist minimal} if $a(\pi)\leq a(\chi\cdot \pi)$ for all characters $\chi\colon F^{\times}\to \C^{\times}.$ 
Note that the principal series representation $\pi=\chi\boxplus\chi^{-1}$ is not twist minimal (as a representation of $G$), since $a(\pi)=2a(\chi)$ and $$a(\chi^{-1}\cdot \pi)=a(1\boxplus \chi^{-2})=a(\chi^{-2}).$$ 
In general, we have $2a(\chi)>a(\chi^{-2}).$ 
On the other hand, if $\pi$ is supercuspidal and has trivial central character, then $\pi$ is twist minimal. 
This is a special case of \cite[Lemma~2.1]{Hu-Nel-Sa}. 

\begin{rem}\label{rem:conductor_class}
For simplicity, we will later restrict to representations $\pi$ with exponent conductor $a(\pi)=4n$. By going through the list of Subsection~\ref{sec:classification} and checking the conductors, which have been computed in \cite{ralf} for example, we see that there are only two possibilities:
\begin{enumerate}
	\item $\pi = \chi\boxplus \chi^{-1}$ with $a(\chi)=2n$;
	\item $\pi$ is supercuspidal with trivial central character and conductor $a(\pi)=4n$.
\end{enumerate} 
\end{rem}

\subsubsection{Whittaker and Kirillov models}\label{wkm}

A linear functional $\Lambda\colon \pi\to \C$ is called a $\boldsymbol{\psi}^{(a)}$-Whittaker functional if it satisfies 
\begin{equation}
	\Lambda(\pi(n)v)=\boldsymbol{\psi}^{(a)}(n)\cdot\Lambda(v) \text{ for all }n\in N(F)\text{ and }v\in \pi.\nonumber
\end{equation}
A representation is called generic if it features a non-trivial $\boldsymbol{\psi}^{(a)}$-Whittaker functional. It is an easy fact that if $\pi$ has a non-trivial $\boldsymbol{\psi}^{(a)}$-Whittaker functional, then it has a non-trivial $\boldsymbol{\psi}^{(a')}$-Whittaker functional for all $a'\in F^{\times}$. In particular, the notion of genericity does not depend on the particular choice of $a\in F^{\times}$. We will thus specialize to $a=1$ for simplicity and drop the prefix $\boldsymbol{\psi}^{(1)}$ from the term Whittaker functional.

It turns out that all the representations listed above are generic. Indeed, all smooth admissible irreducible infinite dimensional representations of $G$ are generic (see \cite[Theorem~4.4.3]{bump}). Furthermore, a Whittaker functional is always unique up to scaling (see \cite[Theorem~4.4.1]{bump}) and gives rise to the unique Whittaker model $\mathcal{W}(\pi,\psi)$. More precisely, if $\Lambda_{\pi}$ denotes the Whittaker functional, then 
\begin{equation}
    \mathcal{W}(\pi,\psi) = \{ W_v(g)=\Lambda_{\pi}(\pi(g)v)\colon v\in \pi\}.\nonumber 
\end{equation}
Note that by construction $W_v(n(x)g)=\psi(x)W_v(g)$. In this model, $G$ acts by right translation. 

Finally, note that the Kirillov model is given by
\begin{equation}
    \mathcal{K}(\pi,\psi)=\{ f(y)=W(a(y)) \colon W\in \mathcal{W}(\pi,\psi)\}.\nonumber
\end{equation}
It is a well-known fact that $\mathcal{C}_c^{\infty}(F^{\times}) \subseteq \mathcal{K}(\pi,\psi)$ for all $\pi$. Furthermore, the action of $DN$ is quite explicitly given by
\begin{equation}
    [g.f](y) = \psi \left(\frac{by}{d}\right) \cdot f\left(\frac{ay}{d}\right) \text{ for }g=\left(\begin{matrix} a & b\\ 0 & d\end{matrix}\right)\in DN \text{ and }f\in \mathcal{K}(\pi,\psi).\label{action_kirillow}
\end{equation}
If $\pi$ is unitary, then the $G$-invariant inner product in the Kirillov model is given by
\begin{equation*}
    \langle f_1,f_2\rangle = \zeta_q(1)\int_{F^{\times}}f_1(y)\overline{f_2}(y) \, \frac{dy}{\vert y\vert},
\end{equation*}
for $f_1,f_2\in \mathcal{K}(\pi,\psi)$. We have the following basis of $\mathcal{C}_c^{\infty}(F^{\times})$. Given a character $\chi\colon \mathcal{O}^{\times}\to S^1$ and $m\in \Z$ we set
\begin{equation*}
    \xi_{\chi}^{(m)}(x) =\begin{cases}
        \chi(x') &\text{ if }x=\varpi^{-m}x' \text{ for }x'\in \mathcal{O}^{\times},\\
        0 &\text{ else.}
    \end{cases} 
\end{equation*}
We see directly that
\begin{equation}
    \langle \xi_{\chi}^{(m)},\xi_{\chi'}^{(m')}\rangle = \delta_{(\chi,m)=(\chi',m')}.\nonumber
\end{equation}
Furthermore, $a(\varpi^s).\xi_{\chi}^{(m)} = \xi_{\chi}^{(m+s)}$. These are convenient elements to use for computations in the Kirillov model. To illustrate this, we note that the new-vector in the Whittaker model (i.e. as a function on $G(F)$) can be very complicated. See for example \cite{Ass19}. On the other hand, it has a very nice description in the Kirillov model:

\begin{lemmy}\label{lm:kirnew}
Let $\pi$ be a generic irreducible smooth admissible representation with trivial central character and $a(\pi) = 2m$, for $m\in \N$. Then the $L^2$-normalized new vector $v_{\textrm{new}}$ is given by $\xi_1^{(0)}=\vvmathbb{1}_{\mathcal{O}^{\times}}$ in the Kirillov model. Furthermore, $v^{\circ}$ is given by $\xi_1^{(a(\chi))}$ in the Kirillov model.
\end{lemmy}  
\begin{proof}
The formula for the new-vector in the Kirillov model can be found in \cite[Section~2]{ralf}. Note that our assumption on the level implies in particular that $a(\pi)\geq 2$. This excludes unramified twists of $\St$. The corresponding formula for $v^{\circ} = \pi(a(\varpi^m))v_{\textrm{new}}$ in the Kirillov model follows directly from \eqref{action_kirillow}.
\end{proof}

\subsubsection{Matrix coefficients of unitary representations}\label{mcur}

In this section, let $\pi$ be a smooth unitary irreducible representation of $G$ with trivial central character. 
Further, let $v\in \pi$ be a $K_0$-finite vector, and we associate the (truncated) matrix coefficient
\begin{equation*}
    \varphi_v(g) =  \frac{\langle v,\pi(g)v\rangle}{\Vert v\Vert^2}\cdot \vvmathbb{1}_{ZK_{0}}(g).
\end{equation*}
Here we recall some important properties of the test functions $\varphi_v(g)$ and the corresponding operators $\pi(\varphi_v)$. The material is standard and is used for example in \cite{marsh_local}, as well as in \cite{Sa}. 

First, recall the orthogonality relations 
\begin{equation}
	\int_{K_0} \langle \tau(k)v_1, v_2\rangle \overline{\langle \tau(k)v_3,v_4  \rangle}dk = \frac{\langle v_1,v_3\rangle\overline{\langle v_2,v_4\rangle}}{\dim(V)},\label{eq:orth}
\end{equation}
where $(\tau,V)$ is an irreducible representation of $K_0$ and $v_1,v_2,v_3,v_4\in V$. These follow from Schur's orthogonality relations for finite fields as stated in \cite[Lemma~4]{marsh_local} for example.

\begin{lemmy}\label{mat_coeff_indep}
Let $(\tau,V)$ denote the $K_0$-representation generated by $v$. Then we have
\begin{equation*}
    \varphi_v = \dim(V)\cdot \varphi_v\ast \varphi_v.\nonumber
\end{equation*}
\end{lemmy}
\begin{proof}
Recall that the convolution is defined by
\begin{align*}
    [\varphi_v\ast\varphi_v](h) &= \int_{Z\backslash G} \varphi_v(g^{-1})\varphi_v(gh)dg \\
    &= \vvmathbb{1}_{ZK_0}(h) \int_{K_0} \varphi_v(k^{-1})\varphi_v(kh)dk.
\end{align*}
Writing the $K_0$-integral out yields
\begin{equation*}
    \int_{K_0} \varphi_v(k^{-1})\varphi_v(kh)dk = \frac{1}{\Vert v\Vert^4}\int_{K_0}\langle \pi(k)v,v\rangle\cdot \langle v,\pi(kh)v\rangle dk.
\end{equation*}
For $h\in K_0$ one uses \eqref{eq:orth} with $v=v_1=v_2=v_4$ and $v_3=\pi(h)v$ to get
\begin{equation*}
    \int_{K_0} \varphi_v(k^{-1})\varphi_v(kh)dk = \frac{1}{\dim(V)\Vert v\Vert^2}\langle v,\pi(h)v\rangle = \frac{1}{\dim(V)}\varphi_v(h).
\end{equation*}
This completes the proof.
\end{proof}

\begin{cor}
The operator $\pi(\varphi_v)$ is self-adjoint and non-negative.
\end{cor}
\begin{proof}
Note that $\varphi_v(g^{-1})=\overline{\varphi_v(g)}$, so that $\pi(\varphi_v)$ is self-adjoint. Non-negativity follows directly from the previous lemma.
\end{proof}

Finally, we observe that $\pi(\varphi_v)$ acts as a projection.
\begin{lemmy}\label{p-adic_proj}
Let $w\in \pi$ be any $K_0$-finite vector. Then 
\begin{equation*}
    \pi(\varphi_v) w = \frac{1}{\dim(V)}\frac{\langle w,v\rangle}{\Vert v\Vert^2} \cdot v,
\end{equation*}
where $(\tau,V)$ is the $K_0$-representation generated by $v$.
\end{lemmy}
\begin{proof}
First note that $(G, K_0)$ is a strong Gelfand pair.\footnote{We say that $(G,H)$ is a strong Gelfand pair if $H\subseteq G$ and we have $\dim\Hom_H(\pi\vert_H,\sigma)\leq 1$ for all irreducible representations $\pi$ of $G$ and $\sigma$ of $H$.} Given any $w'\in \pi$ we can compute
\begin{align}
    \langle\pi(\varphi_v) w,w'\rangle &= \frac{1}{\Vert v\Vert^2}\int_{K_0}\langle v,\pi(k)v\rangle\cdot \langle \pi(k)w,w'\rangle dk \nonumber \\
    &= \frac{\langle w,v\rangle\cdot \langle w',v\rangle}{\Vert v\Vert^2\cdot \dim(V)}.\nonumber
\end{align}
In the last step we have used the orthogonality relation from \eqref{eq:orth}. If $w'$ is orthogonal to $v$, then we find that $\langle\pi(\varphi_v) w,w'\rangle=0$. This allows us to conclude that $\pi(\varphi_v) w = C\cdot v$ for some $C\in \C$. To determine this constant we observe that 
\begin{equation*}
	C= \frac{\langle \pi(\varphi_v) w,v\rangle}{\Vert v\Vert^2} = \frac{\langle w,v\rangle}{\Vert v\Vert^2\dim(V)}. 
\end{equation*}
This completes the proof.
\end{proof}

\begin{rem}
If $\pi$ is supercuspidal, then the matrix coefficients $\langle v,\pi(g)v\rangle$ are compactly supported (modulo centre) and they satisfy typical orthogonality relations. See for example \cite[Theorem~~1.1 and~1.3]{car}. Thus, in principle we do not need to truncate them to have support in $ZK_0$.
\end{rem}

\subsection{Microlocalisation: definition and consequences}\label{mdc}

\begin{defn}\label{def:loc_vec}
Let $\pi$ be an irreducible unitary smooth representation of $G(F)$, let $C\subseteq G(F)$ be a subgroup and let $\tau\in \widehat{C}$ be a character. 
We say that a non-zero vector $v \in \pi$ is \emph{localised at} $(\tau,C)$ if 
\begin{equation}
    \pi(c)v=\tau(c)\cdot v \nonumber
\end{equation}
for all $c\in C$.
\end{defn}

\begin{rem}
We have stated this definition at the group level for convenience. 
Many examples arise from the level of the Lie algebra, as follows. 
Given $C=\exp(L)$ for a suitable lattice $L\subseteq \mathfrak{g}'$ one can construct characters of $C$ by setting 
\begin{equation}
    \tau(\exp(X)) = \Psi(Y)X,
\end{equation}
for suitable $Y\in L^{\wedge}$. 
Using this set-up one can obtain an alternative definition of microlocalised vectors modelled on Definition~\ref{def:microlocal}, where the lattice $L$ plays the role of the localisation parameter $R$.
\end{rem}

\begin{rem}
    Note that if $v$ is microlocalised at $(\tau,C)$ then $v'=\pi(g)v$ is microlocalised at $(\tau^g,gCg^{-1})$. 
\end{rem}

This simple remark allows for an important observation concerning the support of matrix coefficients. Indeed, we deduce that, if there is $h\in C\cap gCg^{-1}$ with $\tau(h)\neq \tau^g(h)$, then $\langle \pi(g)v,v\rangle =0$. See Lemma~\ref{microloc_orth} and its proof below for a related argument. In particular, we obtain that the support of the diagonal matrix coefficient $g\mapsto \langle \pi(g)v,v\rangle$ associated to a vector $v$ microlocalised at $(\tau,C)$ is contained in the set
\begin{equation}
	\{ g\in G(F)\colon \tau\vert_{C\cap gCg^{-1}}=\tau^g\vert_{C\cap gCg^{-1}}\}.\nonumber
\end{equation} 
The group elements $g$ in this set are said to intertwine $\tau$. See \cite[Definition~2.2]{Hu-sup} for a general definition. In practice, we will work with tuples $(\tau,C)$, where we have good control on the set of intertwiners. This will directly translate to important properties of the corresponding microlocalised vectors. See Lemma~\ref{lm:mat_coeff_super} and Lemma~\ref{lm:mat_coeff_princ} below for examples for the support of matrix coefficients.

\subsubsection{Limits to localisation}\label{ltl}

As in Section \ref{sec:limits}, we investigate how large the subgroup $C$ may be taken in terms of the level of $\pi$.
To do this, we assume that $C$ is an open group, which is always the case in the rest of our discussion, and therefore $K(n) \subset C$ for some $n \geq 0$.
We show below that the existence of a vector localised at $C$ implies that $\pi$ has level at most $4n$.

Since the level is connected to subgroups fixing non-zero vectors, we observe that commutators $c_1 c_2 c_1^{-1} c_2^{-1}$, which we denote by $[c_1, c_2]$, $c_1, c_2 \in C$ act trivially on a microlocalised vector $v$.

\begin{lemmy} \label{lm:commutator} 
    For $p >2$ and $n > 0$ (and $n > 1$ if $p = 2$), the commutator subgroup $[K(n), K(n)]$ satisfies that
    \begin{displaymath}
        [K(n), K(n)] \cdot (Z \cap K(2n)) = K(2n).
    \end{displaymath}
\end{lemmy}
\begin{proof}
    The inclusion $[K(n), K(n)] \subset K(2n)$ is easy to see. 
    If $1 + \varpi^n X$ is an element of $K(n)$ with $X \in \Mat_2(\Oc)$, then its inverse is of the form $(1 - \varpi^n X + \varpi^{2n}X')$ for some $X' \in \Mat_2(\Oc)$, by the formula for the geometric series.
    Modulo $\pf^{2n}$, the commutator of two such elements has the form 
    \begin{displaymath}
        (1+\varpi^nX)(1 + \varpi^n Y)(1 - \varpi^nX)(1 - \varpi^n Y) \equiv 1.
    \end{displaymath}

    For the converse, we first observe that $[K(n), K(n)]$ contains certain elementary matrices.
    More precisely, looking at the commutators
    \begin{displaymath}
        \left[ 
        \begin{pmatrix}
            1 & y \\
            & 1
        \end{pmatrix},
        \begin{pmatrix}
            x & \\ 
            & 1
        \end{pmatrix}
        \right]
        \quad \text{and} \quad
                \left[ 
        \begin{pmatrix}
            1 & \\
            y & 1
        \end{pmatrix},
        \begin{pmatrix}
            x & \\ 
            & 1
        \end{pmatrix}
        \right]
    \end{displaymath}
    shows that the subgroups
    \begin{displaymath}
        U^+(2n) := \left\{ 
            \begin{pmatrix}
                1 & \pf^{2n}\\
                & 1
            \end{pmatrix}
        \right\}
        \quad \text{and} \quad
        U^-(2n) := \left\{ 
            \begin{pmatrix}
                1 & \\
                \pf^{2n} & 1
            \end{pmatrix}
        \right\}
    \end{displaymath}
    are contained in~$[K(n), K(n)]$.
    We further compute that
    \begin{equation} \label{eq:commutator-upper-lower}
        \left[ 
        \begin{pmatrix}
            1 & x \\
            & 1
        \end{pmatrix},
        \begin{pmatrix}
            1 & \\ 
            y & 1
        \end{pmatrix}
        \right] =
        \begin{pmatrix}
            1 + xy + x^2 y^2 & -x^2 y \\
            xy^2 & 1 - xy
        \end{pmatrix}.
    \end{equation}

    To proceed, we note the general computation
    \begin{equation} \label{eq:iwahori-decomp-kinda}
        \begin{pmatrix}
            1 & -b(1+d)^{-1} \\
            & 1
        \end{pmatrix}
        \begin{pmatrix}
            1 + a & b \\
            c & 1 + d
        \end{pmatrix}
        \begin{pmatrix}
            1 &  \\
            -c(1+d)^{-1}& 1
        \end{pmatrix}
        = \begin{pmatrix}
            1+a' & \\
            & 1+d
        \end{pmatrix},
    \end{equation}
    where $a' = a - bc(1+d)^{-1}$.
    This implies, first, that the group $K(m)$ is generated by the subgroups $U^+(m)$, $U^-(m)$, and $D \cap K(m)$, for any $m > 0$.
    Second, coupled with the computation \eqref{eq:commutator-upper-lower} and noting that the matrix obtained there has determinant $1$, it also implies that the commutator $[K(n), K(n)]$ additionally contains the elements of the form $\diag((1+z)^{-1}, 1+z)$ for any $z \in \pf^{2n}$.
    We deduce that $(Z \cap K(2n)) [K(n), K(n)]$ contains all elements of the form $\diag(1, (1+z)^2)$ with $z \in \pf^{2n}$.
    Under our conditions, any element in $1 + \pf^{2n}$ has a square root and this shows that $(Z \cap K(2n)) [K(n), K(n)]$ contains $D \cap K(2n)$.
\end{proof}

Lemma \ref{lm:commutator} shows that the existence of a non-zero vector in $\pi$ localised at a subgroup $C$ containing $K(n)$ implies the existence of a non-zero vector fixed by $K(2n)$.
Assuming that the central character is trivial, such a vector is fixed by $K_H(2n, 2n)$.
As in the paragraph containing \eqref{eq:defHecke}, this implies that $\pi$ has level at most $2n + 2n = 4n$.

\subsubsection{Orthogonality}\label{ooo}

\begin{lemmy} \label{microloc_orth}
Let $\pi$ be an irreducible unitary smooth representation of $G(F)$ with trivial central character. Further, for $i=1,2$, let $v_i\in \pi$ be localised at $(\tau_i,C_i)$. If $\tau_1\vert_{C_1\cap C_2} \neq \tau_2\vert_{C_1\cap C_2}$, then
\begin{equation}
    \langle v_1,v_2\rangle_{\pi} =0.\nonumber
\end{equation}
\end{lemmy}
This is a non-archimedean analogue of the orthogonality phenomenon explained in Section~\ref{sec:orth}. 
Note that here we obtain actual orthogonality instead of an approximation.
\begin{proof}
The proof is very simple. Take $c\in C_1\cap C_2$ such that $\tau_1(c)\neq \tau_2(c).$ We then compute
\begin{equation}
    \tau_1(c)\langle v_1,v_2\rangle_{\pi} = \langle \pi(c)v_1,v_2\rangle  = \langle v_1,\pi(c^{-1})v_2\rangle = \tau_2(c)\langle v_1,v_2\rangle_{\pi}.\nonumber
\end{equation}
This implies $\langle v_1,v_2\rangle_{\pi}=0$ as desired.
\end{proof}

Suppose that $\pi$ is supercuspidal, suppose that $C$ is compact modulo centre, and let $\tau\colon C\to S^1$ be a character of $C$. Then the operator
\begin{equation}
    \mathcal{P}_{(\tau,C)}w = \frac{1}{\Vol((C\cap Z)\backslash C)}\int_{(C\cap Z)\backslash C}\tau(c)^{-1} \cdot \pi(c)w\, d\mu(c) \nonumber
\end{equation}
defines a projector onto the space of $(\tau,C)$-microlocalised vectors. On the other hand, we can fix an orthonormal basis $\mathcal{B}$ of this space and define the function 
\begin{equation}
    \varphi_{(\tau,C)}(g) = \sum_{v\in \mathcal{B}}\langle v,\pi(g)v\rangle. \nonumber
\end{equation}
Note that this function is compactly supported, because we are currently assuming that $\pi$ is supercuspidal.\footnote{If $\pi$ is not supercuspidal, then one encounters certain issues due to the lack of compactness of matrix coefficients. We will encounter these below.} Recall from \cite[Theorem~1.1 and~1.3]{car} that for matrix coefficients of supercuspidal representations one has a version of Schur orthogonality relations available. This allows us to write
\begin{equation}
    \mathcal{P}_{(\tau,C)}w = d(\pi)\cdot\int_{Z\backslash G(F)} \varphi_{(\tau,C)}(g)\cdot \pi(g)w \, d\mu(g), \nonumber
\end{equation}
where $d(\pi)$ is the formal degree of $\pi$.

In particular, if $\sharp \mathcal{B}=1$, i.e. the space of $(\tau,C)$-microlocalised vectors is one dimensional, then comparing these two expressions for the projector suggests a very simple formula for matrix coefficients of microlocalised vectors. Indeed, one might expect that
\begin{equation}
	\langle v,\pi(g)v\rangle \overset{?}{=} \begin{cases}
		\frac{1}{d(\pi)\cdot \Vol((C\cap Z)\backslash C)} \tau(c)^{-1} &\text{ if }g=z\cdot c \text{ for }z\in Z, c\in C,\\
		0&\text{ else.}
	\end{cases}\label{eq:matrix_coeff_exp.}
\end{equation}
This feature will be nicely illustrated in our examples below.

\begin{rem}\label{rem:about_f1f2}
Let us stress that we have formulated \eqref{eq:matrix_coeff_exp.} as a heuristic to motivate latter formulae for matrix coefficients of microlocalised vectors. Indeed, let $f_1(g) = d(\pi)\cdot \varphi_{\tau,C}(g)$ and let $f_2(g)$ be the right hand side of \eqref{eq:matrix_coeff_exp.}. Then we have shown that $\pi(f_1)=\pi(f_2)$ as operators on $\pi$. This means when working within $\pi$ we can use $f_1$ and $f_2$ interchangeably but it does \textbf{not} imply that $f_1 = f_2$ (pointwise) as functions in $\mathcal{C}_c^{\infty}(Z\backslash G(F))$. However, for very well chosen tuples $(\tau,C)$, it is the case that $\pi^{(\tau,C)}\neq \{0\}$ determines the isomorphism class of $\pi$. If this is the case, then our argument actually shows that $\sigma(f_1) = \sigma(f_2)$ for all $\sigma$ in the unitary dual of $Z\backslash G(F)$. This is now strong enough to deduce that $f_1(g)=f_2(g)$ holds for all $g\in G(F)$. See \cite[§18.8.1]{Dix} for example.
\end{rem}

\subsection{Examples of microlocalised vectors}\label{emv}

To illustrate the abstract discussion from Section~\ref{mdc} above, we will now give examples for microlocalised vectors. This will be done in a pragmatic manner, keeping our application to the sup-norm problem in mind. Indeed, we will consider the two classes of representations appearing in Remark~\ref{rem:conductor_class}. For each case of $\pi$ that can occur we will construct localised vectors with particularly nice properties. Our plan is roughly as follows:
\begin{itemize}
	\item We start by constructing suitable tuples $(\tau,C)$. If $\pi$ is supercuspidal, this will be a so-called \textit{cuspidal type associated to $\pi$} (see Definition~\ref{def:cusp_type}). For principal series we call this a \textit{general $\chi$-type} (see Definition~\ref{def:chi-type}).
	\item Next we define the corresponding localised vectors. For supercuspidal representations these will be called \textit{minimal vectors} (see Definition~\ref{def:min_vex}). If $\pi$ is principal series, we encounter \textit{microlocal lift vectors} (see Definition~\ref{def:mic_loc_loft}).
	\item For these special localised vectors we can confirm that their matrix coefficients have the form predicted in \eqref{eq:matrix_coeff_exp.}. See Lemma~\ref{lm:mat_coeff_super} if $\pi$ is supercuspidal and Lemma~\ref{lm:mat_coeff_princ} if $\pi$ is principal series.
	\item To provide some intuition for these localised vectors, we also compute expressions for them in the Kirillov model. This allows us to relate them to the generalized new-vector $v^{\circ}$ from Convention~\ref{conv:newvec}. This is Corollary~\ref{cor:decomp_newform_super} for $\pi$ supercuspidal and Corollary~\ref{cor:newform_exp_princ} for principal series.
\end{itemize}

Both examples that are considered here have already appeared in the literature. In the case when $\pi$ is supercuspidal we will use this to our advantage and heavily rely on existing exposition. On the other hand, if $\pi$ is a principal series, we will use the rather accessible induced model and prove many of the results using explicit computations. We hope that this makes the material more accessible and nicely complements the existing literature.

\subsubsection{Minimal vectors}

We start with minimal vectors. These appear in supercuspidal representations $\pi$ and have been studied in \cite{Hu-Nel-Sa,  Hu-sup, Hu-Sh-Yi}, for example. Here we record the most important results following \cite{Ne-Hu_test}. We make the following simplifying assumption.

\begin{ass}\label{ass:super}
Throughout this subsubsection we assume that $\pi$ is a supercuspidal representation of $G(F)$ with trivial central character and exponent conductor $a(\pi)=4r$.
\end{ass}

We first need to introduce some more notation. Let $E/F$ be an unramified quadratic extension of $F$. The valuation ring of $E$ is denoted by $\mathcal{O}_E$, its maximal ideal is $\mathfrak{p}_E$ and the residue field is $\mathfrak{k}_E=\mathcal{O}_E/\mathfrak{p}_E$. We can choose $\varpi$ as uniformiser for $E$. Note that $\# \mathfrak{k}_E=q_E=q^2$. The (normalised) valuation on $E$ is denoted by $v_E$, and we obtain the corresponding absolute value $\vert \cdot \vert_E$ as usual. We define an additive character on $E$ by setting
\begin{equation*}
    \psi_E(x)=\psi(\Tr_{E/F}(x)).
\end{equation*}
Since $E/F$ is unramified, this character has conductor $\mathcal{O}_E$.

If $E=F(\sqrt{D})$, then we define the embedding $E\to \mathfrak{g}'$ by
\begin{equation*}
    a+b\sqrt{D}\mapsto\left(\begin{matrix} a & b \\ bD & a \end{matrix}\right).
\end{equation*}
We now fix a character $\theta\colon E^{\times} \to S^1$ such that
\begin{enumerate}
    \item $\theta$ has even exponent conductor $a(\theta)=2r$; and 
    \item $\theta\vert_{F^{\times}} = 1$.\footnote{This assumption implies that $\theta$ is a minimal character. We refer to \cite[Section~2]{Ne-Hu_test} for general definitions and more information.}
\end{enumerate}
We define an element $b_{\theta}\in \mathfrak{p}_E^{-2r}$ by requiring that
\begin{equation}
	\theta(1+u) = \psi_E(b_{\theta}u) \text{ for all }u\in \mathfrak{p}_E^r. \nonumber
\end{equation}
The existence of such a $b_{\theta}$ is easily verifies, since $u\mapsto \theta(1+u)$ is a character of $\mathfrak{p}_E^{r}/\mathfrak{p}_E^{2r}$. See \cite[Lemma~2.1]{Ne-Hu_test} for a reference. Using this element we define the character $\theta_{b_{\theta}} \colon K(r) \to S^1$ by setting
\begin{equation}
	\theta_{b_{\theta}}(k) = \psi(\tr (b_{\theta}\cdot (k-1))).\label{some_def}
\end{equation}
Here we view $b_{\theta}$ as an element in $\mathfrak{g}'$ using the embedding described above.

We associate the groups\footnote{Usually there are two other groups $H$ and $H^1$, as in Definition 3.12 of \cite{Hu-Nel}. However, since we are assuming that $E$ is unramified and that $a(\theta)$ is even we have $H=J$ and $H^1=J^1$.}
\begin{equation*}
    J = E^{\times}\cdot K(r) \text{ and } J^1=(1+\mathfrak{p}_E)\cdot K(r).
\end{equation*}
One can now extend $\theta$ to $\widetilde{\theta}\colon J^1\to S^1$ by setting
\begin{equation*}
    \widetilde{\theta}(x\cdot k_1) = \theta(x)\cdot \theta_{b_{\theta}}(k_1).
\end{equation*}
We call $(J^1,\widetilde{\theta})$ a simple character. It is instructive to check that this is well-defined. Indeed, suppose that $xk_1=x'k_1'\in J^1$. Then we have $$x^{-1}x' = k_1(k_1')^{-1}\in (1+\mathfrak{p}_E)\cap K(r).$$ Now we can compute
\begin{align}
	\tilde{\theta}(x'k_1') &= \theta(x')\theta_{b_{\theta}}(k_1') = \theta(x) \theta(x^{-1}x')\theta_{b_{\theta}}(k_1') \nonumber\\
	&= \theta(x) \psi_E(b_{\theta}(x^{-1}x'-1))\theta_{b_{\theta}}(k_1') = \theta(x) \theta_{b_{\theta}}(k_1(k_1')^{-1})\theta_{b_{\theta}}(k_1') \nonumber\\
	&= \theta(x)\theta_{b_{\theta}}(k_1) = \tilde{\theta}(xk_1). \nonumber
\end{align}

Next one constructs an extension $(J,\Lambda)$ of $(J^1,\widetilde{\theta})$. This extension is uniquely determined by the properties:
\begin{enumerate}
    \item $\Lambda$ is a character, $\Lambda\vert_{J^1}=\widetilde{\theta}$ and $\Lambda\vert_{Z(F)}=\theta\vert_{F^{\times}}$; and
    \item We have $\Lambda\vert_{E^{\times}} = \theta$.
\end{enumerate}
See \cite[Lemma~3.18]{Ne-Hu_test} for more details as well as for the corresponding general notions.

\begin{defn}\label{def:cusp_type}
We call the triple $(L(0), J,\Lambda)$ a \emph{cuspidal type}. Furthermore, it is said to be associated to $\pi$ if $\pi=\textrm{c-Ind}_J^{G(F)}(\Lambda).$
\end{defn}

\begin{rem}
Note that we have started with a specific character $\theta$ and used it to construct a supercuspidal representations $\pi$ satisfying Assumption~\ref{ass:super}. It is even true that all supercuspidal representations fulfilling these assumptions arise this way.
\end{rem}
 
\begin{defn}\label{def:min_vex}
We call $v \in \pi$ a \emph{minimal vector} if there is a cuspidal type $(L(0),J,\Lambda)$ associated to $\pi$ such that $v$ is an eigenvector of the simple character $(\widetilde{\theta},J^1)$.
We denote such a vector by $v_\MV$.
\end{defn}

In symbols, the defining property of a minimal vector reads
\begin{equation}
	\pi(j)v_{\MV} = \widetilde{\theta}(j)\cdot v_{\MV} \text{ for all }j\in J^1.\nonumber
\end{equation}
Thus, according to Definition~\ref{def:loc_vec}, $v_{\MV}$ is localised at $(\widetilde{\theta},J^1)$. Note that in our situation a minimal vector actually satisfies $\pi(j)v_{\MV} = \Lambda(j)v_{\MV}$ for all $j\in J$. Thus, it is even localised at $(\Lambda ,J)$.\footnote{This is because we are assuming $c(\theta)$ to be even, so that in the language of \cite[Definition~3.21]{Ne-Hu_test} Type 1 and Type 2 vectors agree.} Using the fact that $(L(0),J,\Lambda)$ is a cuspidal type associated to $\pi$, one can easily verify that the dimension of $(\widetilde{\theta},J^1)$ localised vectors is $1$. Furthermore, it is a consequence of type theory, that the existence of $(\widetilde{\theta},J^1)$-localised vectors fully determines the isomorphism class of $\pi$. See for example \cite[15.7]{bushnell}. We conclude that, in view of Remark~\ref{rem:about_f1f2}, we can expect the heuristic \eqref{eq:matrix_coeff_exp.} to be true for the matrix coefficients of minimal vectors. This is one of the most important results concerning minimal vectors, and we formulate it properly.

\begin{lemmy}[Lemma~3.22, \cite{Ne-Hu_test} and Proposition~3.2, \cite{Hu-Nel-Sa}] \label{lm:mat_coeff_super}
Let $v_\MV\in \pi$ be a minimal vector, then the matrix coefficient $\langle \pi(g)v_\MV,v_\MV \rangle$ is supported in $J$. Even more, we have
\begin{equation}
	\langle \pi(g)v_\MV,v_\MV \rangle = \begin{cases}
		\Lambda(j)^{-1}\cdot \Vert v_{\MV}\Vert^2 &\text{ if }g=j\in J,\\
		0 &\text{ else.}
	\end{cases}\nonumber
\end{equation}
\end{lemmy}

We record the following result from \cite[Lemma~2.9]{Hu-Sh-Yi}, which is a version of \cite[Lemma~A.7]{Ne-Hu_test}.
\begin{lemmy}
Up to a constant multiple a minimal vector $v_\MV$ is given in the Kirillov model by
\begin{equation*}
    g_0(y)=\vvmathbb{1}_{\varpi^{-2r}(1+\mathfrak{p}^{r})}(y).
\end{equation*}
\end{lemmy}

We directly obtain the following, which is a version of \cite[Corollary~2.10]{Hu-Sh-Yi}.

\begin{cor}\label{cor:decomp_newform_super}
Let $v_\MV$ denote the minimal vector and let $v^\circ$ be the generalised $K_H(2r,2r)$-invariant new-vector. 
Then we have
\begin{equation*}
    \frac{v^\circ}{\Vert v^\circ\Vert } = \zeta_q(1)^{\frac{1}{2}}q^{-\frac{r}{2}}\sum_{x\in (\mathcal{O}/\mathfrak{p}^r)^{\times}} \frac{\pi(a(x))v_\MV}{\Vert v_\MV\Vert }.
\end{equation*}
\end{cor}

\begin{rem}\label{rem:v-new-dim-supercusp}
The vector $v_\MV$ generates the same $K_0$-type as the generalised new-vector $v^\circ$, see Convention~\ref{conv:newvec}. This $K_0$-type, which we denote by $V_{\textrm{new}}$, has dimension $q^{\frac{a(\pi)}{2}}(1-\frac{1}{q})$. 
This formula for the dimension is obtained for example by using \cite[Lemma~2.2]{Ass} and recalling that $\pi$ is automatically twist minimal in our case.
\end{rem}

\subsubsection{Microlocal lift vectors}

Microlocal lift vectors play a major role in \cite{Nel-QUE-QP, marsh}. 
In \cite{Nel-QUE-QP} the properties of the microlocal lifts are pragmatically obtained via well-known properties of the new-vector in $1\boxplus \chi^{-2}$ and transported via twisting. Here we aim to develop the theory using \textit{microlocal calculations}. From this point of view many ideas can be found in \cite{marsh, Hu-Nel}, and we closely follow \cite{marsh}.

We will assume that $\pi=\chi\boxplus \chi^{-1}$ for $\chi\colon F^{\times}\to S^1$ with $a(\chi) = m = 2r\geq 2$. Recall that in \eqref{eq:chi_BBBB} we lifted $\chi$ to a character $\chi_B\colon B\to S^1$ on $B$. For convenience, we restate the definition:
\begin{equation*}
    \chi_B\left(\left(\begin{matrix} t_1 & x\\ 0 &t_2\end{matrix}\right)\right) = \chi(t_1/t_2).
\end{equation*}
Further, the compact model for $\pi$ was denoted by $I(\chi_B)$. It consists of functions $f\colon K_0\to \C$ satisfying 
\begin{equation*}
    f(bk) = \chi_B(b)f(k) \text{ for }k\in K_0 \text{ and }b\in B(\mathcal{O}).\nonumber
\end{equation*}

Since $a(\chi)=2r$, the map $x\mapsto \chi(1+\varpi^rx)$ is additive on $\mathcal{O}$. Thus, we can write
\begin{equation}
    \chi(1+\varpi^{r}x) = \psi(b_{\chi}\cdot x\cdot \varpi^{-r})\nonumber
\end{equation}
where $b_{\chi}\in \mathcal{O}^{\times}$ is uniquely determined modulo $\mathfrak{p}^r$. We define
\begin{equation*}
    \tau({\chi}) = \varpi^{-m}\left(\begin{matrix} b_{\chi} & 0\\ 0&-b_{\chi}\end{matrix}\right)\in L(-m).
\end{equation*}
More generally, given $\tau \in L(-m)$, we define $\theta_{\tau}\colon K(r)\to S^1$ by
\begin{equation}
	\theta_{\tau}(k) = \psi(\tr(\tau(k-1))).\nonumber
\end{equation}
This is analogous to \eqref{some_def}.

We start with an easy but very useful observation.

\begin{lemmy}\label{lm:char_thick_tor}
The following statements are true:
\begin{enumerate}
    \item Let $\tau\in L(-m)$ be such that  $\chi_B\vert_{B(\mathcal{O})\cap K(r)}=\theta_{\tau}\vert_{B(\mathcal{O})\cap K(r)}$. Then we have
    \begin{equation}
        \tau - \left(\begin{matrix} \varpi^{-m}b_{\chi} & \tau_{1,2} \\ 0 & -\varpi^{-m}b_{\chi} \end{matrix}\right) \in L(-r).\nonumber
    \end{equation}
    \item Let $\chi_{D}$ be the character on $D(\mathcal{O})$ given by $\chi_{D}(\diag(t_1,t_2))=\chi(t_1/t_2).$ Then $\widetilde{\chi}(tk)=\chi_{D}(t)\cdot \theta_{\tau(\chi)}(k)$ is the unique extension of $\chi_{D}$ to $\widetilde{D}({r})=D(\mathcal{O})\cdot K(r).$
\end{enumerate}
\end{lemmy}
\begin{proof}
To prove the first statement, we take
\begin{equation*}
    b=\left(\begin{matrix} 1+\varpi^{r}t_1 & \varpi^{r}y \\ 0 & 1+\varpi^{r}t_2\end{matrix}\right)\in B(\mathcal{O})\cap K(r).\nonumber
\end{equation*}
We compute
\begin{equation}
    \chi_B(b) = \psi((t_1-t_2)b_{\chi}\varpi^{-r}).\label{eq:1}
\end{equation}
On the other hand, unpacking the definition of $\theta_{\tau}$ we get
\begin{equation}
    \theta_{\tau}(b) = \psi((\tau_{1,1}t_1+\tau_{2,2}t_2+\tau_{2,1}y)\varpi^{r}). \label{eq:2}
\end{equation}
Note that by assumption \eqref{eq:1} is equal to \eqref{eq:2} for all $t_1,t_2,y\in \mathcal{O}$. This directly gives the first result.

The second claim is essentially \cite[Lemma~2.1]{marsh}. To start, let us note that $\widetilde{\chi}$ is obviously well-defined. It remains to check uniqueness. Assume that $\xi$ is some extension of $\chi_{D}$ and consider $\xi'=\widetilde{\chi}\cdot \xi^{-1}$. By assumption, we have $\xi'\vert_{D(\mathcal{O})}\equiv 1$. On the other hand, $\xi'$ is a character, and so it is trivial on the commutator $[\widetilde{D}({r}),\widetilde{D}({r})]$. 
Observe that 
\begin{equation}
	\left[\left(\begin{matrix} 1 & \mathfrak{p}^r \\ 0 & 1 \end{matrix}\right),\left(\begin{matrix} \mathcal{O}^{\times} & 0 \\ 0 & 1 \end{matrix}\right)\right] = \left(\begin{matrix} 1 & \mathfrak{p}^r \\ 0 & 1\end{matrix}\right) \text{ and } 	\left[\left(\begin{matrix} 1 & 0 \\ \mathfrak{p}^r & 1 \end{matrix}\right),\left(\begin{matrix} 1 & 0 \\ 0 & \mathcal{O}^{\times} \end{matrix}\right)\right] = \left(\begin{matrix} 1 & 0 \\ \mathfrak{p}^r & 1\end{matrix}\right). \nonumber
\end{equation}
In particular, $[\widetilde{D}({r}),\widetilde{D}({r})]$ contains (among others) matrices of the form
\begin{equation*}
    \left(\begin{matrix} 1 & \mathfrak{p}^r\\0&1\end{matrix}\right) \text{ and }\left(\begin{matrix} 1 & 0\\ \mathfrak{p}^{r} &1\end{matrix}\right).
\end{equation*}
We conclude that $\xi'$ must be trivial.
\end{proof}

\begin{defn}[Definition~2.1, \cite{marsh}]\label{def:chi-type}
We call the tuple $(\widetilde{D}(r),\widetilde{\chi})$ the \emph{standard $\chi$-type}. A general $\chi$-type is then a tuple $(J,\lambda)$ obtained by conjugating the standard $\chi$-type by an element in $K_{0}$. (I.e. $J=\Ad_k\widetilde{D}(r)$ and $\lambda=\Ad^{\ast}_k(\widetilde{\chi})$ for $k\in K_{0}$.)
\end{defn}

We can now construct a $(\widetilde{\chi},\widetilde{D}(r))$-localised vector. 

\begin{lemmy}\label{def_and_uniq}
The representation $\pi$ contains an up to scaling unique $(\widetilde{\chi},\widetilde{D}(r))$-localised element denoted by $v_{\ML} \in \pi$. The corresponding element $f_{\ML}\in I(\chi_B)$ is given by 
\begin{equation*}
    f_{\textrm{ML}}(k) = \begin{cases}
        \chi_B(b)\theta_{\tau(\chi)}(k_1) &\text{ if }k=bk_1\in B(\mathcal{O})K(r), \\
        0 &\text{ else.}
    \end{cases}
\end{equation*}
\end{lemmy}
\begin{proof}
For the proof it is sufficient to work in the compact model $I(\chi_B)$. More  precisely, we have to show that $f _{\textrm{ML}}$ is well-defined and essentially unique. The first fact follows directly from Lemma~\ref{lm:char_thick_tor}. 

We turn to uniqueness. Suppose that $f\in I(\chi_B)$ is another non-trivial vector with the same localisation properties. Since we are assuming that $f$ is non-trivial, we can fix $k = \begin{psmallmatrix} a & b \\ c & d \end{psmallmatrix}\in K_0$ with $f(k)\neq 0$. For $b\in B(\mathcal{O})\cap K(r)$ we compute that
\begin{equation}
    \chi_B(b)f(k) = f(bk)=f(k)\theta_{\tau(\chi)}(k^{-1}bk) = f(k)\theta_{\Ad_k(\tau(\chi))}(b).\nonumber
\end{equation}
In particular, we have $\chi_B(b)=\theta_{\Ad_k(\tau(\chi))}(b)$ for $b\in B(\mathcal{O}) \cap K(r).$ In view of the first part of Lemma~\ref{lm:char_thick_tor}, we obtain that
\begin{equation}
	\frac{\varpi^{-m}b_{\chi}}{ad-bc}\cdot \left(\begin{matrix} ad+bc & -2ab \\ 2cd & -(ad+bc) \end{matrix}\right) = \Ad_k(\tau(\chi)) \in \left(\begin{matrix} \varpi^{-m}b_{\chi} & -\frac{b_{\chi}\varpi^{-m}}{ad-bc}ab \\ 0 & -\varpi^{-m}b_{\chi}\end{matrix}\right) + L(-r). \nonumber
\end{equation}
Recall that we are working in odd residual characteristic, so that comparing entries gives the two conditions
\begin{equation}
	bc \in \mathfrak{p}^r \text{ and } cd \in \mathfrak{p}^r.\nonumber
\end{equation}
This implies $k\in B(\mathcal{O})\cdot K(r)$. In particular, using the transformation behaviour of $f$ shows that $f(1)\neq 0$ and this lets us conclude that $f$ is a scalar multiple of $f_{\ML}$.
\end{proof}

\begin{defn}[Microlocal Lift Vector]\label{def:mic_loc_loft}
We call $v\in I(\chi)$ a \emph{microlocal lift vector} if there is a $\chi$-type $(J,\lambda)$ such that $\pi(j)v=\lambda(j)v$ for all $j\in J$. More precisely we say that $v$ is a microlocal lift vector associated to $(J,\lambda)$. The microlocal lift vector $v$ associated to the standard $\chi$-type (i.e. $v_\ML=f_{\textrm{ML}}$) is referred to as standard microlocal lift vector.
\end{defn}

\begin{rem}
This is essentially \cite[Definition~2.3]{marsh}. It should be noted that the standard microlocal lift vector features prominently in \cite{Nel-QUE-QP}.
\end{rem}

\begin{lemmy}[Lemma~2.4, \cite{marsh}]
For each $\chi$-type $(J,\lambda)$, the space of microlocal lift vectors in $I(\chi)$ associated with $(J,\lambda)$ is $1$-dimensional.
\end{lemmy}
\begin{proof}
After conjugation by an element from $K_{0}$, this follows directly from the uniqueness of the standard microlocal lift vector.
\end{proof}

\begin{lemmy}[Lemma~2.5, \cite{marsh}]\label{typ_u}
A nonzero microlocal lift vector is associated to at most one $\chi$-type.
\end{lemmy}
\begin{proof}
Without loss of generality, we can assume that $v$ is the standard microlocal lift vector.  Suppose that there is another $\chi$-type $(J,\lambda)$ such that $\pi(j)v=\lambda(j)v$. By definition, we have $J=k\widetilde{D}(r)k^{-1}$ for $k=\begin{psmallmatrix} \alpha & \beta \\ \gamma & \delta\end{psmallmatrix}\in K_{0}$. Furthermore,
\begin{equation}
    \widetilde{\chi}\vert_{K(r)} = \Ad_k^{\ast}(\widetilde{\chi})\vert_{K(r)}. \label{useful_thing}
\end{equation}
Writing out what this means leads to congruences, which in turn imply $k\in \widetilde{D}(r)$. Indeed, take $g=1+\varpi^r \begin{psmallmatrix} a & b \\ c & d \end{psmallmatrix}\in K(r)$. Then we have
\begin{multline}
	\psi(\varpi^{-r}b_{\chi}(a-d)) = \widetilde{\chi}(g) = \Ad_k^{\ast}(\widetilde{\chi})(g) \\ =  \psi\left(\frac{\varpi^{-r}b_{\chi}}{\alpha\delta-\beta\gamma}\cdot \left[(\alpha\delta+\beta\gamma)(a-d)+2\beta\delta \cdot c-2\alpha\gamma \cdot b\right]\right).\nonumber
\end{multline}
Recall that $g\in K(r)$ can be anything, so that we can vary $(a-d),b,c\in \mathcal{O}$. We obtain the congruence conditions
\begin{equation}
	\frac{\alpha\delta+\beta\gamma}{\alpha\delta-\beta\gamma} \in 1 + \mathfrak{p}^r,\, 2\beta\delta\in \mathfrak{p}^r \text{ and } 2\alpha\gamma\in \mathfrak{p}^r.\nonumber
\end{equation}
Since $\alpha\delta-\beta\gamma =\det(k)\in \mathcal{O}^{\times}$ (and $2\in \mathcal{O}^{\times}$) the only solution to this system is $\beta,\gamma\in \mathfrak{p}^r$. In other words $k\in \widetilde{D}(r)$ as desired.
\end{proof}

\begin{lemmy}[Lemma~2.6, \cite{marsh}]
The stabiliser in $K_0$ of a $\chi$-type $(J,\lambda)$ is equal to $J$.
\end{lemmy}
\begin{proof}
The computations relevant to the proof of this result are very similar to those in the proof of Lemma~\ref{typ_u} and we will only sketch them. First, by conjugation, we reduce to the case $(J,\lambda) = (\widetilde{D}(r),\widetilde{\chi})$. Next, suppose that $k\in K_0$ stabilises $(\widetilde{D}(r),\widetilde{\chi})$. Since $K(r)$ is normal in $K$ we find that $K(r)\subseteq \widetilde{D}(r) \cap k\widetilde{D}(r)k^{-1}$. Thus, we arrive at the identity \eqref{useful_thing} and can proceed as in the proof of Lemma~\ref{typ_u} to show that $k\in \widetilde{D}(r)$ as claimed. 
\end{proof}
 
We now turn towards understanding the matrix coefficients associated to microlocal lift vectors. As in the case of minimal vectors considered above we expect them to behave as predicted in \eqref{eq:matrix_coeff_exp.}.\footnote{Note that this ignores some issues arising from the fact that matrix coefficients of principal series representations are not compactly supported modulo the centre.} Indeed, in lieu of the strong uniqueness properties exhibited by $\chi$-types, this expectation is morally supported by the discussion in Remark~\ref{rem:about_f1f2}. For our purposes we will content ourselves with the following result concerning the support of the matrix coefficients, which we can prove rather directly.

\begin{lemmy}[Lemma~2.7, \cite{marsh}] \label{lm:mat_coeff_princ}
We have $\langle \pi(g)v_\ML,v_\ML\rangle=0$ unless $g\in K(r)D(F)K(r).$
\end{lemmy}
\begin{proof}
Let $\varphi(g) = \langle \pi(g)v_\ML,v_\ML\rangle$. We can express this matrix coefficient in the induced model. To do so we extend $f_{\textrm{ML}}\in I(\chi_B)$ to a function on $G(F)$ via the Iwasawa decomposition and recall that the inner product in the induced model is given by
\begin{equation*}
    \varphi(g)= \langle \pi(g)v_\ML,v_\ML\rangle = \int_{K_0} f_{\textrm{ML}}(kg)\overline{f_{\textrm{ML}}(k)}dk.
\end{equation*}
We now compute that
\begin{equation*}
    \langle \pi(g)v_\ML,v_\ML\rangle = [B(\mathcal{O})\colon B(\mathcal{O})\cap K(r)]\int_{K(r)}\theta_{\tau(\chi)}(k)^{-1}f_{\textrm{ML}}(kg)dk.
\end{equation*}
Note that $f_{\textrm{ML}}(kg)=0$ unless $kg\in BK(r).$ In particular, we must have $g\in K(r)BK(r)$. After using the centre and making an easy change of variables we can assume that $$g=\left(\begin{matrix} t & x \\0&1\end{matrix}\right)\in B.$$ With this information gathered we return to the abstract setting. As in Lemma~\ref{microloc_orth} we start by observing that
\begin{equation*}
    \theta_{\tau(\chi)}(k_1)\varphi(g) = \varphi(k_1g)=\theta_{\tau(\chi)}(g^{-1}k_1g)\varphi(g)
\end{equation*}
for $k_1\in K(r)\cap gK(r)g^{-1}$. If we assume that $\varphi(g)\neq 0$, then we must have
\begin{equation}
    \theta_{\tau(\chi)}(k_1) = \theta_{\tau(\chi)}(g^{-1}k_1g) \text{ for }k_1\in K(r)\cap gK(r)g^{-1}. \label{eq:equality}
\end{equation}
We define the lattice $L_g = L(r)\cap gL(r)g^{-1}$ and consider the corresponding dual lattice $L_g^{\ast}$. See \eqref{eq:dual_lat} for the general definition of the dual lattice. We now claim that \eqref{eq:equality} holds if and only if $$\tau(\chi)-g\tau(\chi)g^{-1}\in L_g^{\ast} = L(r)^{\ast}+gL(r)^{\ast}g^{-1}.$$ This is easy to see and stated explicitly in \cite[Lemma~2.1]{How2} for example. We conclude that there are $X,X'\in L(r)$ such that
\begin{equation}
    \left(\begin{matrix} b_{\chi} & 0\\0&-b_{\chi}\end{matrix}\right)+X=g\left(\left(\begin{matrix} b_{\chi} & 0\\ 0 & -b_{\chi}\end{matrix}\right)+X'\right)g^{-1}.\label{another}
\end{equation}
Using the specific shape of $g$ allows us to write this as
\begin{equation*}
    \left(\begin{matrix} 0 & 2xb_{\chi} \\ 0 & 0 \end{matrix}\right) = -X+gX'g^{-1}.
\end{equation*}
Using the three conditions $[gX'g^{-1}]_{ij} \in \mathfrak{p}^r$ for $(i,j) \in \{(1,1),(2,1), (2,2)\}$ allows us to deduce the congruence $$2xb_{\chi} \in \mathfrak{p}^r+t\mathfrak{p}^r+x\mathfrak{p}^r$$ from the top right entry of \eqref{another}. If $x\in \mathfrak{p}^r$, then $g=n(x)a(t)\in K(r)D(F)$ and we are done. Otherwise, we must have $v(t)\leq v(x)-r$. But this implies that $g=a(t)n(x/t)$ and $n(x/t)\in K(r)$. This completes the proof.
\end{proof}

\begin{rem}
The proof given in \cite{marsh} is more conceptual and crucially uses that $\tau(\chi)$ is regular. Similar arguments are very common in the literature. See for example \cite[Lemma~3.3]{How2}. In particular, the crucial technical lemma \cite[Lemma~2.8]{marsh} is a slight refinement of \cite[Lemma~3.1]{How2}. Here we took a more direct approach, because for $\GL(2)$ this seems to be more transparent.
\end{rem}

We turn our attention to the Kirillov model and let $g_{\textrm{ML}}\in \mathcal{K}(\pi,\psi)$ be the element in the Kirillov model corresponding to $v_\ML$. In \cite[Lemma~48]{Nel-QUE-QP} it is shown that
\begin{equation}
    \frac{g_{\textrm{ML}}(y)}{\Vert g_{\textrm{ML}}\Vert} = \eta\cdot \zeta_q(1)^{-\frac{1}{2}}q^{-\frac{r}{2}}\chi(y)\vert y\vert^{\frac{1}{2}}\vvmathbb{1}_{\mathfrak{p}^{-r}}(y),\label{eq:ML_vec_kir}
\end{equation}
for some $\eta\in S^1$. An easy computation shows that this is orthogonal to $\xi_{\mathbf{1}}^{(m)}\in \mathcal{K}(\pi,\psi)$. As we have seen in Lemma~\ref{lm:kirnew} the latter function belongs to the generalised new-vector $v^{\circ}$. We conclude that the two vectors $v_{\ML}$ and $v^{\circ}$ are orthogonal. 

However, for our global application we will need to construct a localised vector that is not orthogonal to the generalised new-vector (see more in Lemma \ref{lm:summary_Q}). This motivates us to study translates of $v_\ML$. It turns out that certain $K_0$-translates are sufficient. Put $$\gamma_{a} =\left(\begin{matrix} 0&-1\\1& a\end{matrix}\right) \text{ and }\widetilde{\gamma}_{a} =\left(\begin{matrix} 1&0\\ a&1\end{matrix}\right).$$ Note that
\begin{equation*}
    \widetilde{\gamma}_{a\varpi^i} = \diag(1,a)\widetilde{\gamma}_{\varpi^i}\diag(1,a^{-1}) \text{ and }\gamma_{a\varpi^i} =\diag(1,a)\gamma_{\varpi^i}\diag(a^{-1},1). 
\end{equation*}
Furthermore, for $a\in \mathcal{O}^{\times}$ we have
\begin{equation*}
    \left(\begin{matrix} a^{-1} & 1 \\ 0 & a\end{matrix}\right)\gamma_{a^{-1}} = \widetilde{\gamma}_a.
\end{equation*}
In particular, we only need to consider $\widetilde{\gamma}_a$ for $a\in \mathfrak{p}.$ More generally, we are motivated  by the double coset decomposition
\begin{equation*}
    K_0 = \bigsqcup_{a\in \mathcal{O}/\mathfrak{p}^{r}} B(\mathcal{O})\left(\begin{matrix} 0 & -1 \\ 1 & a \end{matrix}\right)K(r) \sqcup \bigsqcup_{a\in \mathfrak{p}/\mathfrak{p}^{r}} B(\mathcal{O})\left(\begin{matrix} 1 &0\\a & 1\end{matrix}\right)K(r).
\end{equation*}
This is a mild extension of \cite[Lemma~3.4]{Ass}. Alternatively it can be checked by considering the right-action of $K_0$ on the vector $(0,1) \in \mathcal{O}^2$ and identifying the quotient $B(\mathcal{O}) \backslash K_0$ with the set of primitive vectors in $\mathcal{O}^2$ quotiented by the action of $\mathcal{O}^\times$.

In principle one can use Lemma~\ref{lm:newvec_princ_ind} and work completely in the induced model. However, it seems to be clearer to work in the Kirillov model. To translate between the two models we need the following lemma.

\begin{lemmy}\label{lm:kirr_princ}
For $a\in \mathcal{O}$ let $v_a = \pi(\gamma_{a})v_\ML$ and $\widetilde{v}_a=\pi(\widetilde{\gamma}_{a})v_\ML$. Denote the corresponding elements in the Kirillov model by $g_a,\widetilde{g}_a\in \mathcal{K}(\pi,\psi)$. We have
\begin{equation*}
    \frac{g_a(y)}{\Vert g_a\Vert} = \zeta_q(1)^{-\frac{1}{2}}q^{-r/2}\chi(y)^{-1} \vert y\vert^{\frac{1}{2}}\vvmathbb{1}_{2b_{\chi}a\varpi^{-m}+\mathfrak{p}^{-r}}(y).
\end{equation*}
Furthermore, if $a=a_0\cdot\varpi^i$ for $0< i<r$ and $a_0\in \mathcal{O}^{\times}$, then
\begin{equation*}
    \frac{\widetilde{g}_a(y)}{\Vert \widetilde{g}_a\Vert} = \epsilon(\frac{1}{2},\chi^2)\zeta_q(1)^{-\frac{1}{2}}\chi(y\varpi^{2i})q^{-\frac{r}{2}}\vvmathbb{1}_{(2a_0b_{\chi}+\mathfrak{p}^{r-i})\varpi^{i-m}}(y)
\end{equation*}
Finally, \eqref{eq:ML_vec_kir} holds with $\eta=\epsilon(\frac{1}{2},\chi^2)$.
\end{lemmy}
\begin{proof}
Our goal is to compute $W_v(a(y)) = \Lambda_{\pi}(\pi(a(y)v)$. To do so we need a workable expression of $\Lambda_{\pi}$. For principal series this can be written down explicitly. Suppose $v$ is given by $f_v\in I(\chi)$. Then, using the Iwasawa decomposition, we can extend $f_v$ to $G(F)$ by setting
\begin{equation*}
    f_v(bk)=\chi(t_1/t_2)\vert t_1/t_2\vert^{\frac{1}{2}}f(k) \text{ for }b=\left(\begin{matrix}
        t_1 & y \\ 0 & t_2
    \end{matrix}\right).
\end{equation*}
Let $w=\left(\begin{matrix} 0 &1 \\ -1 & 0\end{matrix}\right)$ be the long Weyl element. We can now define 
\begin{equation*}
    \Lambda_{\pi}(f_v) = \int_{F}^{\textrm{st}} f_v\left( w\left(\begin{matrix} 1 & x \\ 0&1\end{matrix}\right)\right)\psi(x)^{-1}dx.\nonumber
\end{equation*}
Here $\int_{F}^{\textrm{st}} = \lim_{r\to\infty}\int_{\mathfrak{p}^{-r}}$ is the stabilized integral. This allows us to write
\begin{equation*}
    W_v(a(y)) =  \chi(y)^{-1}\vert y \vert^{\frac{1}{2}}\underbrace{\int_F^{\textrm{st}} f_v\left(\left(\begin{matrix} 0& 1 \\ -1 & -x\end{matrix}\right)\right) \psi(yx)^{-1}dx}_{=I(v,y)}.
\end{equation*}
We expand the integral as 
\begin{equation}
    I(v,y)=I_0(v,y) + \sum_{n=1}^{\infty} \chi(\varpi^{2n}) I_n(v,y) \nonumber
\end{equation}
for
\begin{align*}
    I_0(v,y) &= \int_{\mathcal{O}} f_v\left(\left(\begin{matrix} 0& 1 \\ -1 & -x\end{matrix}\right)\right) \psi(yx)^{-1}dx \text{ and }\\
    I_n(v,y) &= \int_{\mathcal{O}^{\times}} \chi(x)^{2} f_v\left(\left(\begin{matrix} 1& 0 \\ \varpi^nx & 1\end{matrix}\right)\right) \psi(yx^{-1}\varpi^{-n})^{-1}dx.
\end{align*}
Here we have used the common identity
\begin{equation*}
    \begin{pmatrix} 0& 1 \\ -1 & -x\end{pmatrix}= \begin{pmatrix} -x^{-1} & 1 \\ 0 & -x\end{pmatrix} \begin{pmatrix} 1 & 0 \\ x^{-1} & 1\end{pmatrix}.
\end{equation*}

At this point we will start to use the fact that we are dealing with specific vectors $v=v_a$. In this case $f_{v_a} = \pi(\gamma_a)f_{\textrm{ML}}$ is supported in $B(\mathcal{O})\gamma_aK(r)$. We start by computing $I_0(v_a,y)$. To do so we observe that
\begin{equation*}
    f_{v_a}\left(\left(\begin{matrix} 0& 1 \\ -1 & -x\end{matrix}\right)\right) = f_{\textrm{ML}}\left(\left(\begin{matrix} 0&1\\-1&-x\end{matrix}\right)\left(\begin{matrix}0 & -1\\1 & a\end{matrix}\right)\right) =f_{\textrm{ML}}\left(\left(\begin{matrix} 1 & a\\-x & 1-ax \end{matrix}\right)\right).
\end{equation*}
Note that on the support of $f_{\ML}$ we have $x\in \mathfrak{p}^r$, so that in particular $1-ax\in \mathcal{O}^{\times}$. We can write
\begin{multline*}
    f_{\textrm{ML}}\left(\left(\begin{matrix} 1 & a\\-x & 1-ax \end{matrix}\right)\right) = f_{\textrm{ML}}\left(\left(\begin{matrix} (1-ax)^{-1} & a\\0 & 1-ax \end{matrix}\right)\left(\begin{matrix} 1& 0\\ \frac{-x}{1-ax} & 1\end{matrix}\right)\right) \\ =\chi(1-ax)^{-2} f_{\textrm{ML}}\left(\left(\begin{matrix} 1& 0\\ \frac{-x}{1-ax} & 1\end{matrix}\right)\right).
\end{multline*}
In particular, we find that 
\begin{equation*}
    f_{\textrm{ML}}\left(\left(\begin{matrix} 1 & a\\-x & 1-ax \end{matrix}\right)\right) = \chi(1-ax)^{-2}\cdot \vvmathbb{1}_{\mathfrak{p}^{r}}(x).\nonumber
\end{equation*}
We conclude that
\begin{align*}
    I_0(v_a,y) &=\int_{\mathfrak{p}^{r}}\chi(1-ax)^{-2}\psi(yx)^{-1}dx \\
    &= q^{-r}\int_{\mathcal{O}}\chi(1-ax\varpi^{r})^{-2}\psi(yx\varpi^{r})^{-1}dx \\
    &= q^{-r}\int_{\mathcal{O}}\psi(2b_{\chi}ax\varpi^{-r}-yx\varpi^{r})dx\\
    &= q^{-r}\vvmathbb{1}_{2b_{\chi}a\varpi^{-m}+\mathfrak{p}^{-r}}(y).
\end{align*}
Finally, we observe that the integrals $I_n(v_a,y)$ must all vanish due to the support of $f_{v_a}$. Putting everything together yields
\begin{equation*}
    W_{v_a}(a(y))=q^{-r}\chi(y)^{-1}\vert y\vert^{\frac{1}{2}}\vvmathbb{1}_{2b_{\chi}a\varpi^{-m}+\mathfrak{p}^{-r}}(y).\nonumber
\end{equation*}
Using the definition of the inner product in the Kirillov model, one easily verifies that
\begin{equation*}
    \langle g_a,g_a\rangle = \zeta_q(1)q^{-r}.
\end{equation*}
This gives the desired result.

We now turn towards $v=\widetilde{v}_a$. By an easy support argument, which uses $a\in \mathfrak{p}$, we see that $I_0(\widetilde{v}_a,y)=0$. Thus, let us consider $I_n(\widetilde{v}_a,y).$  To do so we first compute
\begin{equation*}
    f_{\widetilde{v}_a}\left(\left(\begin{matrix} 1 & 0\\ \varpi^nx & 1\end{matrix}\right)\right) = f_{\textrm{ML}}\left(\left(\begin{matrix} 1 & 0\\a+\varpi^nx & 1\end{matrix}\right)\right) = \vvmathbb{1}_{\varpi^{-n}(-a+\mathfrak{p}^{r})}(x).\nonumber
\end{equation*}
We put $i=v(a)$, write $a=a_0\varpi^i$ and consider two cases. First suppose that $i<r$. Then $I_n(\widetilde{v}_a,y)=0$ unless $n=i$ and we have
\begin{align*}
    I_i(\widetilde{v}_a,y) &= \int_{-a_0^{-1}+\mathfrak{p}^{r-i}}\chi(x)^{-2}\psi(yx\varpi^{-i})^{-1}dx \\
    &= \chi(a_0)^{2}\int_{1+\mathfrak{p}^{r-i}}\chi(x)^{-2}\psi(a_0^{-1}yx\varpi^{-i})dx 
\end{align*}
This integral is sometimes referred to as incomplete Gauß sum and is well understood. See \cite[Remark~5.9]{Ass19} for the correct size and \cite[Remark~3.3.7]{Ass_thesis} for a precise evaluation. Indeed, we have
\begin{equation*}
     I_i(\widetilde{v}_a,y)=\epsilon(\frac{1}{2},\chi^2)\chi(ya_0^{-1}\varpi^{i})^{2}q^{-\frac{m}{2}}\vvmathbb{1}_{(2b_{\chi}+\mathfrak{p}^{r-i})\varpi^{i-m}}(a_0y).
\end{equation*}
We conclude that
\begin{equation*}
    W_{\widetilde{v}_a}(a(y)) = \epsilon(\frac{1}{2},\chi^2)\chi(y\varpi^{2i})q^{-\frac{m}{2}}\vvmathbb{1}_{(2a_0b_{\chi}+\mathfrak{p}^{r-i})\varpi^{i-m}}(y).
\end{equation*}
Checking the normalisation is now straight forward. 

Finally, suppose $a\in \mathfrak{p}^{r}$. In this case we see that $I_n(\widetilde{v}_a,y)=0$ unless $n\geq r$. In the latter case we have
\begin{align*}
    I_n(\widetilde{v}_a,y) &=\int_{\mathcal{O}^{\times}} \chi(x)^{-2}\psi(-yx\varpi^{-n})dx \\
    &= q^{-\frac{m}{2}}\epsilon(\frac{1}{2},\chi^2)\chi^2(y\varpi^{-n})\vvmathbb{1}_{\varpi^{n-m}\mathcal{O}^{\times}}(y),\nonumber
\end{align*}
We conclude that
\begin{align*}
    \widetilde{g}_a(y) &= \chi(y)^{-1}\vert y\vert^{\frac{1}{2}}\sum_{n\geq r} \chi(\varpi^{2n})I_n(\widetilde{v}_a,y) \\
    &= \epsilon(\frac{1}{2},\chi^2)\chi(y)\vert y\vert^{\frac{1}{2}}q^{-\frac{m}{2}}\sum_{n\geq r}\vvmathbb{1}_{\varpi^{n-m}\mathcal{O}^{\times}}(y) = \epsilon(\frac{1}{2},\chi^2)\chi(y)\vert y\vert^{\frac{1}{2}}q^{-\frac{m}{2}}\vvmathbb{1}_{\mathfrak{p}^{-r}}(y). 
\end{align*}
The $L^2$-normalization is now easily computed.
\end{proof}

\begin{rem}
One can compute $g_a$ following \cite{Nel-QUE-QP} using the non-equivariant twisting isomorphism. The relevant formulae for the corresponding new-vector translates are well known in this case. 
\end{rem}

We can now compute the corresponding decomposition of the generalised new-vector $v^\circ$.
This then gives clear candidates for which localised vectors to use in our treatment of the sup-norm problem.
We recall the sketch in the archimedean case given in Section \ref{sec:sketch-K-period} and Section \ref{sec:intro-rel-trace-formula}, where we want a certain period of the chosen vector to be large.
We will apply an analogous idea at the place $p$.
In this case, the period is given by an integral over $K_H(m, m)$, which translates to the inner product with the new-vector $v^\circ$ (see the next section).
Though this is more complicated in the archimedean world (due to approximations), in our $p$-adic setting we can directly compute inner products of localised vectors with the new-vector using the following orthogonal decomposition.

\begin{cor}\label{cor:newform_exp_princ}
Define the matrix 
\begin{equation} \label{eq:matrix-b-chi}
    c_t=\left(\begin{matrix} b_{\chi}/t & -\frac{1}{2} \\ 1 & \frac{t}{2b_{\chi}}\end{matrix}\right).
\end{equation}
Then we have
\begin{equation*}
    \frac{v^\circ}{\Vert v^\circ\Vert} = \sum_{t\in (\mathcal{O}/\mathfrak{p}^r)^{\times}} \frac{a_t }{\Vert v_\ML\Vert}\cdot \pi(c_t)v_\ML,\nonumber
\end{equation*}
with
\begin{equation*}
    a_t=\zeta_q(1)^{\frac{1}{2}}q^{-r/2}\chi(t\varpi^{-m})\psi(-b_{\chi}\varpi^{-m}).
\end{equation*}
\end{cor}
\begin{proof}
We compute in the Kirillov model. To do so we let $w_t\in \mathcal{K}(\pi,\psi)$ be the function associated to the vector $\frac{\pi(c_t)v_\ML}{\Vert g_{\textrm{ML}}\Vert}.$ Note that
\begin{equation*}
    c_t=n(b_{\chi}/t)\cdot \gamma_{\frac{t}{2b_{\chi}}}.\nonumber
\end{equation*}
In particular, we have
\begin{align*}
    w_t(y) &= \frac{1}{\Vert g_{\frac{t}{2b_{\chi}}}\Vert}\cdot [n(b_{\chi}/t).g_{\frac{t}{2b_{\chi}}}](y) = \psi(yb_{\chi}/t) \frac{g_{\frac{t}{2b_{\chi}}}(y)}{\Vert g_{\frac{t}{2b_{\chi}}}\Vert} \\
    &= \zeta_q(1)^{-\frac{1}{2}}q^{-r/2}\psi(yb_{\chi}/t)\chi(y)^{-1} \vert y\vert^{\frac{1}{2}}\vvmathbb{1}_{t\varpi^{-m}+\mathfrak{p}^{-r}}(y)
\end{align*}
We can modify $b_{\chi}\in \mathcal{O}^{\times}$ if necessary so that
\begin{equation*}
        \chi(1+z\varpi^{r}) = \psi(b_{\chi}z\varpi^{-r})\text{ for }z\in \mathcal{O}.\nonumber
\end{equation*}
We observe that
\begin{align*}
    w_t(t(1+z\varpi^r)\varpi^{-m}) &= \zeta_q(1)^{-\frac{1}{2}}q^{r/2}\chi(t\varpi^{-m})^{-1}\chi(1+z\varpi^r)^{-1}\psi(\varpi^{-m}(1+z\varpi^r)b_{\chi}) \\
    &= \zeta_q(1)^{-\frac{1}{2}}q^{r/2}\chi(t\varpi^{-m})^{-1}\psi(b_{\chi}\varpi^{-m}). 
\end{align*}
We can now compute the inner product between $g^\circ$ and $w_t$:
\begin{align*}
    \langle g^\circ,w_t\rangle &= \zeta_q(1)\int_{\varpi^{-m}\mathcal{O}^{\times}}\overline{w_t(y)}\frac{dy}{\vert y\vert} \\
    &= \zeta_q(1)\int_{1+\mathfrak{p}^r}\overline{w_t(t\varpi^{-m}y)}dy \\
    &= \zeta_q(1)q^{-r}\int_{\mathcal{O}} \overline{w_t(t(1+z\varpi^r)\varpi^{-m})}dz \\
    &= \zeta_q(1)^{\frac{1}{2}}q^{-r/2}\chi(t\varpi^{-m})\psi(-b_{\chi}\varpi^{-m}) =a_t.\nonumber
\end{align*}
This gives us the desired coefficients. Finally, we compute
\begin{equation*}
    \sum_{t\in (\mathcal{O}/\mathfrak{p}^r)^{\times}}\vert a_t\vert^2 = 1.\nonumber
\end{equation*}
Since the supports of the $w_t$'s are disjoint, they are orthogonal, and we are done.
\end{proof}

We end our discussion of microlocal lift vectors by making two remarks that will turn out to be useful later on.

\begin{rem}
    Lemma \ref{lm:mat_coeff_princ} shows that the matrix coefficient of $v_\ML$ localises around the torus $D(F)$. 
    It is easy now to compute that the matrix coefficient of $\pi(c_1)v_\ML$ localises around
    \begin{equation} \label{eq:torus-princ}
        c_1 D(F) c_1^{-1} = \left\{
            \frac12 
            \begin{pmatrix}
                \alpha+\beta & b_\chi (\alpha-\beta) \\
                b_\chi^{-1} (\alpha-\beta) & \alpha+\beta
            \end{pmatrix}
            : \alpha, \beta \in F^\times
        \right\}
    \end{equation}
\end{rem}

\begin{rem}\label{rem:v-new-dim-princ}
	It was shown in \cite[Lemma~2.17]{Sa} (see also \cite{marsh_local}) that the $K_0$-representation generated by $f_{\circ}$ is irreducible. It turns out that the standard microlocal lift vector $v_{\ML}$ generates the same representation. In the induced model we can write this representation down explicitly:
	\begin{equation*}
		V_{\textrm{new}} = \{ f\in I(\chi_B)\colon f(kc)=f(k) \text{ for }k\in K_0 \text{ and }c\in K(m)\}.
	\end{equation*}
	In particular, the dimension of this space is given by the cardinality of $B(\mathcal{O})\backslash K_0/K(m)$. Computing this cardinality reduces to computing $B(\mathcal{O}/\mathfrak{p}^m)\backslash G(\mathcal{O}/\mathfrak{p}^m)$. Since $\sharp G(\mathcal{O}/\mathfrak{p}^m) = q^{4m}(1-q^{-1})(1-q^{-2})$ and $\sharp B(\mathcal{O}/\mathfrak{p}^m) = q^{3m}(1-q^{-1})^{2}$ we conclude that
	\begin{equation}
		\dim_{\C} V_{\textrm{new}} =\sharp B(\mathcal{O})\backslash K_0/K(m) = \sharp B(\mathcal{O}/\mathfrak{p}^m)\backslash G(\mathcal{O}/\mathfrak{p}^m) = q^{m}(1+1/q). \nonumber
	\end{equation}
	This is in line with the proof of Proposition~2.13 for non-supercuspidal representations in \cite{Sa}. 
\end{rem}

\subsection{An example of a local period}\label{sec:loc_perio_fin}

We will now probe how the microlocalised vectors from our examples interact with vectors fixed under the (generalised) Hecke congruence subgroups $K_H(m,m)$, i.e. generalised new-vectors as in Convention~\ref{conv:newvec}. To capture this, we define\footnote{The notational clash with the archimedean section is intentional, and it is resolved in the next global section.}
\begin{equation}
    Q(v)=\frac{1}{\Vol(K_H(m,m))}\int_{K_H(m,m)}\langle \pi(k)v,v\rangle_{\pi} \, dk. \label{eq:loc_per_def}
\end{equation}
Given an irreducible unitary smooth representation $\pi$, we write $\mathcal{P}_{m,m}$ for the orthogonal projection on $\pi^{K_H(m,m)}$.

First, we observe that
\begin{equation}
    Q(v) = \langle \mathcal{P}_{m,m}(v),v\rangle = \Vert \mathcal{P}_{m,m}(v)\Vert^2. \label{eq:Q_to_proj}
\end{equation}
If we further assume that $\dim \pi^{K_H(m,m)} = 1$, then we can evaluate $Q(v)$ as follows. Let $v^{\circ}\in \pi^{K_H(m,m)}$ be the unit vector, so that the orthogonal projection $\mathcal{P}_{m,m}\colon \pi\to \pi$ is given by
\begin{equation*}
    \mathcal{P}_{m,m}(v)=\langle v,v^{\circ}\rangle\cdot v^{\circ}.
\end{equation*}
We can now write 
\begin{equation}
    Q(v) = \langle \mathcal{P}_{m,m}(v),v\rangle = \vert \langle v,v^{\circ}\rangle\vert^2.\nonumber
\end{equation}
In particular, if $v^{\circ} = a\cdot v+v^{\bot}$ with $\langle v,v^{\bot}\rangle=0$, then
\begin{equation*}
    Q(v) = \vert a \vert^2.
\end{equation*}

One of the main conclusions of our computations in the previous sections is the following evaluation of $Q$ at microlocalised vectors appearing in the decomposition of the new-vector.
These are the vectors for which $Q$ is large, at least among the class of microlocalised vectors.

\begin{lemmy}\label{lm:summary_Q}
If $\pi$ is a smooth unitary irreducible representation of $G(F)$ with trivial central character and exponent conductor $a(\pi)=4r$, then we have
\begin{equation*}
    Q(v) = \zeta_q(1)q^{-r},
\end{equation*}
where $v=v_\MV/\Vert v_\MV\Vert$  is the standard minimal vector if $\pi$ is supercuspidal and $v=\pi(c_t)v_\ML/\Vert v_\ML\Vert$ with $t\in (\mathcal{O}/\mathfrak{p}^r)^{\times}$ if $\pi$ is principal series.
\end{lemmy}
\begin{proof}
This follows directly from the discussion above using Corollary~\ref{cor:decomp_newform_super} if $\pi$ is supercuspidal and Corollary~\ref{cor:newform_exp_princ} if $\pi$ is principal series.
\end{proof}

Heuristically one can explain the size of $Q(v)$ as in the archimedean case. Indeed, in analogy to the reasoning from Section~\ref{sec:arch_mat_scetch}, one finds that
\begin{equation}
	Q(v) \approx \frac{\Vol(J\cap K_H(m,m))}{\Vol(K_H(m,m))} \approx q^{-r}\nonumber
\end{equation}
in the supercuspidal case, for example. Here $J$ is the group carrying the support of the matrix coefficient (see Lemma~\ref{lm:mat_coeff_super}).

\subsection{A non-archimedean volume estimate}\label{anavs}

In this section we prove non-archimedean volume estimates similar to those in Theorem \ref{thm:archi-main-estimate}.
The test function we consider is a truncated matrix coefficient of a microlocalised vector with large period $Q$, as in Lemma \ref{lm:summary_Q}.
By Lemma \ref{lm:mat_coeff_super} and Lemma \ref{lm:mat_coeff_princ}, these functions localise around a torus $T$.
Thus, such a test function takes the desired form \eqref{eq:intro-archi-test-fn} mentioned in the introduction.

In the case of supercuspidals, we take the minimal vector $v_\MV$, for which $T$ is the torus $E^\times$, where $E/F$ is an unramified extension.  
In the case of principal series, we consider the vector $\pi(c_1) v_\ML$, with $c_1$ defined in Corollary \ref{cor:newform_exp_princ}, for which $T$ is a conjugate of the diagonal torus given in \eqref{eq:torus-princ}.

More precisely, let $n \in \N$,  and let $\pi$ be an irreducible representation of conductor $a(\pi)=4n$ and trivial central character. 
Write $H(2n) = K_H(2n, 2n)$ for simplicity.
We fix the $L^2$-normalized vector $v\in \pi$ to be $v=\frac{\pi(c_1)v_\ML}{\Vert v_{\ML}\Vert}$ if $\pi$ is principal series and $v=\frac{v_\MV}{\Vert v_{\MV}\Vert}$ if $\pi$ is supercuspidal. 
Our goal is to understand the integral
\begin{equation}
	I_{\gamma}(f\ast f) = \frac{1}{\Vol(H(2n))^2}\int_{H(2n)}\int_{H(2n)} [f\ast f](x^{-1} \gamma y)dxdy,\label{eq:def_Igamma}
\end{equation}
for 
\begin{equation}
	f(g) = \dim(V_{\textrm{new}})\cdot \langle v,\pi(g)v\rangle \cdot \vvmathbb{1}_{K_0}(g).\nonumber
\end{equation}
Note that, by Lemma~\ref{mat_coeff_indep}, this function satisfies $f\ast f=f$.

For $\gamma\in H(2n)$, this integral $I_{\gamma}(f)$ is easily evaluated. Indeed, after a simple change of variables we have
\begin{multline}
	I_{\gamma}(f) = \Vol(H(2n))^{-1}\cdot \int_{H(2n)} f(x) dx \\ = \dim(V_{\textrm{new}}) \langle v ,\mathcal{P}_{2n,2n}(v)\rangle = \dim(V_{\textrm{new}}) \cdot Q(v). \nonumber
\end{multline}
Recalling Lemma \ref{lm:summary_Q}, Remark \ref{rem:v-new-dim-supercusp}, and Remark \ref{rem:v-new-dim-princ}, we record the estimate
\begin{equation}
	I_{\gamma}(f\ast f)\ll q^n \text{ for }\gamma\in H(2n),\label{eq:benchmark_1}
\end{equation}
which is essentially sharp.

For general $\gamma$, we will ignore any oscillation of the matrix coefficients $f$ and simply bound
\begin{equation}
	I_{\gamma}(f\ast f) \ll I_{\gamma}(\vert f\vert). \nonumber
\end{equation}
According to the support properties of the matrix coefficients $f$ we have
\begin{equation}
	\vert f\vert \ll \Vol(\mathcal{U})^{-1} \vvmathbb{1}_{\mathcal{U}}, \nonumber
\end{equation}
where $\Uc = T(\Oc) \cdot K(n)$ and $T(\Oc) = T \cap G(\Oc)$.

We follow Section \ref{sec:intro-transv}, where all steps are directly applicable and no approximations are involved.
Thus, we perform the $x$-integral for each $y$ such that $\gamma y \in H(2n) \Uc$ and obtain that
\begin{multline}
    \frac{1}{\Vol(H(2n))^2}\int_{H(2n)} \int_{H(2n)} \vert f\vert (x^{-1} \gamma y) \, dx \, dy \\ \ll \frac{\Vol\{ y \in H(2n) \mid \gamma y \in H(2n) \Uc \}}{\Vol(H(2n))}\cdot \frac{\Vol(\Uc\cap H(2n))}{\Vol(H(2n))\Vol(\Uc)}. \label{eq:123es_kommt_die_sonne}
\end{multline}
In this section we will prove a bound for the volume 
\begin{equation}
	V_{\gamma}(\mathcal{U}) = \Vol\{ y \in H(2n) \mid \gamma y \in H(2n) \Uc \}.\nonumber
\end{equation}

\begin{rem}\label{rm:trivial_bound}
Note that if $\gamma\in H(2n)K(n)$, then we have
\begin{equation}
	V_{\gamma}(\mathcal{U}) = \Vol(H(2n)). \nonumber
\end{equation}
Indeed, here we can reduce to the case that $\gamma = 1$, since $K(n)$ is normal and contained in $\Uc$, so $\gamma y \in H(2n) \Uc$ for all $y \in H(2n)$. We obtain the bound\footnote{We use that $\Vol(H(2n)) \asymp q^{-4n}$, $\Vol(\Uc)\asymp q^{-2n}$ and $\Vol(\Uc\cap H(2n))\asymp q^{-5n}$.}
\begin{equation}
	I_{\gamma}(f\ast f) \ll \frac{\Vol(\mathcal{U}\cap H(2n))}{\Vol(H(2n))\Vol(\Uc)} \ll q^{n}, \text{ for }\gamma\in H(2n)K(n).\nonumber
\end{equation}
For $\gamma\in H(2n)$ we recover \eqref{eq:benchmark_1}, but it is valid in a larger region.
\end{rem}

\subsubsection{Distances}

The goal is to bound the volumes above in terms of a notion of distance between $\gamma$ and $H(2n)$.

\begin{defn}
    For a subgroup $H \leq K_0$ and $g \in K_0$, we set $d_H(g) = 0$ if~$g \in H$. Otherwise, we define $d_H(g) = q^{-l}$ for~$l \in \Z_{\geq 0}$ such that $g \in H \cdot K(l) \setminus H \cdot K(l+1)$. 
\end{defn}

As in Section \ref{sec:intro-transv}, we may reduce $\gamma$ to a certain shape, although we do not have a Cartan decomposition.
We always assume $\gamma\in K_0$ and write $\gamma=(\gamma_{ij})$. We start with the case where $\gamma_{22}\in \Oc^{\times}$. Since $H(2n)$ contains $D(\Oc)$, we can reduce to the case that $\gamma_{22} = 1$. Observe also that $d_{H(2n)}(\gamma) = q^{-l}$ implies that at least $\gamma_{12}$ or $\gamma_{21}$ has valuation precisely $l$.

\begin{rem}
If we assume that $d_{H(2n)}(\gamma)<1$, so that $\gamma \in K(1) \cdot H(2n)$, then we always have $\gamma_{11},\gamma_{22}\in \Oc^{\times}$. This is due to the invertibility of the determinant.
\end{rem}

Recall that the cases where $\gamma \in H(2n) K(n)$ are handled in Remark~\ref{rm:trivial_bound}. Thus, we assume from now on that $\gamma \in H(2n) K(l)\setminus H(2n)K(l+1)$ for $l < n$.

\begin{rem}\label{rm:ich_muss_das_schwarze_vogel_schiessen}
The reduction $\gamma_{22}=1$ only works if $\gamma_{22}\in \Oc^{\times}$. If this is not the case, then we must have $d_{H(2n)}(\gamma)=1$. Furthermore, we also deduce that 

Note that if $\gamma_{22}\in \pf$ and $\gamma\in K_0$, then we must have $d_{H(2n)}(\gamma)=1$. In particular, $\gamma_{12},\gamma_{21}\in \Oc^{\times}$. We write $\gamma=w\cdot \gamma'$ and, since $w$ normalises $H(2n)$, we observe that
\begin{equation}
	V_{\gamma}(\mathcal{U}) = \{ y \in H(2n) \mid \gamma' y \in H(2n) w\Uc \}.\nonumber
\end{equation}
At this point, we observe that $H(2n)w\Uc = H(2n)\Uc$.\footnote{For example if $\mathcal{U}$ is the thickening of the standard torus conjugated by $c_1$, then we simply compute $$ wc_1 = \left(\begin{matrix} 1 & 1/2b_{\chi} \\ -b_{\chi} & \frac{1}{2}\end{matrix} \right) =  \left(\begin{matrix} -b_{\chi}^{-1} & 0\\0&b_{\chi} \end{matrix}\right)c_1\left(\begin{matrix} -1 & 0 \\0&1\end{matrix}\right).$$} We conclude that
\begin{equation}
		V_{\gamma}(\mathcal{U}) =  	V_{\gamma'}(\mathcal{U})\nonumber
\end{equation}
in this case.

In particular, for $\gamma=w$, we have
\begin{equation}
	V_{w}(\mathcal{U}) = V_1(\mathcal{U}) = \Vol(H(2n)). \nonumber	
\end{equation}
This is larger than what we expect, given that $d(\gamma) = 1$. This phenomenon is related to the failure of transversality for rank two as described in \cite[Remark~5.2]{Hu-Nel}.
\end{rem}

\subsubsection{Further reductions}

Since $K(n) \subset \Uc$, we can reduce the computation of $V_\gamma(\mathcal{U})$ modulo $\pf^n$, as is done in \cite[Sec. 9.4]{Hu-Nel}.
More precisely, we consider cosets $H(2n) / H(2n) \cap K(n)$, which can be represented by diagonal matrices with elements in $\Oc/\pf^n$.
Indeed, this follows by a computation similar to \eqref{eq:iwahori-decomp-kinda}. 
Of course, we will always pick up the relevant volume in the central direction, which is contained in $D$.
We then have
\begin{multline}
   V_{\gamma}(\mathcal{U})
    = \Vol(H(2n) \cap K(n)) \cdot \sharp(\Oc/\pf^n)^\times \\ \cdot \sharp\{ \alpha \in (\Oc/\pf^n)^\times \mid \gamma \diag(\alpha, 1) \in D(\Oc/\pf^n) T \}. \label{eq:reduction}
\end{multline}
Assume now that
\begin{displaymath}
    \gamma \diag(\alpha, 1) = d \cdot t,
\end{displaymath}
with $d \in D(\Oc/\pf^n)$ and $t \in T$.
Again, since $T$ contains $Z(\Oc)$, we may assume that $d = \diag(\beta, 1)$.

Recall that in both cases, supercuspidal and principal series, the torus element $t$ can be written as
\begin{displaymath}
    t = \begin{pmatrix}
        a & u_1 b \\
        u_2 b & a
    \end{pmatrix},
\end{displaymath}
where $a, b \in \Oc$ and $u_1, u_2 \in \Oc^\times$.
More precisely, for the minimal vector $u_1 = 1$ and $u_2 = D$, whereas for the shifted microlocal lift vector $u_1 = \beta_\chi$ and $u_2 = \beta_\chi^{-1}$.

\subsubsection{Polynomials with large coefficients and bounds}
To bound the last factor in the \eqref{eq:reduction} using the distance function, we reduce to counting solutions to polynomials with $p$-adically large coefficients.

Note now that our assumptions (i.e. $\gamma_{22}=1$) imply that $a = 1$ by looking at the lower-right entry in the matrix equation.
Explicitly, we have the condition
\begin{displaymath}
    \begin{pmatrix}
        \alpha \gamma_{11} & \gamma_{12} \\
        \alpha \gamma_{21} & 1 
    \end{pmatrix}
    =
    \begin{pmatrix}
        \beta & \beta u_1 b\\
        u_2 b & 1
    \end{pmatrix},
\end{displaymath}
modulo $\pf^n$.
Observe that, since $\beta, u_1, u_2 \in \Oc^\times$, we must have $\nu(\gamma_{12}) = \nu(\gamma_{21}) = \nu(b)$.
In particular, if $d_{H(2n)}(\gamma) = q^{-l}$, then $\nu(\gamma_{12})=\nu(\gamma_{21})=l$.
Let us also assume that $\gamma_{11} \in \Oc^\times$. (This happens for example if $d_{H(2n)}(\gamma)<1$.)

By computing the ratio of elements on the first row, we obtain the conditions
\begin{displaymath}
    \frac{\beta u_1 b}{\beta} = u_1 b = \frac{\gamma_{12}}{\alpha \gamma_{11}} \quad \text{and} \quad u_2 b = \alpha \gamma_{21}.
\end{displaymath}
Putting these together, $\alpha$ must satisfy
\begin{displaymath}
    \gamma_{21} \alpha^2 - c \equiv 0\, (\text{mod } \pf^n),
\end{displaymath}
for some $c \in \Oc$.
Assuming that $\alpha$ is close to the identity, we can produce a non-trivial bound from this condition as in \cite[Lemma 6.4]{marsh}.

\begin{lemmy}
    The number of $\alpha = 1 + z \in (\Oc/\pf^n)^\times$ with $z \in \pf$ satisfying $c_2 \alpha^2 - c_0 \equiv 0$ with $\nu(c_2) = l \leq n$ and $c_0 \in \Oc$ is bounded by $q^{(n+l)/2}$.
\end{lemmy}
\begin{proof}
    The polynomial equation can be rewritten as $c_2 z^2 + c_1 z + c_0 = 0$ with $c_1 \in \Oc$ and same conditions.
    Assuming the existence of any solution, we can shift $z$ by that solution and reduce to 
    the case where $c_0 = 0$.
    Furthermore, if $\nu(c_1) \leq l$, then $\nu(c_2z^2 + c_1 z) = \nu(c_1 z) \geq n$ implies that $\nu(z) \geq n-l$.
    Otherwise, we have that
    \begin{displaymath}
        \nu(c_2 z (z + c_1/c_2)) \geq n
    \end{displaymath}
    and, putting $\nu(c_2)=l$ on the other side, this implies that $\nu(z) \geq (n-l)/2$ or $\nu(z + c_1/c_2) \geq (n-l)/2$.
    The number of such $z$ is bounded in all cases by $q^{n-(n-l)/2}=q^{(n+l)/2}$.
\end{proof}

Recalling the definition of $d_{H(2n)}(\gamma)$, the previous lemma implies that 
\begin{equation}
    \sharp\{ \alpha \in (\Oc/\pf^n)^\times \mid \gamma \diag(\alpha, 1) \in D(\Oc/\pf^n) T \} \leq q \cdot \frac{q^{n/2}}{d_{H(2n)}(\gamma)^{1/2}} \label{eq:verdamm_mich_noch_eins}
\end{equation}
for $d_{H(2n)}(\gamma) \geq q^{-n}$ and $\gamma_{11},\gamma_{22}\in \Oc^{\times}$. The factor of $q$ accounts for the assumption that $\alpha$ is close to the identity, i.e. reducing to $\alpha \equiv 1$ modulo $\pf$.

\subsubsection{The main volume bound}\label{sec:p_adic_vol}

The results of this subsection imply the following important result.

\begin{theorem}\label{th:local_vol_bound}
For $\gamma\in K_0$ we have
\begin{equation}
		V_{\gamma}(\mathcal{U}) \ll  \Vol(H(2n)) \cdot q\max_{\gamma'\in\{\gamma,w\gamma\}}\min\left(\frac{q^{-n/2}}{d_{H(2n)}(\gamma')^{\frac{1}{2}}},1\right). \nonumber
\end{equation}
\end{theorem}
\begin{proof}
First, we note that if $d_{H(2n)}(\gamma)<q^{-n}$, then we have
\begin{equation}
	V_{\gamma}(\Uc) =\Vol(H(2n)) \nonumber 
\end{equation}
by Remark~\ref{rm:trivial_bound}. Thus, we assume $d_{H(2n)}(\gamma) \geq q^{-n}$. Using \eqref{eq:reduction} with $\Vol(H(2n)\cap K(n))\ll \Vol(H(2n))q^{-2n}$ and $\sharp (\Oc/\pf^n)\ll q^n$ yields
\begin{equation}
	V_{\gamma}(\Uc) \ll \Vol(H(2n)) q^{-n} \cdot \sharp \{ \alpha \in (\Oc/\pf^n)^\times \mid \gamma \diag(\alpha, 1) \in D(\Oc/\pf^n) T \}. \nonumber
\end{equation}
If $\gamma_{11},\gamma_{22}\in \Oc^{\times}$, then we can apply \eqref{eq:verdamm_mich_noch_eins} and obtain
\begin{equation}
		V_{\gamma}(\Uc) \ll q^{1-n}\cdot \Vol(H(2n)) \frac{q^{n/2}}{d_{H(2n)}(\gamma)^{\frac{1}{2}}}.\nonumber
\end{equation}
(At $d_{H(n)}(\gamma) = q^{-n}$ this agrees with the trivial bound up to a factor of $q$ and it is strictly worse from then on.)

Finally, if $\gamma_{11}\in \pf$ or $\gamma_{22}\in \pf$, then we have $\gamma_{12},\gamma_{21}\in \Oc^{\times}$. Thus, we can apply the argument above to $\gamma'=w\gamma$. As in Remark~\ref{rm:ich_muss_das_schwarze_vogel_schiessen} we obtain
\begin{equation}
	V_{\gamma}(\Uc) = V_{\gamma'}(\Uc) \ll q\cdot \Vol(H(2n))\min\left(\frac{q^{-n/2}}{d_{H(2n)}(\gamma')^{\frac{1}{2}}},1\right).\nonumber
\end{equation}
This completes the proof.
\end{proof}

We obtain the following important result concerning the size of our integral:

\begin{cor}\label{cor:vol_bound_p-adig}
We have
\begin{equation}
	I_{\gamma}(\vert f\vert) \ll \frac{q^{1-n}}{\Vol(\Uc)}\max_{\gamma'\in \{\gamma,w\gamma\}} \min\left(\frac{q^{-n/2}}{d_{H(2n)}(\gamma')^{1/2}},1\right).\nonumber
\end{equation}
\end{cor}
\begin{proof}
This follows directly from the previous theorem and \eqref{eq:123es_kommt_die_sonne}.
\end{proof}

\subsubsection{On the relation to matrix coefficients}\label{sec:p-adic_relmat}

In this subsection, we will discuss how much information we have lost when trivially estimating $I_{\gamma}(f)\ll I_{\gamma}(\vert f\vert)$. To do so, we will make a connection to estimates for matrix coefficients of new-vectors.

Recall that $\mathcal{P}_{2n,2n}$ is the projection on the (normalised) generalised new-vector in $\pi$ as above. We compute
\begin{align}
	 I_{\gamma}(f\ast f) = \dim(V_{\textrm{new}})\langle \mathcal{P}_{2n,2n}(v),\pi(\gamma)\mathcal{P}_{2n,2n}(v)\rangle, \nonumber
\end{align}
for $\gamma\in K_0$. We further have
\begin{equation}
	\langle \mathcal{P}_{2n,2n}(v),\pi(\gamma)\mathcal{P}_{2n,2n}(v)\rangle = Q(v)\langle v^{\circ},\pi(\gamma)v^\circ\rangle,\nonumber
\end{equation}
where $v^{\circ}$ is as in Convention~\ref{conv:newvec}. After recalling that $\dim(V_{\textrm{new}}) \asymp\Vol(\Uc)^{-1}$ and $Q(v)=\zeta_q(1)q^{-n}$ we obtain
\begin{equation}
	\langle v^{\circ},\pi(\gamma)v^{\circ}\rangle \ll q\max_{\gamma'\in\{\gamma,w\gamma\}}\min\left( \frac{q^{-n/2}}{d_{H(2n)}(\gamma')^{\frac{1}{2}}},1\right). \label{eq:Ich_mache_bescheuerte_labels_bis_du_Lammbock_schaust}
\end{equation}
Note that, by Cauchy--Schwarz, the trivial bound for the matrix coefficient is $1$. 
Thus, up to the factor $q$ our bound is always as least as good as the trivial bound.

Our method is quite soft and reduces essentially to a counting problem. 
More elaborate methods, more precisely the method of stationary phase applied to suitable integral representations of matrix coefficients, were applied in \cite{Hu-Sa}. 
In \cite[Proposition~5.2]{Hu-Sa} the following bounds are established:
\footnote{Since we are assuming $a(\pi)=4n$ we have $n_1=2n$ and because we are working directly with the normalized generalised new-vector (i.e. the unique vector fixed by $K_{H}(2n,2n)$) our vector $v^{\circ}$ corresponds to $v_{\pi}'/\Vert v_{\pi}\Vert$ in \cite{Hu-Sa}. 
Also note that our subgroup $H(2n)$ is denoted by $K^{\ast}(2n)$ in loc.cit. 
With these translations made we have $	\Phi'_{\pi}(g) = \langle v^{ \circ},\pi(g) v^{\circ}\rangle$.}
\begin{equation}
	\langle v^{\circ},\pi(\gamma) v^{\circ}\rangle \ll q^{\frac{j-2n}{2}+O(1)}, \text{ for }\gamma \in H(j)\setminus H(j+1) \text{ with }1\leq j\leq 2n. \label{eq:Mann_Mann_Achim}
\end{equation}
After recalling the definition of $d_{H(2n)}(\gamma)$, we observe that 
\begin{equation}
	d_{H(2n)}(\gamma) = q^{-l} \text{ if and only if }\gamma\in H(l)\setminus H(l+1), \nonumber
\end{equation}
for $0<l<2n$. Up to a term $q^{O(1)}$, these bounds are sharp and they are superior to our volume bounds. A full comparison is displayed in Table~\ref{table:Whatnot}.

\begin{table} 
    \label{table:Whatnot}
	\centering
	
	\begin{tabular}{l || c | c | c | c }
		 $l$ & $l=\infty$ & $n\leq l\leq 2n$ & $0<l <n$ & $l=0$   \\
		 \hline\hline 
		 $A_{\textrm{Vol}}(\gamma)$ & $1$  & $1$ & $\frac{l-n}{2}$ & $\frac{l'-n}{2}$  \\
		 \hline 
		 $A_{\textrm{Hu-Sa}}(\gamma)$ & $1$ & $\frac{l}{2}-n$ & $\frac{l}{2}-n$ & - 
 	\end{tabular}
    \caption{\emph{Bounds for matrix coefficients:} Here we have written $d_{H(2n)}(\gamma)=q^{-l}$ and $l' = \max\left(-\log(d_{H(2n)}(w\gamma))/\log(q),1\right)$. 
    Then we can write the bound for the matrix coefficient as $
        \langle v^{\circ},\pi(\gamma)v^{\circ}\rangle \ll q^{-A_{\ast}(\gamma)+O(1)},$
    where $A_{\textrm{Vol}}(\gamma)$ is the exponent derived from \eqref{eq:Ich_mache_bescheuerte_labels_bis_du_Lammbock_schaust} and $A_{\textrm{Hu-Sa}}(\gamma)$ is the exponent coming from \eqref{eq:Mann_Mann_Achim}.}
\end{table}

\section{Global framework}

So far we have abstractly discussed aspects of \textit{microlocalisation} and the \textit{quantitative orbit method}. We have done so in the realm of the local representation theory of $\PGL(2)$ over $\R$ and $\Q_p$. We will now switch gears and turn towards the global setting. Our goal for the rest of this text is to showcase how the tools developed earlier can be applied to concrete problems in the analytic theory of automorphic forms. More precisely, the application we have in mind is the \textit{sup-norm problem} for Hecke--Maa\ss\  newforms and we will ultimately prove in Section~\ref{sec:glob_perII} below the new bound stated in Theorem~\ref{th:Intro}. This requires some preparations.

Recall that we are ultimately interested in bounding classical Hecke--Maaß newforms $f\colon\Gamma_H(N)\backslash \mathbb{H}\to \C$, for $N=p^{4n}$. In order to effectively piece together the local theory developed above, we work in the adelic framework. In this section, we introduce the necessary terminology. Most importantly, we will define an adelic automorphic form $\varphi'$, see \eqref{eq:shifted_new}. This form is related to the classical form $f$, and in particular bounding the sup-norm of $\varphi'$ is equivalent to bounding the sup-norm of $f$.

\subsection{Adelic groups and measures}

We let $\mathbb{A}$ and $\mathbb{A}^{\times}$ denote the adeles and ideles, respectively, of $\Q$. 
The adeles come equipped with the modulus
\begin{equation*}
    \vert\cdot \vert_{\mathbb{A}}\colon \mathbb{A}\to \R_{\geq 0},\quad (x_{\infty},x_2,x_3,\ldots) \mapsto \vert x_{\infty} \vert_{\infty}\cdot \prod_p \vert x_p\vert_p.
\end{equation*}

We embed $\Q$ and $\Q^{\times}$ diagonally into $\mathbb{A}$ and $\mathbb{A}^{\times}$, respectively. 
The absolute values on $\Q_{\infty}=\R$ and $\Q_p$ are then normalised so that $\vert q\vert_{\mathbb{A}}=1$ for all $q\in \Q$. 
This matches the normalisation in Section~\ref{sec:p-adic} where we take $F=\Q_p$ and $\varpi=p$.

If $\psi_p$ is the standard additive character of $\Q_p$ and $\psi_{\infty}(x)=\exp(2 \pi i x)$, then we set
\begin{equation*}
    \psi_{\A}(x) = \psi_{\infty}(x_{\infty})\cdot \prod_p\psi_p(x_p).
\end{equation*}
Again we note that $\psi_{\A}(q)=1$ for $q\in \Q$, so that we can view $\psi_{\A}$ as a character on $\Q\backslash \A$. 

Finally, recall the notation
\begin{equation*}
    \widehat{\Z} = \prod_p\Z_p \quad \text{and} \quad \widehat{\Z}^{\times} = \prod_p\Z_p^{\times}.
\end{equation*}

We now turn towards global aspects of the general linear group. At the archimedean place, we define 
\begin{displaymath}
    G_{\infty} = \GL_2(\R), \quad G_{\infty}'=\PGL_2(\R), \quad G_{\infty}^+=\GL_2^+(\R).
\end{displaymath}
Similarly, we put 
\begin{displaymath}
    G_p=\GL_2(\Q_p).    
\end{displaymath}
All other subgroups from previous sections, such as $N, D, B$ and $Z$, are analogously equipped with subscripts to indicate the field over which consider them.

We set $K_{\infty}=\SO_2(\mathbb{R})$, $K_{\infty}'=\textrm{PSO}_2(\R)$ and $K_p=\GL_2(\Z_p)$. 
Put these together to we obtain a compact subgroup
\begin{equation*}
    K = K_{\infty} \times \prod_p K_p.
\end{equation*}
Recall that we equipped the local groups $K_{\infty}$ and $K_p$ with appropriately normalised Haar measures, so that we can define the product measure $\mu_K$ on $K$. 
In our normalisation we have $\Vol(K)=2\pi$. Let 
\begin{equation}
    N=p^{4n}  \text{ and }N_0=p^{2n}, \text{ for }n\in \N \text{ and }p\text{ odd prime.} \nonumber
\end{equation}
We will encounter the subgroup
\begin{equation*}
    K^{\ast}(N_0) = K_{\infty}\times K_{H,p}(2n,2n)\times \prod_{l\neq p} K_l,
\end{equation*}
where $K_{H,p}(m,m')$ denotes the subgroup of $K_p$ defined in \eqref{eq:defHecke} with $F=\Q_p$. It is easy to see that
\begin{align*}
    \Vol(K^{\ast}(N_0)) &= \frac{2\pi}{[K\colon K^{\ast}(N_0)]} \\ 
    &= \frac{2\pi}{N(1+p^{-1})} = 2\pi N^{-1+o(1)}.
\end{align*}
By $K_{\textrm{fin}}$ and $K^{\ast}(N_0)_{\textrm{fin}}$ we denote the finite parts of these groups.

The adele group $G(\A)$ is now defined as the restricted product (with respect to $\{K_p\}_{p}$)
\begin{equation*}
    G(\A) =G_{\infty}\times {\prod_p}' G_p.
\end{equation*}
The adelic points of all other groups mentioned are defined analogously. 
In particular, we have
\begin{equation*}
    N(\A)\cong \A,\quad D(\A)  \cong (\A^{\times})^2, \quad Z(\A)\cong \A^{\times}.
\end{equation*}
We identify $Z(\A)\backslash D(\A)$ with
\begin{equation*}
    \{ a(y) = \diag(y, 1) \colon y \in \A^{\times}\} \cong \A^\times.
\end{equation*}
These groups can thus be endowed with the measures coming from $\A$, $\A^{\times}$ and $(\A^{\times})^2$, respectively. 
By construction, these agree with the measures that are inherited from the restricted product structure. Similarly, the local Haar measures define product measures on $G(\A)$ and $G(\A)$, normalised so that
\begin{equation*}
    \int_{G(\A)} f(g) \, dg = \int_{Z(\A)} \int_{\A} \int_{\A^{\times}} \int_K f(zn(x)a(y)k) \, dk \, \frac{dy}{\vert y\vert_{\A}^2} \, dx \, dz.
\end{equation*}

\subsubsection{Adelic quotients}
We recall the approximation theorems
\begin{equation*}
    \A = \R\times \widehat{\Z}+\Q \quad \text{and} \quad \A^{\times} = \Q^{\times} \cdot \R^{\times} \cdot \widehat{\Z}^{\times}.
\end{equation*}
Note that this allows us to compute that
\begin{equation*}
    \Vol(\Q\backslash \A) = \Vol([0,1]\times \widehat{\Z})=1.
\end{equation*}

The strong approximation theorem (see \cite[Theorem~3.3.1]{bump}) for $\GL(2)$ implies that
\begin{equation*}
    G(\A) = G(\Q) \cdot G_{\infty}^+ \cdot K^{\ast}(N_0)_{\textrm{fin}},
\end{equation*}
using the fact that the determinant from $K^{\ast}(N_0)_{\textrm{fin}}$ to $\widehat{\Z}^{\times}$ is surjective. 
We can check that
\begin{displaymath}
    G(\Q)\cap G_{\infty}^+ \cdot K^{\ast}(N_0) = \Gamma_H(N_0,N_0),
\end{displaymath}
where 
\begin{displaymath}
     \Gamma_H(N_1,N_0) = \left[\begin{matrix} \Z & N_0\Z \\ N_1\Z & \Z\end{matrix} \right]\cap \SL_2(\Z).
\end{displaymath}
We therefore find that the quotient 
\begin{equation*}
    X_{N_0} = Z(\A)G(\Q)\backslash G(\A)/K^{\ast}(N_0) 
\end{equation*}
is isomorphic to $\Gamma_H(N_0,N_0)\backslash \mathbb{H}$. Finally, we write
\begin{displaymath}
    [G] = Z(\A) G(\Q) \backslash G(\A).
\end{displaymath}

\subsection{Automorphic forms and representations}\label{subsec:aut_fomrs}
A cuspidal Hecke--Maaß newform $\varphi$ on $\Gamma_H(N)\backslash \mathbb{H}$ can be translated to be $\Gamma_H(N_0,N_0)$ invariant and then be interpreted as a function on the adelic quotient $X_{N_0}$. This is due to the isomorphism with $\Gamma_H(N_0,N_0)\backslash \mathbb{H}$.
As in Chapter 5, Part C, of \cite{gelbart}, or Section 3 of \cite{cass}, or Section 3.6 of \cite{bump}, it generates an irreducible automorphic representation, say, $\pi \subset L^2([G])$.

The representation $\pi$ is isomorphic to $\pi_\infty \otimes \bigotimes_l \pi_l$ and $\pi_\infty \cong \Pc(s)$, where $\frac14 - s^2$ is the Laplace eigenvalue of $\varphi$. In the following, we will often make the simplifying assumption that $\pi_{\infty}$ is tempered (i.e. $\pi$ is tempered at infinity and $s\in i\R$). We write $s=iT$ and assume $T$ to be large. The local representations $\pi_l$ with $l\neq  p$ are class 1 and uniquely determined by the Hecke eigenvalue $\lambda_{\varphi}(l)$ of $\varphi$ at $l$. Pinning down the ramified local representation $\pi_p$ is harder. We recall from Remark~\ref{rem:conductor_class} that due to our restriction on central character and level we find the following options:

\begin{itemize}
    \item There is a character $\chi\colon \Q_p^{\times}\to S^1$ of (exponent) conductor $a(\chi)=2n$ and $\pi_p=\chi\boxplus \chi^{-1}$; or
    \item The representation $\pi_p$ is supercuspidal and associated to an unramified quadratic extension $E/\Q_p$ and a character $\theta\colon E^{\times}\to S^1$ of (exponent) conductor $a(\theta)=2n$.
\end{itemize}

\subsubsection{Uniqueness theorems} \label{sec:uniqueness}
Note that $\varphi$ is $K_0(N)$-invariant. On the tensor product side, this latter vector therefore lies in the tensor product of
\begin{displaymath}
    W_l = \{ w \in \pi_l \mid \pi_l(k) w = w \, \text{for all} \, k \in K_{H,l}(v_l(N),0) \}
\end{displaymath}
over all places $l$, tensored with $\pi_{\infty}^{K_{\infty}}$.

A crucial input in our approach to the sup-norm problem are the following uniqueness theorems.\footnote{This is not quite standard terminology. Depending on the situation, one speaks more generally of Gelfand pairs or multiplicity-one statements.}
If $v$ is the archimedean place or an unramified non-archimedean place, that is, a prime $l\neq p$, then $W_v$ is (at most) one-dimensional by the uniqueness of spherical vectors. We refer here to Theorem 2.4.2 and Theorem 4.6.2 in \cite{bump}.

At the ramified place $p$, we have the theory of newforms or new-vectors, as described in Section \ref{sec:new-vectors}. If $a(\pi_p)=4n$ the exponent conductor of $\pi_p$, then
\begin{displaymath}
    \{ w \in \pi_p \mid \pi_p(k) w = w \, \text{for all} \, k \in K_{H,p}({a(\pi_p)},0) \}
\end{displaymath}
is one-dimensional. This is due to newform theory, in the sense of Atkin and Lehner, as explained in Section 3 of \cite{cass}. We conclude that $W_p$ is also one-dimensional.

All in all, we see that the space of $K_0(N)$-invariant vectors in the representation $\pi$ is spanned by $\varphi$.
Furthermore, since $N=N_0^2$, then the space of $K^{\ast}(N_0)$-invariant vectors is spanned by 
\begin{equation}
	\varphi'=\pi(a(N_0))\varphi. \label{eq:shifted_new}
\end{equation}
This follows directly because $a(N_0)^{-1}K^{\ast}(N_0)_{\textrm{fin}}a(N_0) = K_0(N)$. 
As observed in \cite{marsh_local, Sa}, it is advantageous to work with the form $\varphi'$ instead of $\varphi$, more generally for some well-chosen factorisation of $N$, when bounding the sup-norm in the level aspect. We follow this approach below.

\section{The global period I}

We keep the notation of the previous section, but to shorten notation we set $H=K^{\ast}(N_0)$ and $Z_H = Z\cap H$. Note that $Z_H$ is compact and we have normalised measures so that $\Vol(Z_H)=2$. Define also $\Gamma=\PSL_2(\Z)$. 

Furthermore, we let $\pi$ be a cuspidal automorphic representation of $\GL_2(\Q)$ with trivial central character and level $N=p^{4n}$. We assume that $\pi$ is tempered and spherical at infinity. Its spectral parameter is denoted by $T$ and is assumed to be large. Let $\varphi'$ be the normalised shifted newform defined in \eqref{eq:shifted_new}.

In this section we make the first steps towards estimating the global sup-norm of $\varphi'$ (and thus also of our classical Hecke--Maa\ss\  newform $f$) using the local ideas developed earlier. The first step, carried out in Section~\ref{sec:integrals-matrix-coeffs} below, is to replace the rather singular period $\varphi'(g)$, given by evaluation at a single point, by a more flexible one, following the sketch in Section~\ref{sec:sketch-K-period}. 

Indeed, using uniqueness of the $K^{\ast}(N_0)$-invariant vector in $\pi$, we are able to write
\begin{equation}
	\mathcal{Q}_g(\Psi)\cdot \vert \varphi'(g)\vert^2 = \left\vert \frac{1}{\Vol(Z_H\backslash H)}\int_{Z_H\backslash gHg^{-1}} \Psi(h)dh \right\vert^2. \label{dummy}
\end{equation}
Here, $\mathcal{Q}_g(\Psi)$ is a weight that can be expressed as an integral of matrix coefficients --- see \eqref{eq:defQg} below for the precise definition. The key idea is now to choose $\Psi\in \pi$ to be a microlocalised vector. If the microlocal support of $\Psi$ is chosen carefully, one has good hope of obtaining satisfying lower bounds for $\mathcal{Q}_g(\Psi)$ and simultaneously good control on the $gHg^{-1}$-periods of $\Psi$. In practice, the right hand side of \eqref{dummy} is bounded by embedding $\Psi$ in a suitable spectral family and applying a relative trace formula. This is recorded in Proposition~\ref{pr:pre_amplified_guy} below, which is our main global work horse. 

We then use the local theory developed in Section~\ref{sec:archi-microlocalisation-big-chap} at the archimedean place and in Section~\ref{sec:p-adic} at the ramified place to choose suitable test functions and to gather important estimates. As a first result we establish a \emph{hybrid baseline (or local) bound}: see Theorem~\ref{th:local_bound} below.
Thus, in this section we rigorously execute all the steps sketched out in Section~\ref{sec:sketch_baseline}, including an analogous treatment of the level aspect.

\subsection{Integrals of global matrix coefficients} \label{sec:integrals-matrix-coeffs}

We now start as in Section \ref{sec:sketch-K-period} and roughly follow Section 3.5 in \cite{Nel-Un}.

First, the value $\varphi'(1)$ can be related to an integral of matrix coefficients. 
For this we define the (densely defined) linear functional $l: L^2([G]) \longrightarrow \C$ by
\begin{displaymath}
    l(\Psi) = \frac{1}{\Vol(H)}\int_{H} \Psi(h) \, dh = \frac{2}{\Vol(H)}\int_{Z_H\backslash H} \Psi(h) \, dh,
\end{displaymath}
for smooth elements $\Psi \in \mathcal{C}^{\infty}([G])$. 
It is straight-forward to check that
\begin{displaymath}
    l(\pi(h) \Psi) = l(\Psi) \text{ for }h\in H.
\end{displaymath}
Since $K_{\infty}\subseteq G_{\infty}$ has measure $0$, we can not define $l$ on the full Hilbert space $L^2([G])$ directly. However, for our purposes it is sufficient to extend $l\vert_{\pi}\colon \pi\cap \mathcal{C}^{\infty}([G])\to \C$ to the Hilbert (sub)space $\pi\subseteq L^2([G])$. To do so, we observe that $l$ can be written as the evaluation at the identity of a vector, that is
\begin{equation}
	l\vert_{\pi}(\Psi) = \left(\int_{H} \pi(k) \Psi dk\right)(1).\label{eq:observe}
\end{equation}
At this point we observe that the operator
\begin{equation}
	\Psi \mapsto \int_{H} \pi(k) \Psi dk \nonumber
\end{equation}
is a projection on subspace of $H$-invariant elements in $\pi$. Since the latter space is spanned by the $L^2$-normalized element $\varphi'$, see Section~\ref{sec:uniqueness}, we find that
\begin{equation}
	\int_{H} \pi(k) \Psi dk = \langle \Psi,\varphi'\rangle\cdot \varphi'. \nonumber
\end{equation}
Finally, note $\varphi'$ is smooth so that we can evaluate it at the identity.  In view of \eqref{eq:observe} we obtain
\begin{equation}
    l(\Psi) = \varphi'(1) \cdot \langle \Psi, \varphi' \rangle. \label{l_as_inner_prod}
\end{equation}

Next, decomposing into $H$-types and using orthogonality, we compute that
\begin{align}
    \mathcal{Q}(\Psi) = \frac{1}{\Vol(H)}\int_{H} \langle \pi(h) \Psi, \Psi \rangle \, dh 
    = \abs{\langle \Psi, \varphi' \rangle}^2. \label{Q_as_inner_prod}
\end{align}
One can compare this to the notation in Section 3.5 of \cite{Nel-Un}, where we could consider $\mathcal{Q}(\Psi \otimes 1)$, taking the trivial representation of $H$.\footnote{One could also take, for instance, other $\SO(2)$-types, generalising the present discussion to weight-$k$ Maaß forms.} 
We combine \eqref{l_as_inner_prod} and \eqref{Q_as_inner_prod} to obtain 
\begin{equation*}
    |l(\Psi)|^2 = |\varphi'(1)|^2 \cdot \mathcal{Q}(\Psi).
\end{equation*}
This is our replacement for \cite[(3.1)]{Nel-Un}.

The advantage of the integral expression of $\mathcal{Q}$ is that it factorises over all places and is not dependent on the global realisation of the vector $\Psi$.
Interpreting $\pi$ now as a tensor product and assuming $\Psi$ is a factorisable vector, we can write
\begin{equation*}
    \mathcal{Q}(\Psi)=Q_{\infty}(\Psi_{\infty}) \cdot \prod_{p} Q_p(\Psi_p),
\end{equation*}
where 
\begin{displaymath}
    Q_v(\Psi_v) = \frac{1}{\Vol(H_v)} \int_{H_v} \langle \pi(h_v) \Psi_v, \Psi_v \rangle \, dh_v.
\end{displaymath}
These local integrals of matrix coefficients are defined the same way in Section \ref{sec:rel-char-estimates}, at the real place, and in Section~\ref{sec:loc_perio_fin}, at finite places. 

To obtain more general values of $\varphi$, let $g \in G$ and define $H^g = gHg^{-1}$. %
As above, the vector $\pi(g) \varphi$ spans the space of $H^g$-invariant vectors in $\pi$. 
Defining 
\begin{displaymath}
    l_g(\Psi) = \frac{1}{\Vol(\bar{H}^g)}\int_{\bar{H}^g} \Psi(h) \, dh
\end{displaymath}
and
\begin{displaymath}
    \mathcal{Q}_g(\Psi) 
    = \frac{1}{\Vol(\bar{H}^g)}\int_{\bar{H}^g} \langle \pi(h) \Psi, \Psi \rangle \, dh, \label{eq:defQg}
\end{displaymath}
the same computations as above now show that
\begin{equation} \label{eq:lg}
    |l_g(\Psi)|^2 = |(\pi(g)\varphi')(1)|^2 \cdot \mathcal{Q}_g(\Psi) = |\varphi'(g)|^2 \cdot \mathcal{Q}_g(\Psi).
\end{equation}

\subsection{The pretrace inequality}

Let $\mathcal{B}(\pi)$ be any orthonormal basis for the representation $\pi$, viewed as a subspace of $L^2([G])$. For any test function $f \in C_c^\infty (G)$, we have a pretrace inequality
\begin{equation} \label{eq:pretrace-proto}
    \sum_{\Psi \in \Bc(\pi)} \left| \int_{H} \pi(f) \Psi \right|^2 
    \leq \int_{x,y \in H} \sum_{\gamma \in G'(\Q)} (f \ast f^\ast)(x^{-1} \gamma y) \, dx \, dy
\end{equation}
as given in \cite[Lemma 5.5]{Nel-Un} or \cite[Lemma 3.1]{marsh}. A nice way to derive this inequality is by an adaption of the argument presented in \cite[Section~3]{BHMM}. We apply this inequality with a different subgroup and test function for each $g \in G$.

Recall that we put $H^g = gHg^{-1}$ and define
\begin{displaymath}
    f_g(x) = f(g^{-1} x g)
\end{displaymath}
for any function $f$ on $G$.
By the definition of convolution operators and a change of variables, we see that
\begin{displaymath}
    \pi(f_g) = \pi(g) \pi(f) \pi(g^{-1}).
\end{displaymath}
Another simple application of the definition and a change of variables shows that 
\begin{displaymath}
    (f_g \ast f_g^\ast)(x) = (f \ast f^\ast)(g^{-1} x g).
\end{displaymath}

We apply now the inequality \eqref{eq:pretrace-proto} to $H^g$ and $f_g$ to obtain
\begin{align*}
    \sum_{\Psi \in \Bc(\pi)} |l_g(\pi(f_g) \cdot \Psi)|^2 
    & \leq \frac{1}{\Vol(\bar{H}^g)^2}\int_{x, y \in \bar{H}^g} \sum_{\gamma \in G'(\Q)} (f \ast f^\ast)(g^{-1} x^{-1} \gamma y g) \, dx \, dy \\
    & = \frac{1}{\Vol(\bar{H})^2}\int_{x, y \in \bar{H}} \sum_{\gamma \in G'(\Q)} (f \ast f^\ast)(x^{-1} g^{-1} \gamma g y) \, dx \, dy.
\end{align*}
To treat the left-hand side above, we change the basis $\Bc(\pi)$ by applying the unitary operator $\pi(g)$.
Thus, by independence of the choice of basis and the relation \eqref{eq:lg},
\begin{align*}
    \sum_{\Psi \in \Bc(\pi)} |l_g(\pi(f_g) \cdot \Psi)|^2 
    & = \sum_{\Psi \in \Bc(\pi)} |l_g(\pi(f_g)\pi(g) \cdot \Psi)|^2 \\
    & = |\varphi'(g)|^2 \cdot \sum_{\Psi \in \Bc(\pi)} \mathcal{Q}_g(\pi(f_g)\pi(g) \Psi).
\end{align*}
Unpacking the definition of $\mathcal{Q}_g$ and using the factorisation of $\pi(f_g)$ quickly shows that
\begin{displaymath}
    \mathcal{Q}_g(\pi(f^g)\pi(g) \Psi) = \mathcal{Q}(\pi(f) \Psi),
\end{displaymath}
where we recall that $\mathcal{Q} = \mathcal{Q}_1$. 
All in all, we deduce the main pretrace inequality
\begin{equation} \label{eq:pretrace-g}
    |\varphi'(g)|^2 \cdot \sum_{\Psi \in \Bc(\pi)} \mathcal{Q}(\pi(f) \Psi) 
    \leq\frac{1}{\Vol(H)^2} \int_{x, y \in H} \sum_{\gamma \in G'(\Q)} (f \ast f^\ast)(x^{-1} g^{-1} \gamma g y) \, dx \, dy.
\end{equation}

\subsection{The test function}\label{sec:test_fct}

We now fix our choice of test function $f$. 
This is done place by place, and we write
\begin{equation*}
    f=f_{\infty}\otimes f_p \otimes \bigotimes_{l\neq p} f_l.
\end{equation*}
We choose $f_{\infty}$ as in Theorem \ref{thm:archi-main-estimate}. Recall that $f_{\infty}$ is supported in a small neighbourhood of the identity, which we denote by $\mathfrak{U}$.
At the unramified places $l\neq p$, we simply put $f_l = \vvmathbb{1}_{ZK_l}$ 
(these places will be used later when we discuss amplification). 
Finally, at the ramified place $p$, we have to consider two scenarios:
\begin{itemize}
    \item Suppose $\pi_p=\chi\boxplus \chi^{-1}$ for $\chi\colon \Q_p^{\times}\to S^{1}$ with $a(\chi)=2n$. Then we put
    \begin{equation}
        f_p(g) = \dim(V_{\textrm{new},p})\cdot \frac{\langle  v_{p},\pi(g)v_p\rangle}{\Vert v_p\Vert^2}\cdot \vvmathbb{1}_{ZK_p}(g),\nonumber
    \end{equation}
    where $v_p=\pi_p(c_1)v_\ML$ with $c_1$ as in \eqref{eq:matrix-b-chi} and $v_{ML}$ being the standard microlocal lift vector.
    \item Suppose $\pi_p$ is supercuspidal with $a(\pi)=4n$. Then we define
    \begin{equation*}
        f_p(g) = \dim(V_{\textrm{new},p})\cdot \frac{\langle  v_{p},\pi(g)v_p\rangle}{\Vert v_p\Vert^2}\cdot \vvmathbb{1}_{ZK_p}(g)
    \end{equation*}
    where $v_p=v_{MV}$ is the standard minimal vector. Note that here the matrix coefficient is automatically compactly supported, so that no cut-off is necessary.
\end{itemize}

With this choice we can compute that
\begin{equation*}
    \sum_{\phi \in \Bc(\pi)} \mathcal{Q}(\pi(f) \phi) = \left(\sum_{v\in \Bc(\pi_{\infty})} \mathcal{Q}_{\infty}(\pi(f_{\infty}) v )\right)\cdot \mathcal{Q}_{p}( v_p)
\end{equation*}
This follows from Lemma~\ref{p-adic_proj} and a convenient choice of (factorisable) basis $\Bc(\pi)$. Now we can use our local results to evaluate 
\begin{equation*}
    \mathcal{Q}_{p}(v_p)= \zeta_p(1)p^{-n}.
\end{equation*}
Here we used Lemma~\ref{lm:summary_Q}. %
Finally, recall the definition of the integral $I_{\gamma}(f_p\ast f_p)$ given in \eqref{eq:def_Igamma}. 

Next, we observe that the support of $f$ is in $Z(\A)\cdot (\mathfrak{U}\times K_{\fin})$. Very conveniently, we have $$Z(\A)\backslash [Z(\A)\cdot (\mathfrak{U}\times K_{\fin})] \cap G'(\Q) = \Gamma\cap \mathfrak{U}.$$ This is a consequence of strong approximation.x

We summarise this discussion in the following proposition:
\begin{prop}\label{pr:pre_amplified_guy}
Suppose that $g_{\fin} \in K_{\fin}$, Then, with the choice of test function $f$ exhibited above, we have
\begin{multline*}
    |\varphi'(g)|^2 \cdot \frac{\zeta_p(1)}{p^n}\sum_{v \in \Bc(\pi_{\infty})} \mathcal{Q}(\pi(f_{\infty}) v)  \\ \ll  \sum_{\gamma \in \Gamma\cap g\mathfrak{U}g^{-1}}I_{g^{-1}\gamma g}(f_p\times f_p)\cdot  \int_{x,y\in K_{\infty}'}(f \ast f^\ast)(x^{-1} g^{-1} \gamma g y) \, dx \, dy.
\end{multline*}
\end{prop}

\begin{rem}\label{rm:support_shrink}
It will later be useful to shrink the support of the test function $f_p$ a little bit. Indeed, we define
\begin{equation}
	f_p'=f_p\cdot \vvmathbb{1}_{ZK_H(1)}.\nonumber 
\end{equation}
Note that we still have 
\begin{equation}
	f_p'\ast f_p' = \Vol(K_H(1))\cdot f_p'.\nonumber
\end{equation}
On the other hand we can choose the orthonormal basis $\mathcal{B}(\pi_p)$ of $\pi_p$ such that it contains a finite subset  $\mathcal{B}_p$ so that
\begin{equation}
	\pi(f_p)w = \Vol(K_H(1)) w.
\end{equation}
See \cite[Corollary~2.16]{Sa} for a similar argument. In particular $v_p\in \mathcal{B}_p$.
\end{rem}

\subsection{The local bound}

At this point we can already discuss how to establish the local bound, i.e. the baseline bound. 

\begin{theorem}[Local bound]\label{th:local_bound}
Let $\varphi'$ be as in \eqref{eq:shifted_new} and fix a compact set $\Omega$. Then, for $\epsilon>0$, we have
\begin{equation}
    \Vert \varphi'\vert_{\Omega}\Vert_{\infty}^2 \ll_{\Omega,p,\epsilon}  (N_0T)^{\frac{1}{2}+\epsilon}.\nonumber
\end{equation}
\end{theorem}
\begin{proof}
By strong approximation we can assume that $\Omega\subset G_{\infty}^+$. Our starting point is Proposition~\ref{pr:pre_amplified_guy}. Note that if $\gamma=1$, then 
\begin{equation}
	\int_{x,y\in K_{\infty}'}(f \ast f^\ast)(x^{-1} g^{-1} \gamma g y) \, dx \, dy = \Vol(K_{\infty}')\cdot \int_{K_{\infty}'} (f \ast f^\ast)(y) dy \ll T^{\frac{1}{2}+\epsilon}, \nonumber
\end{equation}
where we have used the second estimate in Theorem~\ref{thm:archi-main-estimate}. Making use of the remaining two estimates of this theorem yields
\begin{equation}
	|\varphi'(g)|^2 \cdot \frac{\zeta_p(1)}{p^nT^{\frac{1}{2}}}  \ll_{\epsilon} T^{\frac{1}{2}+\epsilon}I_1(f_p\times f_p) + T^{\epsilon}\sum_{1\neq \gamma \in \Gamma\cap \mathfrak{U}}I_{g^{-1}\gamma g}(f_p\times f_p)\cdot  d(g^{-1}\gamma g)^{-1}. \nonumber
\end{equation}
Using \eqref{eq:123es_kommt_die_sonne} together with the volume bounds from Theorem~\ref{th:local_vol_bound} yields the estimate
\begin{equation}
	|\varphi'(g)|^2 \ll_{p,\epsilon} T^{1+\epsilon}p^{2n} + T^{\frac{1}{2}+\epsilon}p^{2n}\sum_{1\neq \gamma \in \Gamma\cap \mathfrak{U}}  d(g^{-1}\gamma g)^{-1}. \nonumber
\end{equation} 
Finally, we trivially estimate
\begin{equation}
	\sum_{1\neq \gamma \in \Gamma\cap \mathfrak{U}}  d(g^{-1}\gamma g)^{-1} \ll_{\mathfrak{U},\Omega} 1 \label{eq:trivial_count}
\end{equation}
and recall that $N_0=p^{2n}$ to obtain the bound
\begin{equation}
	|\varphi'(g)|^2 \ll_{p,\epsilon} T^{1+\epsilon}p^{2n}. \nonumber
\end{equation} 
Taking the square root gives the desired result.
\end{proof}

\begin{rem}
This local bound can be achieved in a variety of ways. For example, using the generating set for $X_{N_0}$ from \cite{Sa}, one gets a global (i.e. without restriction to $\Omega$) version of Theorem~\ref{th:local_bound} by estimating the Fourier expansion. See also Proposition~\ref{prop:Fourier} below. To the best of our knowledge, a local bound as above first appeared in the work of Marshall \cite{marsh_local}.
\end{rem}

Note that our proof only uses the most basic volume estimates. In particular, at the ramified place $p$ we have not used anything but the trivial bound. Doing so we have lost a lot of valuable information on the geometric side. In the next sections we will upgrade the simple pretrace inequality used here to an amplified pretrace inequality. This will allow us to exploit the full strength of our volume bounds. Beside this, we will choose a suitable generating domain for the quotient $[G]$ which allows us to address the dependence on the compact set $\Omega$.

\section{Amplification and counting}\label{sec:amp}

This section is an intermezzo in which we prepare some further tools that are necessary before we can return to the study of the global period in Section~\ref{sec:glob_perII}. Our goals here are twofold.

First, in Section~\ref{se:amplifier} we construct a refined test function at the unramified places, which will complement the choice made at the places $\infty$ and $p$ in Section~\ref{sec:test_fct}. This unramified test function is essentially a cleverly chosen combination of Hecke operators and often referred to as an \emph{amplifier}. 

Second, recall that in order to establish the baseline bound in the proof of Theorem~\ref{th:local_bound}, we have estimated the geometric side using the essentially trivial count \eqref{eq:trivial_count}. Such an naive approach is not sufficient if one wants to improve on the baseline bound. We will provide sufficiently uniform counting results in Section~\ref{sec:count} below. See Corollary~\ref{cor:summary_count} for a packaging of these that will be used later on.

As in the previous section, we are working with a cuspidal automorphic representation $\pi=\bigotimes_{v}\pi_v$. Let $N$ denote the level of $\pi$ and let $\varphi'$ be as in \eqref{eq:shifted_new}.

\subsection{Hecke operators and the amplifier}\label{se:amplifier}

For each $m\in \N$ with $(m,N)=1$ we can define a test function $\kappa_m$ such that $R(\kappa_m)$ is a Hecke operator. In particular, we have
\begin{equation}
	R(\kappa_m)\varphi' = \lambda_{\pi}(m)\varphi'.\nonumber
\end{equation}   
Here $\lambda_{\pi}(m)$ are Hecke eigenvalues of $\pi$ (resp. of $\varphi$ or $\varphi'$), which we normalise such that $\vert \lambda_{\pi}(m)\vert \leq \tau(m)$ under the Ramanujan conjecture.

We define the functions $\kappa_m$ following \cite[Section~3.5]{Sa}. First we factor them as
\begin{equation}
	\kappa_m = \prod_{(l,N)=1}\kappa_{l,v_l(m)} \text{ for }\kappa_{l,v_l(m)}\in \mathcal{C}_c^{\infty}(Z\backslash G(\Q_p)).\nonumber
\end{equation}
The local test functions are now defined by 
\begin{equation}
	\kappa_{l,r}(g) = \begin{cases}
		l^{-\frac{r}{2}} & \text{ if }g\in ZK_la(l^r)K_l, \\
		0 &\text{ else.}
	\end{cases}\nonumber
\end{equation}
Note that these are bi-$K_l$-invariant and $\kappa_{l,0}$ is simply the characteristic function on $ZK_l$. For later use we record the Hecke relations
\begin{equation}
	\kappa_r\ast \kappa_s = \sum_{t\mid (r,s)} \kappa_{\frac{rs}{t^2}} \text{ for }(rs,N)=1.\label{eq:hecke_rel}
\end{equation}
See for example \cite[Proposition~4.6.4]{bump} for details. The Hecke eigenvalues inherit this recursion property. In particular, they satisfy $\lambda(l)^2=\lambda(l^2)+1$ for all primes $l\nmid N$. This allows us to deduce the important inequality
\begin{equation}
	\vert \lambda(l)\vert + \vert \lambda(l^2)\vert > \frac{1}{2}. \label{eq:key_amp}
\end{equation}

Define the weights
\begin{equation}
	x_m = \frac{\overline{\lambda_{\pi}}(m)}{\vert \lambda_{\pi}(m)\vert}.\nonumber
\end{equation}
For a parameter $X$, fix the finite set of primes $\mathcal{P}_X \subseteq \{l\in [X,2X]\colon (l,N)=1\}$ and define the test function
\begin{equation}
	f_{\textrm{ur}}^X = \left(\sum_{l\in \mathcal{P}_X} x_l\kappa_l\right)\ast \left(\sum_{l\in \mathcal{P}_X} x_l\kappa_l\right) + \left(\sum_{l\in \mathcal{P}_X} x_{l^2}\kappa_{l^2}\right)\ast \left(\sum_{l\in \mathcal{P}_X} x_{l^2}\kappa_{l^2}\right). \label{eq:def_test_ur}
\end{equation}
This version of the amplifier is taken from \cite[(9.14)]{BHMM}.\footnote{This amplifier is state of the art. For our purposes more rudimentary amplifiers would be sufficient, but this one is still very convenient to work with.} This function is constructed as a sum of self-convolutions, so that $R(f_{\textrm{ur}}^X)$ is a positive operator.
It will be useful to write
\begin{equation}
	f_{\textrm{ur}}^X = f_{\textrm{ur},1}^X\ast f_{\textrm{ur},1}^X+ f_{\textrm{ur},2}^X\ast f_{\textrm{ur},2}^X
\end{equation}
for the obvious functions $f_{\textrm{ur},i}^X$.

We find that
\begin{equation}
	R(f_{\textrm{ur}}^X)\varphi' = C_X(\pi)\cdot \varphi'. \nonumber
\end{equation}
with
\begin{align}
	C_X(\pi) &= \left(\sum_{l\in \mathcal{P}_X} \vert \lambda_{\pi}(l)\vert\right)^2+\left(\sum_{l\in \mathcal{P}_X} \vert \lambda_{\pi}(l^2)\vert\right)^2 \nonumber\\
	&\geq \frac{1}{2} \left(\sum_{l\in \mathcal{P}_X}\vert \lambda(l)\vert+\vert \lambda(l^2)\vert \right)^2> \frac{1}{8}(\sharp \mathcal{P}_X)^2. \label{eq:basic_ampli_inequality}
\end{align}
In the last step we used \eqref{eq:key_amp}. See also \cite[(9.17)]{BHMM} for the same chain of inequalities. 
On the other hand, we can write
\begin{equation}
	f_{\textrm{ur}}^X = \sum_{m} a_m \cdot  \kappa_m \label{eq:expansion_hecke_rel}
\end{equation}
using the Hecke relations \eqref{eq:hecke_rel}. More precisely we only need to compute
\begin{align}
	f_{\textrm{ur},1}^X\ast f_{\textrm{ur},1}^X &= \sum_{l_1,l_2\in \mathcal{P}_X} x_{l_1}x_{l_2}\cdot \kappa_{l_1}\ast \kappa_{l_2} = \sum_{l_1,l_2\in \mathcal{P}_X}x_{l_1}x_{l_2}\cdot (\kappa_{l_1l_2}+\delta_{l_1=l_2}\kappa_1) \nonumber
\end{align}
and 
\begin{align}
	f_{\textrm{ur},2}^X\ast f_{\textrm{ur},2}^X &= \sum_{l_1,l_2\in \mathcal{P}_X} x_{l_1^2}x_{l_2^2}\cdot\kappa_{l_1^2}\ast \kappa_{l_2^2} = \sum_{l_1,l_2\in \mathcal{P}_X}x_{l_1^2}x_{l_2^2}\cdot(\kappa_{l_1^2l_2^2}+\delta_{l_1=l_2}(\kappa_{l_1^2}+\kappa_1)). \nonumber
\end{align}
Adding these two expressions yields that the coefficients $a_m$ in \eqref{eq:expansion_hecke_rel} are given by
\begin{equation}
	a_m=\begin{cases}
		x_{l_1^2}x_{l_2^2} & \text{ if }m=l_1^2l_2^2 \text{ for }l_1,l_2\in \mathcal{P}_X,\\
		x_{l_1}x_{l_2}\cdot (1+\delta_{m=\square}) &\text{ if }m=l_1l_2 \text{ for }l_1,l_2\in \mathcal{P}_X,\\
		2\cdot\sharp\mathcal{P}_X &\text{ if }m=1,\\
		0 &\text{ else.}
	\end{cases} \label{bounds_am}
\end{equation} 
For our purposes, it is sufficient to keep in mind that $a_1 \asymp \sharp \mathcal{P}_X$, $\vert a_m\vert \asymp 1$ if $m=l_1^il_2^i$ for $l_1,l_2\in \mathcal{P}_X$ and $i=1,2$.

Note that any one of the functions $\kappa_m$ are bi-$K_l$ invariant for all finite places $l\nmid N$. In particular, for any $\Psi\in \pi$, the elements $R(\kappa_m)\Psi$ are spherical at all finite places $(l,N)=1$. Thus, when working with test functions that are combinations of the $\kappa_m$, it is sufficient to work with the corresponding subspace of $\pi$ consisting of elements that are unramified away from $N$. We will denote this subspace by $\pi_N$. (Note that we can identify $\pi_N$ with $\pi_{\infty} \otimes \bigotimes_{p\mid N}\pi_p$). Let $\mathcal{B}^N(\pi)$ be a basis of $\pi_N$. We observe that for all $\Psi\in \mathcal{B}^N(\pi)$ we have 
\begin{equation}
	R(\kappa_m)\Psi = \lambda_{\pi}(m)\Psi.\nonumber
\end{equation}
This justifies the notation $ \lambda_{\pi}(m)$. We will always assume that $\mathcal{B}^N(\pi)\subseteq \mathcal{B}(\pi)$.

\begin{rem}
We have only discussed the \textit{amplifying} effect that $f_{\textrm{ur}}^X$ has on $\varphi'$. This is sufficient because we can drop the contribution of all other terms on the spectral side by positivity. It is however instructive to consider other functions. Indeed, suppose $\pi'$ is some other cuspidal automorphic representation containing a non-trivial $H$-invariant element $\psi'$. Then we obtain
\begin{equation}
	R(f_{\textrm{ur}}^X)\psi' = C_X(\pi')\cdot \psi' \nonumber
\end{equation}
for 
\begin{equation}
	C_X(\pi') = \left(\sum_{l\in \mathcal{P}_X} \sgn(\lambda_{\pi}(l))\cdot \lambda_{\pi'}(l)\right)^2+\left(\sum_{l\in \mathcal{P}_X} \sgn(\lambda_{\pi}(l^2))\cdot \lambda_{\pi'}(l^2)\right)^2.\nonumber
\end{equation}
Since $\pi$ and $\pi'$ are distinct their Hecke eigenvalues are expected to be uncorrelated and we expect that 
\begin{equation}
	C_X(\pi') = o((\sharp \mathcal{P}_X)^2).\nonumber
\end{equation} 
While showing this in practice is beyond our capabilities, it shows that at least morally $R(f_{\textrm{ur}}^X)$ dampens the contribution of representations besides $\pi$ in the relative trace formula.
\end{rem}

\subsection{Analysing the support condition} \label{sec:analysing-supp-condition}

Throughout this section we let $N=p^{4n}$ and let $g_{\infty}\in G(\R)$ and $g_p\in K_p$. We will need to control the size of the sets
\begin{multline}
	\mathfrak{M}(g_{\infty},g_p;m,\delta,l) = \{ \gamma\in G'(\Q) \colon \kappa_m(\gamma)\neq 0,\, d_{H(2n)}(g_p^{-1}\gamma g_p) \leq p^{-l} \\ \text{ and }d(g_{\infty}^{-1}\gamma g_{\infty})\leq \delta\}.\label{eq:def_mathfrakM}
\end{multline}
We want to find representatives for this set which allow us to run more classical counting arguments.

Recall the action of $\SL_2(\R)$ on the upper half-plane $\mathbb{H} \subset \C$ (see e.g. \cite[Chap.~1]{iwaniec-spectral}). 
The latter inherits the invariant distance function by setting $d(z,i) = d(g_z)$ for $z \in \mathbb{H}$ and a matrix $g_z \in \SL_2(\R)$ such that $g_z.i = z$, where $i$ is the imaginary unit.
By invariance, we set $d(z, w) = d(g_z^{-1}w, i)$.
We then define the point-pair invariant 
\begin{displaymath}
    u(z,w) = \frac{\abs{z-w}^2}{\Im z \cdot \Im w}.
\end{displaymath}
As in (1.3) in \cite[Chap.~1]{iwaniec-spectral}, we have
\begin{displaymath}
    \cosh d(z,w) = 1 + \frac12 \cdot u(z,w).
\end{displaymath}
This implies the bound
\begin{equation}\label{eq:distances-translation}
    d(z,w)^2 \asymp u(z,w) 
\end{equation}
for small $d(z,w)$.

The following lemma is essentially \cite[Lemma~3.19]{Sa} and translates between the adelic and classical setting.
\begin{lemmy}\label{lm:ana_support}
Let $l\geq 0$ and let $\delta>0$ be sufficiently small.
Let $g_p\in \SL_2(\Z)$, $z=g_{\infty}.i$ and, for a positive integer $L$, define
\begin{align*}
	M(z,g_p;\delta', L, m) = \sharp \{ \gamma\in \Mat_2(\Z) \colon & u(\gamma.z,z)\leq \delta',\, L \mid [g_p^{-1}\gamma g_p]_{12}, \\
    & L \mid [g_p^{-1}\gamma g_p]_{21} \text{ and }\det(\gamma)= m\}.
\end{align*}
Then we have
\begin{equation*}
	\sharp \mathfrak{M}(g_{\infty},g_p;m,\delta,l) \leq M(z, g_p; \delta', p^l, m),
\end{equation*}
where $\delta^{2}\asymp \delta'$.
\end{lemmy}
\begin{proof}
Suppose $\gamma\in \mathfrak{M}(g_{\infty},g_p;m,\delta,l)$. We first observe that
\begin{equation}
	\gamma \in ZK_p.\nonumber
\end{equation}
This is because $l\geq 0$ and $g_p\in K_p$, so that the condition $d_{H(2n)}(g_p^{-1}\gamma g_p)$ in particular implies that $\gamma\in ZK_p$. At places $l\neq p$ we use the support of $\gamma_m$ to obtain that
\begin{equation}
	\gamma\in  Z K_l a(l^{v_l(m)})K_l.\nonumber
\end{equation}
Thus we find $\xi \in Z(\Q)$ such that $\gamma = \xi \gamma'$ with $\gamma'\in \textrm{Mat}_{2\times 2}(\Z)$ and $\det(\gamma')=m$. 
Next, the distance condition follows now from \eqref{eq:distances-translation}.

Finally, we have to take the condition $d_{H(2n)}(g_p^{-1}\gamma g_p)=p^{-l}$ into account. We relax this condition to 
\begin{equation}
	g_p^{-1}\gamma'g_p \in K_{H,p}(l,l). \nonumber
\end{equation}
Under the assumption $g_p\in \textrm{SL}_2(\Z)$, this reduces to 
\begin{equation}
	g_p^{-1}\gamma' g_p = \left( \begin{matrix} a & b \\ c & d \end{matrix} \right) \in\textrm{Mat}_2(\Z) \nonumber
\end{equation}
with $p^l\mid (c,b)$.
\end{proof}

\begin{rem}\label{rm:intro_Mop}
In the volume bounds Corollary~\ref{cor:vol_bound_p-adig} we cannot distinguish between $d_{H(2n)}(g_p^{-1}\gamma g_p)=p^l$ and $d_{H(2n)}(wg_p^{-1}\gamma g_p)=p^l$, where $w$ is the long Weyl element. Thus, a priori one would also encounter the following sets of matrices:
\begin{multline}
	\mathfrak{M}^{\textrm{op}}(g_{\infty},g_p;m,\delta,l) = \{ \gamma\in G'(\Q) \colon \kappa_m(\gamma)\neq 0,\, d_{H(2n)}(w g_p^{-1}\gamma g_p) \leq  p^{-l} \\ \text{ and }d(g_{\infty}^{-1}\gamma g_{\infty})\leq \delta\}.\label{eq:def_mathfrakMop}
\end{multline}
Essentially re-producing the arguments given above one finds that
\begin{equation}
	\sharp \mathfrak{M}^{\textrm{op}}(g_{\infty},g_p;m,\delta,l) \leq 	M^{\textrm{op}}(z,g_p;\delta', L, m),\nonumber
\end{equation}
where $\delta^{2}\asymp \delta'$ sufficiently small and 
\begin{multline}
	M^{\textrm{op}}(z,g_p;\delta', L, m) = \sharp \{ \gamma\in \textrm{Mat}_2(\Z)\colon u(\gamma.z,z)\leq \delta', \\ 
    L\mid [g_p^{-1}\gamma g_p]_{ii} \text{ for } i=1,2 \text{ and }\det(\gamma)= m\}.\nonumber
\end{multline}
These matrices seem hard to count. 
To circumvent this issue, we later on employ the support shrinking technique, as in Remark \ref{rm:support_shrink}, to avoid the sets $\mathfrak{M}^{\textrm{op}}(g_{\infty},g_p;m,\delta,l)$ altogether.
\end{rem}

\subsubsection{Useful inequalities and congruences}\label{sec:useful_ineq}
To aid in counting the matrices $\gamma$, we observe the following bounds on matrix entries.
These can also be found in \cite[Sec.~10.1]{BHMM}, but we give here slightly more conceptual arguments.

Firstly, assume that $\gamma \in \Mat_2(\Z)$ with $\det \gamma = m > 0$ and set $\tilde{\gamma} = m^{-1/2} \gamma \in \SL_2(\R)$. 
Then the hyperbolic distance $d(g_\infty^{-1} \gamma g_\infty) = \delta$ can be given using the Cartan decomposition (or geodesic polar coordinates, see \cite[Sec.~1.3]{iwaniec-spectral}).
More precisely, there are $k',k''\in \SO(2)$ such that $g_\infty^{-1} \tilde{\gamma} g_\infty$ is equal to $k' \diag(e^{-\delta/2}, e^{\delta/2}) k''$.
Expanding the exponential into a series, if $\delta < 1$, we find there is $k \in \SO(2)$ such that 
\begin{equation} \label{eq:condition-at-infty}
    g_\infty^{-1} \gamma g_\infty = \sqrt{m}(k + O(\delta)).
\end{equation}

Write now
\begin{displaymath}
    \gamma = \begin{pmatrix}
        a & b \\ c & d
    \end{pmatrix}.
\end{displaymath}
Applying the conjugation invariant trace to equation \ref{eq:condition-at-infty}, the former being conjugation invariant, we obtain that 
\begin{equation}
    a+d \ll \sqrt{m}, \label{eq:trace_ineq}
\end{equation} 
using that orthogonal matrices are bounded.
Next, using the Iwasawa decomposition, we can write
\begin{displaymath}
    g_\infty = \begin{pmatrix}
        y & x \\ 0 & 1
    \end{pmatrix}.
\end{displaymath}
Performing the conjugation shows that the lower-left entry of $g_\infty^{-1} \gamma g_\infty$ is simply $c y$.
Thus, the condition implies that 
\begin{equation}
    cy \ll \sqrt{m}. \label{eq:bound_c}
\end{equation}

Finally, the condition $u(\gamma z, z) < \delta'$ gives the bound
\begin{equation} \label{eq:u-condition}
    \abs{-c z^2 + (a-d) z + b}^2 \ll m \delta' y^2,
\end{equation}
where $z = g_\infty.i \in \mathbb{H}$, as in \cite[(1.10)]{Te}, for example.

On the non-archimedean side of the medal, these inequalities are replaced by congruences coming from the condition
\begin{equation}
	p^l\mid [g_p^{-1}\gamma g_p]_{ij} \text{ for }\{i,j\} = \{ 1,2\}, \nonumber
\end{equation}
which are satisfied for all $\gamma$ that contribute to the count $M(z,g_p;\delta';p^l,m)$.  Let\footnote{Observe that the lower right entry of $g_p$ is denoted by $w$, which is the symbol usually reserved for the long Weyl element. We hope this double usage of the letter $w$ does not cause any confusion.}
\begin{equation}
	g_p=\left(\begin{matrix} t & u \\ v & w\end{matrix}\right)\in \SL_2(\Z).\nonumber
\end{equation}
Thus we have $tw-uv=1$. We note that, since $wK_{H,p}(l,l)w^{-1}=K_{H,p}(l,l)$ we can always replace $g_p$ by $g_pw$.
In particular, we can assume that $(v,p)=1$. 
Put $s=v_p(w)$.\footnote{We have $g_p^{-1}\gamma g_p\in K_{H,p}(l,l)$ if and only if $(g_pw)^{-1}\gamma (g_pw)\in K_{H,p}(l,l)$. 
However, $$g_pw = \left(\begin{matrix}-u & t \\ -w & v\end{matrix}\right).$$ Since $tw-uv=1$ we can not have $v_p(w)>0$ and $v_p(v)>0$ simultaneously.}\textsuperscript{,\,}\footnote{We remark in passing that this reduction step is reminiscent of the Atkin-Lehner operators applied in the treatment of the square-free level case, e.g. \cite{HT}. Indeed, the normalising Weyl element $w$ corresponds to the Fricke involution in our balanced group.}  

Recall that we want to count $M(z,g_p;\delta',p^l,m)$ with $\delta'<1$ and $l>0$. (The case $l=0$ is straight forward and follows directly from the counting results of \cite{Te}.) We write a matrix $\gamma$ contributing to this count as
\begin{equation}
	\gamma= \left( \begin{matrix} a & b \\ c & d\end{matrix}\right)\in \Mat_2(\Z). \nonumber
\end{equation}
We compute
\begin{equation} \label{eq:p-adic-conjugation-condition}
	g_p^{-1}\gamma g_p = \left(\begin{matrix} \ast & (a-d)uw+bw^2-cu^2 \\ -(a-d)vt-bv^2+ct^2 & \ast \end{matrix}\right).
\end{equation}
Setting $$A=a-d,$$ we obtain the congruence conditions
\begin{align}
	Auw+bw^2-cu^2 &\equiv 0 \text{ mod }p^l \nonumber\\
	-Avt-bv^2+ct^2 &\equiv 0 \text{ mod }p^l.\nonumber
\end{align}
We put
\begin{equation}
	T_g = \left(\begin{matrix} u & w \\ -t & -v \end{matrix}\right) \text{ and } v_g= \left( \begin{matrix} u^2 \\ -t^2\end{matrix}\right).\nonumber
\end{equation}
We can write the congruence conditions as
\begin{equation}
	\left(\begin{matrix} w & 0 \\ 0 & v\end{matrix}\right)T_g\left(\begin{matrix} A \\ b \end{matrix}\right) \equiv c\cdot v_g\text{ mod }p^l.\nonumber
\end{equation}
We note that $\det(T_g)=1$ and 
\begin{equation}
	T_g^{-1} = \left(\begin{matrix} -v & -w \\ t & u\end{matrix}\right). \nonumber %
\end{equation}
We can further reduce this by writing $w=p^s w_0$ and $v=p^rv_0$ with $\gcd(p,v_0 w_0)=1$. Let $\overline{w_0}$ (resp. $\overline{v}_0$) be a representative of the modular inverse of $w_0$ (resp. $v_0$) modulo $p^n$. (Here we anticipate that later we will only need $l\leq n$.) Then we further reduce to
\begin{equation} \label{eq:reduced_congruence}
	\left(\begin{matrix} p^s & 0 \\ 0 & p^r\end{matrix}\right)T_g\left(\begin{matrix} A \\ b\end{matrix}\right)= c\left(\begin{matrix} \overline{w_0} & 0 \\ 0& \overline{v_0}\end{matrix}\right)v_g + p^l\Z^2 = \left(\begin{matrix} c\cdot u^2\overline{w_0} \\ -c\cdot t^2 \overline{v_0}\end{matrix}\right)+p^l\Z^2.
\end{equation}
As mentioned above we can later reduce to the case where $r=0$, so that $v=v_0$ is co-prime to $p$.

\subsection{Counting results}\label{sec:count}

It is now key to have uniform estimates for $M(z,g_p;\delta',L,m)$. 
The easiest results can be obtained when only focusing on $\delta'$. 
Here we can drop the constraints coming from $g_p$ and $L$ altogether. 
In particular, it will be sufficient to work with $z$ in the standard fundamental domain for $\SL_2(\Z)$, which we denote by $\mathcal{F}$. 
We import the following counting result from \cite[Appendix~1]{IS95}.

\begin{lemmy} \label{lm:IS_counting}
Suppose $z\in \mathcal{F}$, $g_p\in \SL_2(\Z)$ and assume that $\delta'\leq 1$. Then we have
\begin{equation}
	M(z,g_p;\delta', L, m) \ll m^{\epsilon}(1+m(\delta')^{\frac{1}{4}}+\sqrt{m\delta'}y).\nonumber
\end{equation}
The same bound holds for $M^{\textrm{op}}(z,g_p;\delta', L, m)$.
\end{lemmy}

Note that the implicit constant is independent of $L$ and $g_p$, since we are simply forgetting the corresponding conditions. 
However, we would expect and, indeed, require some saving in $L$. 
To see this saving, we present a slight adaptation of the Harcos-Templier \cite{HT} counting argument. 
It is based on the following simple, yet powerful lattice point estimate.

\begin{lemmy}
    Let $B_w(R)\subseteq \R^2$ be a ball of radius $R$ centred at $w \in \R^2$.
    Next, let $\Lambda \subset \R^2$ be a lattice of full rank, with first successive minimum $\lambda_1$ and covolume $V$.
    Then we have 
    \begin{equation} \label{eq:lat-pt-count}
        \sharp \Lambda \cap B_w(R) \ll 1 + \frac{R}{\lambda_1}+ \frac{R^2}{V}.
    \end{equation}
    Furthermore, if $\Lambda^\ast$ is the set of primitive lattice points, i.e. not of the form $k \cdot v$ for $k \in \Z_{\geq 2}$ and $v \in \Lambda$, then
    \begin{equation} \label{eq:lat-pt-count-primitive}
        \sharp \Lambda^\ast \cap B_0(R) \ll 1 + \frac{R^2}{V}.
    \end{equation}
\end{lemmy}
\begin{proof}
    This is well-known and, for convenience, we refer to \cite[Lemma 1]{BHM}, which covers the first assertion.
    To prove the second assertion we easily modify the argument as follows.

    Assume there exists $0 \neq v \in \Lambda^\ast \cap B_0(R)$.
    Then the $\R$-span of $\Lambda^\ast \cap B_0(R)$ is one- or two-dimensional.
    In the first case, there are only 2 primitive points.
    In the latter, we follow the argument of \cite{BHM} using translates of a fundamental parallelipiped, whose volume is $V$.
\end{proof}

Given $z = x + iy \in \Hb$, our analysis revolves around the lattice
\begin{equation}
	\Lambda_{z} = z\cdot \Z+\Z \subseteq \C\nonumber
\end{equation}
and its sublattices given by congruence conditions. 
Note that the covolume of  $\Lambda_{z}$ is given by
\begin{equation}
	\det(\Lambda_{z}) = y.\nonumber
\end{equation}
The first successive minimum, i.e. the shortest non-zero length, of $\Lambda_{z}$ is denoted by $\lambda_{z}$.
Reduction theory is classically done so that, if $z$ lies in the standard fundamental domain $\mathcal{F}$ for $\SL_2(\Z)$, then the shortest vector in $\Lambda_z$ is $1 \in \C$ and thus $\lambda_z = 1$.
More generally, it is true that $\lambda_z \geq \min(y, 1)$.

We can now turn towards the counting results. 
Adapting the Harcos-Templier method, more precisely inspired by Lemma 2 in Section 2.3 of \cite{HT}, we obtain the following.

\begin{lemmy}
	Let $z=x+iy\in \Hb$ with $y\gg 1$, $g_p\in \SL_2(\Z)$, $0<\delta'<1$ and $0\leq l\leq n$ then we have
	\begin{equation}
		\sum_{\substack{m=l_1l_2 \\ l_1,l_2\in \mathcal{P}_X}} M(z,g_p;\delta',p^l,m) \ll 
		X^2 + \frac{X^3\sqrt{\delta'}}{p^l}y+\frac{X^4\delta'}{p^{l}}.\nonumber
	\end{equation}
\end{lemmy}
\begin{proof}
	We assume $l>0$, since otherwise the result follows directly from the counting results in \cite{Te}. We use the notation for $g_p$ and $\gamma$ introduced in Section~\ref{sec:useful_ineq}. In particular, we reduce to the case where $(v,p)=1$ and $w=p^sw_0$.
	
	We start by estimating the subset of matrices where $c\neq 0$. 
    For clarity, we consider several cases. 
    First, if $s\geq l$, then the congruence coming from \eqref{eq:reduced_congruence} degenerates to
	\begin{equation}
		c\cdot\underbrace{u^2\overline{w_0}}_{\in (\Z/p^l\Z)^{\times}}\equiv 0\text{ mod }p^l. \nonumber
	\end{equation}
	In particular, we must have $p^l\mid c$. Forgetting the remaining congruences,\footnote{This will cost us a power of $p^l$ in the end, but this is irrelevant in our application.} we can simply count pairs $(A,b)$ using \eqref{eq:u-condition}, which we write as
	\begin{equation}
		\vert -cz^2+Az+b\vert^2\ll X^2\delta'y^2.\label{eq:key_to_lattice1}
	\end{equation}
	In view of \eqref{eq:lat-pt-count} we get
	\begin{equation}
		\sharp\{ (a-d,b)\} \ll 1 + \frac{Xy\sqrt{\delta'}}{\lambda_z} + X^2\delta'y.\nonumber
	\end{equation}
	Note that, since $y\gg 1$ we have $\lambda_z^{-1}\ll 1$ and that by \eqref{eq:bound_c} the bound  $c\ll X/y$ holds. This allows us to simply sum over all such $c$ and obtain
	\begin{equation}
		\sharp \{ (c,b,a-d)\} \ll \frac{X}{p^l} + \frac{X^2\sqrt{\delta'}}{p^l}+ \frac{X^3\delta'}{p^l}.\nonumber
	\end{equation} 
	Here we also used that $p^l\mid c$ and $c\neq 0$, as well as $y^{-1}\ll 1$ in the first term. 
    Finally, recall that \eqref{eq:trace_ineq} states that $a+d\ll X$. This allows for an easy count of possible values for $a+d \in \Z$. Altogether we have 
	\begin{equation}
		\sharp \{ (a-d,a+d,b,c) \colon c\neq 0\} \ll \frac{X^2}{p^l} + \frac{X^3\sqrt{\delta'}}{p^l} + \frac{X^4\delta'}{p^{l}}\nonumber
	\end{equation}
	in the case $s\geq l$.
    Each quadruple obviously fixes $\gamma$.
	
	Next, we consider the case $0 \leq s<l$. Here the argument is similar. We first look at \eqref{eq:reduced_congruence} but consider it modulo $p^s$. We infer that $p^s\mid c$, which we keep in mind (this is an empty condition if $s=0$, but this is allowed). 
	
	For what follows it will be useful to identify $\C$ with $\R^2$ via $$\C\ni x+iy \mapsto \left(\begin{matrix} y \\ x\end{matrix}\right)\in  \R^2.$$ In particular, $\mathbb{Z}^2$ interpreted as column vectors corresponds to the Gau\ss ian integers in $\C$. Returning to the argument we fix $c$ and combine \eqref{eq:key_to_lattice1} with \eqref{eq:reduced_congruence}. 
    The latter says that $(A, b)^T$ shifted by a fixed vector lies in the lattice
    \begin{equation}
		\Lambda_{z,g_p} = p^l\cdot \left(\begin{matrix} y & 0\\ x & 1\end{matrix}\right)T_g^{-1}a(p^{-s})\cdot \Z^2. \nonumber
	\end{equation}
    Shifting the centre of the ball in \eqref{eq:key_to_lattice1}, we see that
	\begin{equation}
		\sharp \{ (a-d,b)\} \ll \sharp \Lambda_{z,g_p}\cap B_{\mathfrak{v}_c}(X\sqrt{\delta'}y), \nonumber
	\end{equation} 
	where $B_{\mathfrak{v}_c}(R)$ is the ball or radius $R$ around
	\begin{equation}
		\mathfrak{v}_c = c\left(\begin{matrix} 2xy\\ x^2-y^2\end{matrix}\right)-c\left(\begin{matrix} y & 0 \\ x & 1\end{matrix}\right)T_g^{-1} a(p^{-s})\cdot v_g, \nonumber
	\end{equation}
    which we identify with a point in $\C$, as usual. 
    Let $\lambda_{z,g_p}$ denote the minimum of $\Lambda_{z,g_p}$, and note that $\det(\Lambda_{z,g_p}) = yp^{2l-s}$. Applying \eqref{eq:lat-pt-count} yields
	\begin{equation}
		\sharp \{ (a-d,b)\} \ll 1 + \frac{X\sqrt{\delta'}y}{\lambda_{z,g_p}}+\frac{X^2\delta'y}{p^{2l-s}}.\nonumber
	\end{equation}

	At this point we observe that, since $p^la(p^{-s})\Z^2\subseteq p^{l-s}\cdot \Z^2$ for $s<l$ and $T_g\in \SL_2(\Z)$, we have $\Lambda_{z,g_p}\subseteq p^{l-s}\Lambda_z$. In particular, we have $\lambda_{z,g_p}\geq p^{l-s}\lambda_z$. Inserting this and $\lambda_z\gg 1$ above and summing over the possibilities for $p^s\mid c$ with $0\neq c\ll X/y$ yields
	\begin{equation}
		\sharp \{ (c,a-d,b)\colon c\neq 0\} \ll X + \frac{X^2\sqrt{\delta'}}{p^l}+\frac{X^3\delta'}{p^{2l}}.\nonumber
	\end{equation}
	Finally, we account for $a+d$ trivially using $a+d\ll X$ as above. The result is
	\begin{equation}
		\sharp \{ (c,a-d,a+d,b)\colon c\neq 0\} \ll X^2 + \frac{X^3\sqrt{\delta'}}{p^l}+\frac{X^4\delta'}{p^{2l}}.\nonumber
	\end{equation}
	
    Next, we consider the case of $c=0$.
    The congruence condition from \eqref{eq:reduced_congruence} on the second component implies that
	\begin{equation}
		tA+vb\equiv 0 \text{ mod }p^l\Z.\nonumber
	\end{equation}
	Thus, we can write
	\begin{equation}
		b=-\overline{v_0}t\cdot A+p^lb',\nonumber
	\end{equation}
    for some $b'\in \Z$, where we recall that we can assume $v = v_0$ and $\gcd(v_0, p) = 1$.
	At this point we modify the standard lattice point counting argument and first sort the lattice points by their greatest common divisor. This enables us to input a primitivity condition and improve on the trivial bound from \eqref{eq:lat-pt-count}. Therefore, let $b_0=\gcd(A,b')$ and write $b' = b_0b''$ and $A=b_0A'$.
    Keep in mind that we have achieved $(A',b'')=1$.
	
	Fixing a value of $b_0 \in \N$, condition \eqref{eq:u-condition} gives the problem of counting primitive lattice points $(A',b'')$ such that
	\begin{equation*}
		\vert A'(z-\overline{v_0}t)+p^lb''\vert^2\ll \frac{X^2 \delta' y^2}{b_0^2}.
	\end{equation*}
	By \eqref{eq:lat-pt-count-primitive}, the number of such points is bounded by
	\begin{equation*}
		1+\frac{X^2\delta'}{p^lb_0^2}y.
	\end{equation*}
    Here, the lattice $\Z (z - \overline{v_0} t) + p^l \Z$ is a sublattice of index $p^l$ inside $\Z (z - \overline{v_0} t) + \Z = \Lambda_z$. In particuar we find that the covolume of $\Z (z - \overline{v_0} t) + p^l \Z$ is $p^l\cdot \det(\Lambda_z) = p^ly$.
	
	Next, we bound the number of possibilities for $b_0$. For this, we observe that $ad \ll X^2$, by applying the determinant, and that $a+d \ll X$, by applying the trace. Since $(a-d)^2 = (a+d)^2 - 4ad$, we obtain that $a-d \ll X$. This implies that $b_0 \ll X$, because $b_0 \mid (a-d)$, showing that there are $\ll X$ possibilities for $b_0$. We now sum the bound for the number of primitive points above over $b_0 \ll X$ and find that the number of pairs $(A, b)$ is
	\begin{equation*}
		\ll X+\frac{X^2\delta'}{p^l}y.
	\end{equation*}
	
	Since $ad = m = l_1 l_2$, we conclude that each such pair gives rise to at most $X$ matrices contributing to the count. We get\footnote{The standard argument gives $$\sharp \{ (a-d,ad,b) \} \ll  X+\frac{\sqrt{\delta'}X^2}{\lambda_{z,l}}y+\frac{X^3\delta'}{p^l}y$$ instead, where $\lambda_{z,l}$ is the minimum of the lattice $\Z\cdot z+p^l\Z$.}
	\begin{equation}
		\sharp \{ (a-d,ad,b)\} \ll  X^2+\frac{X^3\delta'}{p^l}y.\nonumber
	\end{equation}
	Combining the two cases completes the counting for $M$ and the result as stated follows by simplifying some of the terms.
\end{proof}

For counting matrices with square determinant $m = l_1^2 l_2^2$, we employ a standard strategy to detect the sparseness in the average over the determinants.
We write
\begin{equation}
	M(z,1;\delta',p^l,m) = M_{\ast}(z,1;\delta',p^l,m)+M_\parab(z,1;\delta',p^l,m),\nonumber
\end{equation}
where the subscript $\parab$ stands for the sub-count of parabolic matrices with non-vanishing bottom-left entry.
Recall that parabolic matrices satisfy the equation $\tr^2 = 4 \det$.
Now, for the first term, we proceed with a similar argument.

\begin{lemmy}
Let $z=x+iy\in \Hb$ with $y\gg 1$, $g_p\in \SL_2(\Z)$, $0<\delta'<1$ and $0\leq l$.
Then, we have
\begin{equation}
	\sum_{\substack{m=l_1^2l_2^2 \\ l_1,l_2\in \mathcal{P}_X}} M_{\ast}(z,g_p;\delta',p^l,m) \ll_{\epsilon} X^{2+\epsilon}\cdot \left(1+\frac{\sqrt{\delta'}X^2}{p^l}y+\frac{\delta'X^4}{p^l}\right).\nonumber
\end{equation}
\end{lemmy}
\begin{proof}
Note that every matrix contributing to the count has square determinant $\det(\gamma) = m \asymp X^4$. We first assume $c\neq 0$ and count admissible entries $a-d,b,c$ as above. This gives
\begin{equation}
		\sharp\{ (a-d,b,c)\} \ll X^2 + \frac{\sqrt{\delta'}X^4}{p^l} + \frac{X^6\delta'}{p^{l}}.\nonumber
\end{equation}
Having fixed one such triple, instead of counting possibilities for $a+d$ trivially using an upper bound for the trace, we recall that
\begin{equation}
	(a-d)^2+4bc = (a+d)^2-4m = (a+d-2l_1l_2)(a+d+2l_1l_2).\nonumber
\end{equation}
Since we are excluding parabolic matrices, we obtain that the left-hand side is non-zero. Arguing as in the proof of \cite[Lemma~2]{HT} we find at most $X^{\epsilon}$ possibilities for tuples $(a+d,m)$. This gives
\begin{equation}
	\sharp\{ (a-d,a+d,b,c)\colon c\neq 0,\, (a+d)^2-4m\neq 0\} \ll X^{\epsilon}\left(X^2 + \frac{\sqrt{\delta'}X^4}{p^l} + \frac{X^6\delta'}{p^{l}}\right).\nonumber
\end{equation}

Finally, we account for $c=0$ by following the argument in \cite[Lemma~4]{HT}. Indeed, as in the previous lemma we use lattice point counting to estimate
\begin{equation}
	\sharp \{a-d,b\} \ll X^2+\frac{X^4\delta'}{p^l}y. \nonumber
\end{equation}
Then we use $ad=l_1^2l_2^2$ for primes $l_1l_2$ to find for each (fixed) $a-d$ there are $\ll X^{\epsilon}$ possibilities for pairs $(a,d)$. Combining the contributions completes the proof.
\end{proof}

For bounding the parabolic contribution, we present a new strategy that produces a result similar in shape to \cite[Lemma~4.4]{Te}, where the proof is left out, but with a saving in $p^l$.
\begin{lemmy}
	Let $z=x+iy\in \Hb$ with $y\gg 1$, $0<\delta'<1$ and $0\leq l\leq n$ then we have
	\begin{equation}
		\sum_{\substack{m = (l_1 l_2)^2 \\ l_1,l_2 \in \mathcal{P}_X}} M_{\parab}(z,g_p;\delta',p^l,m) \ll 
		X^{2+\epsilon} + \frac{X^{4+\epsilon} y \sqrt{\delta'}}{p^l}.\nonumber
	\end{equation}
\end{lemmy}
\begin{proof}
	Let $\gamma$ be parabolic with $\det \gamma = m = a^2$, where $a = l_1 l_2$.
	Then there exists $\xi \in \SL_2(\Z)$ such that
	\begin{displaymath}
		\xi \gamma \xi^{-1} = \gamma' = \begin{pmatrix}
			a & b \\ 0 & a
		\end{pmatrix}.
	\end{displaymath}
	Thus, it suffices to count the possibilities for $\xi \in U(\Z) \backslash \SL_2(\Z)$, where $U(\Z)$ are upper triangular matrices in $\SL_2(\Z)$, and $\gamma'$ of the form above.
	We note that the classes $U(\Z) \backslash \SL_2(\Z)$ are parametrised by the second row of their representative.
	
	To start, the number of $\gamma = \diag(a, a)$ is simply bounded by the number of terms in the sum, that is $X^2$.
	From now on, we assume that $\gamma$ is not a central matrix.
	This implies that $\gamma'$ as above must satisfy that $b \neq 0$.
	
	At the place $p$, write
	\begin{displaymath}
		\xi g_p = \begin{pmatrix}
			t & u \\ v & w
		\end{pmatrix}.
	\end{displaymath}
	Then we can rewrite the condition as
	\begin{displaymath}
		(\xi g_p)^{-1} \gamma' (\xi g_p) \in K_{H,p}(l, l).
	\end{displaymath}
	Using the computation \eqref{eq:p-adic-conjugation-condition}, we deduce from the special shape of $\gamma'$ that
	\begin{displaymath}
		p^l \mid v^2 b \quad \text{and} \quad p^l \mid w^2 b.
	\end{displaymath}
	Since $\gcd(w,v)=1$, it follows that $p^l \mid b$.
	
	Moving to the archimedean place, put $z'= \xi z = x'+ iy'$.
	Then the condition can be rewritten as $u(\gamma' z', z') \leq \delta'$, by invariance of the distance function.
	In terms of matrix entries as in \eqref{eq:u-condition}, this implies that $b \ll a y' \sqrt{\delta'}$.
	
	Since $p^l \mid b$ and $b \neq 0$ by assumption, this gives a bound of $X^2 y' \sqrt{\delta'} p^{-l}$ on the number of $b$ for every fixed $\xi$.
	There is no such $\xi$ unless $1/y' \ll  a \sqrt{\delta'}$.

	Now we sum up all possibilities, sorting by $\xi$, to get that
	\begin{displaymath}
		\sum_{\substack{\abs{cz + d}^2 \leq y X^2 \sqrt{\delta'}, \\ (c,d)=1}} \frac{X^{2+\epsilon} y \sqrt{\delta'}}{p^l \abs{cz + d}^2} \ll p^{-l} X^2 y \sqrt{\delta'}. 
	\end{displaymath}
	Here we have used the formula for $y'$ and that
	\begin{displaymath}
		\sum_{\substack{\abs{cz + d}^2 \leq Y,\\(c,d)=1}} \frac{1}{\abs{cz + d}^2} \leq Y^{\epsilon} \sum_{(c,d)=1} \frac{1}{\abs{cz + d}^{2+\epsilon}} \ll_{\epsilon} Y^{\epsilon}, \nonumber
	\end{displaymath}
	for any $Y$. In the last step we have used absolute convergence of the sum, see \cite[Section~3.2]{iwaniec-spectral}, as well as $y\gg 1$. We now finish the count by considering all possibilities for the determinant.
\end{proof}

Finally, we also need the following result.

\begin{lemmy}
Let $z=x+iy\in \Hb$ with $y\gg 1$, $g_p\in \SL_2(\Z)$, $0<\delta'<1$ and $0\leq l$ then we have
\begin{equation}
	M(z,g_p;\delta',p^l,1) \ll 1+\frac{\sqrt{\delta'}}{p^l}y.\nonumber
\end{equation}
\end{lemmy}
\begin{proof}
Counting $M(z,g_p,\delta',p^l,1)$ is straight forward following the argument from \cite[Lemma~11]{BHMM} for example. We omit the details.
\end{proof}

We now present our counting results in a compact form, which is convenient for the application below. %

\begin{cor}\label{cor:summary_count}
Let $z = x+iy \in \Hb$ with $y\gg 1$, $g_p\in \SL_2(\Z)$, $0<\delta'<1$ and $0\leq l$. Then we have
\begin{equation}
	M(z,g_p,\delta',p^l,1) \ll 1+\frac{\sqrt{\delta'}}{p^{l}}y. \nonumber
\end{equation}
Furthermore, we have
\begin{equation}
	\sum_{\substack{m=l_1^il_2^i \\ l_1,l_2\in \mathcal{P}_X\\ i=1,2}}m^{-\frac{1}{2}} M(z,g_p;\delta',p^l,m) \ll X^{1+\epsilon}+ X^{2+\epsilon}\frac{\sqrt{\delta'}}{p^l}y+\frac{X^{4+\epsilon}\delta'}{p^{l}}.\nonumber
\end{equation}
\end{cor}

\begin{rem}\label{rem:com_counting_IS}
Our counting is based on a quite elaborate lattice point counting scheme and has its origin in a series of papers by Templier and Harcos-Templier. Historically, the development of the counting method  (and also of the amplifier) was driven by optimizing the sup-norm bounds in the (square-free) level aspect. Indeed, while there have been many improvements in the level aspect over the years, the spectral bound from \cite{IS95} has not been surpassed at the time of writing the present article. 

This being said, efficient counting turns out to be crucial for us. 
We will illustrate this in Section~\ref{sec:arch_naive} below, where Corollary~\ref{cor:first_sub_conv_bound} shows that our soft analysis and the old counting techniques of Iwaniec and Sarnak are not enough to recover the best known exponent in the spectral aspect.
See the discussion after Theorem \ref{thm:archi-main-estimate} for the comparison of the analytic results.

For later reference, we record the following consequence of Lemma~\ref{lm:IS_counting}. For $z\in \mathcal{F}$, $g_p\in \textrm{SL}_2(\Z)$ and $\delta'\leq 1$ we have
\begin{equation}
	\sum_{\substack{m=l_1^il_2^i \\ l_1,l_2\in \mathcal{P}_X\\ i=1,2}}m^{-\frac{1}{2}} M(z,g_p;\delta',p^l,m) \ll X^{\epsilon}(X + X^4(\delta')^{\frac{1}{4}}+X^2\sqrt{\delta'} y).\nonumber
\end{equation}
Observe that here we encounter a term featuring $X^4(\delta')^{\frac{1}{4}}$, while in the refined counting result this is improved to $X^4\delta'p^{-2l}$. 
It is this improvement which will be responsible for recovering the Iwaniec--Sarnak spectral exponent even with slightly weaker bounds on the archimedean test functions.
\end{rem}

\section{The global period II}\label{sec:glob_perII}

In this section, we combine all the pieces and finally prove our main Theorem~\ref{th:Intro} (see also Theorem~\ref{th:main_theorem_impro} below). The strategy behind the proof was sketched in Section~\ref{sec:sketch_amp} and we carry out the details here. We first upgrade the relative trace fromula from Proposition~\ref{pr:pre_amplified_guy} by including the amplifier discussed in Section~\ref{se:amplifier}. This results in Proposition~\ref{pr:amp_trace}, which is the starting point of our analysis. For pedagogical purposes we will present three sub-baseline sup-norm bounds gradually increasing the technical difficulty of the argument:
\begin{itemize}
	\item In Theorem~\ref{th:naive_bound}, we present the most naive argument, which is closest to the sketch given in Section~\ref{sec:sketch_amp}. Here we estimate the archimedean test function using the geometric ideas developed in Section~\ref{sec:archi-microlocalisation-big-chap} and handle the finite ramified place trivially. We also directly import the original counting result from Iwaniec and Sarnak \cite{IS95}. As a result we obtain a bound that is sub-baseline in the spectral aspect and matches the baseline bound in the level aspect.
	\item The next result is Proposition~\ref{prop:amplified_bound}. Here, we refine our argument by improving the analysis at the ramified finite place. We do so by bounding the orbital integrals (i.e. newform matrix coefficients) on the geometric side using the soft geometric tools developed in Sections~\ref{sec:p_adic_vol} and~\ref{sec:p-adic_relmat}. This refinement also requires the use of the new counting results obtained in Section~\ref{sec:count}. As a result we obtain a  hybrid sub-baseline bound for the sup-norm, which matches the Iwaniec--Sarnak exponent in the spectral aspect.
	\item Finally, the best bounds are obtained in Theorem~\ref{th:main_theorem_impro} below. Here, we replace the soft treatment of the $p$-adic test-function by sharp bounds for newform matrix coefficients obtained by Hu and Saha. This allows for an improvement of the exponent in the level aspect. Indeed, the achieved hybrid sup-norm bound finally has the Iwaniec--Sarnak exponents in all aspects.
\end{itemize}

Throughout this section, we let $\pi$ be a cuspidal automorphic representation for $\GL_2(\Q)$ with trivial central character and level $N = p^{4n}$. We write $N_0=p^{2n}$. We assume that $\pi$ is spherical at infinity and denote its spectral parameter by $T$. For convenience, we assume that $T$ is positive and sufficiently large. In particular, $\pi$ is assumed to be tempered at infinity.

\subsection{An amplified relative pretrace inequality}

\begin{prop}\label{pr:amp_trace}
Let $f_{\infty}\in \mathcal{C}_c^{\infty}(Z\backslash G(\R))$ and $f_p\in \mathcal{C}_c^{\infty}(Z\backslash G(\Q_p))$ be two test functions, which we combine to $$f=f_{\infty}\otimes f_p\otimes \bigotimes_{l\neq p} \vvmathbb{1}_{ZK_l}.$$ Then, for $g\in G(\R)\times G(\Q_p)$ we have
\begin{multline}
	\vert \varphi'(g)\vert^2\cdot \sum_{\Psi\in \mathcal{B}^N(\pi)} \mathcal{Q}(\pi(f)\Psi) \\ \leq \frac{8}{(\sharp \mathcal{P}_X)^2} \sum_{m=1}^{\infty}a_m\int_{x,y \in H} \sum_{\gamma \in G'(\Q)} (f \ast f^\ast)(x^{-1} g^{-1}\gamma g y)\cdot \kappa_m(\gamma) \, dx \, dy.\nonumber 
\end{multline}
with coefficients $a_m$ as in \eqref{eq:expansion_hecke_rel}.
\end{prop}
\begin{proof}
We first define the test-function 
\begin{equation}
	\tilde{f}_i = f_{\infty} \otimes f_p \otimes f_{\textrm{ur},i}^X \text{ for }i=1,2.\nonumber 
\end{equation}
If $\Psi\in \mathcal{B}^N(\pi)$, then 
\begin{equation}
	\pi(\tilde{f}_i)\Psi = \mathcal{S}_{\pi,i}(X)\cdot \pi(f)\Psi
\end{equation}
for $$\mathcal{S}_{\pi,i}(X) = \sum_{l\in \mathcal{P}_X} \vert \lambda_{\pi}(l^i)\vert.$$ In particular $C_X(\pi)=\mathcal{S}_{\pi,1}(X)^2+\mathcal{S}_{\pi,2}(X)^2$. Of course, $\pi(\tilde{f}_i)\Psi=0$ if $\Psi \in \mathcal{B}(\pi)\setminus \mathcal{B}^N(\pi)$. This allows us to write
\begin{equation}
	C_X(\pi)\cdot  \sum_{\Psi\in \mathcal{B}^N(\pi)} \mathcal{Q}(\pi(f)\Psi) = \sum_{i=1,2}\sum_{\Psi\in \mathcal{B}(\pi)} \mathcal{Q}(\pi(\tilde{f}_i)\Psi).\nonumber
\end{equation}
Applying \eqref{eq:pretrace-proto} to each piece and expanding $f_{\textrm{ur}}^X = f_{\textrm{ur},1}^X\ast f_{\textrm{ur},1}^X+f_{\textrm{ur},2}^X\ast f_{\textrm{ur},2}^X$ using \eqref{eq:expansion_hecke_rel} yields
\begin{multline}
		C_X(\pi)\cdot \vert \varphi'(1)\vert^2\cdot  \sum_{\Psi\in \mathcal{B}^N(\pi)} \mathcal{Q}(\pi(f)\Psi) \\ \leq  \sum_{m=1}^{\infty}a_m\int_{x,y \in H} \sum_{\gamma \in G'(\Q)} (f \ast f^\ast)(x^{-1} \gamma y) \, dx \, dy \cdot \kappa_m(\gamma).\nonumber
\end{multline}
After dividing by $C_{X}(\pi)$ and applying the lower bound \eqref{eq:basic_ampli_inequality} the proof of the initial case is complete. To deal with general $g$ one follows the arguments leading to \eqref{eq:pretrace-g}.
\end{proof}

\subsection{The naive archimedean bound}\label{sec:arch_naive}

Let us start by using the formula from Proposition~\ref{prop:classif-rep} to prove a sub-baseline bound in the spectral aspect.

\begin{theorem}\label{th:naive_bound}
Let $\pi$ be a cuspidal automorphic representation for $\textrm{GL}_2(\Q)$ with trivial central character, level $N=p^{4n}$ and spectral parameter $T$. We assume that $T$ is sufficiently large, so that in particular $\pi$ is tempered and spherical at infinity. Let $g=g_{\infty} g_p$ with $g_{\infty}.i=x+iy\in \mathcal{F}$ and $g_p\in \SL_2(\Z)$. For the unique $L^2$-normalized $K_H(2n,2n)$-invariant form $\varphi'$ in $\pi$ defined in \eqref{eq:shifted_new}. We have
\begin{equation}
	\vert \varphi'(g)\vert^2 \ll N^{\frac{1}{2}}T^{\frac{11}{12}+\epsilon}+N^{\frac{1}{2}}T^{\frac{1}{2}+\epsilon}y.\nonumber
\end{equation}
\end{theorem}
\begin{proof}
Our starting point is Proposition~\ref{pr:amp_trace} with the test functions $f_{\infty}$ and $f_p$ as in Section~\ref{sec:test_fct}. The parameter $X$ will be chosen as a power of $T$, so that $\sharp \mathcal{P}_X\gg XT^{-\epsilon}$. Furthermore, recall from Theorem~\ref{thm:archi-main-estimate} and Lemma~\ref{lm:summary_Q} that
\begin{equation}
	\sum_{\Psi\in \mathcal{B}^N(\pi)} \mathcal{Q}( \pi(f)\Psi) \gg p^{-n}T^{-\frac{1}{2}}. \label{eq:lower_boound_spec}
\end{equation}
On the other hand we estimate the $H$-integrals as
\begin{equation}
	\int_{x,y \in H} (f \ast f^\ast)(x^{-1} g^{-1}\gamma gy) \, dx \, dy \ll T^{\epsilon}p^n\cdot \min(d(g_{\infty}^{-1}\gamma g_{\infty})^{-1},T^{\frac{1}{2}}). \label{eq:applied_vol_bound}
\end{equation}
For the  archimedean part we have used the estimate from Theorem~\ref{thm:archi-main-estimate}. On the other hand, we have used only the trivial bound for the non-archimedean part, see Corollary~\ref{cor:vol_bound_p-adig}.

Altogether, using also the support properties of the test functions, we obtain
\begin{equation}
	\vert \varphi'(g)\vert^2\ll \frac{T^{\frac{1}{2}+\epsilon}N^{\frac{1}{2}}}{X^2} \sum_{m=1}^{\infty}a_m \sum_{\substack{\gamma \in G'(\Q),\\ d(g_{\infty}^{-1}\gamma g_{\infty}) \ll 1,\\ d_H(g_p^{-1}\gamma g_p)\leq 1 }}\min(d(g_{\infty}^{-1}\gamma g_{\infty})^{-1},T^{\frac{1}{2}}) \cdot \vert \kappa_m(\gamma)\vert.\nonumber 
\end{equation}
We sort the $\gamma$-sum by $d(g_{\infty}^{-1}\gamma g_{\infty})\in [\delta/2,\delta]$ and sum $\delta$ over a dyadic intervals. This way we get
\begin{multline}
	\vert \varphi'(g)\vert^2\ll \frac{T^{\frac{1}{2}+\epsilon}N^{\frac{1}{2}}}{X^2} \sum_{m=1}^{\infty}\frac{a_m}{m^{\frac{1}{2}}} \bigg[
	\sum_{\substack{T^{-\frac{1}{2}}\leq \delta \ll 1,\\ \text{dyadic}}}\delta^{-1} \cdot \sharp \mathfrak{M}(g_{\infty},g_p;m,\delta,0) \\
	+ T^{\frac{1}{2}}\cdot \sharp \mathfrak{M}(g_{\infty},g_p;m,T^{-\frac{1}{2}},0)\bigg].\nonumber 
\end{multline}
We first discuss the contribution from $d(g_{\infty}^{-1}\gamma g_{\infty})<T^{-\frac{1}{2}}$. Indeed, by Lemma~\ref{lm:ana_support} and Lemma~\ref{lm:IS_counting} we obtain
\begin{equation}
	\sharp \mathfrak{M}(g_{\infty},g_p;m,T^{-\frac{1}{2}},0) \leq M(z,g_p,T^{-1},1,m) \ll_{\epsilon} m^{\epsilon}(1+mT^{-\frac{1}{4}}+m^{\frac{1}{2}}T^{-\frac{1}{2}}y).\nonumber
\end{equation}
Inserting this estimate above leads to a contribution
\begin{align}
	&\frac{T^{\frac{1}{2}+\epsilon}N^{\frac{1}{2}}}{X^2} \sum_{m=1}^{\infty}\frac{a_m}{m^{\frac{1}{2}}} T^{\frac{1}{2}}\cdot \sharp \mathfrak{M}(g_{\infty},g_p;m,T^{-\frac{1}{2}},0)\nonumber\\
	&\qquad \ll \frac{T^{1+\epsilon}N^{\frac{1}{2}}}{X^2} \sum_{m=1}^{\infty} a_m(m^{\epsilon-\frac{1}{2}}+m^{\frac{1}{2}+\epsilon}T^{-\frac{1}{4}}+m^{\epsilon}T^{-\frac{1}{2}}y) \nonumber\\
	&\qquad \ll \frac{T^{1+\epsilon}N^{\frac{1}{2}}}{X^2} ( X + X^4T^{-\frac{1}{4}}+X^2T^{-\frac{1}{2}}y).\nonumber
\end{align}
In the last step we have estimated the $m$-sum by recalling that the sequence $a_m$ is supported on $m=l_1^il_2^i$ for $i=1,2$ and $l_1,l_2\in \mathcal{P}(X)$. We have also used the bound $a_m\ll 1+\delta_{m=1}\cdot X$. In summary, we see that the $\gamma$'s with $d(g_{\infty}^{-1}\gamma g_{\infty})<T^{-\frac{1}{2}}$ contribute at most
\begin{equation}
	\frac{T^{\frac{1}{2}+\epsilon}N^{\frac{1}{2}}}{X^2} \left( XT^{\frac{1}{2}}+X^4T^{\frac{1}{4}} + X^2y\right) \nonumber
\end{equation}
to our bound for $\vert \varphi'(g)\vert^2$. It turns out that this will match the contribution of $\delta \geq T^{-\frac{1}{2}}$, which we will discuss next.

As before, we start by using Lemma~\ref{lm:ana_support}. This allows us to write
\begin{equation}
	\vert \varphi'(g)\vert^2\ll \frac{T^{\frac{1}{2}+\epsilon}N^{\frac{1}{2}}}{X^2} 	\sum_{\substack{T^{-1}\leq \delta' \ll 1,\\ \text{dyadic}}}(\delta')^{-\frac{1}{2}}\sum_{m=1}^{\infty}\frac{\vert a_m\vert}{m^{\frac{1}{2}}} M(z,g_p;\delta',1,m).\nonumber 
\end{equation}
We now recall that $a_1\ll X$ and that $M(z,g_p;\delta',1,m) \ll 1+\sqrt{\delta'}y$ by Lemma~\ref{lm:IS_counting}. On the other hand, we have $a_m\ll1$ if $m=l_1^il_2^i$ for $i=1,2$ and $l_1,l_2\in \mathcal{P}_X$. Thus, using the bound from Remark~\ref{rem:com_counting_IS} for the remaining $m$-sum, we obtain
\begin{equation}
	\vert \varphi'(g)\vert^2\ll \frac{T^{\frac{1}{2}+\epsilon}N^{\frac{1}{2}}}{X^2} 	\sum_{\substack{T^{-1}\leq \delta' \ll 1,\\ \text{dyadic}}}\left( X(\delta')^{-\frac{1}{2}}+X^4(\delta')^{-\frac{1}{4}} + X^2y\right).\nonumber
\end{equation}
The length of the dyadic $\delta'$-sum is $\ll T^{\epsilon}$. Thus executing it trivially yields
\begin{equation}
	\vert \varphi'(g)\vert^2\ll \frac{T^{\frac{1}{2}+\epsilon}N^{\frac{1}{2}}}{X^2} \left( XT^{\frac{1}{2}}+X^4T^{\frac{1}{4}} + X^2y\right).\nonumber
\end{equation}
Finally, choosing $X=T^{\frac{1}{12}}$ gives the desired result.
\end{proof}

\begin{cor}\label{cor:first_sub_conv_bound}
We keep the notation from Theorem~\ref{th:naive_bound} and let $\Omega\subset\mathcal{F}$ be a fixed compact set. Put $\Omega'=\SL_2(\Z)\Omega$. Then we have
\begin{equation}
	\Vert \varphi'\vert_{\Omega'}\Vert_{\infty} \ll_{\Omega,\epsilon} N^{\frac{1}{4}}T^{\frac{1}{2}-\frac{1}{24}+\epsilon}.\nonumber
\end{equation} 
\end{cor}

Note that we match here the local bound (i.e. Theorem~\ref{th:local_bound}) in the level aspect. However, we have improved the exponent in the spectral aspect by $1/24$. Note that this is worse than the improvement of $1/12$ that is established in \cite{IS95}. This is an important observation, because we have used the same counting (and even an improved amplifier). The only difference is that the volume bound we use to handle the geometric side is slightly weaker than the Iwaniec--Sarnak analogue. See Theorem~\ref{eq:comp_IS_bounds} and the surrounding discussion.

It turns out that we can still match the saving obtained by Iwaniec and Sarnak in the spectral aspect by using the improved counting results developed above. This will be implemented in the next section. Doing so we will also perform a non-trivial analysis at the place $p$ and obtain a hybrid sub-local bound, which is our final result.

\subsection{The general bound}

We will use $g_p\in \SL_2(\Z)$ and allow $g_{\infty}$ to vary such that $z= x+iy = g_{\infty}.i$ satisfies $y\gg 1$. We will start estimating the amplified (relative) trace inequality from Proposition~\ref{pr:amp_trace} essentially as in the previous subsection. However, we will replace $f_p$ by $f_p' = f_p\cdot \vvmathbb{1}_{ZK_H(1)}$ as defined in Remark~\ref{rm:support_shrink}. Note that the lower bound from \eqref{eq:lower_boound_spec} still holds, up to possibly changing the implicit constant, which is now allowed to depend on $p$. The point of restricting the support of the $p$-adic test function is to slightly simplify the counting problem. This trick was already alluded to in Remark~\ref{rm:intro_Mop} above. Indeed, suppose $f_p'(x^{-1}g_p^{-1}\gamma g_py)\neq 0$ for $x,y\in H(2n)$, then $g_p^{-1}\gamma g_p\in K_H(1)$. However, this implies that $p\mid [g_p^{-1}\gamma g_p]_{21}$ and we infer that 
\begin{equation}
	d_{H(2n)}(wg_p^{-1}\gamma g_p) = 0. \nonumber
\end{equation} 
As a result we can drop the maximum over $\gamma'$ in our volume bound Corollary~\ref{cor:vol_bound_p-adig}. This allows us to replace the estimate \eqref{eq:applied_vol_bound}, which was trivial in the level aspect, by the finer estimate
\begin{multline}
	\int_{x,y \in H} (f \ast f^\ast)(x^{-1} g^{-1}\gamma gy) \, dx \, dy \ll_p T^{\epsilon}\cdot \min(d(\gamma)^{-1},T^{\frac{1}{2}})\\ \cdot \min\left(p^{n/2}d_{H(2n)}(g_p^{-1}\gamma g_p)^{-\frac{1}{2}},p^n\right). \label{eq:gybrid_vol_bound_geo_side} 
\end{multline}
We now drop any additional support property coming from $f_p'$ and continue the argument as above. We sort the $\gamma$-sum according to $d(\gamma)\in [\delta/2,\delta]$ and $d_{H(2n)}(\gamma)=p^{-l}$. In view of Lemma~\ref{lm:ana_support} we obtain
\begin{equation}
	\vert \varphi'(g)\vert^2 \ll_{p,\epsilon} \frac{T^{\frac{1}{2}+\epsilon}p^n}{X^2} \sum_{m=1}^{\infty}\frac{\vert a_m\vert}{m^{\frac{1}{2}}} \sum_{\substack{T^{-1}\leq \delta'\leq 1\\ \text{dyadic}}}\sum_{l=0}^{n} M(g_{\infty},g_p;\delta',p^l,m) \cdot \frac{p^{\frac{n+l}{2}}}{(\delta')^{\frac{1}{2}}}.\nonumber 
\end{equation}
Recall the bounds $a_1\ll X$ and $a_m\ll 1$ for $m>1$ as well as the support of the weights $a_m$. Inserting the counting results from Corollary~\ref{cor:summary_count} yields 
\begin{equation}
	\vert \varphi'(g)\vert^2 \ll_{p,\epsilon} \frac{T^{\frac{1}{2}+\epsilon}p^n}{X^2} \sum_{\substack{T^{-1}\leq \delta'\leq 1\\ \text{dyadic}}}\sum_{l=0}^{n} \frac{p^{\frac{n+l}{2}}}{(\delta')^{\frac{1}{2}}}\cdot \left(X+X^2\frac{\sqrt{\delta'}}{p^l}y+\frac{X^4\delta'}{p^{l}}\right).\nonumber 
\end{equation}
Here we have ignored the contribution of $\delta'<T^{-1}$, which can be handled similarly and will give a comparable contribution. Executing the $\delta'$ and $l$ sums trivially yields
\begin{equation}
	\vert \varphi'(g)\vert^2 \ll_{p,\epsilon} \frac{T^{1+\epsilon}p^{(2+\epsilon)n}}{X} + T^{\frac{1}{2}+\epsilon}p^{(\frac{3}{2}+\epsilon)n}(X^2+y).\nonumber
\end{equation}
At this point we choose $X=(Tp^n)^{\frac{1}{6}}$. After recalling that $N=p^{4n}$ we obtain
\begin{equation}
		\vert \varphi'(g)\vert^2 \ll_{p,\epsilon}T^{\frac{5}{6}+\epsilon}N^{\frac{11}{24}+\epsilon} + T^{\frac{1}{2}+\epsilon}N^{\frac{3}{8}+\epsilon} y.\nonumber
\end{equation}
We record this as a proposition.

\begin{prop}\label{prop:amplified_bound}
Let $g_{\infty}\in G(\R)$ and $g_p\in \SL_2(\Z)$ and assume that $g_{\infty}.i = x+iy\in \Hb$ satisfies $y\gg 1$. Then we have
\begin{equation}
 	\vert \varphi'(g_{\infty})\vert \ll_{p,\epsilon}T^{\frac{1}{2}-\frac{1}{12}}N^{\frac{1}{4}-\frac{1}{48}} + T^{\frac{1}{4}}N^{\frac{3}{16}}\sqrt{y}.\nonumber
\end{equation}
\end{prop}

Note that this immediately provides a numerical improvement over the result from Corollary~\ref{cor:first_sub_conv_bound}. For our main result we want to make this estimate global.

\begin{theorem}\label{th:main_theorem}
Let $\pi$ be a cuspidal automorphic representation for $\textrm{GL}_2(\Q)$ with trivial central character, level $N=p^{4n}$ and spectral parameter $T$. We assume that $T$ is sufficiently large, so that in particular $\pi$ is tempered and spherical at infinity. For the unique $L^2$-normalized newform $\varphi$ in $\pi$ we have
\begin{equation}
	\Vert \varphi\Vert_{\infty} \ll_{p,\epsilon} (TN)^{\epsilon}\cdot \left(T^{\frac{1}{2}-\frac{1}{12}+\epsilon}N^{\frac{1}{4}-\frac{1}{48}+\epsilon}\right).\nonumber
\end{equation}
\end{theorem}
\begin{proof}
Recall the shifted newform $\varphi'$ from \eqref{eq:shifted_new} and note that $\Vert \varphi\Vert_{\infty}=\Vert \varphi'\Vert_{\infty}$. By strong approximation it is sufficient to estimate $\varphi'(g)$ where $g=g_{\infty}g_p$ where $g_{\infty}.i=z=x+iy\in \Hb$ lies in the standard fundamental domain for $\SL_2(\Z)$ and $g_p\in \SL_2(\Z)$. We first observe that, for $\sqrt{y}>T^{\frac{1}{6}}N^{\frac{1}{24}}$, the Fourier bound from Proposition~\ref{prop:Fourier} yields
\begin{equation}
	\vert \varphi'(g_{\infty})\vert \ll_{p,\epsilon} T^{\frac{1}{2}-\frac{1}{6}+\epsilon}N^{\frac{1}{4}-\frac{1}{24}+\epsilon}.\nonumber 
\end{equation}
For $\sqrt{y}\leq T^{\frac{1}{6}}N^{\frac{1}{24}}$ we use the pretrace bound and obtain
\begin{equation}
	\vert \varphi'(g_{\infty})\vert \ll_{p,\epsilon}T^{\frac{1}{2}-\frac{1}{12}+\epsilon}N^{\frac{1}{4}-\frac{1}{48}+\epsilon}.\nonumber
\end{equation}
This is the desired bound.
\end{proof}

\begin{rem}
We have stated this global sup-norm bound for the classical newform $\varphi$ (i.e. the unique $L^2$-normalized $K_H(N)$-invariant vector.) This is simply a cosmetic choice because, as already observed in the proof, we have $\Vert \varphi\Vert_{\infty}=\Vert \varphi'\Vert_{\infty}$. However, in the bounds for the restricted sup-norm, for example Theorem~\ref{th:naive_bound}, it is essential that we work with $\varphi'$. This is because we fix the compact set $\Omega$ independent of $N$ and $T$. But then $\Vert \varphi\vert_{a(N_0)\Omega}\Vert_{\infty} = \Vert \varphi'\vert_{\Omega}\Vert_{\infty}$ and the volume of $\Omega$ and $a(N_0)\Omega$ differ significantly. Classically speaking, if we bound $\varphi'$ at the point $i\in \Hb$, this corresponds to bounding $\varphi$ at $\frac{i}{N_0}$.
\end{rem}

Note that arguably the most important analytic ingredients for our proof are the volume bounds given in Theorem~\ref{thm:archi-main-estimate} (resp. Corollary~\ref{cor:vol_bound_p-adig}) in the archimedean (resp. $p$-adic aspect). In contrast to earlier works, these bounds are somehow of geometric nature and do not rely on direct analysis of the Selberg--Harish-Chandra transform (or other elaborate constructions of test functions using Paley--Wiener theory) as in \cite{IS95} or \cite{Te}. 
In the $p$-adic case, we completely avoid any stationary phase analysis of new-vector matrix coefficients as for example in \cite{Hu-Sa}. While our bound matches the long-standing Iwaniec--Sarnak barrier in the spectral aspect we have a rather weak exponent in the level. This can be improved using matrix coefficients bounds. One obtains the following very harmonic result:

\begin{theorem}\label{th:main_theorem_impro}
Let $\pi$ be a cuspidal automorphic representation for $\textrm{GL}_2(\Q)$ with trivial central character, level $N=p^{4n}$ and spectral parameter $T$. We assume that $T$ is sufficiently large, so that in particular $\pi$ is tempered and spherical at infinity. For the unique $L^2$-normalized newform $\varphi$ in $\pi$ we have
\begin{equation}
	\Vert \varphi\Vert_{\infty} \ll_{p,\epsilon} (TN_0)^{\frac{1}{2}-\frac{1}{12}+\epsilon}.\nonumber
\end{equation}
\end{theorem}
\begin{proof}
Again we have to estimate $\varphi'(g)$ (as in \eqref{eq:shifted_new}) for $g=g_{\infty}g_p$ with $g_p\in \SL_2(\Z)$ and $g_{\infty}.i=x+iy$ with $y\gg 1$. Note that for $\sqrt{y}>T^{\frac{1}{6}}N^{\frac{1}{12}}$ we can simply use Proposition~\ref{prop:Fourier}. This will be complemented using the amplification method.

We want to run the same argument as above. However, we can use \eqref{eq:Mann_Mann_Achim} to improve \eqref{eq:gybrid_vol_bound_geo_side} and obtain
\begin{multline}
	\int_{x,y \in H} (f \ast f^\ast)(x^{-1} g^{-1}\gamma gy) \, dx \, dy \ll_{p,\epsilon} T^{\epsilon}\cdot  \min(d(\gamma)^{-1},T^{\frac{1}{2}})\\ \cdot \min\left(d_{H(2n)}(\gamma)^{-\frac{1}{2}},p^n\right).
\end{multline}
With this at hand, we obtain
\begin{equation}
	\vert \varphi'(g)\vert^2 \ll_{p,\epsilon} \frac{T^{\frac{1}{2}+\epsilon}p^n}{X^2} \sum_{m=1}^{\infty}\frac{\vert a_m\vert}{m^{\frac{1}{2}}} \sum_{\substack{T^{-1}\leq \delta'\leq 1\\ \text{dyadic}}}\sum_{l=0}^{2n} M(g_{\infty},g_p;\delta',p^l,m) \cdot \frac{p^{\frac{l}{2}}}{(\delta')^{\frac{1}{2}}}.\nonumber 
\end{equation}
From here we follow the same argument as above. Namely, we recall support and size of the amplification weights $a_m$, insert our counting results and execute the $\delta'$- and $l$-sums trivially. We arrive at
\begin{equation}
	\vert \varphi'(g)\vert^2 \ll_{p,\epsilon} \frac{T^{1+\epsilon}p^{(2+\epsilon)n}}{X} + T^{\frac{1}{2}+\epsilon}p^{(1+\epsilon)n}(X^2+y).\nonumber
\end{equation}
Here we choose $X=(TN_0)^{\frac{1}{6}}$ and obtain
\begin{equation}
	\vert \varphi'(g)\vert^2 \ll_{p,\epsilon} (TN_0)^{\frac{5}{6}+\epsilon} + (TN_0)^{\frac{1}{2}+\epsilon}y.\nonumber
\end{equation}
In the range $\sqrt{y}\leq T^{\frac{1}{6}}N^{\frac{1}{12}}$ this bound is sufficient.
\end{proof}

\section*{Acknowledgements}\label{sec:ack}
We are indebted to Paul Nelson, Akshay Venkatesh, and Farrell Brumley, whose lectures and papers were the direct inspiration for the present article, and who kindly took the time to discuss and share their knowledge with us.

We thank the organisers of Seminar 2322a: Analysis of Automorphic Forms and L-Functions in Higher Rank, at Mathematisches Forschungsinstitut Oberwolfach in 2023 for giving us the opportunity to learn and discuss many of the topics we included in the present article.
The lecture notes \cite{Ven}, \cite{Brum}, \cite{Nel-MFO} were taken there.
We thank the hosting institute for optimal working conditions.
We also thank the organisers of the Automorphic Forms Summer School at the Erdős centre, Budapest, in 2022, where we benefited from Paul Nelson's lectures on the orbit method \cite{Nel-budapest}.
We are glad to have taken part in the Bonn--Paris Learning Seminar, held in autumn 2021 and spring 2022, where we studied \cite{NV} and met Louis Ioos, who shared his expertise in more general aspects of the orbit method with us.

Finally, we thank the referees for a very thorough reading and for many suggestions that significantly improved the quality of our paper.

The first author is supported by the Germany Excellence Strategy grant EXC-2047/1-390685813
and also partially funded by the Deutsche Forschungsgemeinschaft (DFG, German Research
Foundation) – Project-ID 491392403 – TRR 358.

The second author acknowledges support through the European Research Council Advanced Grant 101054336 and the European Union's Horizon 2020 research and innovation programme under the Marie Skłodowska-Curie
grant agreement No 101034255.

\appendix

\section{Bounds via the Fourier expansion}

In this appendix, we prove the following proposition.

\begin{prop}\label{prop:Fourier}
Let $\pi$ be a cuspidal automorphic representation with conductor $N=p^{4n}$, trivial central character and spectral parameter $T$. We assume that $T$ is positive and sufficiently large.\footnote{This assumption can be dropped if one replaces $T$ be $1+\vert T\vert$ in the upper bounds.} Let $\varphi$ be the $L^2$-normalized newform and write $\varphi' = \pi(a(N_0))\varphi$ as in \eqref{eq:shifted_new}, for $N_0=N^{\frac{1}{2}}$. Let $g_{\infty}\in G(\R)$ with $g_{\infty}.i=x+iy\in \mathcal{F}$ and let $g_p\in K_p$ be arbitrary. Then we have
\begin{equation}
	\varphi'(g_{\infty}g_p) \ll_{p,\epsilon} (TN)^{\epsilon} \cdot \left((N_0T)^{\frac{1}{6}}+\frac{(N_0T)^{\frac{1}{2}}}{y^{\frac{1}{2}}}\right).\nonumber
\end{equation}	
\end{prop}

This will be proven using the Fourier/Whittaker expansion of $\varphi$. We will need some preparations. We keep the notation from  Proposition~\ref{prop:Fourier} throughout this section and refer to Section~\ref{subsec:aut_fomrs} for details concerning $\pi$ and $\varphi$. 

Let $g_{\infty}\in G(\R)$ such that
\begin{equation}
	g_{\infty}.i = x+iy\in \mathcal{F} \nonumber
\end{equation}
and write $g_p'=g_p\cdot a(N_0)$. Our goal is to estimate
\begin{equation}
	\varphi'(g_{\infty}g_p) = \varphi(g_{\infty}g_p').\nonumber
\end{equation}

Recall that the global Whittaker function 
\begin{equation}
	W_{\varphi}(g) = \int_{\Q\backslash \A} \varphi(n(x)g)\psi(x)dx
\end{equation}
of $\varphi$ associated to the standard additive character $\psi$ of $\Q\backslash \A$ factors as
\begin{equation}
	W_{\varphi} = c_{\pi}\prod_v W_v.\nonumber
\end{equation}
Here we normalise $W_l(1)=1$ for all finite places $l$ (including $l=p$) and
\begin{equation}
	\vert W_{\infty}(a(q)g_{\infty})\vert = \vert Tqy\vert^{\frac{1}{2}} \dot e^{\frac{\pi t}{2}} \vert K_{iT}(2\pi \vert qy\vert)\vert.\nonumber
\end{equation}
Note that by Rankin-Selberg theory we have
\begin{equation}
	c_{\pi} \ll (TN)^{\epsilon}T^{-\frac{1}{2}}.\nonumber
\end{equation}
(Note that we have re-scaled by $T^{-\frac{1}{2}}$ in order to appropriately normalise the $K$-Bessel function.) Furthermore,
\begin{equation}
	\prod_{l\neq p} W_l(a(q)) = \lambda_{\pi}(m)\cdot m^{-\frac{1}{2}} \text{ for }q=m\cdot p^r \nonumber
\end{equation}
with $m\in \Z_{\neq 0}$ and $r\in \Z$. (If $q$ is not of this shape, then the unramified part vanishes.)

We need to recall some facts concerning the ramified Whittaker new-vector $W_p$. 
Since we are assuming that $\pi$ has level $N=p^{4n}$ and trivial central character, the function $W_p$ is right invariant by the (local) Hecke congruence subgroup $K_{H,p}(4n)$ defined in \eqref{eq:def_KH}. 
Write $K_{H,p}'(4n)\subseteq K_{H,p}(4n)$ for the subgroup consisting of elements $k$ whose bottom right entry $k_{2,2}$ satisfies $k_{2,2}-1\in p^{4n}\Z_p$. By the decomposition given in \cite[(3)]{Sa} every $g\in G(\Q_p)$ can be written as
\begin{equation}
	g\in ZN(\Q_p) g_{t(g),l(g),v(g)}K_{H,p}'(4n),\nonumber
\end{equation}
for $t(g)\in \Z$, $0\leq l(g)\leq 4n$, $v(g)\in \Z_p^{\times}$ and
\begin{equation}
	g_{t,l,v} = \left(\begin{matrix} 0 & p^t \\ -1 & -\frac{v}{p^k}\end{matrix}\right). \nonumber
\end{equation}
Obviously, we have
\begin{equation}
	\vert W_p(g)\vert = \vert W_p(g_{t(g),l(g),v(g)})\vert. \nonumber
\end{equation}

\begin{lemmy}\label{lm:gtl_v}
Let $$g_p=\left(\begin{matrix} t & u \\ v & w\end{matrix}\right)\in K_p$$ and write $g_p'=g_pa(p^{2n})$. Then we have
\begin{equation}
	l(g_p') = \begin{cases} 
		\min(2n+v_p(v),4n) & \text{ if }v_p(w)=0 \\
		\max(0, 2n-v_p(w)) & \text{ if }v_p(w)>0,
	\end{cases} \nonumber
\end{equation}
and $t(g_p') = -2n-2\cdot \min(2n, v_p(v)).$
\end{lemmy}
\begin{proof}
By the Iwasawa decomposition we find $y\in \Q_p^{\times}$ and $k=\begin{psmallmatrix}a & b \\ c & d\end{psmallmatrix}\in K_p$ such that $g_p'\in ZN(\Q_p)a(y)k$. It is the content of \cite[Remark~2.1]{Sa} that $l(g_p')= \min(v_p(c),4n)$ and $t(g_p') = v_p(y)-2l(g_p')$. Therefore, it remains to compute the Iwasawa decomposition of $g_p'$. We first look at the case $w\in \Z_p^{\times}$. Then we can write
\begin{equation}
	\left(\begin{matrix} p^{2n} & \frac{u}{w} \\ 0 & 1 \end{matrix}\right) \left(\begin{matrix} t-\frac{uv}{w} & 0 \\ p^{2n}v& w\end{matrix} \right) = \left(\begin{matrix} p^{2n} t & u \\ p^{2n}v& w\end{matrix}\right) =g_p\cdot a(p^{2n}).\nonumber
\end{equation}
Thus, we have
\begin{equation}
	l(g_p') = \min(2n+v_p(v),4n) \text{ and }t(g_p') = -2n-\min(4n,2v_p(v)). \nonumber
\end{equation}
Next, consider $w=p^sw_0$ for $w_0\in \Z_p^{\times}$ and $s>0$. In this case, we must have $v\in \Z_p^{\times}$. Assume we have an Iwasawa decomposition of the form
\begin{equation}
	\left(\begin{matrix} p^z & 0 \\ 0 & p^z\end{matrix}\right)\left(\begin{matrix} p^y & x \\ 0 & 1 \end{matrix} \right) \left(\begin{matrix} a & b \\ c & d \end{matrix}\right) = g_p'. \nonumber 
\end{equation}
Then the bottom row gives
\begin{equation}
	p^zc = p^{2n}v \text{ and }p^{z}d=w.\nonumber
\end{equation}
Considering the valuations, we obtain
\begin{equation}
	z+v_p(c)=2n \text{ and }z+v_p(d)=s.\nonumber
\end{equation}
But since $ad-bc\in \Z_p^{\times}$, at least one of the two valuations is $0$. Thus, we have two cases. First, if $z=2n$, then $v_p(c)=0$ and $v_p(d) = s-2n\geq 0$. In this case, computing the determinant we find $y=-2n$. Thus,
\begin{equation}
	l(g_p') = 0 \text{ and }t(g_p') = -2n.\nonumber
\end{equation}
On the other hand we can have $z=s$, so that $v_p(d)=0$. Then we have $v_p(c) = 2n-s$ so that we need to assume $2n\geq s$. The determinant forces $y=2n-2s$. Thus,
\begin{equation}
	l(g_p') = 2n-s \text{ and } t(g_p') = -2n. \nonumber
\end{equation}
Combining these concludes the proof.
\end{proof}

To simplify notation we write $l=l(g_p')$ and $t=t(g_p')$. We will assume $g_p$ to be as in Lemma~\ref{lm:gtl_v} and will frequently consider the corresponding cases.

It is easy to check that $$t(a(q)g_{t,l,v}) = t+v_p(q) \text{ and } l(a(q)g_{t,l,v}) = l.$$ We define
\begin{equation}
	s_l=\begin{cases}
		4n &\text{ if }l\leq 2n,\\
		2l &\text{ if }l>2n.
	\end{cases}\nonumber
\end{equation}  
Further write $S_l=p^{s_l}$. All the parameter combinations that appear in the following have been summarised in Table~\ref{table:gtlv_overview} below. Going through each case separately reveals that $-2n+t = -s_l$. 

\begin{table} 
	\centering
	
	\begin{tabular}{l || c | c | c | c }
		& $v_p(v)=0$, & $v_p(v)=0$, & $v_p(v)\leq 2n$,  & $v_p(v)>2n$,  \\
		&  $v_p(w)\geq 2n$ &  $v_p(w)<2n$ &  $v_p(w)=0$ &  $v_p(w)=0$  \\
		\hline\hline 
		$l(g_p')=l$ & $0$   & $2n-v_p(w)$ & $2n+v_p(v)$ & $4n$  \\
		\hline 
		$t(g_p')=t$ & $-2n$ & $-2n$ & $-2n-2v_p(v)$ & $-6n$ \\
		\hline
		$s_l$ & $4n$ & $4n$ & $4n+2v_p(v)$ & $8n$ \\
	\end{tabular}
	\caption{Summary of parameter combinations.}
    \label{table:gtlv_overview}
\end{table}

We return to the study of the local Whittaker function $W_p$. By \cite[Proposition~2.10]{Sa}, we have $W_p(g_{t',l,v}) = 0$ unless $t'\geq -s_l$. In order to control the size of $W_p$ we introduce the quantity $D_l=D_l(\pi_p)\geq 0$ satisfying that, for all $\epsilon>0$, $v\in \Z_p^{\times}$ and all $r\geq 0$, the bound
\begin{equation}
	W_p(g_{-s_l+r,l,v}) \ll p^{(\epsilon-\frac{1}{2})r}\cdot D_l\label{eq:Whittaker_bound_shape}
\end{equation}
holds with an implied constant only depending on $\epsilon$. The existence, i.e. finiteness, of $D_l$ is not obvious, but it follows for example from the results in \cite{Ass19}. The relevant cases will be discussed in more detail below.

Finally, recall the Whittaker expansion
\begin{equation}
	\varphi(g) = \sum_{q\in \Q_{\neq 0}} W_{\varphi}(a(q)g).\nonumber
\end{equation}
We put $g=g_{\infty}g_p'$. In view of the vanishing criteria recalled above, we note that the sum actually runs over $q\in \frac{1}{N_0}\Z_{\neq 0}$.\footnote{This is also contained in \cite[Proposition~2.11]{Sa}. Note that, in the notation of \cite{Sa}, we have $t-q(g_p') = -s_l$, and one can check in the case under consideration here that $q(g_p')=2n$ for all configurations of $l=l(g_p')$. Globally, we refer to \cite[Lemma~3.11]{Sa} and note that $Q^{g_p'}=N_0$.} Inserting all our estimates we obtain
\begin{equation}
	\varphi'(g_{\infty}g_p) \ll (TN)^{\epsilon}D_l\frac{y^{\frac{1}{2}}}{N_0^{\frac{1}{2}}T^{\frac{1}{2}}}\sum_{\substack{m\in \Z_{\neq 0}\\(m,p)=1}} \vert \lambda_{\pi}(m)\vert \sum_{r=0}^{\infty} p^{\epsilon r}\cdot T^{\frac{1}{2}}e^{\frac{\pi T}{2}}\vert K_{iT}(2\pi \frac{\vert m\vert p^ry}{N_0})\vert.\nonumber
\end{equation} 
We start with a slight variation of the standard argument from \cite[Section~3.4]{Sa} anticipated in \cite[Remark~3.9]{Sa}. Define
\begin{equation}
	R = N_0\frac{T+T^{\frac{1}{3}+\epsilon}}{2\pi y} \asymp \frac{N_0T}{y}.\nonumber
\end{equation}
We obtain
\begin{multline}
	\varphi'(g_{\infty}g_p) \ll \\ 
    (TN)^{\epsilon}D_l\frac{y^{\frac{1}{2}}}{N_0^{\frac{1}{2}}T^{\frac{1}{2}}}\sum_{1\leq m\leq R} (m,p^{\infty})^{\epsilon}\vert \lambda_{\pi}(\frac{m}{(m,p^{\infty})})\vert \cdot T^{\frac{1}{2}}e^{\frac{\pi T}{2}}\vert K_{iT}(2\pi \frac{\vert m\vert y}{N_0})\vert + \textrm{Tail}.\nonumber
\end{multline} 
The tail is easily bounded using the exponential decay of the $K$-Bessel function and we will ignore its contribution for now. The remaining part of the sum is handled using Cauchy--Schwarz to separate the Hecke eigenvalues from the Bessel function.

First, we estimate
\begin{equation}
	\sum_{1\leq m\leq R} Te^{\pi T}\vert K_{iT}(2\pi \frac{my}{N_0})\vert^2\ll (NT)^{\epsilon}(T^{\frac{1}{3}}+R).\nonumber
\end{equation}
One can for example follow the proof of \cite[Lemma~3.13]{Sa} to see this. On the other hand, we have
\begin{equation}
	\sum_{1\leq m\leq R} (m,p^{\infty})^{2\epsilon}\vert \lambda_{\pi}(\frac{m}{(m,p^{\infty})})\vert^2 \ll (NT)^{\epsilon}R.\label{eq:lind_on_av}
\end{equation}
This is \cite[Lemma~3.15]{Sa}. Combining everything yields
\begin{equation}
	\varphi'(g_{\infty}g_p)\ll (TN)^{\epsilon}D_l\frac{y^{\frac{1}{2}}R^{\frac{1}{2}}}{N_0^{\frac{1}{2}}T^{\frac{1}{2}}}(T^{\frac{1}{6}}+R^{\frac{1}{2}}) \ll (TN)^{\epsilon}\cdot D_l\cdot (T^{\frac{1}{6}}+\frac{(N_0T)^{\frac{1}{2}}}{y^{\frac{1}{2}}}).\label{eq:fourier_bound_v1}
\end{equation}

We now recall from Section~\ref{subsec:aut_fomrs} that the possible local representations $\pi_p$ come in two different families. They are either supercuspidal or principal series. We consider these cases individually.

\begin{lemmy}\label{lm:A3}
We keep the notation above and assume that $\pi_p$ is supercuspidal. Then we have
\begin{equation}
	\varphi'(g_{\infty}g_p) \ll_{\epsilon} (TN)^{\epsilon}\cdot(T^{\frac{1}{6}}+\frac{(N_0T)^{\frac{1}{2}}}{y^{\frac{1}{2}}}).\nonumber
\end{equation}
\end{lemmy}
\begin{proof}
This will follow from \eqref{eq:fourier_bound_v1} after we establish that $D_l \ll 1$. This can be extracted from the proof of \cite[Lemma~5.1]{Ass19} and we refer to \cite[Lemma~3.4.1]{Ass_thesis} for more details. See also \cite[Section~3.4.4]{Ass_thesis} for a summary. Let us give some more details.

First, for $l\neq 2n$, we always have $D_l\ll 1$. This is established in Case~I and Case~III of the proof of \cite[Lemma~5.1]{Ass19}. We need to take a closer look at $l=2n$, covered in Case~II of \cite[Lemma~5.1]{Ass19}. We recall that since $\pi$ is supercuspidal and we are assuming that it has trivial central character, it follows that $\pi$ is twist minimal.  However, by analysing the critical points appearing in the stationary phase analysis for $W_{p}(g_{-4n,2n,v})$, one finds that these are non-degenerate for twist minimal supercuspidal representations $\pi_p$. See \cite[Remark~3.4.3]{Ass_thesis} for details. Thus, we also obtain $D_{2n}\ll 1$ and the proof is complete.
\end{proof}

\begin{lemmy}\label{lm:A4}
We keep the notation above and assume that $\pi_p$ is principal series. Then we have
\begin{equation}
	\varphi'(g_{\infty}g_p) \ll_{\epsilon} (TN)^{\epsilon}\cdot((TN_0)^{\frac{1}{6}}+\frac{(N_0T)^{\frac{1}{2}}}{y^{\frac{1}{2}}}).\nonumber
\end{equation}
\end{lemmy}
\begin{proof}
We first note that for $l\neq 2n$ we have $D_l\ll 1$, so that the desired bound follows from \eqref{eq:fourier_bound_v1}. This property of $D_l$ can be extracted from the proof of \cite[Lemma~5.12]{Ass19}. See also \cite[Lemma~3.4.16 and Section~3.4.4]{Ass_thesis}.

We have to be more careful with the case $l=2n$. Here we first observe that one has
\begin{equation}
	W_p(g_{-s_l+r,l,v}) \ll 1 \label{eq:Whittaker_bound_improved}
\end{equation}
as long as $r>0$. This follows by combining \cite[Lemma~3.6]{Ass19} with the estimates from \cite[Lemma~5.11]{Ass19}. A compact summary is also contained in \cite[Section~3.4.4]{Ass_thesis}. Thus, we estimate the Whittaker expansion as
\begin{equation}
	\varphi'(g_{\infty}g_p) \ll (TN)^{\epsilon} (\mathcal{S}_0 + \mathcal{S}_{\geq 1}) + \textrm{Tail}, \nonumber
\end{equation} 
where
\begin{equation}
	\mathcal{S}_{0} = \frac{y^{\frac{1}{2}}}{N_0^{\frac{1}{2}}T^{\frac{1}{2}}}\sum_{\substack{1\leq m\leq R\\ (m,p)=1}} \vert \lambda_{\pi}(m)\vert\cdot \vert W_p(a(m)g_{-4n,2n,v(g_p')})\vert \cdot T^{\frac{1}{2}}e^{\frac{\pi T}{2}}\vert K_{iT}(2\pi \frac{\vert m\vert y}{N_0})\vert \nonumber
\end{equation}
and
\begin{equation}
	\mathcal{S}_{\geq 1} = \frac{y^{\frac{1}{2}}}{N_0^{\frac{1}{2}}T^{\frac{1}{2}}}\sum_{\substack{1\leq m\leq R\\ p\mid m}} (m,p^{\infty})^{\frac{1}{2}}\vert \lambda_{\pi}(\frac{m}{(m,p^{\infty})})\vert \cdot T^{\frac{1}{2}}e^{\frac{\pi T}{2}}\vert K_{iT}(2\pi \frac{\vert m\vert y}{N_0})\vert.\nonumber
\end{equation}
The contribution of $S_{\geq 1}$ is easy to handle.  Indeed, we follow exactly the argument above. The only difference is that we replace \eqref{eq:lind_on_av} by the estimate
\begin{equation}
	\sum_{1\leq m\leq R} (m,p^{\infty})\vert \lambda_{\pi}(\frac{m}{(m,p^{\infty})})\vert^2 \ll (NT)^{\epsilon}R,
\end{equation}
which is exactly \cite[Lemma~3.15]{Sa}. One arrives at
\begin{equation}
		\mathcal{S}_{\geq 1} \ll (TN)^{\epsilon}\cdot(T^{\frac{1}{6}}+\frac{(N_0T)^{\frac{1}{2}}}{y^{\frac{1}{2}}}). \nonumber
\end{equation}

We turn towards $\mathcal{S}_0$. There are two cases to consider. First, if $-1$ is not a square in $\Z_p^{\times}$, then it is shown in \cite[Lemma~5.12]{Ass19} that we have the bound $\vert W_p(a(m)g_{-4n,2n,v(g_p')})\vert\ll 1$. In this case we can also argue as above and obtain the strong bound
\begin{equation}
	\mathcal{S}_{0} \ll (TN)^{\epsilon}\cdot(T^{\frac{1}{6}}+\frac{(N_0T)^{\frac{1}{2}}}{y^{\frac{1}{2}}}). \nonumber
\end{equation}
Putting everything together, we see that
\begin{equation}
	\varphi'(g_{\infty}g_p) \ll_{\epsilon} (TN)^{\epsilon}\cdot(T^{\frac{1}{6}}+\frac{(N_0T)^{\frac{1}{2}}}{y^{\frac{1}{2}}}) \nonumber
\end{equation}
if $p\equiv 3\text{ mod 4}$, i.e. if $-1$ is not a square in $\Z_p^{\times}$.

Finally, we assume that $p\equiv 1\text{ mod }4$.  Thus, there is $i_p\in \Z_p^{\times}$ with $i_p^2=-1$. In this case we can not estimate $W_p(a(n)g_{-4n,2n,v(g_p')})$ trivially, since it can take large values.

By \cite[(3)]{Sa} we have
\begin{equation}
	\vert W_{p}(g_{-4n,2n,v})\vert = \vert W_{p}(g_{-4n,2n,v'})\vert \nonumber
\end{equation}
for $v-v'\in p^{2n}\Z_p$. Furthermore, if $(m,p)=1$, then we  view $m\in \Z_p^{\times}$ and one obtains
\begin{equation}
	\vert W_{p}(a(m)g_{-4n,2n,v})\vert = \vert W_{p}(g_{-4n,2n,vm^{-1}})\vert.\nonumber
\end{equation}
Thus we define a function
\begin{equation}
	f\colon (\Z/p^{2n}\Z)^{\times} \to \R_{\geq 0},\, m\mapsto \vert  W_{p}(g_{-4n,2n,v(g_p')m^{-1}})\vert.\nonumber 
\end{equation}
Using Cauchy--Schwarz, \eqref{eq:lind_on_av}, and the definition of $R$ we conclude that
\begin{equation}
	\mathcal{S}_0 \ll (NT)^{\epsilon}\left(\sum_{\substack{1\leq m\leq R\\ (m,p)=1}} f(m)^2 Te^{\pi T}\vert K_{iT}(2\pi \frac{my}{N_0})\vert^2 \right)^{\frac{1}{2}}.\label{eq:need_to_bound}
\end{equation}
In \cite[Lemma~3.12]{Sa} one can find average bounds for $f(m)^2$. However, these will not be sufficient for our purposes. We extract more detailed properties of $f(m)$ from the proof of \cite[Lemma~5.12]{Ass19}. See also the proof of \cite[3.4.16]{Ass_thesis} for more details.

The size of $f$ is dictated by the discriminant
\begin{equation}
	\Delta(m) = 1 + \frac{4v(g_p')^2b_{\chi}^2}{m^2} \in \Z_p.\nonumber
\end{equation}
Here $b_{\chi}\in \Z_{p}^{\times}$ depends only on $\chi$. We first record the bound
\begin{equation}
	f(m) \ll p^{\frac{n}{3}}, \nonumber
\end{equation}
which is stated in \cite[Lemma~5.12]{Ass19} and holds for all $m$. On the other hand, examining the argument in \cite[Case~IV, Lemma~5.11]{Ass19} shows that we also have
\begin{equation}
	f(m)\ll p^{\frac{v_p(\Delta(m))}{2}}.\nonumber
\end{equation}
We set $v_{\pm} = \pm 2i_pv(g_p')b_{\chi}\in \Z_p^{\times}$ and observe that 
\begin{equation}
	\Delta(m) = \frac{1}{m^2}(m-v_+)(m-v_-).\nonumber
\end{equation} 
Thus,  if $v_p(\Delta(m))=\delta$, then 
\begin{equation}
	m\equiv v_+ \text{ mod }p^\delta \text{ or }m\equiv v_- \text{ mod }p^{\delta}.\nonumber
\end{equation}
Put $\alpha_r = r/2$ if $r<\frac{2n}{3}$ and $\alpha_r=n/3$ if $r\geq \frac{2n}{3}$. In view of these observations and \eqref{eq:need_to_bound} we can now estimate
\begin{equation}
	\mathcal{S}_0 \ll  (NT)^{\epsilon} \left( \sum_{\pm} \sum_{r=0}^{\lceil \frac{2n}{3}\rceil} p^{2\alpha_r}\mathcal{S}_r\right)^{\frac{1}{2}}, \nonumber
\end{equation}
for
\begin{equation}
	\mathcal{S}_r = \sum_{\substack{1\leq m\leq R\\ m\equiv v_{\pm} \text{ mod }p^r}} Te^{\pi T}\vert K_{iT}(2\pi \frac{my}{N_0})\vert^2.\nonumber
\end{equation}
Note that for $0\leq r <\frac{2n}{3}$ we could have restricted the sums further to $m\not \equiv v_{\pm} \text{ mod }p^{r+1}$. We have dropped this condition because it is irrelevant for our estimates. The sums $\mathcal{S}_r$ can be estimated by
\begin{equation}
	\mathcal{S}_r \ll  T^{\frac{1}{3}}+\frac{R}{p^r}.\nonumber
\end{equation}
This follows by a slight modification of the argument from \cite[Lemma~3.13]{Sa}. Inserting this above yields
\begin{equation}
	\mathcal{S}_0 \ll  (NT)^{\epsilon} \left( \sum_{\pm} \sum_{r=0}^{\lceil \frac{2n}{3}\rceil} (p^{2\alpha_r}T^{\frac{1}{3}} + Rp^{2\alpha_r-r})\right)^{\frac{1}{2}} \ll (NT)^{\epsilon}((N_0T)^{\frac{1}{6}}+R^{\frac{1}{2}}), \nonumber
\end{equation}
where we have used that $2\alpha_r-r \leq 0$. This was the last sum we had to estimate and it is of the predicted size. Thus, the proof is complete.
\end{proof}

The statement from Proposition~\ref{prop:Fourier} follows from Lemma~\ref{lm:A3} and Lemma~\ref{lm:A4} after recalling from Remark~\ref{rem:conductor_class} that these cover all cases that can occur.

\printbibliography %

\end{document}